\documentclass[authoryear,review,12pt]{elsarticle}
\usepackage{appendix}
\usepackage{comment}
\usepackage{setspace}
\usepackage[margin=0cm]{caption}
\usepackage[font={sf,bf},textfont=md]{caption}
\usepackage[tbtags]{amsmath}
\usepackage{amsthm,amssymb,amsfonts}
\usepackage{natbib}
\usepackage{multirow}
\usepackage{lscape}
\usepackage{graphicx}
\usepackage{subcaption}
\usepackage{makecell}
\usepackage{exscale}
\usepackage{booktabs}
\usepackage{array}
\usepackage{fullpage}
\usepackage{url}
\usepackage{algorithm}
\usepackage{algpseudocode}
\usepackage{bm}
\usepackage{mathtools}
\usepackage{wrapfig}
\usepackage{lipsum}
\usepackage{mathrsfs}
\usepackage{relsize}
\usepackage{dsfont}
\usepackage{multirow}
\usepackage{apalike}
\usepackage[usenames,dvipsnames,svgnames,table]{xcolor}
\usepackage[colorlinks=true,
linkcolor=blue,
urlcolor=blue,
citecolor=blue]{hyperref}

\title{Testing the number of common factors by bootstrapped sample covariance matrix  in high-dimensional factor models and errata}
\date{}

\theoremstyle{plain}
\newtheorem{axiom}{Axiom}
\newtheorem{claim}[axiom]{Claim}
\newtheorem{theorem}{Theorem}[section]
\newtheorem{lemma}[theorem]{Lemma}
\newtheorem{cor}{Corollary}
\newtheorem{assumption}{Assumption}
\newtheorem{remark}{Remark}
\theoremstyle{remark}
\newtheorem{definition}[theorem]{Definition}

\newcommand{\bb}{\bm{b}}

\newcommand{\be}{\bm{e}}
\newcommand{\bbf}{\bm{f}}

\newcommand{\bh}{\bm{h}}

\newcommand{\bt}{\bm{t}}
\newcommand{\bu}{\bm{u}}
\newcommand{\bv}{\bm{v}}
\newcommand{\bw}{\bm{w}}
\newcommand{\bx}{\bm{x}}
\newcommand{\by}{\bm{y}}
\newcommand{\bz}{\bm{z}}

\newcommand{\Ab}{\mathbf{A}}
\newcommand{\Bb}{\mathbf{B}}
\newcommand{\Cb}{\mathbf{C}}
\newcommand{\Db}{\mathbf{D}}
\newcommand{\Eb}{\mathbf{E}}
\newcommand{\Fb}{\mathbf{F}}
\newcommand{\Gb}{\mathbf{G}}
\newcommand{\Hb}{\mathbf{H}}
\newcommand{\Ib}{\mathbf{I}}

\newcommand{\Kb}{\mathbf{K}}
\newcommand{\Lb}{\mathbf{L}}
\newcommand{\Mb}{\mathbf{M}}

\newcommand{\Rb}{\mathbf{R}}
\newcommand{\Sbb}{\mathbf{S}}

\newcommand{\Ub}{\mathbf{U}}
\newcommand{\Vb}{\mathbf{V}}
\newcommand{\Wb}{\mathbf{W}}
\newcommand{\Xb}{\mathbf{X}}
\newcommand{\Yb}{\mathbf{Y}}
\newcommand{\Zb}{\mathbf{Z}}

\newcommand{\bC}{\bm{C}}

\newcommand{\bL}{\bm{L}}

\newcommand{\bbeta}{\bm{\beta}}
\newcommand{\bgamma}{\bm{\gamma}}
\newcommand{\bepsilon}{\bm{\epsilon}}

\newcommand{\bpsi}{\bm{\psi}}

\newcommand{\bGamma}{\bm{\Gamma}}

\newcommand{\bLambda}{\bm{\Lambda}}

\newcommand{\bPi}{\bm{\Pi}}
\newcommand{\bSigma}{\bm{\Sigma}}

\newcommand{\bPhi}{\bm{\Phi}}
\newcommand{\bPsi}{\bm{\Psi}}

\begin{document}
	\begin{frontmatter}

		\author[myfirstaddress]{Yu Long}
		\ead{fduyulong@163.com}
		
				\author[mysecondaddress]{Zhao Peng}
		\ead{zhaop@jsnu.edu.cn}
		
		\author[mythirdaddress]{Zhou Wang}
		\ead{wangzhou@nus.edu.sg}
				\address[myfirstaddress]{School of Statistics and Management, 
					Shanghai University of Finance and Economics, China.}
		\address[mysecondaddress]{School of Mathematics and Statistic 
			and Jiangsu Provincial Key Laboratory of
			Educational Big Data Science and Engineering,
			Jiangsu Normal University, Xuzhou, Jiangsu, China}
				\address[mythirdaddress]{Department of Statistics and Data Science, National University of Singapore,
			Singapore.}

		\begin{abstract}
		This paper studies the impact of bootstrap procedure on the eigenvalue distributions of the sample covariance matrix under a high-dimensional factor structure.
We provide asymptotic distributions for the top  eigenvalues of bootstrapped sample covariance matrix under mild conditions. After bootstrap, the spiked eigenvalues which are driven by common factors will converge weakly to Gaussian limits after proper scaling and centralization. However, the largest non-spiked eigenvalue is mainly determined by the order statistics of the bootstrap resampling weights, and follows extreme value distribution. Based on the disparate behavior of the spiked and non-spiked eigenvalues, we propose innovative methods to test the number of common factors. {Indicated by extensive numerical and empirical studies, the proposed methods perform  reliably and convincingly under the existence of both weak factors and cross-sectionally correlated errors.}
Our technical details  contribute to random matrix theory on spiked covariance model with convexly decaying density and unbounded support, or with general elliptical distributions.
		\end{abstract}
		
		
	\end{frontmatter}

\section{Introduction}
High-dimensional factor models have attracted growing attention in the recent decades with fruitful applications in statistical learning problems such as covariance matrix estimation, forecasting and model selection. A comprehensive overview of some recent advances in factor models is summarized by \cite{fan2021recent}.   A fundamental step in factor analysis is to determine the number of common factors, which is extensively studied in the literature {and still actively debated. For example, in finance and econometrics, it remains an open question to judge whether a new factor adds explanatory power for asset pricing; see \cite{feng2020taming}. Dropping important factors will result in non-negligible estimation error for the factor scores and loading spaces; see \cite{bai2002determining,bai2003inferential}. It also drops information which can be potentially useful in related statistical applications such as detecting structural breaks in \cite{Baltagi2017}.  On the other hand, overestimating the factor number may result in non-negligible errors, too; see \cite{Barigozzi2020}. In high-dimensional settings, overestimation also increases the computational burden. 
	
	In econometrics, most of the existing methods for determining the number of factors are based on different growth rates of the factor and noise eigenvalues. That is, the eigenvalues of the population covariance matrix driven by common factors will diverge to infinity with a significantly faster rate than those driven by idiosyncratic errors. To list a few examples, the information criterion by \cite{bai2002determining} and its improved version by \cite{alessi2010improved}, the eigenvalue ratio approach by \cite{lam2012factor} and \cite{ahn2013eigenvalue} are widely used to estimate factor number in high dimensions. Along this line, sometimes testing procedures can be more preferred than simply providing a point estimation for the number of factors; see the randomized test in \cite{trapani2018randomized} and the random-perturbation-based rank estimator by \cite{kong2020random}. These methods provide  significance level of the corresponding determination by introducing additional randomness into the system.

	In high dimensional statistics, to determine the number of factors, another important line relies on random matrix theory (RMT) on the largest non-spiked eigenvalues of the sample covariance matrix. It has been shown that they follow the Tracy-Widom law asymptotically after proper centralization and scaling, so one can  test the number of common factors based on this property. This line dates back to \cite{onatski2009testing}, which deals with eigenvalues and spectral densities under the generalized dynamic factor model by \cite{forni2000generalized}. More recent extensions are the eigenvalue thresholding approaches in \cite{onatski2010determining}, \cite{cai2020limiting} and \cite{ke2021estimation}. These methods usually propose milder conditions on the strength of factors, but more restrictive assumptions on the dependence structure of noises. The spiked eigenvalues can be specified even if they are not diverging, as long as they exceed the typical BBP phase transition boundary; see \cite{baik2005phase}. Another closely-related direction is the parallel analysis, which can be a special application of RMT; see \cite{dobriban2019deterministic} and references therein.
	
	One major limitation of the RMT-based methods is that the underlying Tracy-Widom distribution is pretty complicated, depending on unknown parameters of the population covariance matrix. The aforementioned approaches usually need to estimate these parameters first, although \cite{onatski2009testing} avoids this problem by transformation. In statistics,  it's well known that bootstrap is a common way for approximating complicated distributions. A natural question arises: \emph{is it possible to approximate the asymptotic distribution of the sample eigenvalues by bootstrap?} In the current paper, we aim to answer this question by studying the impact of bootstrap procedure on the limiting distributions of top sample eigenvalues. Based on the findings, we further propose new test-based methods to determine the number of common factors.
	
}

Bootstrapping the sample covariance matrix is also considered in \cite{karoui2016bootstrap}, under a scaled spiked covariance model. Based on their results, the limiting distributions of spiked sample eigenvalues can be consistently approximated by bootstrap only when the population spiked eigenvalues are well separated from the non-spiked ones. \cite{yao2021rates} further relaxes the technical conditions in \cite{karoui2016bootstrap}, and provides an upper bound for the bootstrap bias in terms of  the tail probability of the eigenvalue distribution before and after bootstrap. In both papers, the bootstrap works only when the spiked eigenvalues are sufficiently large.  Other related works have also considered  bootstrapping the operator norm (\cite{Han2018}) or spectral statistics (\cite{Lopes2019}) of sample covariance matrix, but they are different from the current paper.

The bootstrapped sample covariance matrix considered in this paper is also closely related to the separable covariance model in the literature of RMT; see \cite{bai2019central} and \cite{ding2021spiked}. It can be written as $\hat\Sbb:=n^{-1}\Ab\Zb\Wb\Zb^\top\Ab^\top$, where $n$ is the sample size, $\bSigma=\Ab\Ab^\top$ is the population covariance matrix, $\Zb$ is a  random matrix with independent entries and $\Wb$ is a diagonal  matrix composed of  the bootstrap resampling weights. If $\Wb$ is the identity matrix, it reduces to the traditional sample covariance matrix without bootstrap, which has been extensively studied; see for example \cite{ding2018necessary} and \cite{cai2020limiting}. {In this paper, $\Wb$ is a diagonal random matrix. \cite{johansson2007gumbel} has studied the largest eigenvalue of a random diagonal matrix  plus a scaled Gaussian Unitary Ensemble matrix, but it's different from the case of bootstrap.}

\subsection{Our contributions}
Firstly, in view of bootstrapping the sample covariance matrix, this paper is a valuable supplement to \cite{karoui2016bootstrap} and \cite{yao2021rates}.  For the spiked sample eigenvalues after bootstrap, following \cite{karoui2016bootstrap} and \cite{yao2021rates} it remains unknown what happens when the population eigenvalues are weak. Our Corollary \ref{cor: bootstrap bias} fills this gap by providing asymptotic limit rather than upper bound for the bootstrap bias.  We also provide limiting distributions for the spiked sample eigenvalues with explicit formulas for the scaling and centralization parameters in the paper. Moreover, for the non-spiked sample eigenvalues, it is the first time that the limiting distributions after bootstrap have been revealed. The results not only contribute to factor models, but also to principal component analysis (PCA) or more general spiked covariance models. We also relax some technical conditions in \cite{karoui2016bootstrap} and \cite{yao2021rates}, e.g., we allow the spiked eigenvalues to diverge with different rates.

Secondly, the theoretical framework of the current paper is totally different from that of the separable covariance model in \cite{bai2019central} and \cite{ding2021spiked}. They require that the limiting spectral density of $\hat\Sbb$ exhibits  the ``square root'' type behavior around the edge of its support. Under bootstrap, the spectral density of $\Wb$ is usually convex at the edge and the ``square root'' characteristic does not hold anymore. This is the major reason why the non-spiked eigenvalues after bootstrap converge to extreme value distributions rather than the Tracy-Widom law. Spectral property of sample covariance matrix with convexly decaying density has been studied in \cite{kwak2021extremal}. Unfortunately, they require $\bSigma$ to be the identity matrix and the spectral distribution of $\Wb$ has bounded support, which excludes the case of bootstrap. Up to our knowledge, we are the first to consider general $\bSigma$ with spiked eigenvalues, allowing the spectral distribution of $\Wb$ to have unbounded support. $\hat\Sbb$ is also closely related to  elliptical distributions if the columns of $\Zb$ are from Gaussian  or spherical distribution. Some related results can be found in \cite{hu2019high} and \cite{wen2019tracy}, where the variances of the entries in $\Wb$ are required to be nearly 0. In the current paper, the diagonal entries of $\Wb$ are from non-degenerated distributions. Therefore, our technical  details will also contribute to the spectral analysis of sample covariance matrix with general elliptical distributions. 

{Thirdly, in practice, we provide new direction for testing the number of common factors, which is useful in very general scenarios. The proposed  approaches in this paper are accurate even if the factors are weak. Moreover, we allow the existence of bounded outliers in the spectrum of the idiosyncratic error covariance matrix. This is a major difference of our approaches from the existing RMT-based methods. As a sacrifice, we require the eigenvalues driven by common factors to be diverging, which is more stringent than typical assumptions in the RMT literature.  The reason is that	the phase transition boundary after bootstrap is determined mainly by the order statistics of the bootstrap resampling weights. In other words, it's possible to increase or decrease the typical BBP phase transition boundary by bootstrapping from different distribution families. Therefore, this paper also provides a new direction for documenting the number of factors with different strength.} 

\emph{Conventions}. $c$ denotes a small positive constant varying in different lines. $[a]$ denotes the largest integer not larger than $a$. $a_n\lesssim b_n$ means $a_n\le c^{-1}b_n$ (or $a_n\le O_p(b_n)$ if $a_n$ or $b_n$ is random) while $a_n\asymp b_n$ means that $c\le a_n/b_n\le c^{-1}$ for sufficiently large $n$. For a (Hermitian) matrix $\Ab$, $\text{tr}(\Ab)$ denotes the trace, $\lambda_i(\Ab)$ denotes the $i$-th largest eigenvalue. $\overset{d}{\rightarrow}$ and $\overset{p}{\rightarrow}$ are for convergence in distribution and probability, respectively.   $(n\vee p)=\max\{n,p\}$. $\|\cdot\|$ and $\|\cdot\|_F$ are for spectral and Frobeniuos norms, respectively.

\section{Factor model and bootstrap}\label{sec: factor model}
We consider high-dimensional factor model which can be written in the form of
\begin{equation}\label{factor model}
	x_{ij}=\bL_i^\top\bbf_j+\bpsi_{i}^\top\bepsilon_j, \quad 1\le i\le p, 1\le j\le n,
\end{equation}
where $\bL_i$'s are $r$-dimensional factor loadings, $\bbf_j$'s are $r$-dimensional latent factor scores, $\bpsi_i$'s are  $p$-dimensional deterministic vectors and $\bepsilon_j$'s are $p$-dimensional idiosyncratic errors. $r$ is the number of common factors, which is of the major interest in the current paper. The model can also be written in matrix form as $\Xb=\Lb\Fb^\top+\bPsi\Eb:=\Ab\Zb$,
where $\Xb=(\bx_1,\ldots,\bx_n)=(x_{ij})_{p\times n}$, $\Lb^\top=(\bL_1,\ldots,\bL_p)$, $\Fb^\top=(\bbf_1,\ldots,\bbf_n)$, $\bPsi^\top=(\bpsi_{1},\ldots,\bpsi_{p})$, ${\Eb_{p\times n}}=(\bepsilon_1,\ldots,\bepsilon_n)$, ${\Ab_{p\times (r+p)}}=(\Lb,\bPsi)$ and $\Zb^\top=(\Fb,\Eb^\top)$.  Some assumptions are given as follows.

\begin{assumption}\label{c1}
	There exist a constant $0<c\le 1$  such that:\\
	(a). $\Fb=\Cb\Fb^0$ for some $n\times n$ deterministic matrix $\Cb$ satisfying $\|\Cb\|\le c^{-1}$  and $n^{-1}\|\Cb\|_F^2=1$. The entries of $\Fb^0$ and $\Eb$ are independent (not necessarily identically distributed) real-valued random variables with mean 0, variance 1 and bounded eighth moments. \\
	(b). $r$ is fixed as $\min\{n,p\}\rightarrow \infty$.\\
	(c). $\|\bPsi\|\le c^{-1}$, $\lambda_{[c p]}(\bPsi\bPsi^\top)\ge c$, $[n\lambda_i(\Lb^\top\Lb)]^{-1}(n\vee p)\log n=o(1)$ for any $1\le i\le r$, and $\lambda_i(\Lb^\top\Lb)/\lambda_{i+1}(\Lb^\top\Lb)\ge 1+c$ for any $1\le i\le r-1$.
\end{assumption}

We assume a separable structure in (\ref{factor model}) for the idiosyncratic errors, which is common in the literature especially when the non-spiked eigenvalues are of concern. Similar assumptions are found in \cite{cai2020limiting} and \cite{ke2021estimation}. Assumption \ref{c1}(a)  requires  bounded eighth moments mainly to ensure  we can find proper estimators for the asymptotic variances of the bootstrapped sample eigenvalues.  This assumption can be potentially relaxed to bounded fourth moments using truncation technique as in \cite{cai2020limiting}. We don't pursue this direction in the current paper. We assume {$\Eb$ to have independent entries but allow the factor process to be serially dependent under a separable scheme through the matrix $\Cb$, covering partially the auto-regressive and moving average processes. The condition $n^{-1}\|\Cb\|_F^2=1$ is for identification. }
Assumption \ref{c1}(b) assumes fixed $r$ which is  common in the literature, especially when the target is to determine the number of common factors; see for example \cite{onatski2009testing}, \cite{ahn2013eigenvalue} and \cite{ke2021estimation}. 

Assumption \ref{c1}(c) deserves more explanation.  The spectral norm of $\bPsi$ is bounded so that the idiosyncratic errors are asymptotically negligible compared with the common factors.  The condition  $[n\lambda_i(\Lb^\top\Lb)]^{-1}(n\vee p)\log n=o(1)$  ensures that the common factors dominate in the system, which is critical especially under the high-dimensional settings when $p\gg n$; see also \cite{wang2017asymptotics} and \cite{cai2020limiting}. {Usually $\lambda_i(\Lb^\top\Lb)$ depends on the dimension $p$, thus this condition can also be viewed as a constraint on the growth rates of $n$ and $p$. It also shows how we identify a common factor in this paper. When $p\asymp n$, we believe that a spiked eigenvalue is driven by a common factor only when it's diverging at rate larger than $\log p$, so that it has non-negligible effects on a number of variables in the system. This condition is  slightly more stringent compared with those in the RMT literature. It should be acknowledged that much of the statistical literature does not require growing spikes, while applied studies in econometrics usually make stronger assumptions on the relative growth rate of factor and noise eigenvalues. We follow the latter to ensure that large idiosyncratic noise will not be identified as common factor. Consider a toy example where $x_{ij}=\epsilon_{ij}\sim\mathcal{N}(0,\sigma_i^2)$ independently with $\sigma_1^2=(1+c)(1+\sqrt{p/n})^2$ for a constant $c>0$ while $\sigma_i^2=1$ for $i\ne 1$.  Then, there is no ``common'' factor at all although the leading eigenvalue exceeds the BBP-type phase transition boundary. Assumption \ref{c1}(c) helps avoid such mis-specification.}  The eigenvalues of $\Lb^\top\Lb$ are assumed to be distinct and allowed to diverge with different rates, so that the corresponding eigenvectors are identifiable. 

{Let $\bSigma=\Ab\Ab^\top$. If $\Cb$ is the identity matrix,  $\{\bx_j\}$ will be a stationary process such that $\bSigma$ is the population covariance matrix $\mathbb{E}(\bx_1\bx_1^\top)$.} Under Assumption \ref{c1}, $\bSigma$ has $r$  spiked eigenvalues  significantly larger than  the remaining non-spiked  ones. Similar property holds for the sample covariance matrix, making it possible to estimate or test the number of factors. {The exact limiting behavior of the sample eigenvalues is usually complicated and potentially dependent on unknown parameters.} In this paper, we are interested in bootstrapping the observations $(\bx_1,\ldots,\bx_n)$ and studying the eigenvalues of the bootstrapped sample covariance matrix. A standard bootstrap procedure resamples the columns of $\Xb$ with replacement. Each column is chosen  with probability $n^{-1}$ in each run. We repeat the resampling procedure $n$ times independently to obtain a new $p\times n$ matrix. Then, the bootstrapped sample covariance matrix  can be written as 
\begin{equation}\label{hat S}
	\hat\Sbb =n^{-1}\sum_{j=1}^nw_j\bx_j\bx_j^\top=n^{-1}\Xb\Wb\Xb^\top,
\end{equation}
where {$\Xb$ is the original data and} $\Wb=\text{diag}(w_1,\ldots,w_n)$ is a diagonal matrix with $w_j$'s being the corresponding resampling weights.
We define two types of bootstrap procedures.
\begin{definition}
	We say that $\hat\Sbb$ is from a multiplier bootstrap procedure, if  in (\ref{hat S}) $w_j$'s are independent and identically distributed (i.i.d.) from  exponential distribution $Exp(1)$. We say that $\hat\Sbb$ is from a standard bootstrap procedure if $\bw=(w_1,\ldots,w_n)$ follows $n$-dimensional multinomial distribution with n trials and event probabilities $(n^{-1},\ldots,n^{-1})$.
\end{definition}
Under multiplier bootstrap, $w_j$'s are i.i.d. which simplifies the technical proofs. We use exponential distribution to ensure that $\hat\Sbb$ is semi-positive definite, while our approaches can be extended to more general distribution families {such as Possion.} Under standard bootstrap, $w_j$'s are no longer independent but still identically distributed. The expectation and covariance satisfy
$\mathbb{E}(w_j)=1$, $\text{Var}(w_j)=1-n^{-1}$, and $ \text{Cov}(w_j,w_l)=-n^{-1}\text{ for }j\ne l$.
In the following, we will study the limiting distributions of both the spiked and non-spiked eigenvalues of the bootstrapped sample covariance matrix, and accordingly propose test-based procedures to determine the number of common factors.

\section{Testing with spiked eigenvalues}\label{sec: spiked}
\subsection{Limiting representation}
The non-zero eigenvalues of $\hat\Sbb$ are the same as those of its companion matrix, defined by 
$\hat{\mathcal{S}}=n^{-1}\Wb^{1/2}\Zb^\top\Ab^\top\Ab\Zb\Wb^{1/2}$.
Further define the eigenvalue decomposition 
\begin{equation}\label{A top A}
	\Ab^\top\Ab=\bGamma\bLambda\bGamma^\top=\bGamma_1\bLambda_1\bGamma_1^\top+\bGamma_2\bLambda_2\bGamma_2^\top,
\end{equation}
where $\bGamma=(\bGamma_1,\bGamma_2)=(\bgamma_1,\ldots,\bgamma_{r+p})$ is the eigenvector matrix, $\bLambda=\text{diag}(\lambda_1,\ldots,\lambda_{r+p})$ is composed of the eigenvalues in descending order. $\bGamma_1$ and $\bLambda_1$ are associated with the leading $r$ eigenvectors and eigenvalues, respectively. Let $\hat\lambda_i$ be the $i$th largest eigenvalue of $\hat{\mathcal{S}}$. The next lemma shows some preliminary properties of $\hat\lambda_i$.

\begin{lemma}\label{lem: preliminary}
	Under Assumption \ref{c1}, as $\min\{n,p\}\rightarrow \infty$ we have $(n\lambda_i)^{-1}(n\vee p)\log n=o(1)$ and $\lambda_i/\lambda_{i+1}\ge 1+c$ for any $1\le i\le r$ while $c\le \lambda_{[cp]}\le\cdots\le\lambda_{r+1}\le c^{-1}$ for some $c>0$. Further,
	no matter under the multiplier or standard bootstrap, we always have
	\[
	\hat\lambda_i/\lambda_i-1=O_{p}\bigg(\frac{(n\vee p)\log n}{n\lambda_i}+\frac{1}{\sqrt{n}}\bigg),\quad i\le r,\quad\text{and}\quad\hat\lambda_{r+1}\le O_{p}\bigg(\frac{(n\vee p)\log n}{n}\bigg).
	\]
\end{lemma}

By Lemma \ref{lem: preliminary}, $\hat\lambda_i/\lambda_i$ converges to 1 for $1\le i\le r$. 
However, the convergence rate can be very slow and Lemma \ref{lem: preliminary} is not very helpful in deriving  distributional property.  Motivated by \cite{cai2020limiting}, we define $\theta_i$ as the solution to 
\[
\frac{\theta_i}{\lambda_i}=\bigg[1-\frac{1}{n\theta_i}\sum_{k=1}^p\frac{\lambda_{r+k}}{1-\lambda_i^{-1}\lambda_{r+k}}\bigg]^{-1},\quad \theta_i\in [\lambda_i,2\lambda_i],\quad 1\le i\le r.
\]
Under Assumption \ref{c1}, the existence and uniqueness of $\theta_i$ can be verified easily by the mean value theorem. \cite{cai2020limiting} has shown that $\theta_i$ is a {closer approximation} to the associated eigenvalue of the sample covariance matrix without bootstrap compared with $\lambda_i$. Under our settings, to address the effect of bootstrap, let $\hat\zeta_i$ be the solution to 
\[
\hat\zeta_i=\frac{1}{n}\sum_{j=1}^nw_j\bigg[1-\frac{w_j}{n\theta_i}\sum_{k=1}^p\frac{\lambda_{r+k}}{1-\theta_i^{-1}\lambda_{r+k}\hat\zeta_i}\bigg]^{-1},\quad \hat\zeta_i\in \bigg[\frac{1}{n}\text{tr}\Wb,\frac{2}{n}\text{tr}\Wb\bigg],\quad 1\le i\le r.
\]
We remark that $\hat\zeta_i$ is dependent on the random weights $w_j$'s and claim the next lemma. 
\begin{lemma}\label{lem: hat theta}
	Under Assumption \ref{c1},  for $1\le i\le r$, the solution $\hat\zeta_i$ exists with probability tending to 1 as $\min\{n,p\}\rightarrow \infty$. Moreover, $	\theta_i/\lambda_i=1+O_p((n\lambda_i)^{-1}\text{tr}\bLambda_2)$, and 
	\[
	\begin{split}
		\hat\zeta_i-\frac{\theta_i}{\lambda_i}=&\frac{1}{n}\sum_{j=1}^n(w_j-1)+\bigg(\frac{\text{tr}\bLambda_2}{n\lambda_i}\bigg)^2\times\mathbb{E}[w_1^2(w_1-1)]+o_p(1)\times\bigg(\frac{\text{tr}\bLambda_2}{n\lambda_i}\bigg)^2+o_p(\frac{1}{\sqrt{n}}).\\
	\end{split}
	\]
\end{lemma}

Before moving forward, we need the next assumption.

\begin{assumption}\label{c2}
	Assume that for any $1\le i\le r$, there exists constant $c>0$ such that
	\[
	\xi_i:=\frac{1}{n}\sum_{j=1}^n\bigg\{\sum_{k=1}^{r+p}\gamma_{ik}^4[\nu_{jk}-3(\mathbb{E(}z_{jk}^2))^2]+3[\sum_{k=1}^{r+p}\gamma_{ik}^2\mathbb{E}(z_{jk}^2)]^2\bigg\}-1\ge c,
	\]
	where $(\gamma_{i1},\ldots,\gamma_{i,r+p})^\top=\bgamma_i$ is the eigenvector defined in (\ref{A top A}) and $\nu_{jk}=\mathbb{E}(z_{jk}^4)$. 
\end{assumption}
{Assumption \ref{c2} is a technical condition to ensure that the limiting distributions  of $\hat \lambda_i$ are not degenerate. When $\Cb=\Ib$ and $\bx_j$'s are i.i.d., $\xi_i$ reduces to $\sum_{k=1}^{r+p}\gamma_{ik}^4(\nu_{1k}-3)+2$, which is also in Assumption 4 of  \cite{cai2020limiting}.} It's notable that $\nu_{1k}\ge 1$ always holds while  $\sum_{k=1}^{r+p}\gamma_{ik}^4\le (\sum_{k=1}^{r+p}\gamma_{ik}^2)^2\le 1$. Therefore, under such cases, Assumption \ref{c2} holds as long as $\nu_{1k}\ge 1+c$ or $\max_k|\gamma_{ik}|\le1-c$ for some $c>0$. We have the next theorem.

\begin{theorem}[Limiting representation]\label{thm: representation}
	Under Assumptions \ref{c1} and \ref{c2}, as $\min\{n,p\}\rightarrow \infty$, no matter under the standard or  multiplier bootstrap, we always have for any $1\le i\le r$,
	\begin{equation}\label{representation}
		\frac{\hat\lambda_i}{\theta_i}-1=\frac{1}{n}\bgamma_i^\top\Zb\Wb\Zb^\top\bgamma_i-\frac{1}{n}\text{tr}\Wb-\frac{\theta_i}{\lambda_i}+\hat\zeta_i+O_p\bigg(\frac{1}{\sqrt{n}}\frac{(n\vee p)\log n}{n\lambda_i}+\frac{1}{n}\bigg)+o_p\bigg(\frac{p^2}{(n\lambda_i)^2}\bigg).
	\end{equation}
\end{theorem}
Based on Lemma \ref{lem: hat theta} and Theorem \ref{thm: representation}, one can verify that $\hat\lambda_i/\theta_i-1=O_p\{[(n\lambda_i)^{-1}(n\vee p)]^2+n^{-1/2}\}$ for any $1\le i\le r$, which is a faster rate compared with that in Lemma \ref{lem: preliminary}. Moreover, the asymptotic distribution of $\hat\lambda_1$ is mainly determined by the right hand side (RHS) of (\ref{representation}), which depends on both the sample matrix $\Xb$ and the random weights $w_j$'s. However, the calculations of $\theta_i$ and $\hat\zeta_i$ rely on the  population eigenvalues $\{\lambda_k\}_{k=1}^{r+p}$, which are unknown. In real applications, it will be more preferred to study the limiting distribution of $\hat\lambda_i$ conditional on the sample matrix $\Xb$.

\subsection{Conditional on samples}
Since $\theta_i$ is unknown but close to $\lambda_i$, a natural idea is to replace it with the $i$th largest eigenvalue of sample covariance matrix before bootstrap, i.e., $n^{-1}\Xb\Xb^\top$. 
Using the decomposition (\ref{A top A}), the non-zero eigenvalues of $n^{-1}\Xb\Xb^\top$ are the same as those of $
\tilde\Sbb=n^{-1}\bGamma\bLambda^{1/2}\bGamma^\top\Zb\Zb^\top\bGamma\bLambda^{1/2}\bGamma^\top$, 
or its companion matrix $\tilde{\mathcal{S}}=n^{-1}\Xb^\top\Xb$. 
In the current paper, quantities marked by ``hat'' always stand for ``after bootstrap'', while those marked by ``tilde'' stand for ``before bootstrap''. We denote the eigenvalues and eigenvectors of $\tilde{\mathcal{S}}$ as $\tilde\lambda_i$ (descending)  and $\tilde\bu_i$ respectively, while the eigenvectors of $\tilde\Sbb$ are  $\tilde\bgamma_i$. The first step is to investigate the limiting properties of the quantities $\tilde\lambda_i$, $\tilde\bgamma_i$, and $\tilde\bu_i$ for $1\le i\le r$. 

\begin{lemma}[Without bootstrap]\label{lem: without bootstrap}
	Under Assumptions \ref{c1} and \ref{c2}, as $\min\{n,p\}\rightarrow \infty$, for the eigenvalues and eigenvectors of $\tilde\Sbb$ and $\tilde{\mathcal{S}}$, we have:
	
	(a).  $\tilde\lambda_i/\lambda_i=1+o_p(1)$ for $1\le i\le r$ while  $\tilde\lambda_{r+1}=O_p(1)$. Moreover, 
	\[
	\sqrt{n}\bigg(\frac{\tilde\lambda_i}{\theta_i}-1\bigg)=\frac{1}{\sqrt{n}}\bigg(\bgamma_i^\top\Zb\Zb^\top\bgamma_i-n\bigg)+O_p\bigg(\frac{(n\vee p)}{n\lambda_i}+\frac{1}{\sqrt{n}}\bigg),\quad 1\le i\le r.
	\]
	
	(b).  For any $1\le i,j\le r$ and $i\ne j$, we have
	\[
	\begin{split}
		\bgamma_i^\top\tilde\bgamma_j\le& O_p\bigg(\frac{p}{n\sqrt{\lambda_i\lambda_j}}\frac{\lambda_i}{\max\{\lambda_i,\lambda_j\}}+\frac{\sqrt{\lambda_i\lambda_j}}{\sqrt{n}\times \max\{\lambda_i,\lambda_j\}}\bigg),\quad
		(\tilde\bgamma_j^\top\bgamma_j)^2=1+O_p\bigg(\frac{p}{n\lambda_j}+\frac{1}{n}\bigg).
	\end{split}
	\]
	
	(c). Write $\tilde\bu_i^\top=(\tilde u_{i1},\ldots,\tilde u_{in})$ and $\tilde\sigma_i^2=\sum_{j=1}^n\tilde u_{ij}^4$.  Then, 
	$
	\tilde\sigma_i^2=n^{-1}(\xi_i+1)+o_p(n^{-1})
	$
	for any $1\le i\le r$, where $\xi_i$ is defined in Assumption \ref{c2}.
\end{lemma}

Lemma \ref{lem: without bootstrap} provides comprehensive results on the asymptotic behavior of the spiked eigenvalues $\tilde\lambda_i$ for $1\le i\le r$ and the corresponding eigenvectors $\tilde\bgamma_i$, $\tilde\bu_i$. {{The limiting representation of  $\tilde\lambda_i$ and the convergence of $\tilde\bgamma_i$ are also shown in \cite{cai2020limiting}, but they haven't provided the convergence rates in (b).} 
	Denote $\mathbb{P}^*$ as the probability measure conditional on the sample $\Xb$.  Then, we can define  $\overset{d^*}{\rightarrow}$, $\overset{p^*}{\rightarrow}$, $o_{p^*}(1)$ and $O_{p^*}(1)$ accordingly under $\mathbb{P}^*$. Now we present the limiting distribution of $\hat\lambda_i/\tilde\lambda_i$ for $i\le r$ conditional on $\Xb$.
	
	\begin{theorem}[Conditional on sample]\label{thm: conditional on sample}
		Suppose that Assumptions \ref{c1} and \ref{c2} hold as $\min\{n,p\}\rightarrow \infty$. Conditional on $\Xb$, if $(n\lambda_i)^{-1}p=o(n^{-1/4})$ for some $1\le i\le r$, under the multiplier bootstrap,  	with probability tending to one we have
		\begin{equation}\label{conditional distribution}
			\tilde\sigma_i^{-1}(\hat\lambda_i/\tilde\lambda_i-1)\overset{d^*}{\longrightarrow}\mathcal{N}(0,1).
		\end{equation}
		On the other hand, if $n^{-1/4}=o[(n\lambda_i)^{-1}p]$, with probability tending to 1 we have 
		\begin{equation}\label{conditional power}
			\mathbb{P}^*(|\tilde\sigma_i^{-1}(\hat\lambda_i/\tilde\lambda_i-1)|\le s)\rightarrow 0,
		\end{equation}
		for any constant $s\in\mathbb{R}$. 
		Similar results to (\ref{conditional distribution}) and (\ref{conditional power}) hold under the standard bootstrap by replacing $\tilde\sigma_i$ with $\sqrt{\tilde\sigma_i^2-n^{-1}}$.
	\end{theorem} 
	
	\begin{remark}\label{remark:probability}
		By ``probability tending to one'' hereafter, we mean that there exist a series of events $\{\Xi_n\}$ holding with probability tending to one under the measure $\{\Xb_n\}$, such that (\ref{conditional distribution}) and (\ref{conditional power}) hold conditional on these events.
	\end{remark}
	
	The condition $(n\lambda_i)^{-1}p=o(n^{-1/4})$ in Theorem \ref{thm: conditional on sample} is satisfied  when $\lambda_i$ is sufficiently large. Under such cases, the spiked sample eigenvalues after bootstrap always converge weakly to Gaussian limits after proper  scaling and centralization, no matter under the multiplier ot standard bootstrap. The only difference between the two bootstrap schemes is on the asymptotic variance.  This is because the resampling weights $w_j$'s are weakly dependent under the standard bootstrap. It's worth noting that the scaling and centralization  parameters in (\ref{conditional distribution}) totally depend on the sample matrix $\Xb$, which is observable. Moreover, the condition $(n\lambda_i)^{-1}p=o(n^{-1/4})$ is almost sharp according to (\ref{conditional power}). 
	
	\subsection{Bias of bootstrap} 
	As a byproduct, Theorem \ref{thm: conditional on sample} also helps in understanding why the bootstrap technique may fail to approximate the distribution of sample eigenvalues, shown in  \cite{karoui2016bootstrap} and \cite{yao2021rates}. Following \cite{yao2021rates}, we compare the limiting distributions of $\lambda_i^{-1}(\hat\lambda_i-\tilde\lambda_i)$ and $\lambda_i^{-1}(\tilde\lambda_i-\lambda_i)$. We remark that in their settings, the leading $r$ spiked eigenvalues $\lambda_i$ are  of constant order while the remaining ones are asymptotically vanishing. It's parallel to a spiked covariance model by rescaling the eigenvalues. This is the reason  why we add the scaling coefficient $\lambda_i^{-1}$. The following corollary quantifies  the difference  between the two limiting distributions. 
	
	\begin{cor}[Bias of bootstrap]\label{cor: bootstrap bias}
		Suppose that Assumptions \ref{c1} and \ref{c2} hold as $\min\{n,p\}\rightarrow\infty$.  For any constant $s\in\mathbb{R}$ and $1\le i\le r$, under the standard bootstrap, we have
		\begin{equation}\label{bias hat}
			\begin{split}
				&\mathbb{P}^*\bigg(\sqrt{n}\times\lambda_i^{-1}(\hat\lambda_i-\tilde\lambda_i)\le s\bigg)
				=F_G\bigg(\frac{s}{\sqrt{\xi_i}}-\sqrt{\frac{n}{\xi_i}}\times \frac{\text{tr}^2\bLambda_2}{(n\lambda_i)^2}\times \mathbb{E}[w_1^2(w_1-1)]\bigg)+o_p(1),
			\end{split}
		\end{equation}
		where $F_G(\cdot)$ is the cumulative distribution function (CDF) of standard Gaussian variable. (\ref{bias hat}) also holds under the multiplier bootstrap by replacing $\xi_i$ with $\xi_i+1$. On the other hand, without bootstrap, we have
		\begin{equation}\label{bias tilde}
			\begin{split}
				&\mathbb{P}\bigg(\sqrt{n}\times\lambda_i^{-1}(\tilde\lambda_i-\lambda_i)\le s\bigg)
				=F_G\bigg(\frac{s}{\sqrt{\xi_i}}-\sqrt{\frac{n}{\xi_i}}\frac{\text{tr}\bLambda_2}{n\lambda_i}\bigg)+o(1).
			\end{split}
		\end{equation}
		
	\end{cor}
	
	By Corollary \ref{cor: bootstrap bias},  the standard bootstrap procedure is asymptotically consistent as long as $(n\lambda_i)^{-1}\text{tr}\bLambda_2=o(n^{-1/2})$, because the two tail probabilities are asymptotically equal. This condition is slightly sharper than that in \cite{yao2021rates}, where they require $(n\lambda_1)^{-1}\text{tr}(\bSigma)=o(n^{-1/2})$. More importantly, Corollary \ref{cor: bootstrap bias}  provides asymptotic bias for the bootstrap  procedures when the factors are weak. For instance, when $n^{-3/4}p\ll\lambda_i\ll n^{-1/2}p$, 
	\begin{equation}\label{difference tail}
		\begin{split}
			&\mathbb{P}^*(\sqrt{n}\times\lambda_i^{-1}(\hat\lambda_i-\tilde\lambda_i)\le s)-\mathbb{P}(\sqrt{n}\times\lambda_i^{-1}(\tilde\lambda_i-\lambda_i)\le s)=F_G(s\xi_i^{-1/2})+o_p(1).
		\end{split}
	\end{equation}
	If the common factors are extremely weak, i.e., $\lambda_i\ll n^{-3/4}p$, Corollary \ref{cor: bootstrap bias} indicates that the two tail probabilities will both converge to 0 for any $s\in\mathbb{R}$. However, the coherence of the two tail probabilities under such cases doesn't mean that the bootstrap can accurately approximate the limiting distribution of $\tilde\lambda_i$. 
	Another interesting finding  is that the multiplier bootstrap is always biased, mainly because the asymptotic variances of the two limiting distributions don't match.

	\subsection{Testing procedure}\label{sec:testing 2}
	We now provide the testing procedure to determine the number of common factors. By Theorem \ref{thm: conditional on sample}, the asymptotic distribution in (\ref{conditional distribution})  holds when the factors are strong and $r\ge i$.  Therefore, we consider the null hypothesis and the alternative one as
	\begin{equation}\label{H01}
		H_{0i}: r\ge i, \quad v.s.\quad H_{1i}: r<i, \quad \text{for some } i\ge 1.
	\end{equation}
	We reject the null hypothesis $H_{0i}$ as long as  $\tilde\sigma_i^{-1}|\hat\lambda_i/\tilde\lambda_i-1|\ge F_G^{-1}(1-\alpha/2)$ under a predetermined significance level $\alpha$, where $F_G^{-1}(\cdot)$ is the quantile function of the standard normal distribution.  By letting $i=1$, we can test the existence of common factors.
	
	Furthermore, in order to determine the number of factors  $r$, we implement the testing procedure sequentially as in \cite{onatski2009testing}. Specifically, for $i=1,\ldots, r_{\max}$ where $r_{\max}$ is a predetermined upper bound (fixed), we sequentially  calculate $\tilde\sigma_i^{-1}|\hat\lambda_i/\tilde\lambda_i-1|$  until  $H_{0i}$ is rejected at some $i=k$. Then,  $\hat r=k-1$ is the estimated number of common factors.  The significance level $\alpha$ is usually small in order to control the type one error. {Our simulation studies show that the results are not sensitive to the value of $\alpha$ for $\alpha\in[0.01,0.1]$. However, in finite samples, the above procedure tends to overestimate the number of factors when $p/n$ is small. To improve the performance, we propose to slightly modify the criterion. Specifically, if  $p/n<0.5$, we reject $H_{0i}$ when $\tilde\sigma_i^{-1}|(\hat\lambda_i+c_n)/\tilde\lambda_i-1|\ge F_G^{-1}(1-\alpha/2)$ for some $c_n=O_p(n^{-1/2})$, to enhance the power of the tests in (\ref{H01}) for $i>r$. In this paper, $c_n=2\tilde\sigma_0^2(1+\sqrt{p/n})^2n^{-1/2}$, where $\tilde\sigma_0^2=(np)^{-1}\|\Xb\|_F^2$.}

	\section{Testing with non-spiked eigenvalues}\label{sec:non-spiked}
	\subsection{Limiting behavior of largest non-spiked eigenvalue}
	{Testing with non-spiked eigenvalues is another important direction for determining the number of common factors; see \cite{onatski2009testing}, \cite{cai2020limiting} and \cite{ke2021estimation}. Without bootstrap,  the largest non-spiked eigenvalues of the sample covariance matrix have been shown to follow the Tracy-Widom law. }However, after bootstrap, the limiting behavior of $\hat\lambda_{i}$ for $i>r$ remains an open problem. In this section, we  fill this gap under the  multiplier bootstrap. For the standard bootstrap, it's more challenging because $w_j$'s are dependent and the marginal distribution is discrete. We leave it as future work.
	
	Theoretical analysis of the non-spiked sample eigenvalues is much more challenging because there is no clear gap between $\lambda_{r+1},\ldots,\lambda_{r+p}$. Similarly to \cite{cai2020limiting}, we need more assumptions. Let $\bSigma_2=\bGamma_2\bLambda_2\bGamma_2^\top$ and  $m_{\bSigma_2}(z)$ be the unique solution in $\mathbb{C}^+$ to 
	\[
	m_{\bSigma_2}(z)=-\frac{1}{z-n^{-1}\text{tr}[\{\Ib+m_{\bSigma_2}(z)\bSigma_2\}^{-1}\bSigma_2]},\quad z\in\mathbb{C}^+.
	\]
	{Then, $m_{\bSigma_2}(z)$ is the limit of the Stieltjes transform associated with $n^{-1}\Zb^\top\bSigma_2\Zb$, and corresponds to a probability function $F_{m}(\cdot)$. Let $
		\lambda_+=\inf\{x\in\mathbb{R}:F_{m}(x)=1\}$ and $d_+=-\lim_{z\in\mathbb{C}^+\rightarrow \lambda_+}m_{\bSigma_2}(z)$. See Assumption 8 in \cite{cai2020limiting} for more details on $F_m(\cdot)$, $\lambda_+$ and $d_+$, which motivates us to propose the next assumption.}
	\begin{assumption}\label{c3}
		Further assume that:
		
		(a) The empirical spectral distribution of $\bSigma_2$ converges to some probability distribution $F_{\bSigma_2}$ not degenerate at $0$. 
		
		(b) There are at most finite number of eigenvalues $\lambda_i$ satisfying $\lim\sup_{n\rightarrow\infty}\lambda_i d_+\ge1$.
		
		(c) $p/n=\phi_n\rightarrow \phi\in(0,\infty)$ as {$\min\{n,p\}\rightarrow \infty$} for some constant $\phi$ while the moments $\sup_{i,j}\mathbb{E}|z_{ij}|^q<\infty$ for any integer $q>0$.
	\end{assumption}
	
	Assumption \ref{c3}(a) ensures the existence of non-degenerate $m_{\bSigma_2}(z)$ and $F_{\bSigma_2}$. Assumption \ref{c3}(b)   is actually more general than Assumption 8 in \cite{cai2020limiting}, where $\lim\sup_{n\rightarrow \infty}\lambda_{r+1}d_+<1$. Assumption \ref{c3}(b) is equivalent to allowing a finite number of eigenvalues of $\bSigma_2$ to be separated from the support of $F_{\bSigma_2}$, as long as they are still bounded { as required in Assumption \ref{c1}(c).}  In other words, we allow the existence of outliers in the spectrum of idiosyncratic error covariance matrix. 
	In econometrics, these outliers may exist due to some large marginal variances or the cross-sectional correlations of the idiosyncratic errors. As claimed in the introduction, this is also a major difference between our bootstrapped method and those based on traditional RMT,  such as  \cite{onatski2009testing}, \cite{cai2020limiting} and \cite{ke2021estimation}.  Such a refinement mainly benefits from the largest resampling weight, which is of order $\log n$ thus reducing the effects of bounded outliers in $\bSigma_2$. Assumption \ref{c3}(c) requires that $p$ and $n$ are of the same order, which is a common assumption in the RMT literature. The moment condition can be potentially relaxed, which is not the major concern of the current paper.
	
	Like in Section \ref{sec: spiked}, we need to find a proper approximation to  $\hat\lambda_{r+1}$.
	Let $\mathcal{T}_n$ be the set of all permutations of $\{1,\ldots,n\}$. Then, the orders of $\{w_1,\ldots,w_n\}$ follow uniform distribution on $\mathcal{T}_n$. We use $\{t_1,\ldots,t_n\}\in\mathcal{T}_n$ to denote the orders such that  $w_{t_1}\ge \cdots\ge w_{t_n}$. We  define $\lambda_0$ as the unique solution to  the following equation:
	\begin{equation}\label{lambda_0}
		\frac{1}{w_{t_1}}=\frac{1}{n}\sum_{i=1}^p\lambda_{r+i}\bigg[\lambda_0-\frac{\lambda_{r+i}}{n}\sum_{j=2}^n\frac{w_{t_j}}{1-w_{t_j}/w_{t_1}}\bigg]^{-1},\quad \lambda_0\in\bigg[\frac{\lambda_{r+1}}{n}\sum_{j=2}^n\frac{w_{t_j}}{1-w_{t_j}/w_{t_1}},\infty\bigg).
	\end{equation}
	The definition of $\lambda_0$ is motivated by Theorem 1.1 in \cite{couillet2014analysis} and (2.10) in \cite{yang2019edge}. 
	We start with the simple case where $r=0$. That is, $\Xb=\bPsi\Eb$ and $\hat\lambda_{r+1}$ is equal to the largest eigenvalue of 
	$
	n^{-1}\bSigma_2^{1/2}\Eb\Wb\Eb^\top\bSigma_2^{1/2}$, where $ \bSigma_2=\bGamma_2\bLambda_2\bGamma_2^\top=\bPsi^\top\bPsi$.
	Lemma \ref{lemma: ratio} below will indicate that $\lambda_0$ is a good approximation to $\hat\lambda_{r+1}$.
	\begin{lemma}\label{lemma: ratio}
		Under Assumptions \ref{c1} and \ref{c3}, if $r=0$, we have
		\[
		\sqrt{n}(\hat\lambda_{r+1}/\lambda_0-1)=n^{-1/2}(\phi_n\bar\lambda)^{-1}(\bx_{t_1}^\top\bx_{t_1}-p\bar\lambda)+o_p(1),
		\]
		as $n\rightarrow \infty$ under the multiplier bootstrap, 
		where $\phi_n=p/n$ and  $\bar\lambda=p^{-1}\sum_{i=1}^p\lambda_{r+i}$.
	\end{lemma}
	Lemma \ref{lemma: ratio} indicates that the ratio $\hat\lambda_{r+1}/\lambda_0$ converges to 1 with rate $n^{-1/2}$, while $\lambda_0$ is random and  dependent on the resampling weights $w_j$'s. To conclude the asymptotic distribution of $\hat\lambda_{r+1}$, it suffices to discuss the fluctuation of $\lambda_0$. See the next lemma.
	
	\begin{lemma}\label{lemma: lambda_0}
		Under the same assumptions as in Lemma \ref{lemma: ratio}, as $n\rightarrow \infty$ we have
		\[
		\begin{split}
			\lambda_0=&\phi_n\bar\lambda w_{t_1}+(n\phi_n\bar\lambda)^{-1}\sum_{i=1}^p\lambda_{r+i}^2+O_p(\frac{1}{\log n}),\\
			\Longrightarrow \hat\lambda_{r+1}=&\phi_n\bar\lambda w_{t_1}+(n\phi_n\bar\lambda)^{-1}\sum_{i=1}^p\lambda_{r+i}^2+O_p\bigg(\frac{1}{\log n}+\frac{\log n}{\sqrt{n}}\bigg).
		\end{split}
		\]
		Therefore, as $n\rightarrow\infty$,
		\begin{equation}\label{probability}
			\mathbb{P}\bigg(\frac{\hat\lambda_{r+1}-(n\phi_n\bar\lambda)^{-1}\sum_{i=1}^p\lambda_{r+i}^2}{\phi_n\bar\lambda}-\log n>x\bigg)\rightarrow \exp(-\exp(-x)),\quad x\in\mathbb{R}. 
		\end{equation}
		(\ref{probability}) also holds with probability tending to one if replacing $\mathbb{P}$ with $\mathbb{P}^*$.
	\end{lemma}
	
	By Lemmas \ref{lemma: ratio} and \ref{lemma: lambda_0},  the limiting distribution of the largest non-spiked eigenvalue after bootstrap is determined by the order statistics of resampling weights. More precisely, it  depends on the largest weight $w_{t_1}$, which converges weakly to the Gumbel distribution after centralization and scaling. This is consistent with the conclusion in \cite{kwak2021extremal}, although the limiting distributions are not the same. The Tracy-Widom law in traditional RMT doesn't hold anymore, because the exponential distribution is convex at the edge with unbounded support.  However, Lemmas \ref{lemma: ratio} and \ref{lemma: lambda_0} are only for the special case $r=0$.
	Theorem \ref{thm:non-spike} below provides the results for general cases where $r\ge 0$.
	
	\begin{theorem}\label{thm:non-spike}
		Under Assumptions \ref{c1} and \ref{c3} and the multiplier bootstrap, as $n\rightarrow \infty$, for any fixed $r\ge 0$ and small constant $c>0$,   it holds that $
		\hat\lambda_{r+1}-\hat\varphi_1= O_p(n^{-2/3+c})$,
		where $\hat\varphi_i$ is the $i$th largest  eigenvalue of $n^{-1}\bSigma_2^{1/2}\Zb\Wb\Zb^\top\bSigma_2^{1/2}$. 
	\end{theorem}
	$\hat\varphi_1$ is from the model without common factors. {When the entries of $\Fb$ are serially independent, the results in Lemmas \ref{lemma: ratio} and \ref{lemma: lambda_0} hold directly for $\hat\varphi_1$. If $\Fb=\Cb\Fb^0$ and Assumption \ref{c1} holds, this will only generate an error of rate $O_p(n^{-1/2+c})$ to $\hat\varphi_1$ for arbitrary small $c>0$, which has negligible effect on its asymptotic distribution. }Then, we conclude that Lemmas \ref{lemma: ratio} and \ref{lemma: lambda_0} also hold for the general factor models with $r\ge 0$.

	\subsection{Testing procedure}
	One can not test the number of common factors directly based on Lemma \ref{lemma: ratio}, Lemma \ref{lemma: lambda_0} or Theorem \ref{thm:non-spike}. Firstly, the centralization and scaling parameters in (\ref{probability}) rely on the unknown population eigenvalues $\lambda_i$'s. Secondly, the convergence rates in Lemmas \ref{lemma: ratio} and \ref{lemma: lambda_0} are actually very slow.  In finite samples, the error is non-negligible and the theoretical critical values are not reliable. In the following, we fix the two problems by repeating the bootstrap procedure and using approximated critical values. 
	
	When testing with non-spiked eigenvalues, the null hypothesis will be different. Motivated by \cite{onatski2009testing}, consider the null hypothesis and the alternative 
	one as 
	\begin{equation}\label{H tilde}
		H_{0i}^*:r= i-1,\quad H^*_{1i}:r\ge i, \text{ for some }i\ge 1.
	\end{equation}
	Then, under $H^*_{0i}$, Lemma \ref{lemma: lambda_0} will hold for $\hat\lambda_i$, which is at most of order $\log n$. On the contrary, under $H^*_{1i}$, Lemma \ref{lem: preliminary} shows that $\hat\lambda_i\asymp \lambda_i\gg \log n$. Therefore, one may reject the null hypothesis $H^*_{0i}$ as long as $\hat\lambda_i>c_{\alpha}$ for some critical value $c_{\alpha}$. Our target is to provide a reasonable approximation to $c_{\alpha}$ given significance level $\alpha\in (0,1)$.

	By Theorem \ref{thm:non-spike}, under $H^*_{0i}$, the asymptotic distribution of $\hat\lambda_i$ will be exactly the same as that of $\hat\varphi_1$. If $\bSigma_2$ is given, the limiting distribution of $\hat\varphi_1$ can be approximated by a standard Monte Carlo method. In fact, we don't need to know  $\bSigma_2$ exactly. Given any fixed integer $k\ge 0$, define $\lambda_0^{(k)}$ as the solution to 
	\[
	\frac{1}{w_{t_1}}=\frac{1}{n}\sum_{i=k+1}^p\lambda_{r+i}\bigg[\lambda_0^{(k)}-\frac{\lambda_{r+i}}{n}\sum_{j=2}^n\frac{w_{t_j}}{1-w_{t_j}/w_{t_1}}\bigg]^{-1},\quad \lambda_0^{(k)}\in\bigg[\frac{\lambda_{r+k+1}}{n}\sum_{j=2}^n\frac{w_{t_j}}{1-w_{t_j}/w_{t_1}},\infty\bigg).
	\]
	That is, we remove the leading $k$ eigenvalues $\lambda_{r+1},\ldots,\lambda_{r+k}$ from $\bSigma_2$. It's not hard to verify $\lambda_0^{(k)}-\lambda_0=O_p(n^{-1+c})$ for any constant $c>0$. Therefore, if we define
	$
	\bSigma_2^{(k)}=\sum_{i=k+1}^p\lambda_{r+i}\bgamma_{r+i}\bgamma_{r+i}^\top$,
	the largest eigenvalue of $n^{-1}(\bSigma_2^{(k)})^{1/2}\Zb\Wb\Zb^\top(\bSigma_2^{(k)})^{1/2}$ will have exactly the same limiting distribution as that of $\hat\varphi_1$. In other words,  if we remove the leading $k$ eigenvalues from the population covariance matrix $\bSigma$ for some constant $k\ge r$,  the asymptotic distribution of the largest non-spiked eigenvalue will not change. Empirically, the population covariance matrix is unknown, so we implement this step on the sample covariance matrix $\tilde\Sbb$ or directly on $\Xb$,  illustrated in  Algorithm \ref{alg1}. 
	
	\begin{algorithm}[H]
		\caption{Find limiting distribution for $\hat\varphi_1$}\label{alg1}
		\begin{algorithmic}[1]
			\Require data matrix $\Xb_{p\times n}$, a fixed  $r_{\max}$, number of Monte Carlo experiments $R$.
			\Ensure  approximation to the limiting distribution of $\hat\varphi_1$.
			
			\State Do Singular-Value-Decomposition for $\Xb$ as $\Xb=\Ub_x\Db_x\Vb_x=\sum_{i=1}^{\min\{n,p\}}d_{x,i}\bu_{x,i}\bv_{x,i}^\top$. Remove the leading $r_{\max}$ singular values and define $\tilde\Xb=\sum_{i=r_{\max}+1}^{\min\{n,p\}}d_{x,i}\bu_{x,i}\bv_{x,i}^\top$.
			
			\State For $j=1,\ldots,R$, generate $w_1,\ldots,w_n$ i.i.d. from exponential distribution $Exp(1)$. Calculate the largest eigenvalue of $n^{-1}\tilde\Xb\Wb\tilde\Xb^\top$, denoted as $\hat\varphi_1^{j}$.
			
			\State Output the empirical distribution of $\{\hat\varphi_1^j\}_{1\le j\le R}$.
		\end{algorithmic}
	\end{algorithm}
	
	With Algorithm \ref{alg1}, we reject the null hypothesis $H^*_{0i}$ as long as $\hat\lambda_i>\hat c_{1-\alpha}$, where $\hat c_{1-\alpha}$ is the $(1-\alpha)$th sample quantile of $\{\hat\varphi_1^j\}_{1\le j\le R}$. By letting $i=1$, we can test the existence of  factors.
	On the other hand,
	in order to determine the number of common factors, it's not necessary to use the sequential tests. Since the eigenvalues are in descending order, we can use the same thresholding approach as in \cite{onatski2010determining}, \cite{cai2020limiting} and \cite{ke2021estimation} by simply defining
	$
	\hat r=\sum_{i=1}^{r_{\max}}I(\hat\lambda_i>\hat c_{1-\alpha})$, where $I(\cdot)$ is the indicator function. In finite samples,  the outputted $\hat\varphi_1^j$'s from Algorithm \ref{alg1} tend to be smaller than the real $\hat\varphi_1$, because we have deleted more singular values than needed. { In real applications, the procedure can be  implemented recursively by updating $r_{\max}$ with $\hat r$ obtained from the last step until convergence. The reported results in the numerical and empirical studies are from the recursive procedure.}

	\section{Decision rule}\label{sec: decision}
	The proposed approaches to determining $r$ are based on the bootstrapped  sample eigenvalues $\hat\lambda_i$'s, whose limiting distributions are  mainly dependent on the randomness of $w_j$'s. To obtain $\hat\lambda_i$, so far we have only considered implementing the bootstrap procedure once. As a randomized approach, the result can be unstable especially when the sample size $n$ is small. Due to the randomness of $w_j$'s, people may report different conclusions even if they are implementing exactly the same procedure on the same data. In fact, if infinite users conduct the test on the same data, the reported p-values will be uniformly distributed on $[0,1]$. See also the criticism in \cite{geyer2005fuzzy} and \cite{he2021vector}. In statistics, a common method to  improve the stability is to conduct the bootstrap more times.

	We start with the sequential test procedure with spiked eigenvalues in Section \ref{sec: spiked}.  Recall the  hypotheses in (\ref{H01}). Given $i$, we independently repeat the bootstrap resampling procedure $B$ times, obtain $B$ sample eigenvalues $\{\hat\lambda_i^b\}_{b=1}^B$ and define
	\[
	D_i^s(\alpha,B)=B^{-1}\sum_{b=1}^B I\{|\tilde\sigma_i^{-1}(\hat\lambda_i^b/\tilde\lambda_i-1)|\le F_G^{-1}(1-\alpha/2)\}.
	\]
	Then, conditional on samples, $D_i^s(\alpha,B)$ can be regarded as the average of some i.i.d. Bernoulli random variables. Obviously, as $B\rightarrow \infty$,
	\begin{equation}\label{Dis}
		D_i^s(\alpha,B)\rightarrow \mathbb{P}^*\{|\tilde\sigma_i^{-1}(\hat\lambda_i^1/\tilde\lambda_i-1)|\le F_G^{-1}(1-\alpha/2)\}.
	\end{equation}
	Based on Theorems \ref{thm: conditional on sample}, under $H_{0i}$ we have $D_i^s(\alpha,B)\rightarrow 1-\alpha$ {if $(n\lambda_i)^{-1}p=o(n^{-1/4})$} while $D_i^s(\alpha,B)\rightarrow 0$ under $H_{1i}$  for any $\alpha\in(0,1)$ as $(n,B)\rightarrow \infty$ with a fast rate of $\sqrt{B}$. Therefore, any constant in $(0,1-\alpha)$ can be a good threshold to distinguish $H_{0i}$ from $H_{1i}$. We reject $H_{0i}$ as long as $D_i^s(\alpha,B)\le C_{th}$ for some predetermined $C_{th}\in(0,1-\alpha)$. This is no longer a regular test. Instead, we call it a decision rule for the number of common factors. As stated by \cite{geyer2005fuzzy}, it's sufficient to report $D_i^s(\alpha,B)$ in real applications. The selection of $C_{th}$ depends on the users' tolerance of under/over estimation errors, quite similar to the determination of significance level in a standard test. Large $C_{th}$ is in favor of $H_{1i}$ but with higher risk of underestimating the number of common factors. Small $C_{th}$ will lead to the opposite result. In this paper, we always use the middle value  $C_{th}=(1-\alpha)/2$.
	The decision rule is also applicable to the thresholding method  in Section \ref{sec:non-spiked} by defining
	\[
	D_i^{ns}(\alpha,B)=B^{-1}\sum_{b=1}^BI(\hat\lambda_i^b<\hat c_{1-\alpha}),\quad \hat r=\sum_{i=1}^{r_{\max}}I\{D_i^{ns}(\alpha,B)<C_{th}\},
	\]
	where $\hat c_{1-\alpha}$ is given by Algorithm \ref{alg1}.  
	
	In addition to de-randomizing the results, the decision rule is also very important in improving the accuracy of $\hat r$. {For the tests with spiked eigenvalues,  the condition $(n\lambda_i)^{-1}p=o(n^{-1/4})$ in Theorem \ref{thm: conditional on sample} indicates that the size and power may be unsatisfactory under weak factors or strong noises. }  Equivalently, $D_i^s(\alpha,B)$ may not be very close to $1-\alpha$ and 0 under the two hypotheses. However, the conclusion from the decision rule will not be affected as long as $D_i^s(\alpha,B)$ is still well separated by $C_{th}$. In other words, the decision rule allows more errors in the individual test. On the other hand, for the tests with non-spiked eigenvalues, it's similar. The approximated critical value $\hat c_{1-\alpha}$ may not be accurate in finite samples, especially when $\alpha$ is small or $r_{\max}$ is large. Then, $D_i^{ns}(\alpha,B)$ may not be very close to $1-\alpha$ or 0, but the determined factor number can be still very accurate. The empirical effects of the decision rule will be further studied by  simulating examples in the Supplement.

	Before ending this section, we discuss the computational complexity of the proposed methods.
	For the tests with spiked eigenvalues, computing  $\tilde\lambda_i$'s and $\tilde\sigma_i$'s typically costs $O(npr_{\max})$ operations. The calculation of $\hat\lambda_i$'s will cost additional $O(npr_{\max})$ operations if we ignore the generation of resampling weights $w_j$'s. Therefore, the total  computational cost is $O(npr_{\max})$ if the bootstrap is implemented only once, which is comparable to the methods in the literature. When the decision rule is applied, the computational cost grows to $O(Bnpr_{\max})$. For the approach based on non-spiked eigenvalues, the major computational cost is from Algorithm \ref{alg1}. Step 1 of Algorithm \ref{alg1} needs to specify the leading $r_{\max}$ singular values and vectors of $\Xb$, which requires $O(npr_{\max})$ operations. Step 2 will cost additional $O(npr_{\max})$ operations. Then, the total computational cost is $O(Rnpr_{\max})$. If the decision rule is applied, the total computational cost will be $O((R+B)npr_{\max})$.

	\section{Numerical studies}\label{sec:numerical}
	In the simulation, we generate data 
	according to 
	\begin{equation}\label{data generation}
		\Xb=\vartheta\Lb\bPhi\Fb^\top+\Eb.
	\end{equation}
	In (\ref{data generation}), $\vartheta\ge 0$ is a parameter controlling the strength of  common factors. When $\vartheta=0$, no common factors exist and $\Xb$ is from pure noise process. $\Lb$ is the $p\times r$ loading matrix with $r=3$, whose entries are from i.i.d. $\mathcal{N}(0,1)$.  $\bPhi=\text{diag}(1.5,1.2,p^{-a})$ is a diagonal matrix  to ensure the spiked eigenvalues are not identical so that Assumption \ref{c1}(c) holds, where $0\le a< 0.5$ is a parameter controlling the strength of the third common factor. By letting $a>0$, we allow the existence of weak factor. $\Fb^\top=(\bbf_1,\ldots,\bbf_n)$ is the factor score matrix from auto-regressive (AR) process, { i.e., 
		$
		\bbf_t=\beta_f\bbf_{t-1}+\bh_t$ with $\bh_t\overset{i.i.d.}{\sim}\mathcal{N}(0,\Ib_r)$, so that Assumption \ref{c1}(a) holds. 
		According to our assumptions, the idiosyncratic error matrix $\Eb=(\bepsilon_1,\ldots,\bepsilon_n)$ is generated by
		$
		\bepsilon_t\overset{i.i.d.}{\sim}\mathcal{N}(0,\bSigma_{\epsilon})$,
		where all the diagonal entries of $\bSigma_{\epsilon}$ are equal to 1 while the off-diagonal entries are equal to $\rho/p$ for some constant $\rho>0$.} The parameter $\rho$ controls the cross-sectional correlations of the idiosyncratic errors.  When $\rho$ is large, there will be an outlier in the eigenvalues of idiosyncratic error covariance matrix, { which is not  regarded as being driven by common factors although it may exceed the BBP phase transition boundary.}

	\subsection{Determining the number of common factors }\label{sec: number}
	In total, we have proposed three methods to determine the number of factors: test with Spiked eigenvalues, Multiplier bootstrap and Decision rule ($\hat r_{SMD}$), test with Spiked eigenvalues, Standard bootstrap and Decision rule ($\hat r_{SSD}$), and Eigenvalue Thresholding with Multiplier bootstrap and Decision rule ($\hat r_{ETMD}$). To implement the procedures, we let $r_{\max}=8$, $\alpha=0.05$, $B=200$ and $R=400$. {More numerical studies in the Supplement indicate that the proposed methods are not sensitive to the above tuning parameters.}
	We will compare our methods with some state-of-the-art approaches mentioned in the introduction: {the information criterion in \cite{bai2002determining} with $IC_{p2}$ rule ($\hat r_{IC}$) and its improved version by \cite{alessi2010improved} ($\hat r_{ABC}$), the eigenvalue ratio approach in \cite{ahn2013eigenvalue} ($\hat r_{ER}$), the sequential tests in \cite{trapani2018randomized} ($\hat r_{TRAP}$) and \cite{onatski2009testing} ($\hat r_{ON}$), the eigenvalue thresholding methods in \cite{onatski2010determining} ($\hat r_{ED}$), \cite{cai2020limiting} ($\hat r_{ETC}$), \cite{ke2021estimation} ($\hat r_{ETZ}$), and the deterministic parallel analysis by \cite{dobriban2019deterministic} ($\hat r_{DDPA_+}$).}  For data generating, we set $\beta_f=0.2$ and try different combinations of $(\vartheta,a,\rho)$ as $n=p$ grows. Table \ref{tab:number} reports the  averaged estimations of the factor number over 500 replications by different methods under diversified settings, while in the Supplement we report the corresponding proportions of under and over estimation.

	\begin{table}
		\addtolength{\tabcolsep}{0.8pt}
		\caption{The averaged estimations of factor number by different approaches over 500 replications. The true number of factors is $r=3$ when $\vartheta\ne 0$. \label{tab:number}}
		\renewcommand{\arraystretch}{0.6}
		\scalebox{0.75}{ 
			\begin{tabular}{lllllllllllllllll}
				\toprule[1.2pt]
				$\vartheta$&$\rho$&$a$&$n=p$&$\hat r_{SMD}$&$\hat r_{SSD}$&$\hat r_{ETMD}$&$\hat r_{IC}$&$\hat r_{ABC}$&$\hat r_{ER}$&$\hat r_{TRAP}$&$\hat r_{ON}$&$\hat r_{ED}$&$\hat r_{ETC}$&$\hat r_{ETZ}$&$\hat r_{DDPA_+}$\\\midrule[1.2pt]
				0&0&0&100&0&0&0&0&0&0&0&0.036&0.038&0.002&0.012&0\\
				0&0&0&200&0&0&0&0&0&0&0&0.054&0.05&0&0&0\\
				0&0&0&300&0&0&0&0&0&0&0&0.032&0.024&0&0.006&0\\\hline
				0&3&0&100&0.29&0.266&0&0&0.218&0&0&0.06&1.008&0.778&0.982&0.95\\
				0&3&0&200&0.002&0.002&0&0&0.15&0&0&0.074&1.028&0.894&1.004&0.988\\
				0&3&0&300&0&0&0&0&0.124&0&0&0.112&1.018&0.966&1&0.998\\\hline
				1&0&0&100&3&3&3&3&3&3&2.296&3.028&3.01&3&3&2.96\\
				1&0&0&200&3&3&3&3&3&3&3&3.026&3.014&3.002&3&2.976\\
				1&0&0&300&3&3&3&3&3&3&3&3.052&3.014&3&3&2.996\\\hline
				1&0&0.25&100&3&3&2.992&2.016&3&2&1.98&2.704&3.01&3&3&2.988\\
				1&0&0.25&200&3&3&3&2.024&3&2&2&3.026&3.03&3&3.002&2.988\\
				1&0&0.25&300&3&3&3&2.852&3&2&2&3.04&3.018&3&3&3\\\hline
				1&3&0&100&3.538&3.454&3&3&3.224&3&2.208&3.056&3.972&3.786&3.884&3.93\\
				1&3&0&200&3.002&3&3&3&3.124&3&3&3.092&4.01&3.914&3.996&3.964\\
				1&3&0&300&3&3&3&3&3.086&3&3&3.112&4.006&3.956&4&3.994\\\hline
				1&3&0.25&100&3.6&3.506&2.998&2.008&3.262&2&1.98&1.852&3.96&3.794&3.892&3.934\\
				1&3&0.25&200&3.006&3.002&3&2.024&3.14&2&2&2.32&4.01&3.88&3.998&3.972\\
				1&3&0.25&300&3&3&3&2.848&3.112&2&2&2.788&4.006&3.952&4&4
				\\\bottomrule[1.2pt]
		\end{tabular}}
	\end{table}

	Table \ref{tab:number}  shows clearly the advantage of the proposed approaches over the competitors. Overall speaking, the proposed approaches can accurately determine the number of common factors for a wide range of parameter settings, as long as $n$ is sufficiently large. {$\hat r_{ETMD}$ is slightly more reliable than $\hat r_{SMD}$ and $\hat r_{SSD}$, especially when $\rho$ or $a>0$ and $n,p$ are small. The competitors perform unsatisfactorily. Methods based on very large eigenvalue gap, such as $\hat r_{IC}$, $\hat r_{ER}$ and $\hat r_{TRAP}$,  may lose accuracy significantly when $a>0$, i.e., weak factor exists, leading to underestimation. It should be acknowledged that  $\hat r_{IC}$ and  $\hat r_{TRAP}$ will also work under weak factors with proper tuning parameters. However, it remains a challenge to  select parameters for them, especially in real applications. For fair comparison, we don't pay much attention to tuning.  The improved information criterion $\hat r_{ABC}$ performs  more reliably than $\hat r_{IC}$ when $a>0$, but slightly less accurately when $a=0$ due to the selection bias of an extra tuning parameter.  The performance of $\hat r_{ON}$ is not bad in most scenarios, but there is always a positive proportion of under or over estimation and it requires larger $n$ to converge. This is because the sequential tests in \cite{onatski2009testing} are implemented with a positive significance level  $\alpha=0.01$. Our methods will be less affected by $\alpha$ after applying the decision rule.  The eigenvalue thresholding methods, such as $\hat r_{ED}$, $\hat r_{ETC}$ and $\hat r_{ETZ}$, are more likely to overestimate $r$ when $\rho>0$, i.e., the idiosyncratic error covariance matrix contains an outlier eigenvalue. These methods will regard the outlier as a new factor. Such a property can be appealing in some applications, but in factor models the outlier may not add sufficient explanatory power besides resulting in more computational burden and potential errors. The deterministic parallel analysis $\hat r_{DDPA_+}$ performs similarly to the eigenvalue thresholding ones, mainly because it's also based on the RMT. In conclusion, the proposed methods, especially $\hat r_{ETMD}$, perform convincingly and reliably in very general scenarios.}

	\subsection{Sensitivity to data generating parameters}
	We further compare the performance of the above methods under more scenarios, by considering a wider range of parameter settings for data generating. This is also helpful in understanding when the proposed methods will fail to work. In Table \ref{tab:number}, the benchmark setting will be $n=200$, $p/n=\phi=1$, $\vartheta=1$
	, $\rho=0$ and $a=0$. We will change  one of $(\phi,\vartheta,\rho,a)$ but fixing  the others in  each experiment, to investigate whether the performance of the methods is sensitive to the data generating parameters. The proportions of exact estimation over 500 replications by different methods are plotted in Figure \ref{fig: sense1} under various parameter settings. {For better illustration, we only show the results of $\hat r_{SMD}$, $\hat r_{ETMD}$, $\hat r_{ABC}$, $\hat r_{ON}$ and $\hat r_{DDPA_+}$, because the remaining competitors perform either comparably or worse in Table \ref{tab:number}, while  $\hat r_{SSD}$ is always very close to $\hat r_{SMD}$. 
		
		\begin{figure}[h]
			\centering
			\begin{minipage}{0.24\textwidth}
				\centering
				\includegraphics[width=4cm,height=4cm]{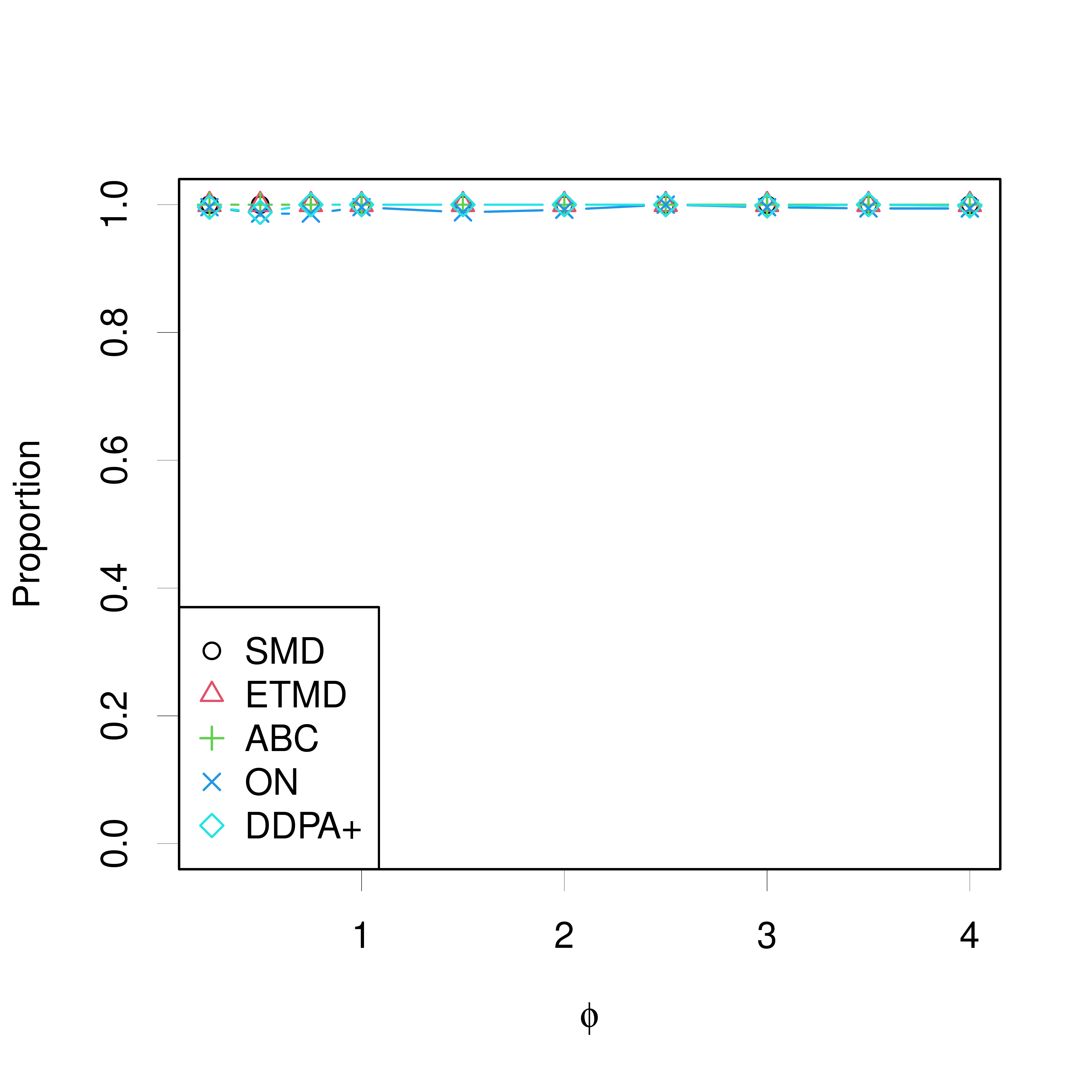}
				(a). $\phi$ changes
			\end{minipage}
			\begin{minipage}{0.24\textwidth}
				\centering
				\includegraphics[width=4cm,height=4cm]{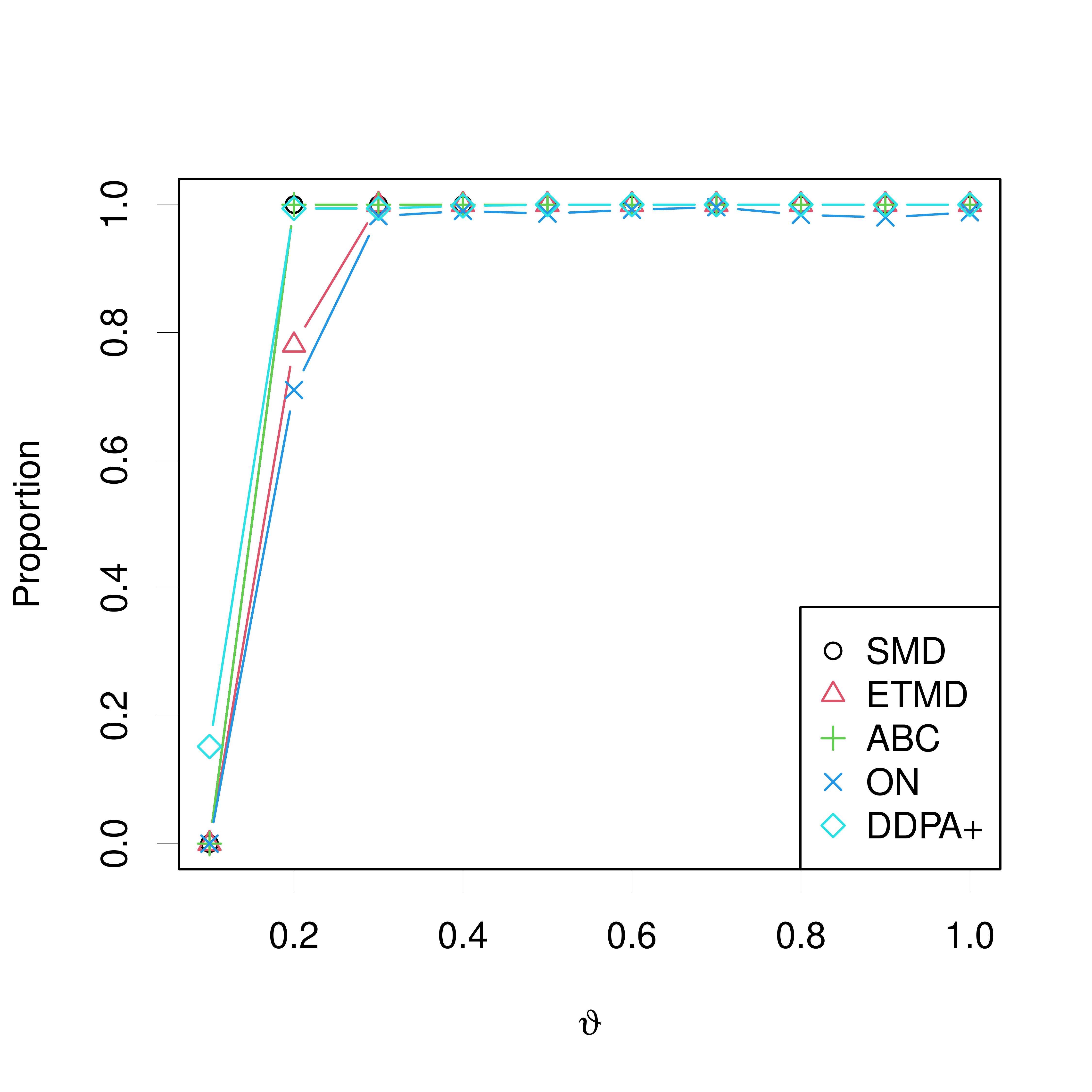}
				(b). $\vartheta$ changes
			\end{minipage}
			\begin{minipage}{0.24\textwidth}
				\centering
				\includegraphics[width=4cm,height=4cm]{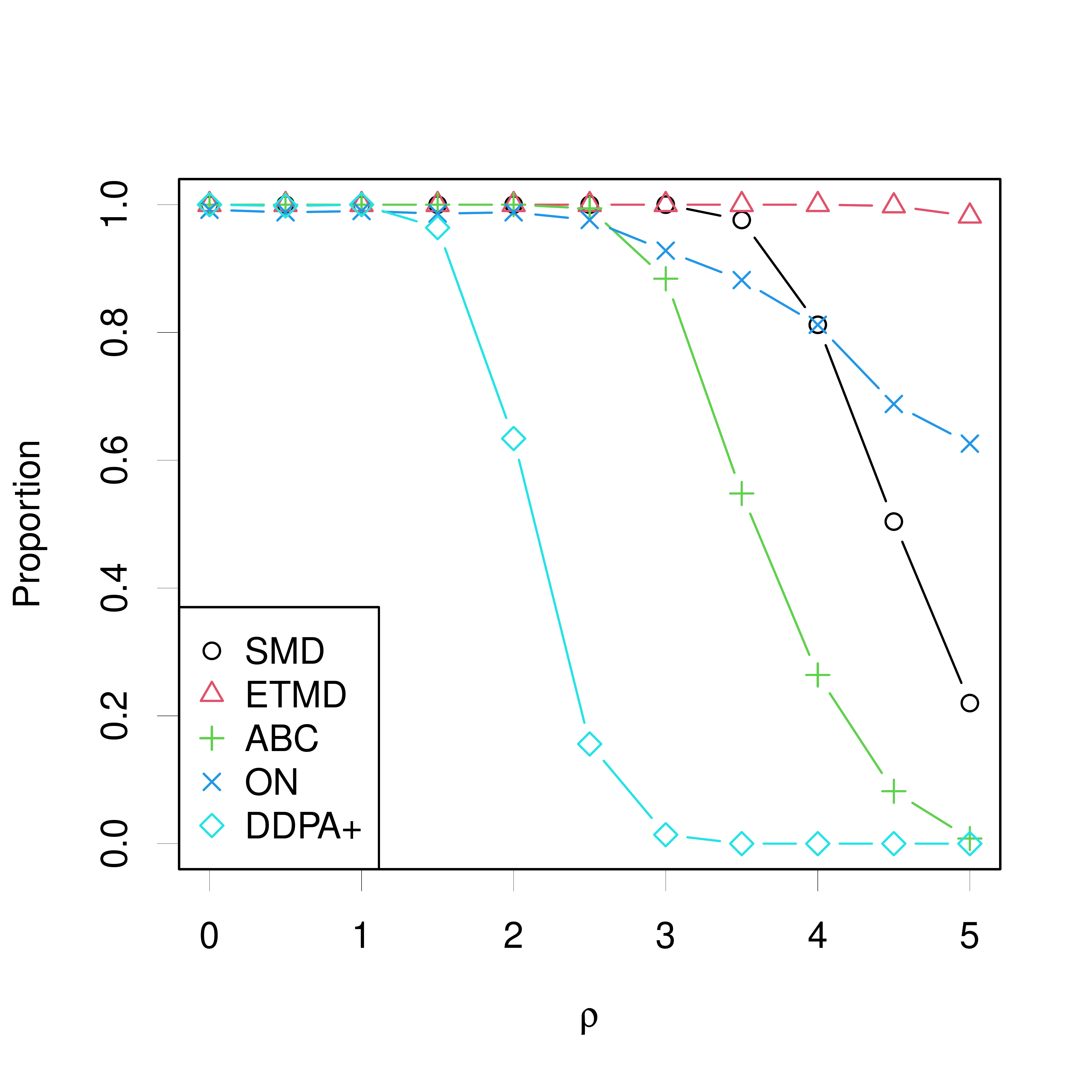}
				(c). $\rho$ changes
			\end{minipage}
			\begin{minipage}{0.24\textwidth}
				\centering
				\includegraphics[width=4cm,height=4cm]{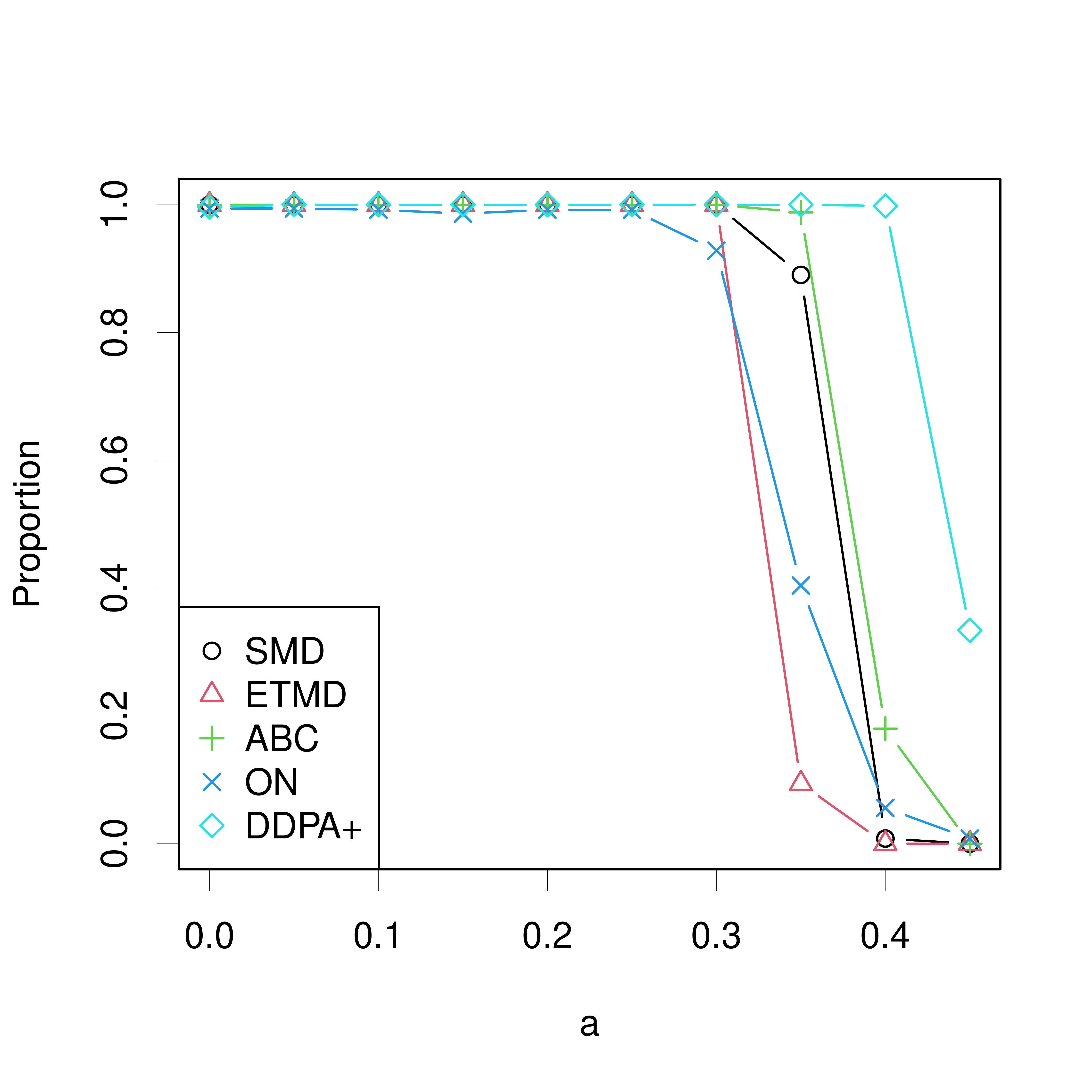}
				(d). $a$ changes
			\end{minipage}
			\caption{The proportions of exact estimation of $r$ over 500 replications by different methods with various data generating parameters. }\label{fig: sense1}
		\end{figure}
		
		In Figure \ref{fig: sense1}(a), all methods are very robust to the parameter $\phi$. $\hat r_{DDPA_+}$ performs the best when the factors are weak (small $\vartheta$ or large $a$), but losing accuracy significantly when $\rho$ is large. On the contrary, $\hat r_{ETMD}$ is the most stable when the noises have strong cross-sectional dependence, but less reliable under weak factors. From the figures, we can see that no method can always outperform the others, which is understandable  in finite samples due to the ambiguity between weak factors and strong idiosyncratic errors.  At least, we can conclude that the proposed methods provide new direction for determining  the number of common factors, and perform comparably and stably under most scenarios.}

	\section{Real data example}\label{sec: real data}
	{We use the proposed approaches to analyze a financial data set, which is an open resource from Kenneth R. French's web page at \url{http://mba.tuck.dartmouth.edu/pages/faculty/ken.french}. It contains monthly returns of 100 portfolios formed on capital size and book-to-market ratio. We focus on the period from January-1964 to December 2022. We standardize  the return series one by one and impute missing values by linear interpolation (missing rate 0.23\%), leading to a data matrix $\Xb_{p\times n}$ with  $p=100$ and $n=708$. 
		
		Figure \ref{fig:ff}(a) shows the eigenvalues of the sample covariance matrix associated with $\Xb$. There  is one extremely large eigenvalue, indicating the existence of at least one powerful common factor. The second and third largest eigenvalues also deviate slightly from the bulk, but they are much smaller than the first one. We use the proposed three approaches to determine the number of common factors. The estimated factor numbers and the computational costs (in seconds) are shown in Table \ref{tab: ff}, including the results of competitors from the literature. The tuning parameters are the same as those in Table \ref{tab:number}. The proposed three methods, $\hat r_{SMD}$, $\hat r_{SSD}$, $\hat r_{ETMD}$ and two competitors $\hat r_{IC}$, $\hat r_{TRAP}$ output the same result that $\hat r =3$. $\hat r_{ER}$ and $\hat r_{ON}$ only report the existence of one factor. $\hat r_{ABC}$, $\hat r_{ED}$ and $\hat r_{DDPA}$ report 4 factors while $\hat r_{ETC}$ and $\hat r_{ETZ}$ report even more factors. The results are consistent with our findings from the numerical studies. In terms of computation, the proposed  $\hat r_{SMD}$ and $\hat r_{SSD}$ are moderately expensive. $\hat r_{ETMD}$ and $\hat r_{ETZ}$ are more costly due to the additional bootstrap step to obtain the approximated critical value.

		\begin{figure}[htbp]
			\centering
			\begin{minipage}{0.24\textwidth}
				\centering
				\includegraphics[width=4cm,height=4cm]{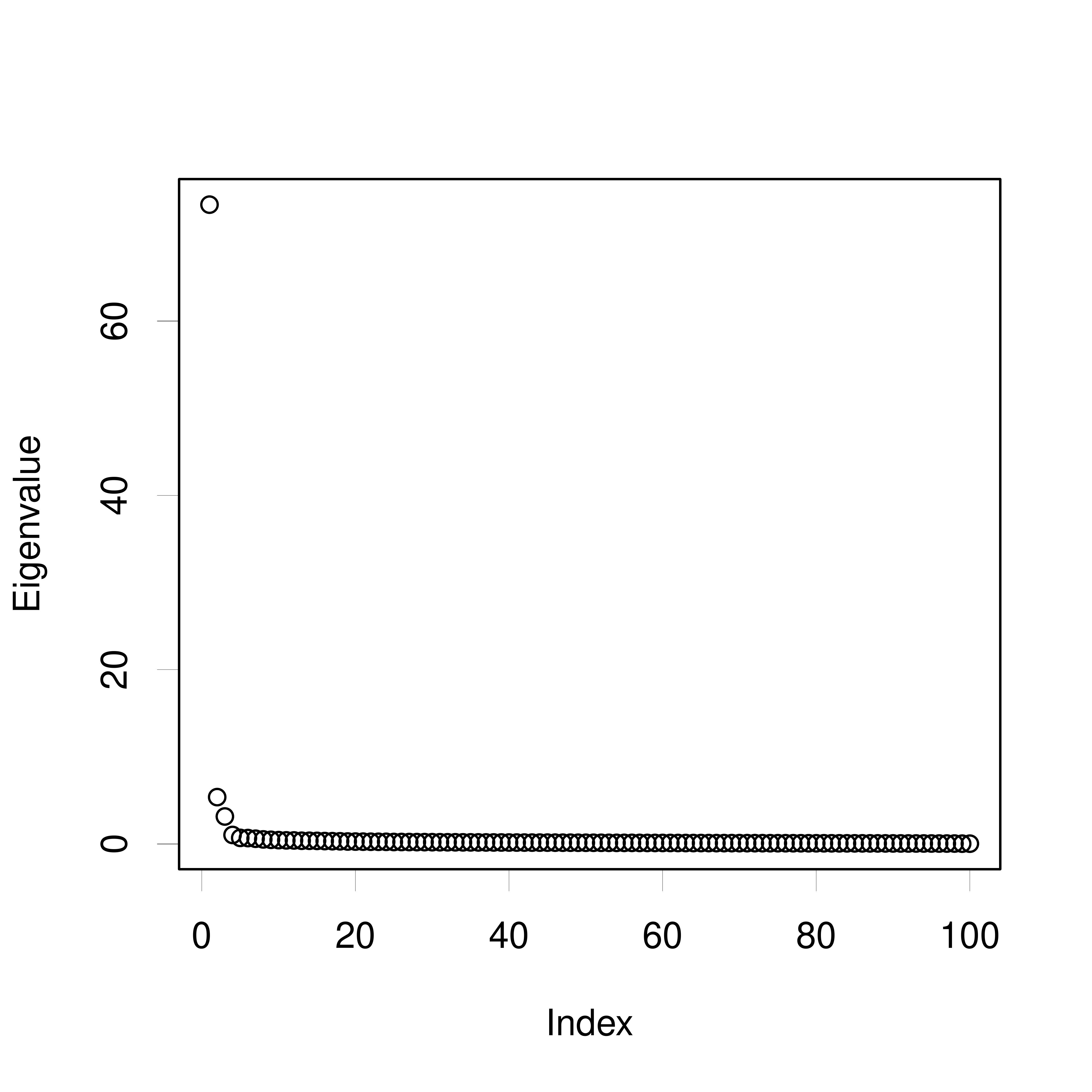}
				(a)
			\end{minipage}
			\begin{minipage}{0.24\textwidth}
				\centering
				\includegraphics[width=4cm,height=4cm]{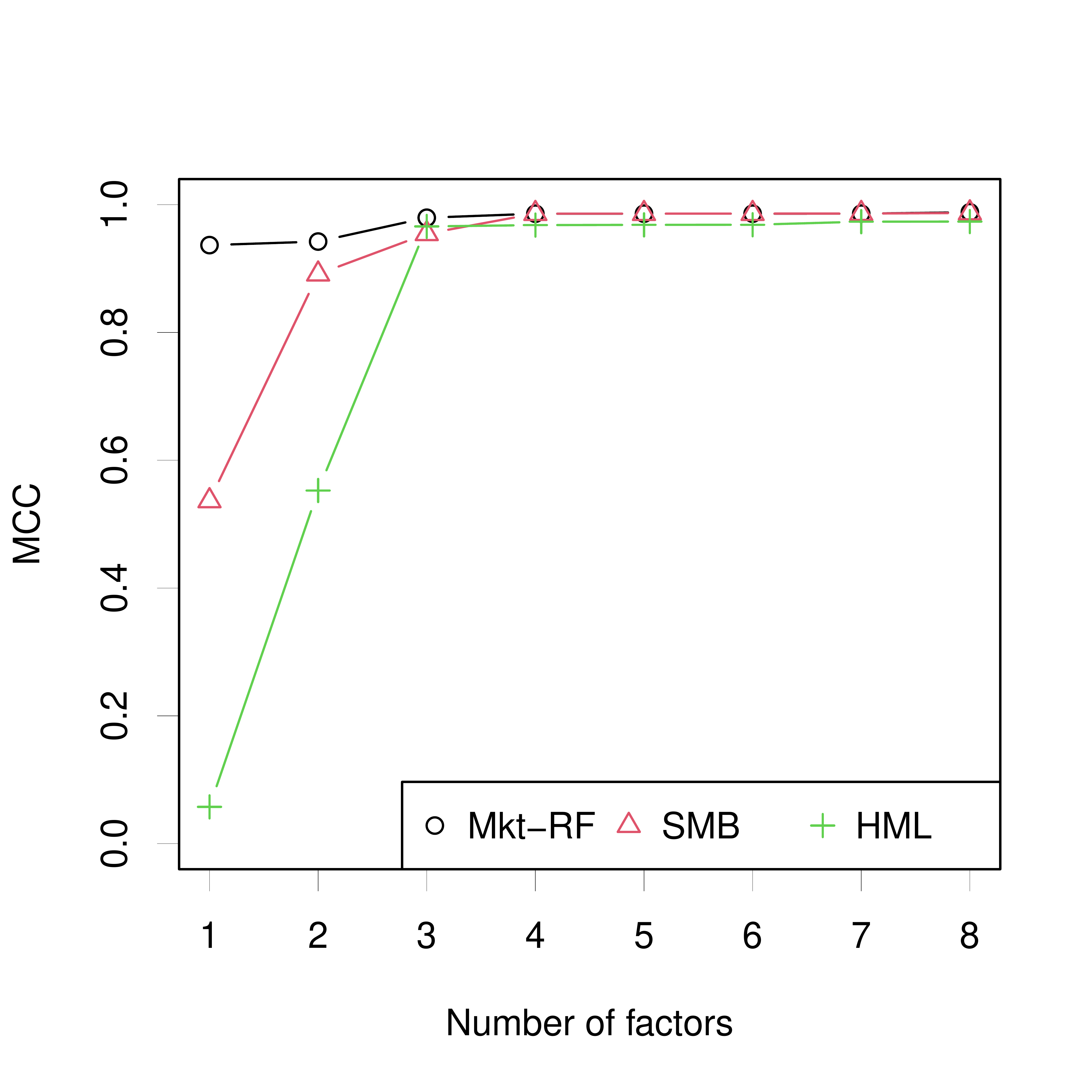}
				(b)
			\end{minipage}
			\begin{minipage}{0.24\textwidth}
				\centering
				\includegraphics[width=4cm,height=4cm]{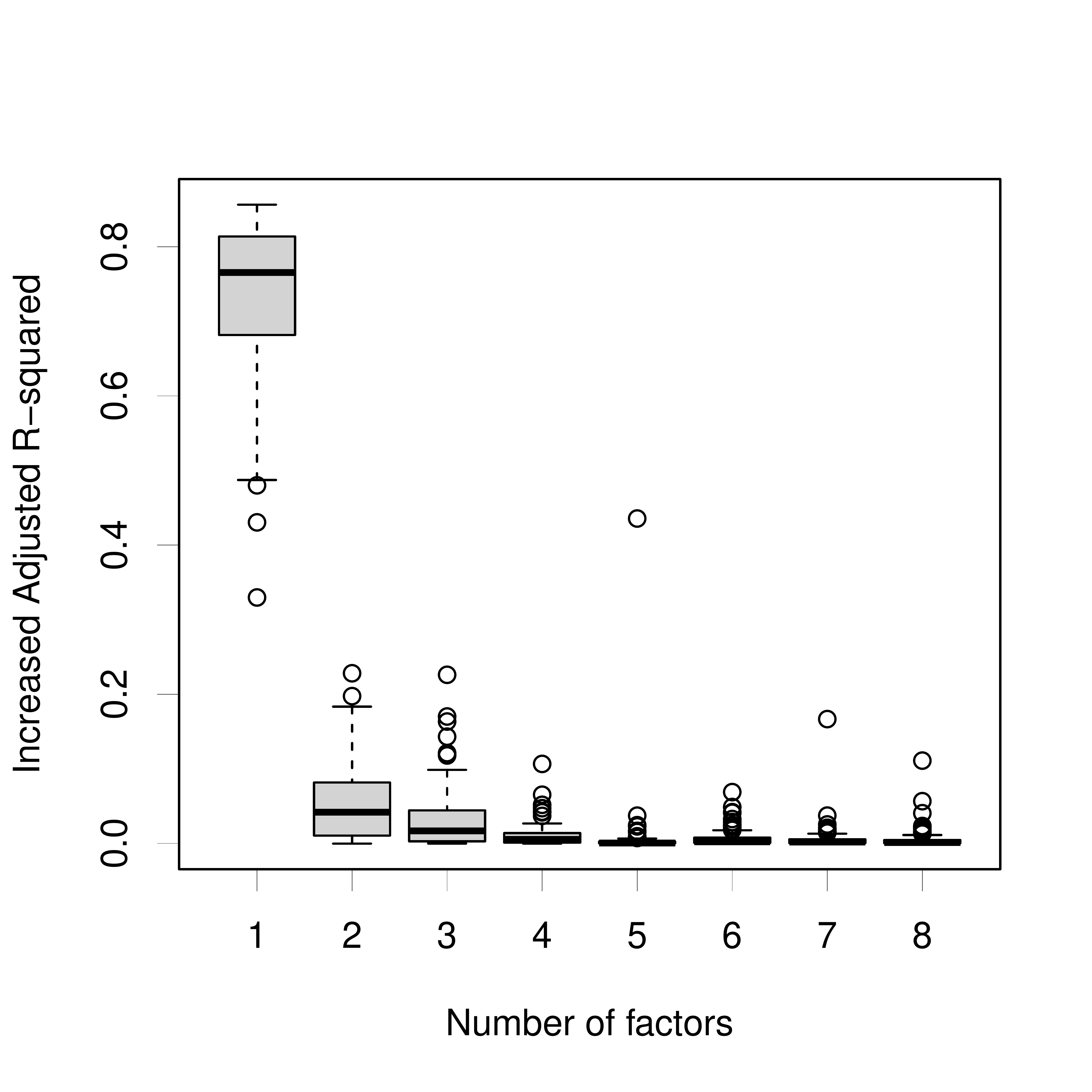}
				(c)
			\end{minipage}
			\begin{minipage}{0.24\textwidth}
				\centering
				\includegraphics[width=4cm,height=4cm]{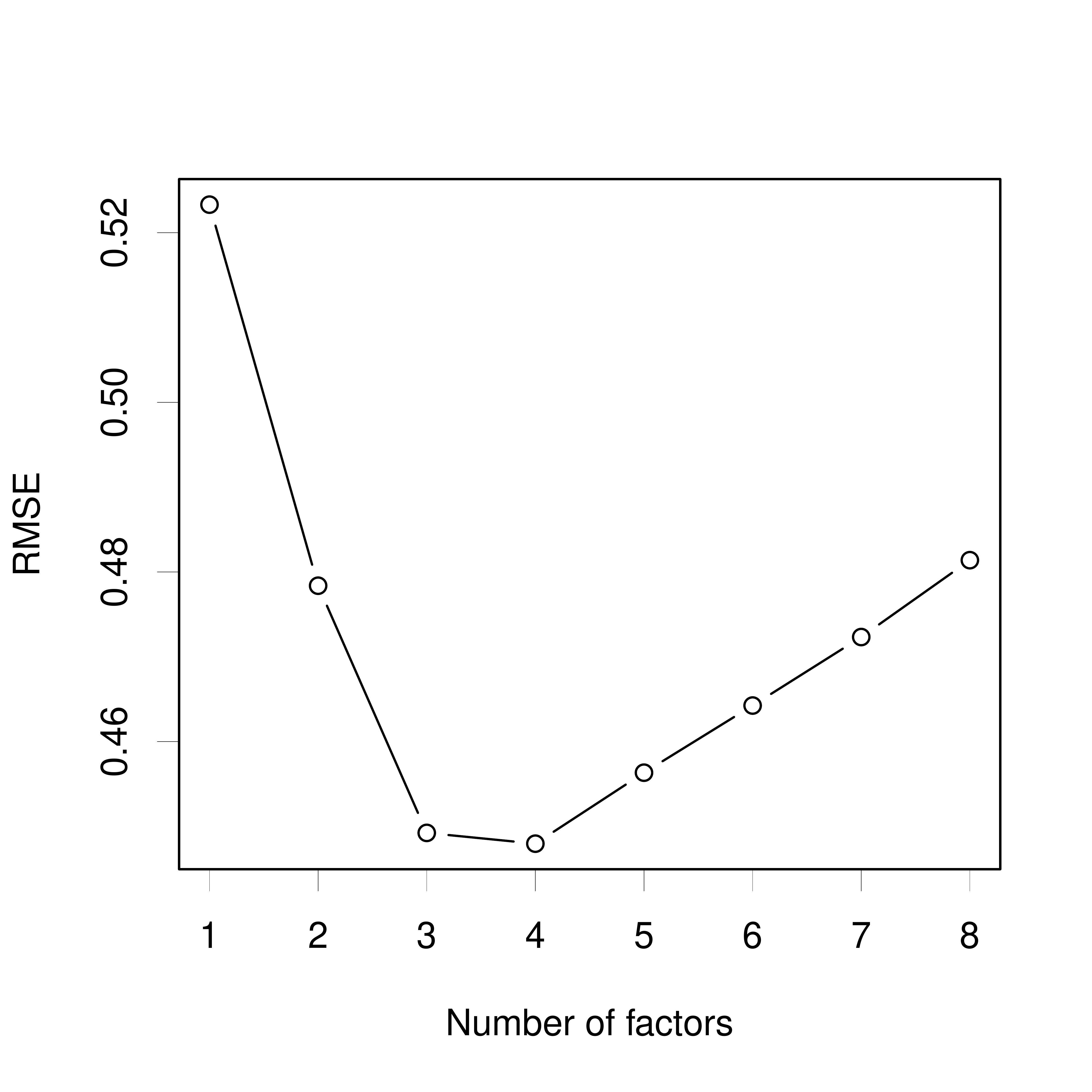}
				(d)
			\end{minipage}
			
			\caption{Figures for the Fama-French portfolio data set: (a) eigenvalues of sample covariance matrix. (b) multiple correlation coefficients between the common factors from the data set and the market risk factor ($\text{Mkt}-\text{RF}$), SMB and HML, respectively as $\hat r$ increases. (c) boxplots of the explanatory power (adjusted R-squared) to the 100 return series when new factors come into system. (d) out-of-sample imputing error (RMSE) as $\hat r$ grows.}\label{fig:ff}
		\end{figure}
		
		\begin{table}[htbp]
			\addtolength{\tabcolsep}{2pt}
			\caption{Estimated number of common factors and the computational cost (in seconds) for the Fama-French portfolio return data set by different methods.\label{tab: ff}}
			\renewcommand{\arraystretch}{0.6}
			\scalebox{0.75}{ 
				\begin{tabular}{lllllllllllllllll}
					\toprule[1.2pt]
					&$\hat r_{SMD}$&$\hat r_{SSD}$&$\hat r_{ETMD}$&$\hat r_{IC}$&$\hat r_{ABC}$&$\hat r_{ER}$&$\hat r_{TRAP}$&$\hat r_{ON}$&$\hat r_{ED}$&$\hat r_{ETC}$&$\hat r_{ETZ}$&$\hat r_{DDPA}$\\\midrule[1.2pt]
					$\hat r$&3&3&3&3&4&1&3&1&4&7&6&4\\
					Cost (s)&1.035&0.961&6.276&0.029&0.261&0.024&0.052&0.335&0.027&0.057&9.961&0.040
					\\\bottomrule[1.2pt]
			\end{tabular}}
		\end{table}
		
		It's well-known in finance that the return of a portfolio is potentially driven by the Fama-French 3 factors, i.e., market risk factor, SMB factor and HML factor, which is consistent with the number of factors estimated by the proposed methods. To check this, we calculate the multiple correlation coefficient (MCC) between the factor score series and each of the Fama-French 3 factors. The factor scores are estimated by PCA given $\hat r$ while the monthly returns of the Fama-French 3 factors are provided by Kenneth's web page. Figure \ref{fig:ff}(b) shows the respective MCCs with $\hat r$ growing. It's seen that the leading three factors from this data set are highly-correlated with the Fama-French 3 factors, while adding the fourth factor only slightly increase the MCC with SMB.

		Next, we investigate how new factors contribute to explaining the variation of the return series. 	For each series, given $\hat r$, we regress the  portfolio return on the estimated factor scores, and use the increased adjusted R-squared to represent the explained variation when more factors are used in the regression. Figure \ref{fig:ff}(c) shows the boxplots of the increased explanatory power to the 100 return series as $\hat r$ grows. It's seen that the explanatory power of the third factor is non-negligible, while the gain from the fourth factor is minor.
		
		Lastly, we verify how new factors help in imputing missing values. We randomly select 50\% of the portfolios, denoted by a set $\mathcal{P}$, and 50\% of the time periods, denoted by a set $\mathcal{N}$. We take $x_{ij}$ as missing when $i\in\mathcal{P}$ and $j\in\mathcal{N}$ . Borrowing the idea from \cite{bai2021matrix}, we first use $\{x_{ij},j\ne \mathcal{N}\}$ to estimate the factor loading space $\hat \Lb$ by PCA given $\hat r$, and then estimate $\hat\Fb$ based on $\{x_{ij},i\ne \mathcal{P}\}$ and $\hat \Lb$. The missing entries are imputed by $\hat x_{ij}=\hat\bL_i^\top\hat\bbf_j$ for $i\in\mathcal{P}$ and $j\in \mathcal{N}$. The out-of-sample imputing error is calculated in terms of Root of Mean Squared Error (RMSE). Because the missing set is selected randomly, we repeat the above procedure 500 times and report the mean of RMSE in Figure \ref{fig:ff}(d) to reduce sampling bias. It's shown that the imputing error is minimized at $\hat r=4$. However, the improvement from $\hat r=3$ to $\hat r=4$ is minor.
		
		In conclusion, we believe that $\hat r=3$ or $\hat r=4$ will be the reasonable  decision for this data set. However, the gain from the fourth common factor is relatively minor in the above experiments. $\hat r =3$ might be a more suitable choice,  which is also consistent with  the asset pricing theory.  In the Supplement, we analyze another real data set in macroeconomics, where $\hat r_{SMD}$, $\hat r_{SSD}$ and $\hat r_{ETMD}$ sill report the same and reasonable result, but $\hat r_{IC}$ and $\hat r_{TRAP}$ lead to underestimation of the factor number.  The proposed three methods are not sensitive to tuning parameters in both examples.
	}

	\section{Conclusion and discussion}\label{sec:discussion}
	The current paper contributes to understanding the effects of bootstrap to the eigenvalue  distribution of sample covariance matrix under high-dimensional factor models or spiked covariance models. It also contributes to the literature of determining the number of common factors or spikes. {In the current paper, we require the spiked eigenvalues driven by common factors to be diverging, which is more stringent than  typical assumptions in the literature of  BBP phase transition; see \cite{bloemendal2016principal}. 
		One reason is that the bootstrap procedure changes the phase transition boundary. A novel and interesting finding is that the exact phase transition boundary seems to mainly depend on the order statistics of the bootstrap resampling weights. By bootstrapping from different distributions, it's possible to increase or decrease the typical BBP phase transition boundary. To verify this, in the Supplement,  we have done more simulation studies on the performance of $\hat r_{ETMD}$ by bootstrapping from more general distribution families such as $Poisson(1)$. It shows that the bootstrap procedure still works, but with different requirement on the strength of common factors. In other words, the bootstrap procedure provides a new direction for documenting the number of factors with different strength. 
		We are interested in studying the exact phase transition boundary under the bootstrap framework with general resampling weights. We would also like to relax the constraint on independent idiosyncratic errors and study more general time series settings.
		It's also of interest to consider the sample correlation matrix instead of covariance matrix after bootstrap using similar techniques in \cite{bao2019tracy}, which usually possesses scale invariant property. We leave these as future works.}

	\section{Acknowledgment}
	The authors would like to thank the editor, associate editor and three anonymous reviewers for their valuable comments and suggestions. We would also like to thank Bao Zhigang  for letting us know the reference \cite{kwak2021extremal}. Long Yu's research is partially supported by the Fundamental Research Funds for the Central Universities, China.

	\section{Supplementary material}

	The supplementary material is composed of five sections. Section \ref{sec: add} provides additional simulation results and real data analysis.  Section \ref{sec: supplement A} proves the theoretical results in Section \ref{sec: spiked} of the main paper, corresponding to the test with spiked eigenvalues. Section \ref{sec: technical} contains some useful technical lemmas used in Section \ref{sec: supplement A}. Section \ref{secc} provides preliminary definitions and technical lemmas for the results in Section \ref{sec:non-spiked} of the main paper, corresponding to the test with non-spiked eigenvalues. Section \ref{sec: proof non spike} completes the proof. In the proof, $\|\Ab\|=\sqrt{\lambda_1(\Ab\Ab^\top)}$ denotes the spectral norm and $\|\Ab\|_F=\sqrt{\text{tr}(\Ab\Ab^\top)}$ denotes the Frobenius norm. $\text{diag}(\Ab)$ is the diagonal matrix whose diagonal elements are the same as those of $\Ab$. $c$ and $C$ indicate some small and large constants which may vary in different lines, respectively.
	
	\begin{appendices}
	\section{Additional simulation results and real data analysis}\label{sec: add}
	\subsection{Simulation: bootstrap from more distribution families}
	Lemma \ref{lemma: lambda_0} in the main paper has indicated that $\hat\lambda_{r+1}$ depends on the order statistics of the bootstrap resmapling weights $w_j$'s under the multiplier bootstrap scheme. When $w_j\overset{i.i.d.}{\sim}\mathcal{N}(0,1)$, $w_{(1)}$ follows the Gumbel distribution and diverges with rate $\log n$. This is also one of the reasons why we require the eigenvalues driven by common factors to be diverging, while this assumption is not proposed in the literature of typical BBP-type phase transition. In other words, the bootstrap procedure will change the phase transition boundary. As claimed in the main paper, the proposed methods can be naturally generalized by bootstrapping from other distributions. We verify this argument in the following experiment. 
	
	We mainly focus on $\hat r_{ETMD}$ to avoid the calculation of asymptotic variances for $\hat r_{SSD}$ and $\hat r_{SMD}$. Instead of bootstrapping from $Exp(1)$, now we will  also try $Poisson(1)$, $Uniform(0.5,1.5)$ and $\mathcal{X}^2(1)$ in Algorithm \ref{alg1}. To see how the bootstrap resampling weights affect the phase transition boundary, we slightly modify the data-generating parameters in Figure \ref{fig: sense1}(b), where the data are from factors plus i.i.d. noises. To be specific, we take the left singular vectors of $\Lb$ as the new loading matrix $\tilde\Lb$, and generate $\Fb$ and $\Eb$ from i.i.d. $\mathcal{N}(0,1)$. Then $\Xb=\vartheta\tilde\Lb\Fb^\top+\Eb$. Under such cases, the leading population eigenvalues will be $\vartheta^2+1$ while the typical BBP phase transition boundary is $(1+\sqrt{p/n})^2=4$ when $n=p$. We will investigate how the performance of $\hat r_{ETMD}$ varies when bootstrapping from different distributions as $\vartheta$
	grows. The parallel analysis method $\hat r_{DDPA_+}$ from \cite{dobriban2019deterministic} is taken as a benchmark to show the typical BBP phase transition. The proportions of exact estimation are shown in Figure \ref{fig: weight} over 500 replications by different methods.
	
	\begin{figure}[h]
		\centering
		\includegraphics[width=7.5cm,height=5cm]{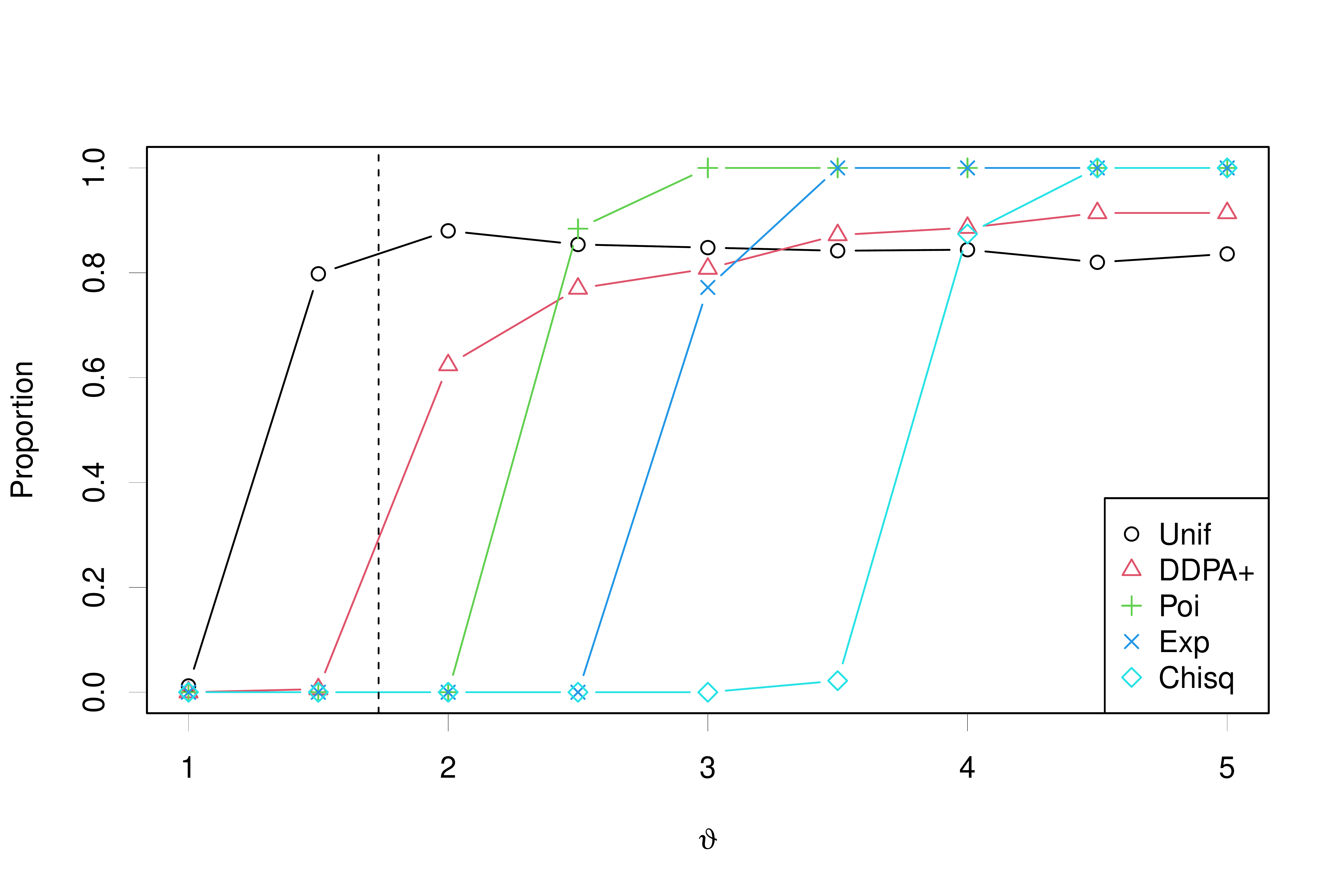}
		\caption{Proportions of exactly estimating $r$ when bootstrapping from different distributions using $\hat r_{ETMD}$, over 500 replications as $\vartheta$ grows. ``Unif'', ``Poi'', ``Exp'' and ``Chisq'' stand for bootstrapping from  $Uniform(0.5,1.5)$, $Poisson(1)$, $Exp(1)$ and $\mathcal{X}^2(1)$, respectively. The dashed vertical line shows the typical phase transition boundary.}\label{fig: weight}
	\end{figure}
	
	It's seen that when bootstrapping from different distributions, $\hat r_{ETMD}$ will still work, but requiring different factor strength. This is understandable because the eigenvalue thresholding method only works when the spiked eigenvalues exceed the phase transition boundary. Motivated by Lemma \ref{lemma: lambda_0}, the phase transition boundary after bootstrap mainly depends on resampling weights. Figure \ref{fig: weight} shows that the transition boundary  grows gradually when bootstrapping from $Uniform(0.5,1.5)$, $Poisson(1)$, $Exp(1)$ and $\mathcal{X}^2(1)$. One potential reason is that the tail of the density becomes thicker and thicker, and the expectation of the associated leading order statistics becomes larger and larger.  In other words, the bootstrap procedure actually provides a flexible way to increase or decrease the phase transition boundary. As a result, we are able to document the number of common factors with different strength, just by resampling from different distributions. Another interesting finding is that when bootstrapping from $Uniform(0.5,1.5)$, $\hat r_{ETMD}$ can still accurately determine $r$ with a large frequency even if $\vartheta$ is below the typical phase transition boundary. 
	
	\subsection{Simulation: robustness to tuning parameters}
	
	We are also interested in how the proposed methods rely on the tuning parameters, i.e., the predetermined upper bound $r_{max}$, the significance level $\alpha$, the number of bootstrap replications $B$ and $R$. 
	For better illustration, we will focus on the most challenging case in Table \ref{tab:number} where $\vartheta=1$, $\rho=3$, $a=0.25$ and $n=p=200$. Similarly to Figure \ref{fig: sense1}, we  plot the proportions of exact estimation over 500 replications for the proposed three methods when one tuning parameter changes but the others are fixed. The benchmark setting for the parameters are $r_{\max}=8$, $\alpha=0.05$, $B=200$ and $R=400$.  The results are reported in Figure \ref{fig: sense2}, where the proposed methods are very accurate under all the considered tuning parameter settings with over 90$\%$ exact  estimations.

	\begin{figure}[h]
		\centering
		\begin{minipage}{0.24\textwidth}
			\centering
			\includegraphics[width=4cm,height=4cm]{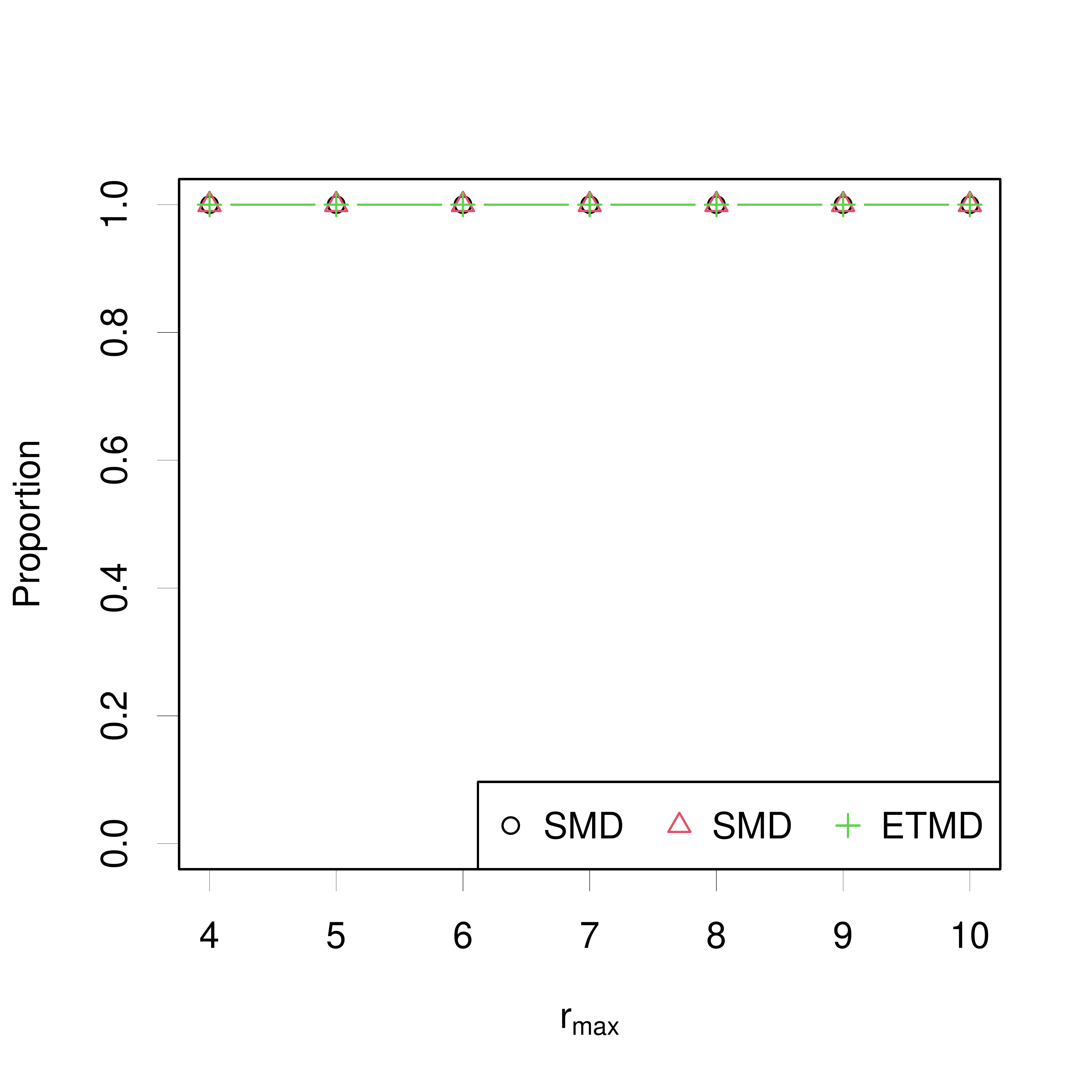}
			(a). $r_{\max}$ changes
		\end{minipage}
		\begin{minipage}{0.24\textwidth}
			\centering
			\includegraphics[width=4cm,height=4cm]{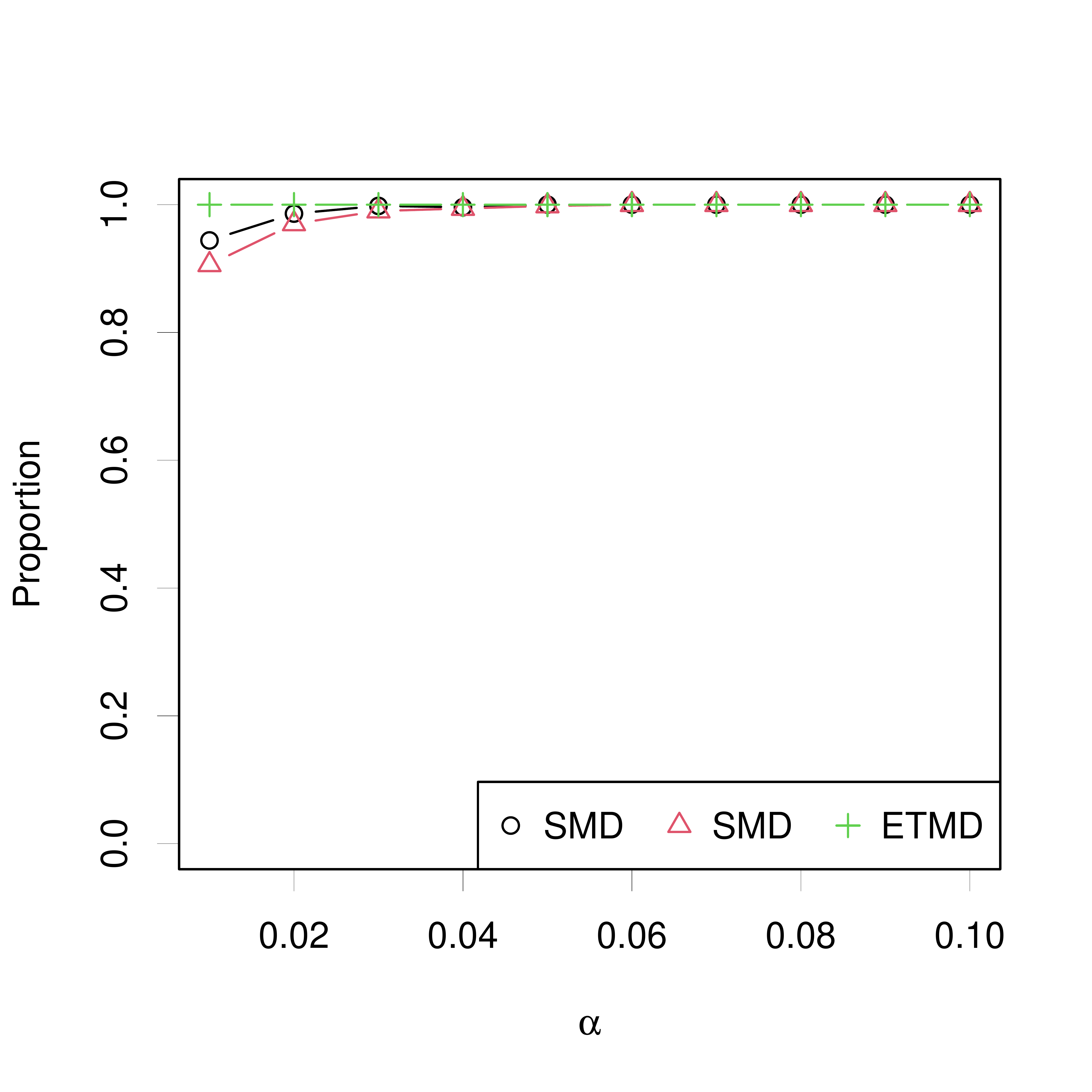}
			(b). $\alpha$ changes
		\end{minipage}
		\begin{minipage}{0.24\textwidth}
			\centering
			\includegraphics[width=4cm,height=4cm]{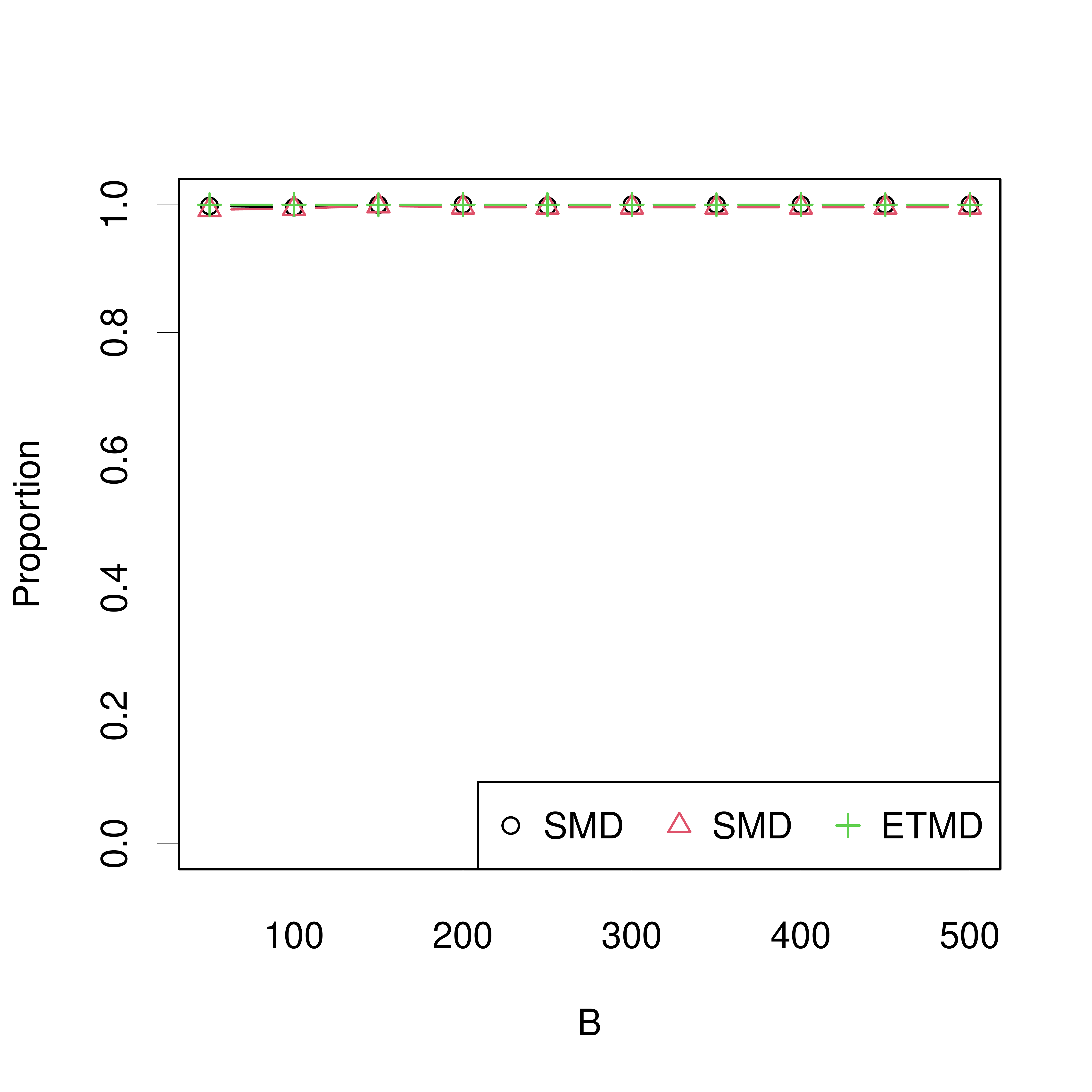}
			(c). $B$ changes
		\end{minipage}
		\begin{minipage}{0.24\textwidth}
			\centering
			\includegraphics[width=4cm,height=4cm]{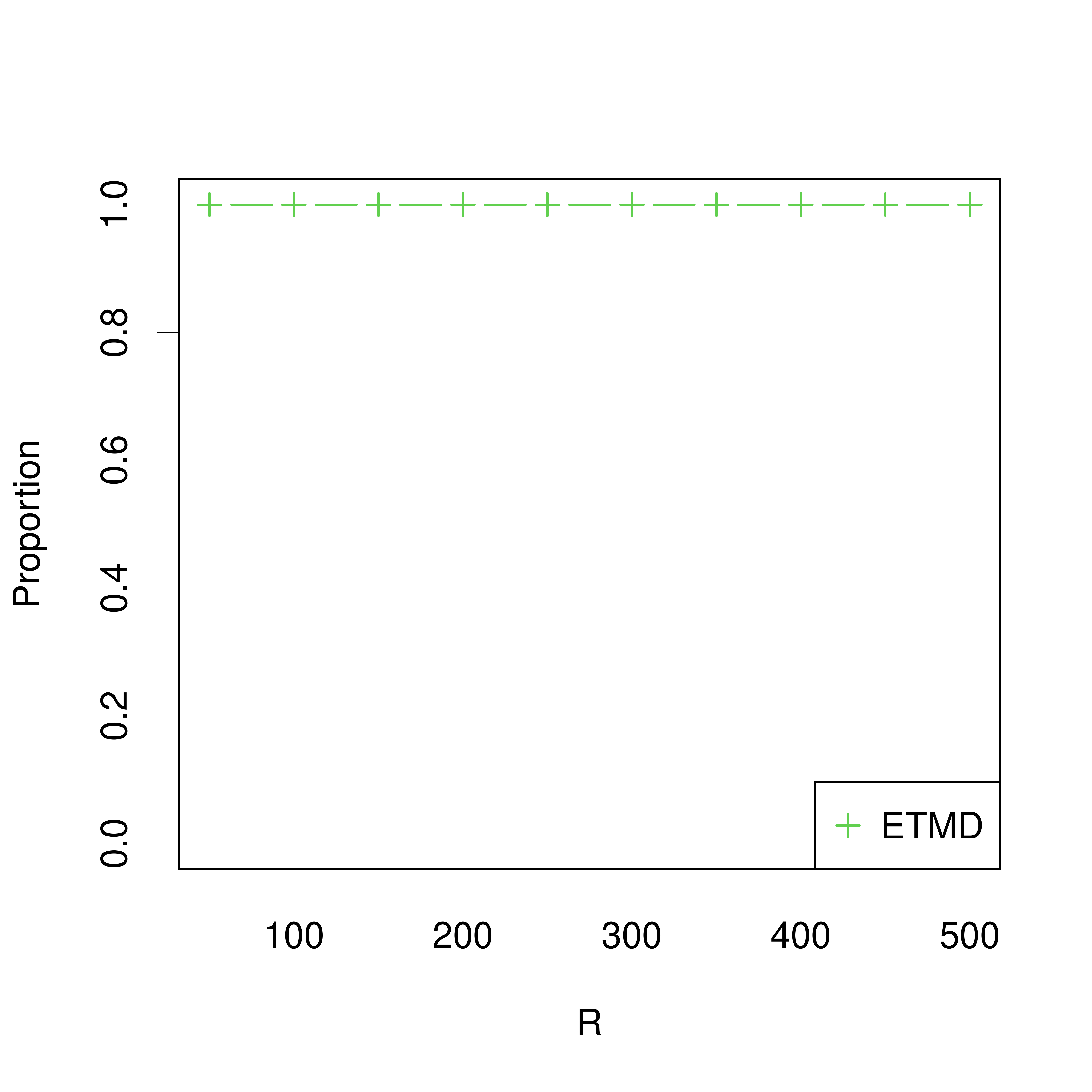}
			(d). $R$ changes
		\end{minipage}
		\caption{The proportions of exactly estimating $r$ over 500 replications by the proposed three methods with different tuning parameters. }\label{fig: sense2}
	\end{figure}

	It's worth mentioning that our approaches do not require $\alpha$ to be asymptotically vanishing by adding the decision rule. This is different from those traditional test-based methods  in \cite{onatski2009testing},\cite{trapani2018randomized}, or the thresholding methods  in \cite{cai2020limiting} and \cite{ke2021estimation}. With the decision rule, we allow more errors for the size and power of the tests. Figure \ref{fig: dr} plots the values of $D_i^s(\alpha,B)$ and $D_i^{ns}(\alpha,B)$ in the decision rule, for the proposed three methods $\hat r_{SMD}$,  $\hat r_{SSD}$ and  $\hat r_{ETMD}$ when $i=1,\ldots,8$. The data generating parameters are the same as those in Table \ref{tab:number} with $\vartheta=1$, $a=0.25$, $\rho=3$ and $n=p=400$. It's seen that in each panel, the three spiked eigenvalues and non-spiked eigenvalues are well separated by the decision rule. This is the reason why we can achieve nearly perfect estimation in most settings. However, if we only implement the bootstrap once, Figure \ref{fig: dr} (a) and (b) show that there is a large probability that the fourth eigenvalue may be identified as a spiked eigenvalue, because $D_4^s(\alpha,B)$ can deviate from the theoretical converging point 0. Fortunately, the decision rule fixes this problem by repeating the bootstrap procedure to stabilize the results.
	
	\begin{figure}[h]
		\centering
		\begin{minipage}{0.31\textwidth}
			\centering
			\includegraphics[width=4.8cm,height=4.8cm]{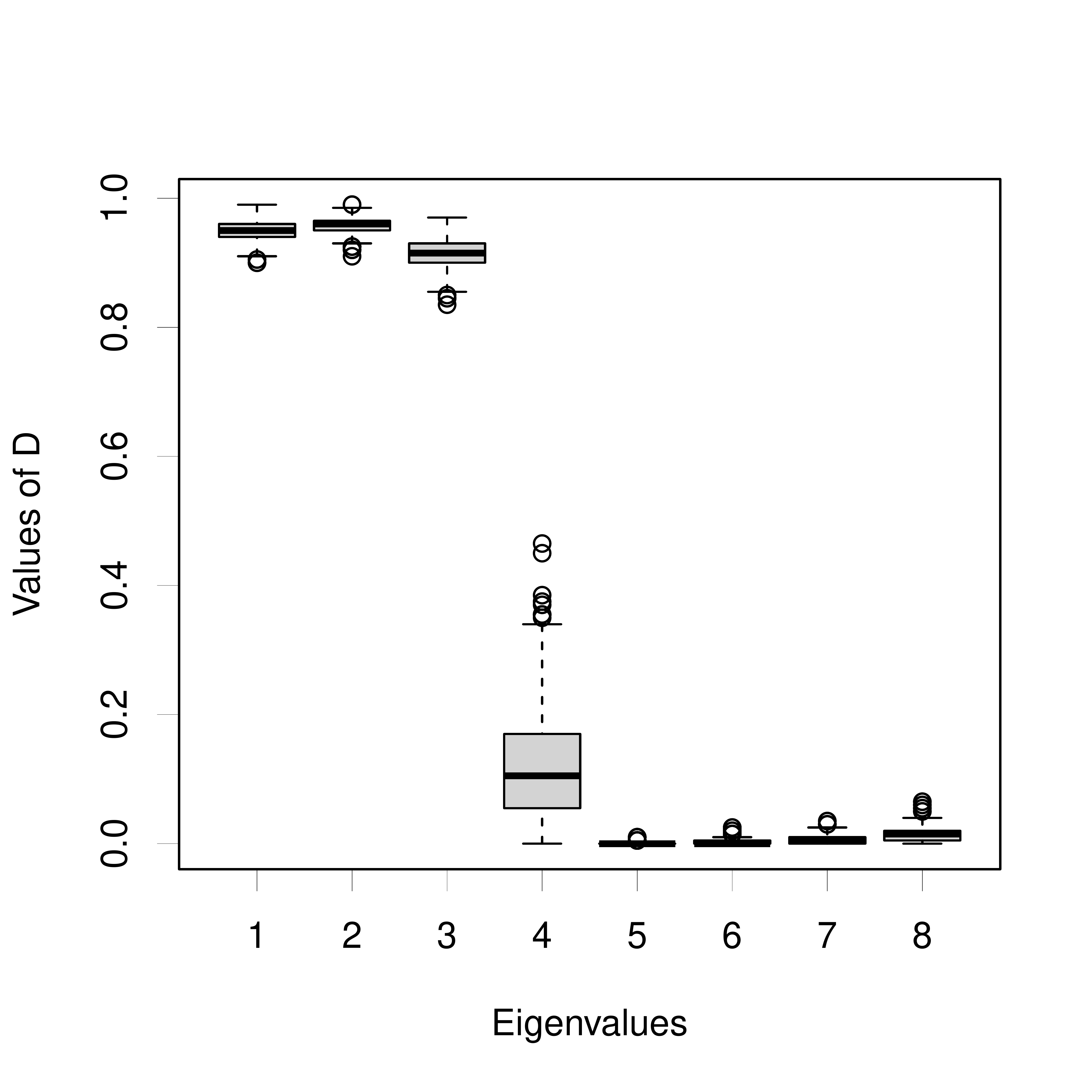}
			(a). $\hat r_{SMD}$
		\end{minipage}
		\begin{minipage}{0.31\textwidth}
			\centering
			\includegraphics[width=4.8cm,height=4.8cm]{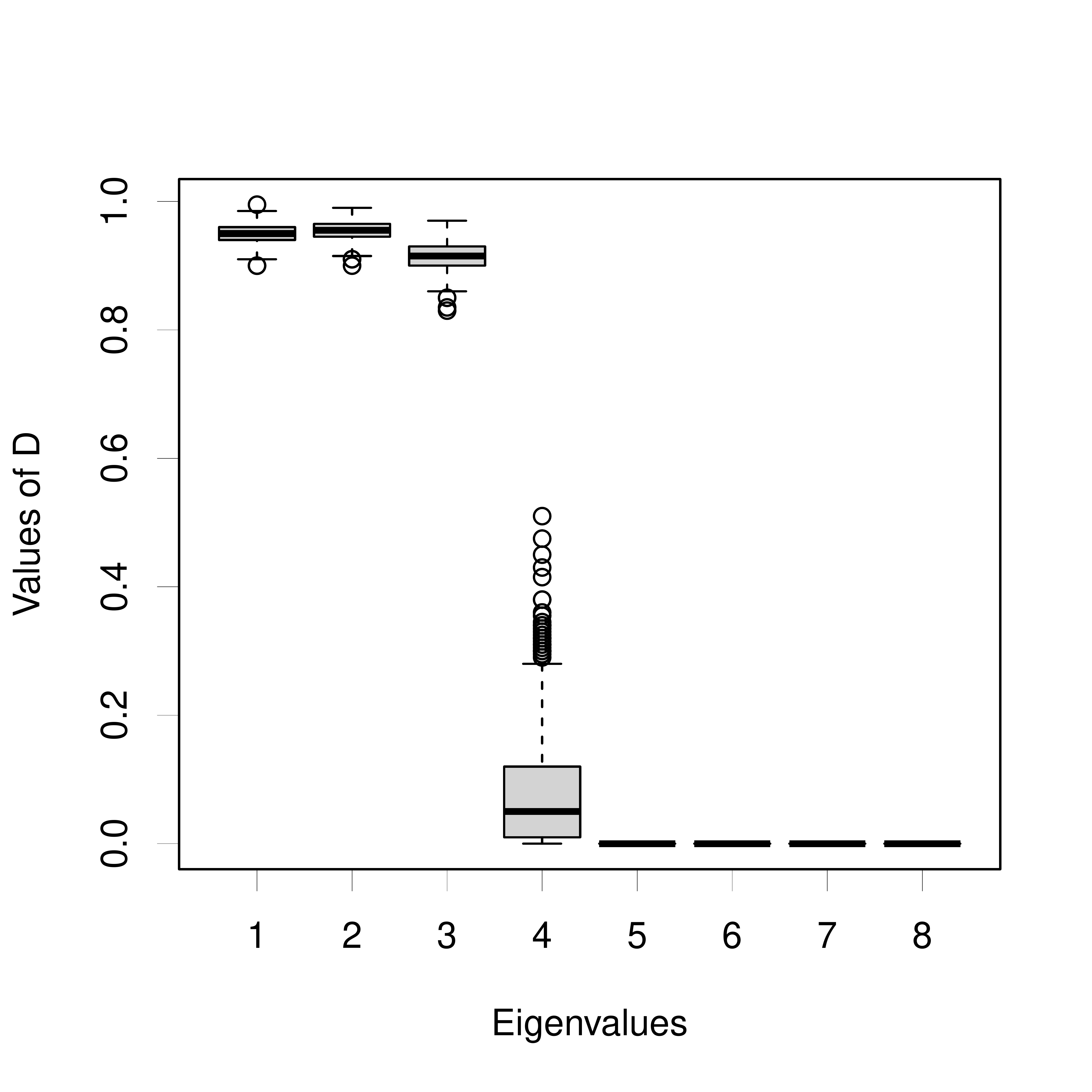}
			(b). $\hat r_{SSD}$
		\end{minipage}
		\begin{minipage}{0.31\textwidth}
			\centering
			\includegraphics[width=4.8cm,height=4.8cm]{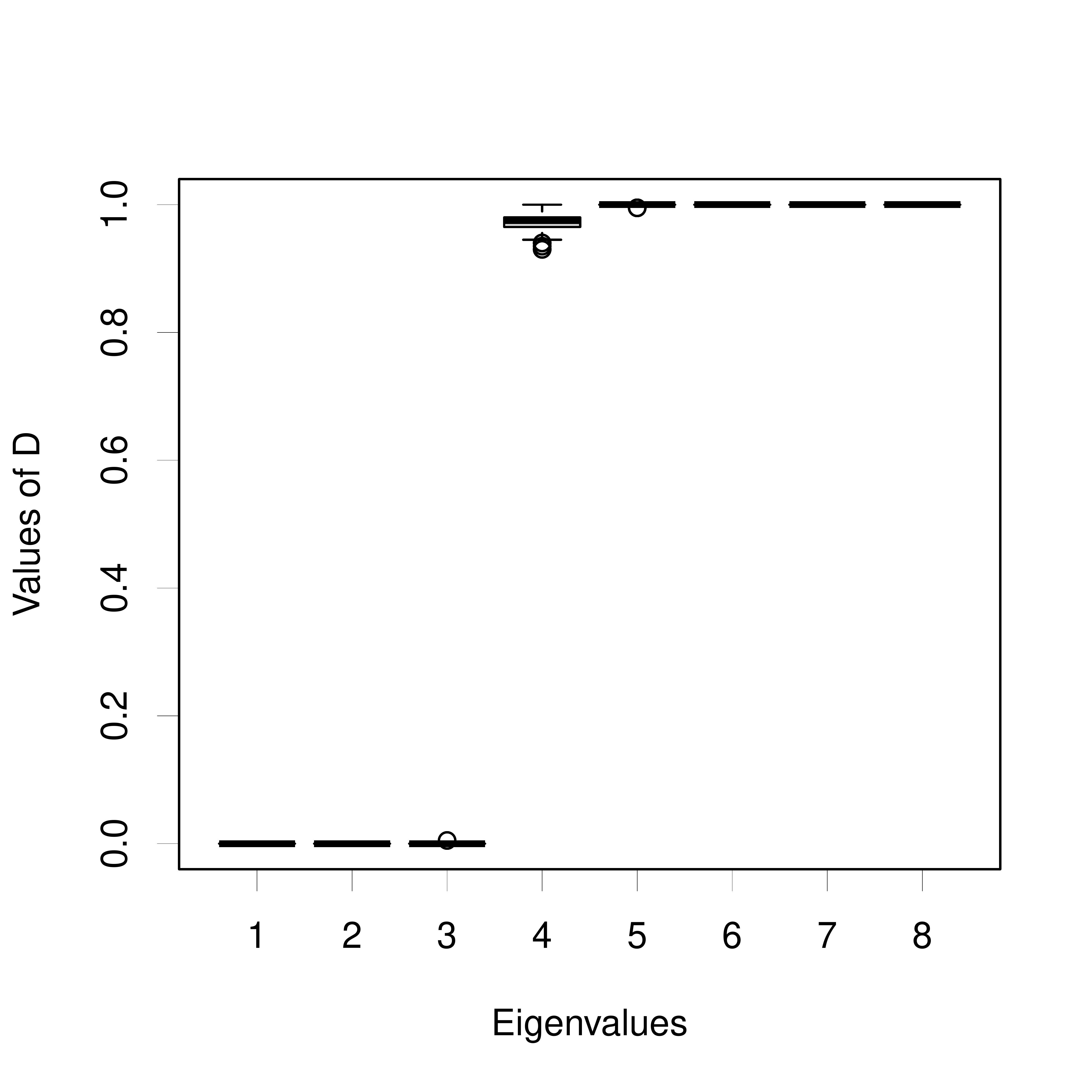}
			(c). $\hat r_{ETMD}$
		\end{minipage}
		\caption{Boxplots of the values $D_i^s(\alpha,B)$ and $D_i^{ns}(\alpha,B)$ from the decision rule over 500 replications, for the leading 8 eigenvalues.}\label{fig: dr}
	\end{figure}

	\subsection{Simulation: proportions of under/over estimation}
	Table \ref{tab: under over} is a supplement to the Table \ref{tab:number} in the main paper, which compares the proportions of under/over estimation of $r$ by different approaches, over 500 replications. The conclusions are the same as those in the main paper.
	\begin{table}
		\addtolength{\tabcolsep}{0pt}
		\caption{The proportions of under estimation (out of the bracket) and over  estimation (in the bracket) for  the factor number by different approaches, over 500 replications. \label{tab: under over}}
		\renewcommand{\arraystretch}{1}
		\scalebox{0.56}{ 
			\begin{tabular}{lllllllllllllllll}
				\toprule[1.2pt]
				$\vartheta$&$\rho$&$a$&$n=p$&$\hat r_{SMD}$&$\hat r_{SSD}$&$\hat r_{ETMD}$&$\hat r_{IC}$&$\hat r_{ABC}$&$\hat r_{ER}$&$\hat r_{TRAP}$&$\hat r_{ON}$&$\hat r_{ED}$&$\hat r_{ETC}$&$\hat r_{ETZ}$&$\hat r_{DDPA_+}$\\\midrule[1.2pt]
				0&0&0&100&0(0)&0(0)&0(0)&0(0)&0(0)&0(0)&0(0)&0(0.012)&0(0.036)&0(0.002)&0(0.012)&0(0)\\
				0&0&0&200&0(0)&0(0)&0(0)&0(0)&0(0)&0(0)&0(0)&0(0.012)&0(0.048)&0(0)&0(0)&0(0)\\
				0&0&0&300&0(0)&0(0)&0(0)&0(0)&0(0)&0(0)&0(0)&0(0.01)&0(0.024)&0(0)&0(0.006)&0(0)\\\hline
				0&3&0&100&0(0.29)&0(0.266)&0(0)&0(0)&0(0.218)&0(0)&0(0)&0(0.03)&0(0.976)&0(0.778)&0(0.98)&0(0.95)\\
				0&3&0&200&0(0.002)&0(0.002)&0(0)&0(0)&0(0.15)&0(0)&0(0)&0(0.056)&0(0.996)&0(0.892)&0(1)&0(0.988)\\
				0&3&0&300&0(0)&0(0)&0(0)&0(0)&0(0.124)&0(0)&0(0)&0(0.098)&0(1)&0(0.966)&0(1)&0(0.998)\\\hline
				1&0&0&100&0(0)&0(0)&0(0)&0(0)&0(0)&0(0)&0.242(0)&0(0.012)&0(0.01)&0(0)&0(0)&0.014(0.002)\\
				1&0&0&200&0(0)&0(0)&0(0)&0(0)&0(0)&0(0)&0(0)&0(0.01)&0(0.014)&0(0.002)&0(0)&0.008(0)\\
				1&0&0&300&0(0)&0(0)&0(0)&0(0)&0(0)&0(0)&0(0)&0(0.018)&0(0.014)&0(0)&0(0)&0.002(0)\\\hline
				1&0&0.25&100&0(0)&0(0)&0.008(0)&0.984(0)&0(0)&1(0)&1(0)&0.238(0.012)&0(0.01)&0(0)&0(0)&0.004(0)\\
				1&0&0.25&200&0(0)&0(0)&0(0)&0.976(0)&0(0)&1(0)&1(0)&0.002(0.01)&0(0.03)&0(0)&0(0.002)&0.004(0)\\
				1&0&0.25&300&0(0)&0(0)&0(0)&0.148(0)&0(0)&1(0)&1(0)&0(0.012)&0(0.018)&0(0)&0(0)&0(0)\\\hline
				1&3&0&100&0(0.538)&0(0.454)&0(0)&0(0)&0(0.224)&0(0)&0.274(0)&0(0.044)&0(0.95)&0(0.786)&0(0.884)&0.002(0.936)\\
				1&3&0&200&0(0.002)&0(0)&0(0)&0(0)&0(0.124)&0(0)&0(0)&0(0.062)&0(0.998)&0(0.914)&0(0.996)&0.006(0.982)\\
				1&3&0&300&0(0)&0(0)&0(0)&0(0)&0(0.086)&0(0)&0(0)&0(0.108)&0(1)&0(0.956)&0(1)&0(0.994)\\\hline
				1&3&0.25&100&0(0.6)&0(0.506)&0.002(0)&0.992(0)&0.002(0.264)&1(0)&1(0)&0.92(0.026)&0(0.946)&0(0.794)&0(0.892)&0.004(0.946)\\
				1&3&0.25&200&0(0.006)&0(0.002)&0(0)&0.976(0)&0(0.14)&1(0)&1(0)&0.776(0.072)&0(0.996)&0(0.88)&0(0.998)&0.004(0.984)\\
				1&3&0.25&300&0(0)&0(0)&0(0)&0.152(0)&0(0.112)&1(0)&1(0)&0.332(0.092)&0(1)&0(0.952)&0(1)&0(1)\\
				\bottomrule[1.2pt]
		\end{tabular}}
	\end{table}

	\subsection{Simulation: verifying theorems}\label{sec: verify}
	Here we verify the major theoretical results in Theorem \ref{thm: conditional on sample}, Corollary \ref{cor: bootstrap bias} and Theorem \ref{thm:non-spike} using simulating data. We start with Theorem \ref{thm: conditional on sample} and set $\vartheta=1$, $n=p=400$, $\beta_f=0.2$, $a=0.4$, $\rho=0$. Therefore, the two leading factors are strong and satisfy the condition for (\ref{conditional distribution})  while the  third factor is pretty weak satisfying the condition for (\ref{conditional power}). 
	
	Similarly to (\ref{Dis}), given the sample matrix $\Xb$ and a constant $s$, the tail probabilities in Theorem \ref{thm: conditional on sample} can be approximated by
	\begin{equation}\label{P hat}
		\hat{\mathbb{P}}_i^*(s):=\hat{\mathbb{P}}^*\bigg(\frac{\hat\lambda_i/\tilde\lambda_i-1}{\tilde\sigma}\le s\bigg)  = \frac{1}{B}\sum_{b=1}^BI\bigg(\frac{\hat\lambda_i^b/\tilde\lambda_i-1}{\tilde\sigma}\le s\bigg),
	\end{equation}
	by repeating the bootstrap procedure $B$ times for some large $B$. In the simulation, we let $B=400$.  We report the averaged $\hat {\mathbb{P}}_i^*(s)$ over 500 replications for $-2\le s\le 2$ and $i=1,2,3$ in Figure \ref{fig: verify} under the two bootstrap schemes, and compare the curves with the CDF of standard normal variable. It's seen from Figure \ref{fig: verify} (a) and (b) that $\hat {\mathbb{P}}_i^*(s)$ is very close to the CDF of standard normal distribution when $i=1,2$. In other words, the first and second largest eigenvalues after bootstrap are asymptotically Gaussian  after proper scaling and centralization, which verifies (\ref{conditional distribution}). However, in Figure \ref{fig: verify}(c),  $\hat {\mathbb{P}}_3^*(s)$ is always close to 0, which verifies (\ref{conditional power}) because the third common factor is very weak when $a=0.4$.

	\begin{figure}[h]
		\centering
		\begin{minipage}{0.31\textwidth}
			\centering
			\includegraphics[width=4.8cm,height=4.8cm]{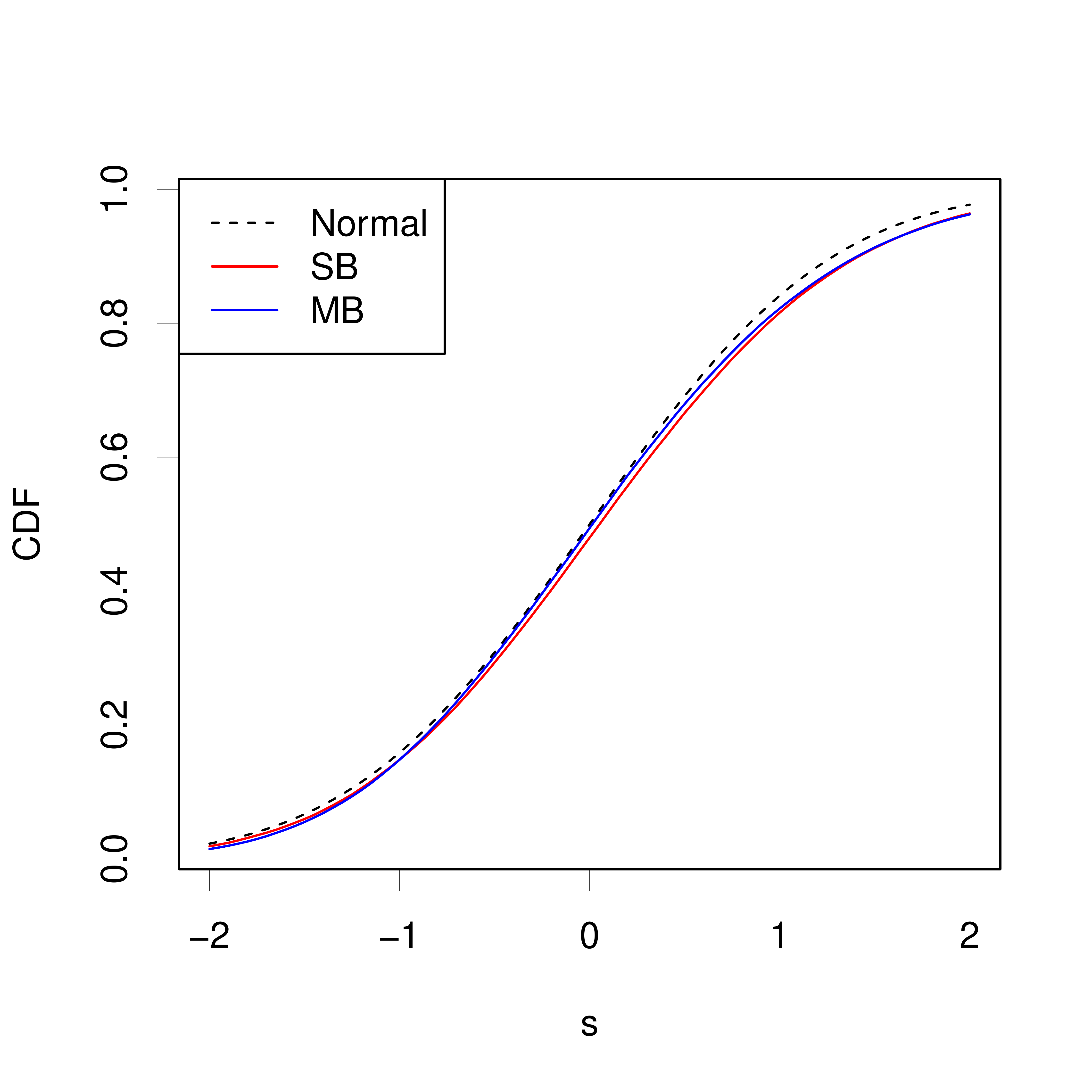}
			(a). $i=1$
		\end{minipage}
		\begin{minipage}{0.31\textwidth}
			\centering
			\includegraphics[width=4.8cm,height=4.8cm]{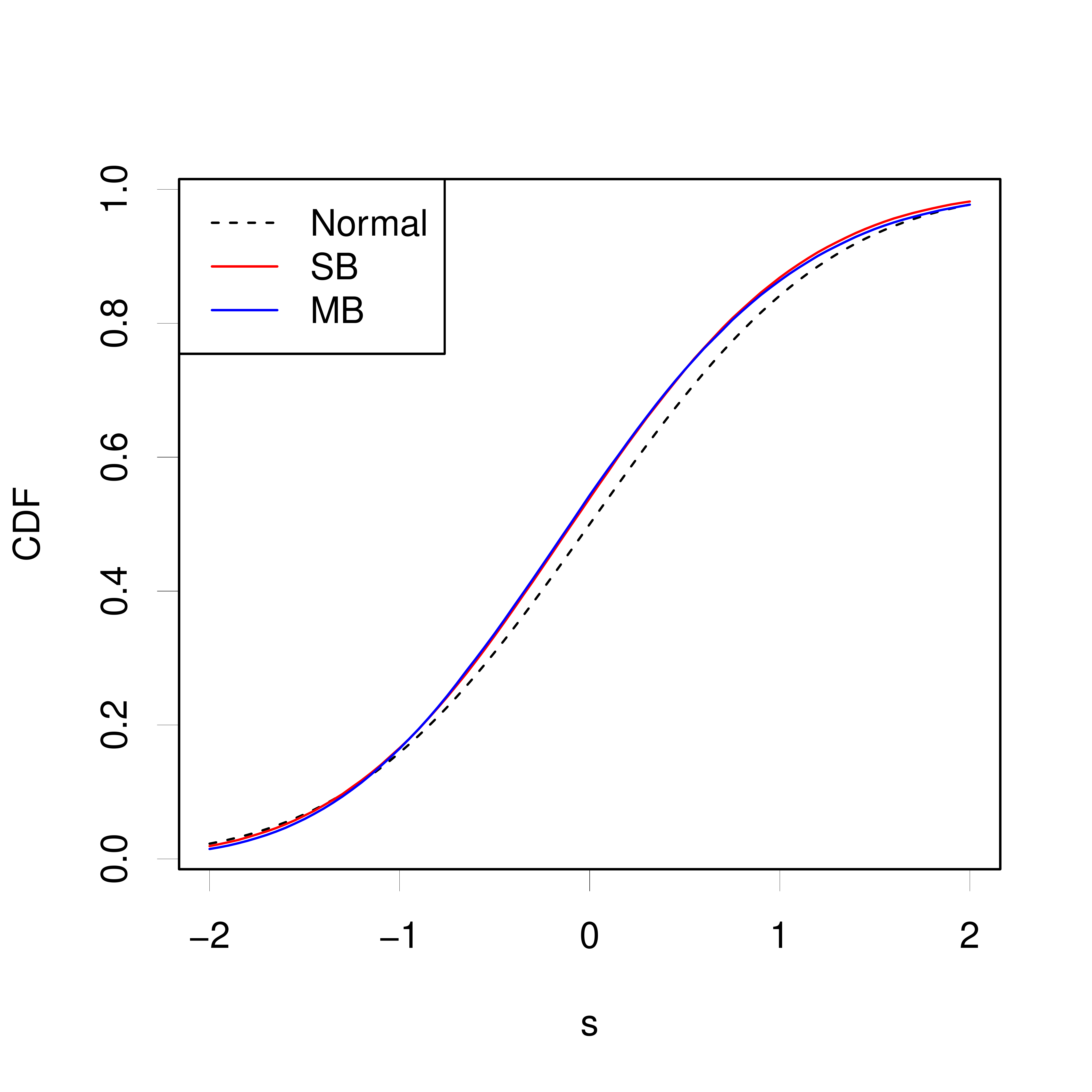}
			(b). $i=2$
		\end{minipage}
		\begin{minipage}{0.31\textwidth}
			\centering
			\includegraphics[width=4.8cm,height=4.8cm]{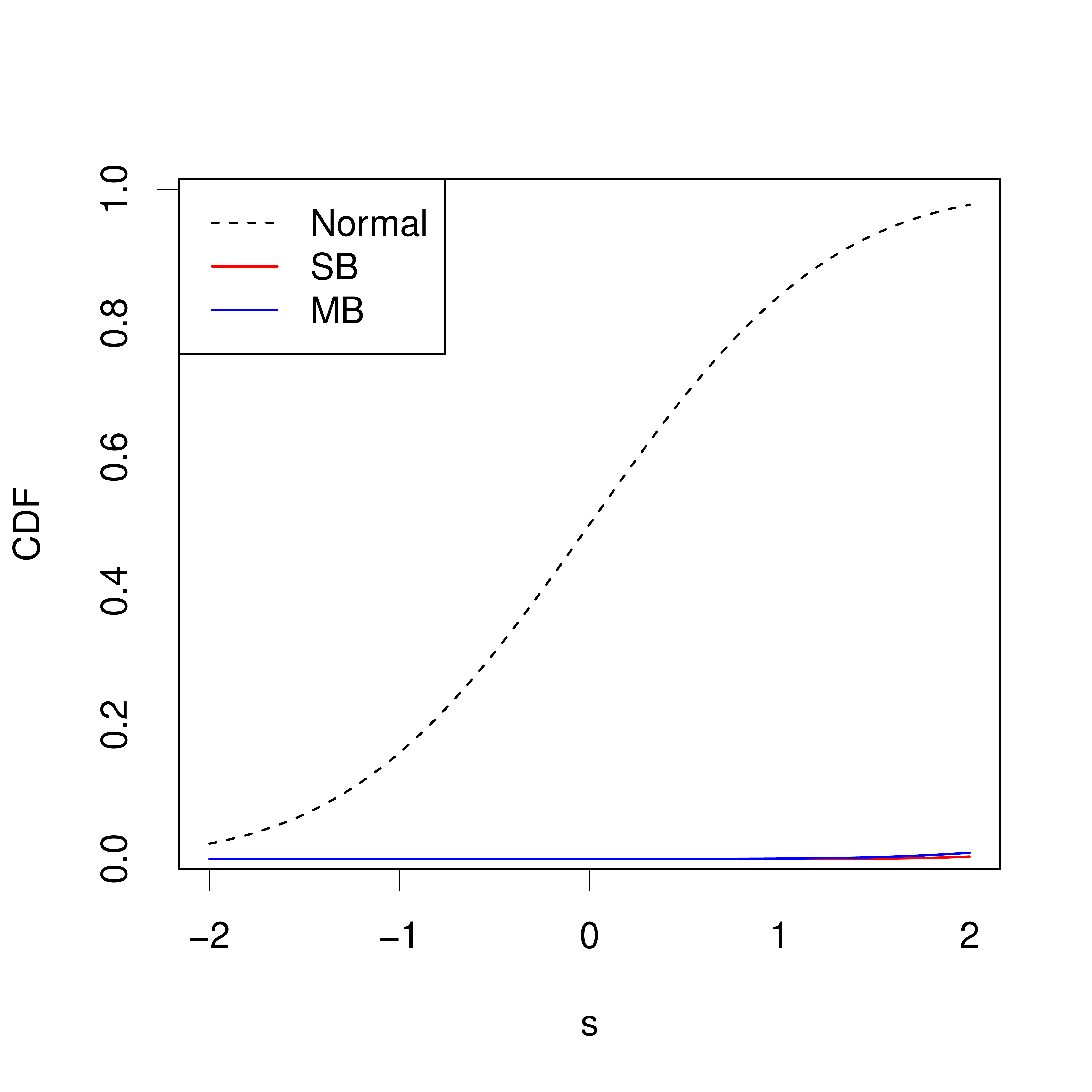}
			(c). $i=3$
		\end{minipage}
		\caption{Averaged $\hat{\mathbb{P}}_i^*(s)$ over 500 replications for verifying  Theorem \ref{thm: conditional on sample}. ``SB'' stands for standard bootstrap while ``MB'' stands for multiplier bootstrap. The dashed line is the CDF of standard normal distribution. }\label{fig: verify}
	\end{figure}
	
	Next, we verify bootstrap bias in Corollary \ref{cor: bootstrap bias}. To simplify the calculation of population covariance matrix, we let $\beta_f=0$. We will focus on the third eigenvalue to see how the two tail probabilities change when $a$ grows. The other data generating parameters are set the same as in Figure \ref{fig: verify}. Similarly to (\ref{P hat}) and slightly abusing the notation, we  let
	\[
	\hat{\mathbb{P}}_3^*(s):= \frac{1}{B}\sum_{b=1}^BI\bigg(\sqrt{n}\times \lambda_3^{-1}(\hat\lambda_3^b-\tilde\lambda_3)\le s\bigg),
	\]
	and report the averaged $\hat{\mathbb{P}}_3^*(s)$ over 500 replications in Figure \ref{fig: bias} for different values of $s$ and $a$ under the two bootstrap schemes. We compare the results with the benchmark tail probability in (\ref{bias tilde}), which is approximated by the frequency of the event $\{\sqrt{n}\times \lambda_3^{-1}(\tilde\lambda_3-\lambda_3)\le s\}$ happening over the 500 replications.
	Figure \ref{fig: bias}(a) indicates that when $a=0$, the standard bootstrap is roughly unbiased because the tail probability curve is very close to that without bootstrap. However, the tail probability curve of the multiplier bootstrap always has a small bias, as expected. When $a$ grows to $0.25$ so that the third factor becomes weaker, the tail probabilities start to deviate from each other. As claimed in (\ref{difference tail}),  $\hat{\mathbb{P}}_3^*(s)$ under the  standard bootstrap tends to be slightly larger than the benchmark tail probability under such cases, especially when $s$ is close to 0 so that the bias terms $\sqrt{n/\xi_3}\times (n\lambda_3)^{-1}\text{tr}\bLambda_2$ on the RHS of (\ref{bias tilde}) can dominate in finite samples. When $a$ further increases, all the tail probabilities tend to 0 in Figure \ref{fig: bias}.

	\begin{figure}[h]
		\centering
		\begin{minipage}{0.31\textwidth}
			\centering
			\includegraphics[width=4.8cm,height=4.8cm]{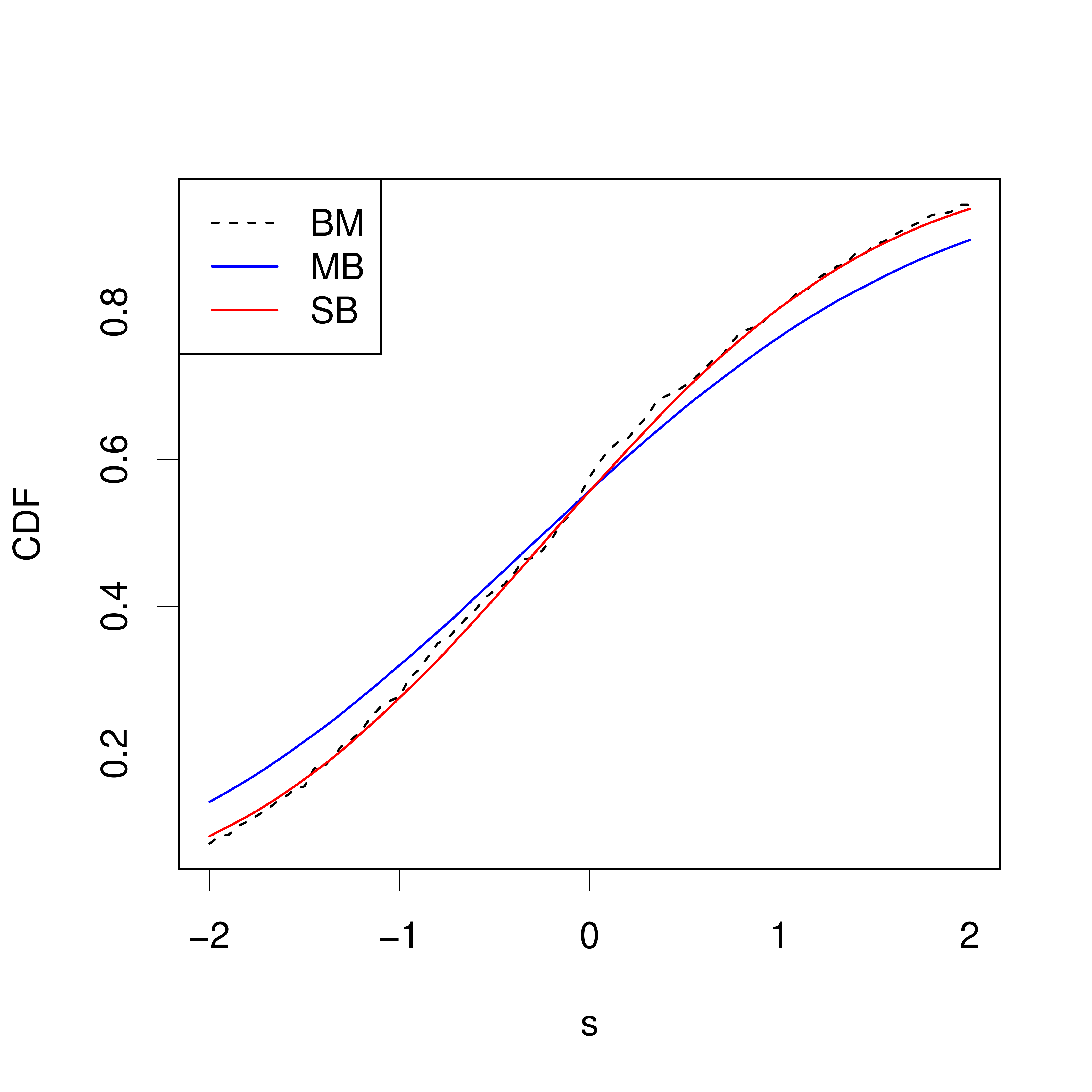}
			(a). $a=0$
		\end{minipage}
		\begin{minipage}{0.31\textwidth}
			\centering
			\includegraphics[width=4.8cm,height=4.8cm]{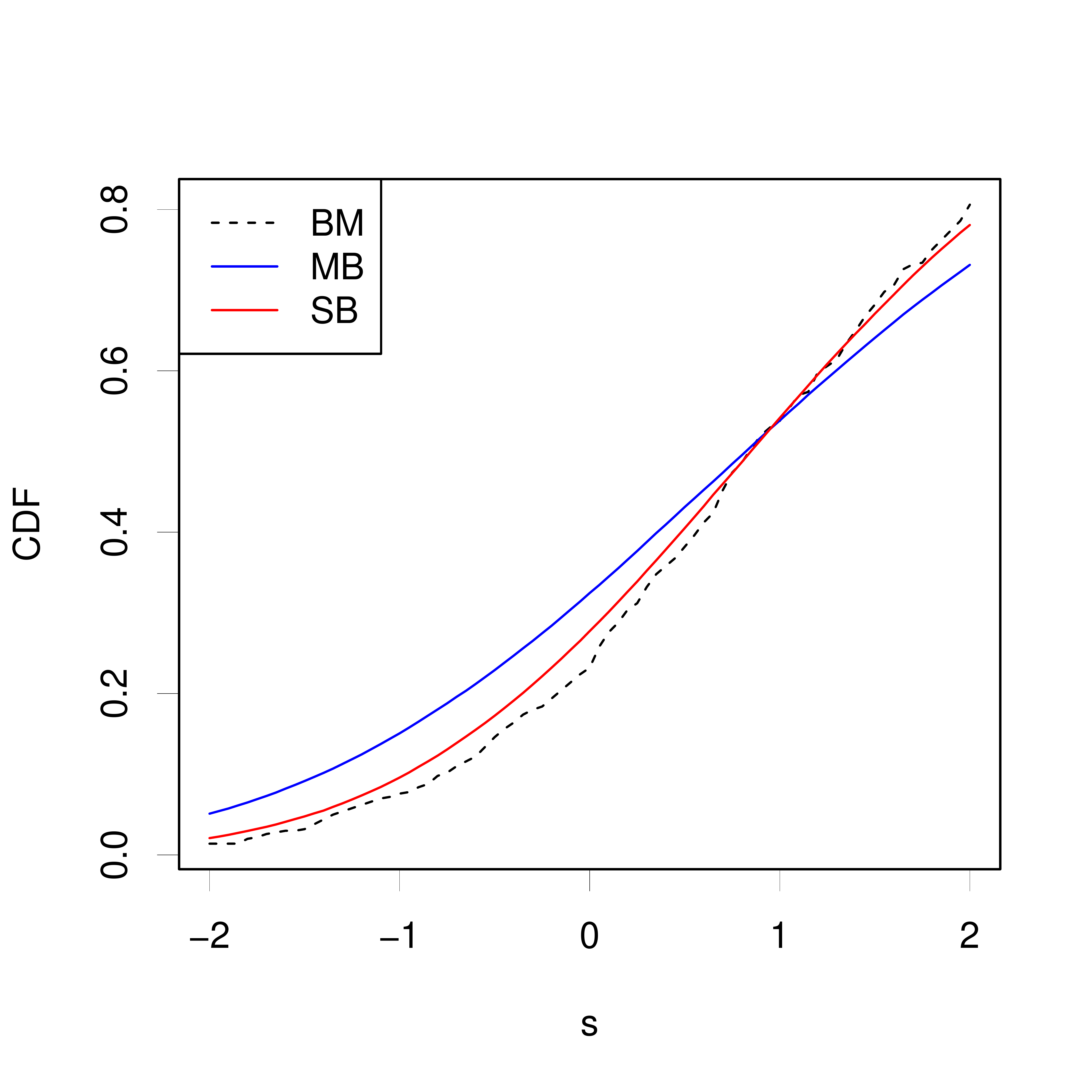}
			(b). $a=0.25$
		\end{minipage}
		\begin{minipage}{0.31\textwidth}
			\centering
			\includegraphics[width=4.8cm,height=4.8cm]{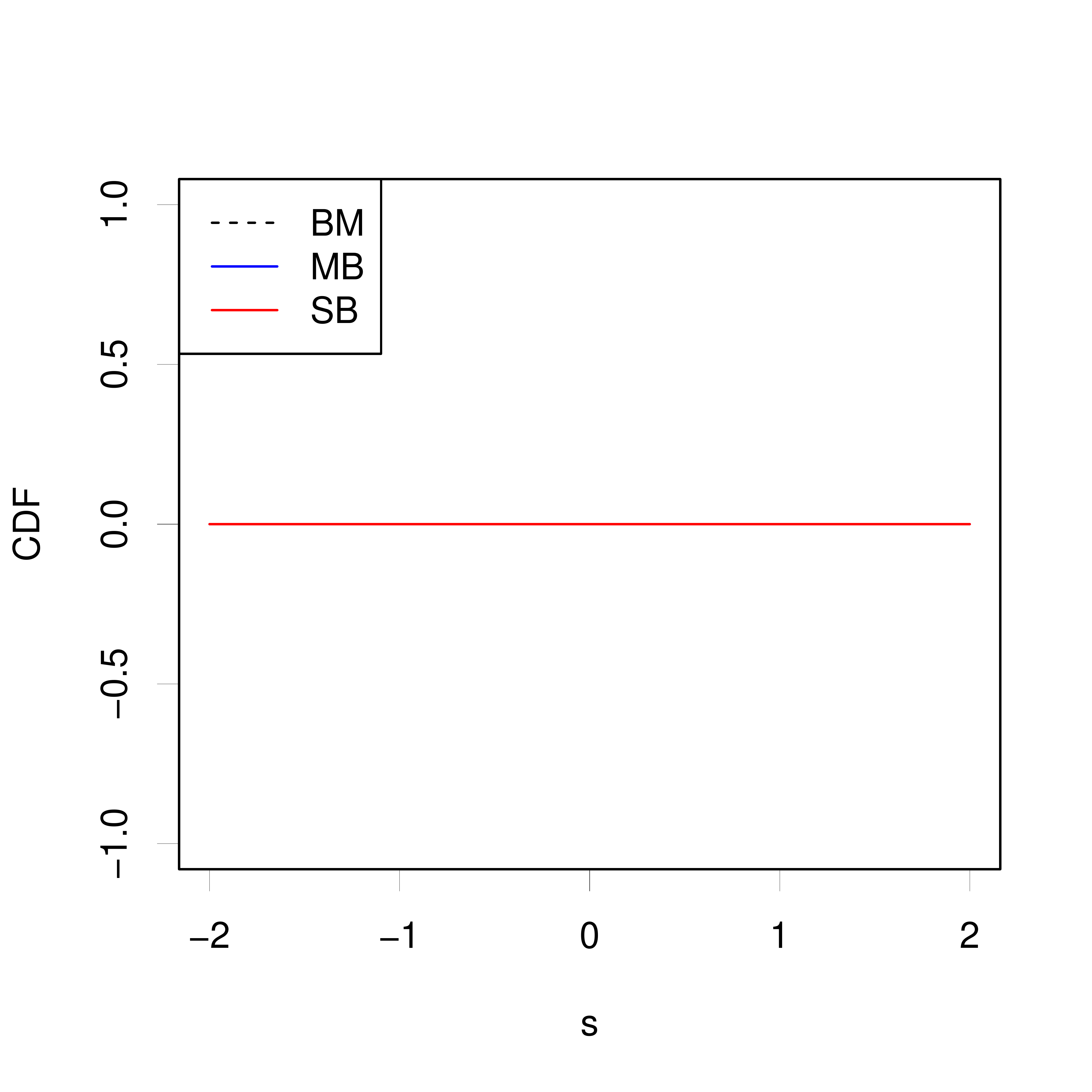}
			(c). $a=0.4$
		\end{minipage}
		\caption{The averaged tail probabilities for the third sample eigenvalue after and before bootstrap over 500 replications for verifying  Corollary \ref{cor: bootstrap bias}, as the third common factor gradually becomes weak. ``SB'' and ``MB'' stand for the standard and multiplier bootstrap respectively,  while ``BM' stands for the benchmark tail probability in (\ref{bias tilde}).   }\label{fig: bias}
	\end{figure}

	\begin{figure}[h]
		\centering
		\begin{minipage}{0.31\textwidth}
			\centering
			\includegraphics[width=4.8cm,height=4.8cm]{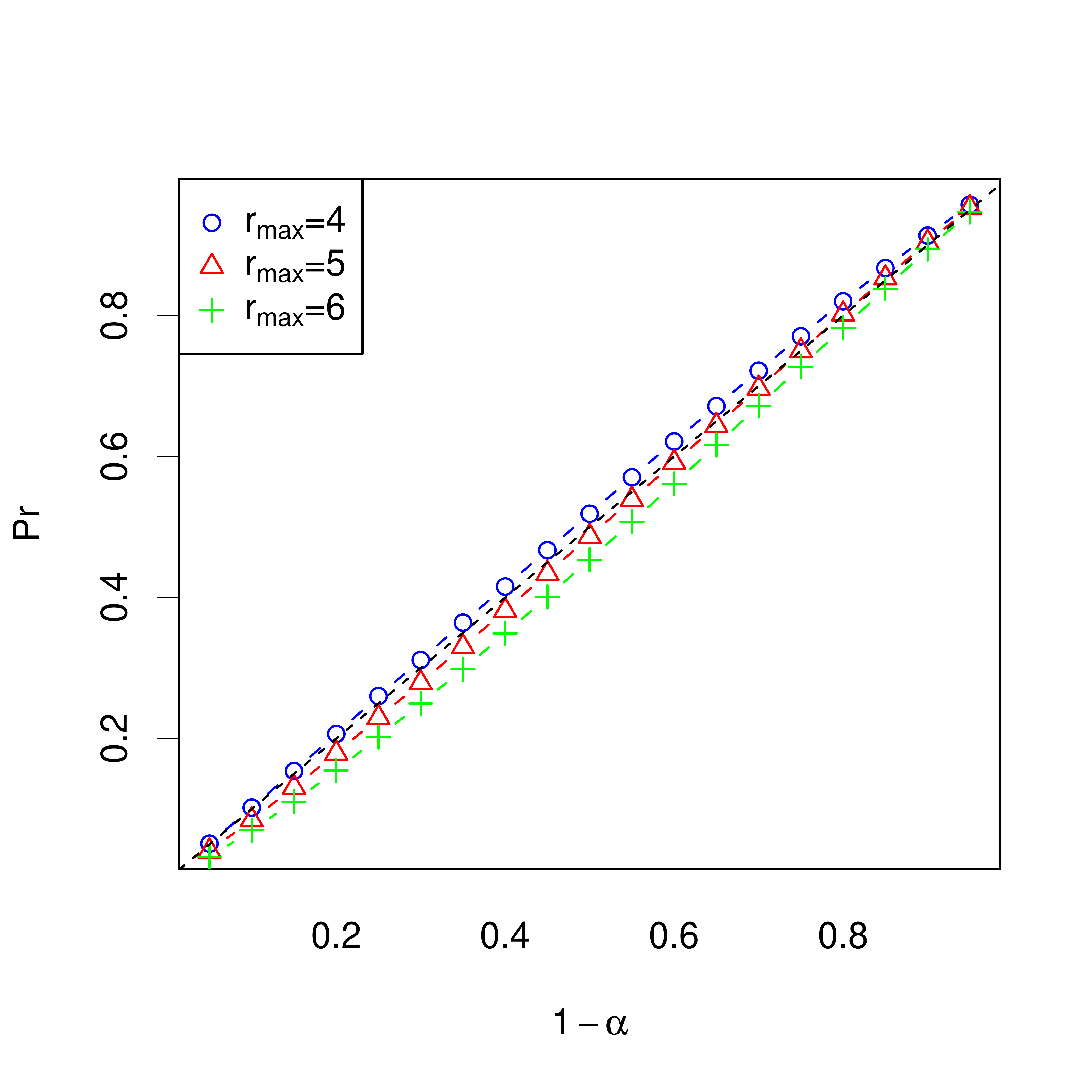}
			(a). $a=0$
		\end{minipage}
		\begin{minipage}{0.31\textwidth}
			\centering
			\includegraphics[width=4.8cm,height=4.8cm]{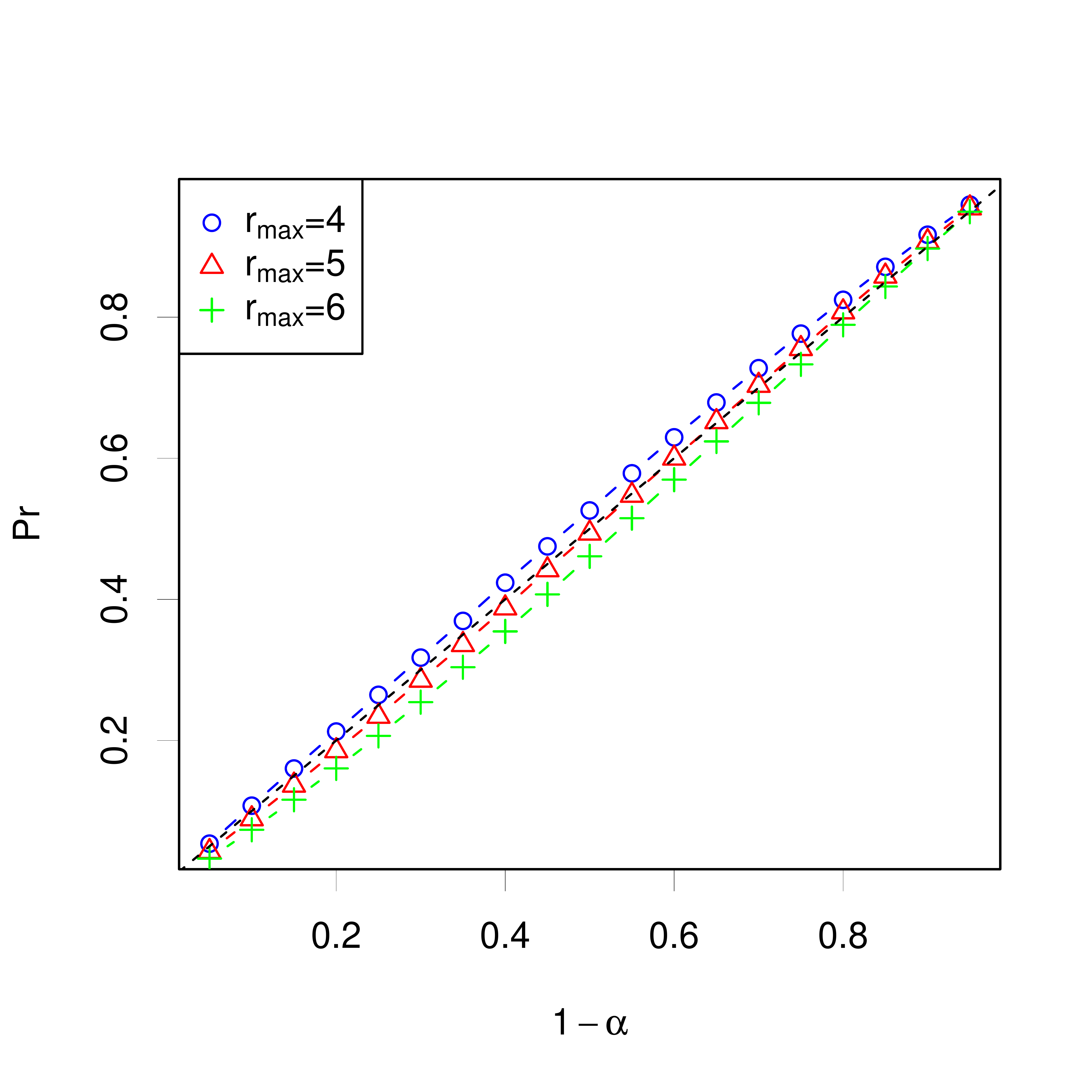}
			(b). $a=0.25$
		\end{minipage}
		\begin{minipage}{0.31\textwidth}
			\centering
			\includegraphics[width=4.8cm,height=4.8cm]{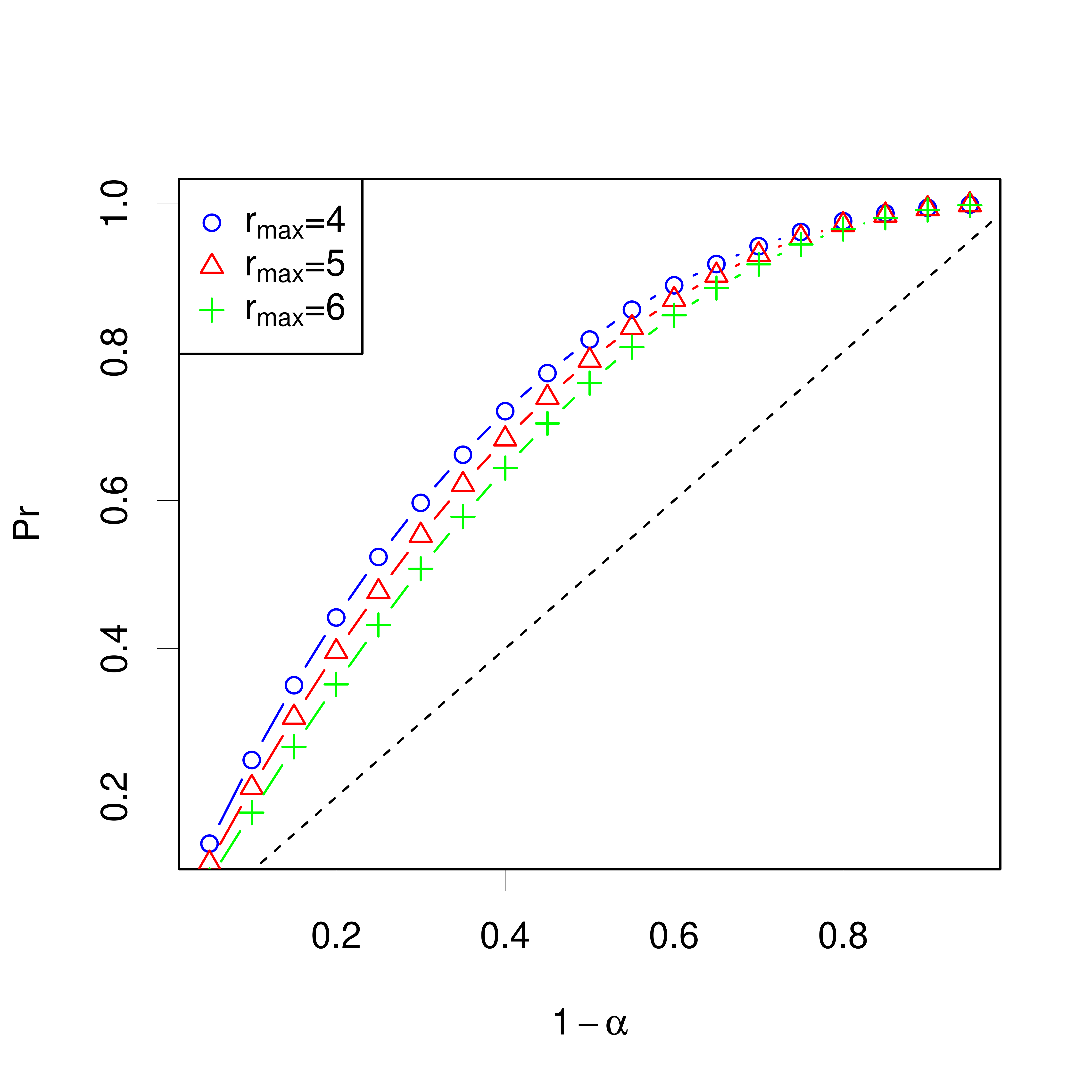}
			(c). $a=0.4$
		\end{minipage}
		\caption{The averaged $\hat{\mathbb{P}}_{r+1}^*(\alpha)$ for the largest non-spiked eigenvalue over 500 replications with different values of $\alpha$ and $r_{\max}$. }\label{fig: non-spike}
	\end{figure}

	Lastly, we verify Theorem \ref{thm:non-spike}. To be more consistent with the proposed approach, we will directly compare the empirical distributions of $\hat\lambda_{r+1}$ and $\{\hat\varphi_1^j\}_{j=1}^R$ from Algorithm \ref{alg1}. To be more specific, given the sample matrix $\Xb$, we run Algorithm \ref{alg1} and obtain a series of sample quantiles $\hat c_{1-\alpha}$ from $\{\hat\varphi_1^j\}_{j=1}^R$ with different values of $\alpha$ and $R=400$. Then, similarly to (\ref{P hat}), we further repeat the bootstrap procedure $B$ times and let 
	\[
	\hat{\mathbb{P}}_{r+1}^*(\alpha)=B^{-1}\sum_{b=1}^BI(\hat\lambda_{r+1}^b<\hat c_{1-\alpha}). 
	\]
	Figure \ref{fig: non-spike} plots the averaged $\hat{\mathbb{P}}_{r+1}^*(\alpha)$ over 500 replications with different values of $\alpha$ and $r_{\max}$. The data generating parameters are the same as in Figure \ref{fig: bias}, except that we let $\beta_f=0.2$. By Figure \ref{fig: non-spike} (a) and (b), the empirical probabilities $\hat{\mathbb{P}}_{r+1}^*(\alpha)$ are very close to $1-\alpha$ as long as $a$ is not too large, indicating that the Algorithm \ref{alg1} can accurately approximate the distribution of $\hat\lambda_{r+1}$. The accuracy decreases when we use larger $r_{\max}$, as expected. When $a=0.4$ so that the third common factor is very weak, Algorithm \ref{alg1} loses accuracy. This is because we are in finite samples and the condition $(n\lambda_i)^{-1}p\log n=o(1)$ will not hold anymore under this case.
	
	\subsection{Real example 2: macroeconomic indices}
	In the second real example, we analyze a macroeconomic data set, namely the FRED-MD data set, which was introduced by \cite{McCracken2015FRED}. It's an open resource from \url{https://research.stlouisfed.org/econ/mccracken/fred-databases/}, containing monthly series of 127 macroeconomic variables since January 1959. This data set is generally regarded as the standard case of stronger factor structures among all common empirical applications in the related literature. 
	We refer to the original paper for more details. Following the code in \cite{McCracken2015FRED},  we transform the data to stationary series, drop 5 variables with largest missing rates, and remove all the outliers which deviate from the sample medians by more than 10 interquartile ranges. We focus on the period from January 1961 to December 2021, covering 732 months. The series are standardized while the missing entries are imputed by linear interpolation. Eventually, a data matrix $\Xb_{p\times n}$ is obtained, with $p=122$ and $n=732$.
	
	Similarly to the financial example in the main paper, we plot the sample eigenvalues in Figure \ref{fig:fred}(a). 
	It's more visible that 7 eigenvalues deviate from the bulk, and the gaps of the leading 7 eigenvalues are not as significant as those in the financial example. The estimated numbers of factors by different methods and the computational costs are reported in Table \ref{tab: fred}, with the same tuning parameters as in Table \ref{tab: ff} except that $r_{\max}=12$. The proposed three methods and $\hat r_{ON}$, $\hat r_{DDPA_+}$ lead to an estimate of $\hat r=7$, which is the majority vote. 
	
	\begin{figure}[h]
		\centering
		\begin{minipage}{0.24\textwidth}
			\centering
			\includegraphics[width=4cm,height=4cm]{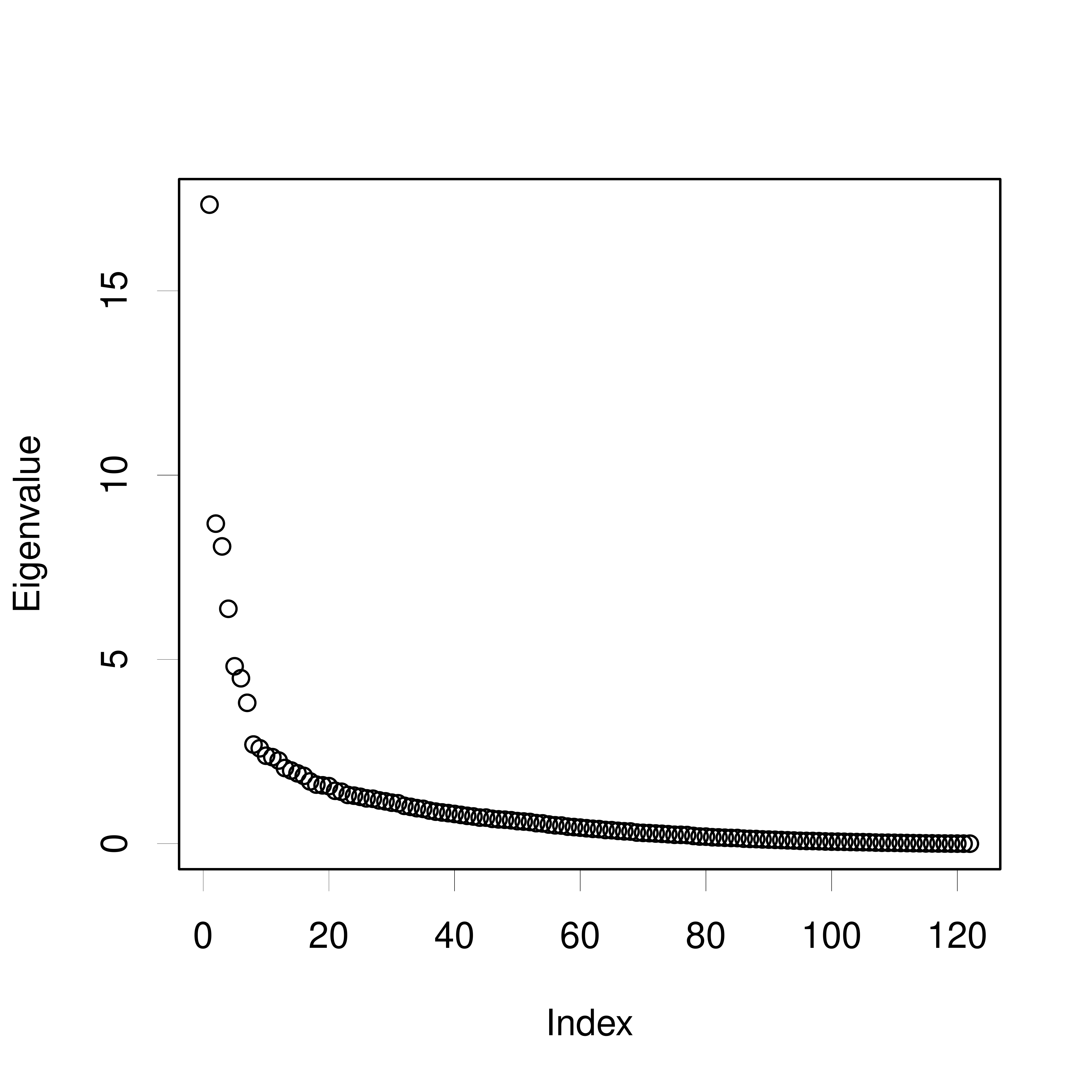}
			(a)
		\end{minipage}
		\begin{minipage}{0.24\textwidth}
			\centering
			\includegraphics[width=4cm,height=4cm]{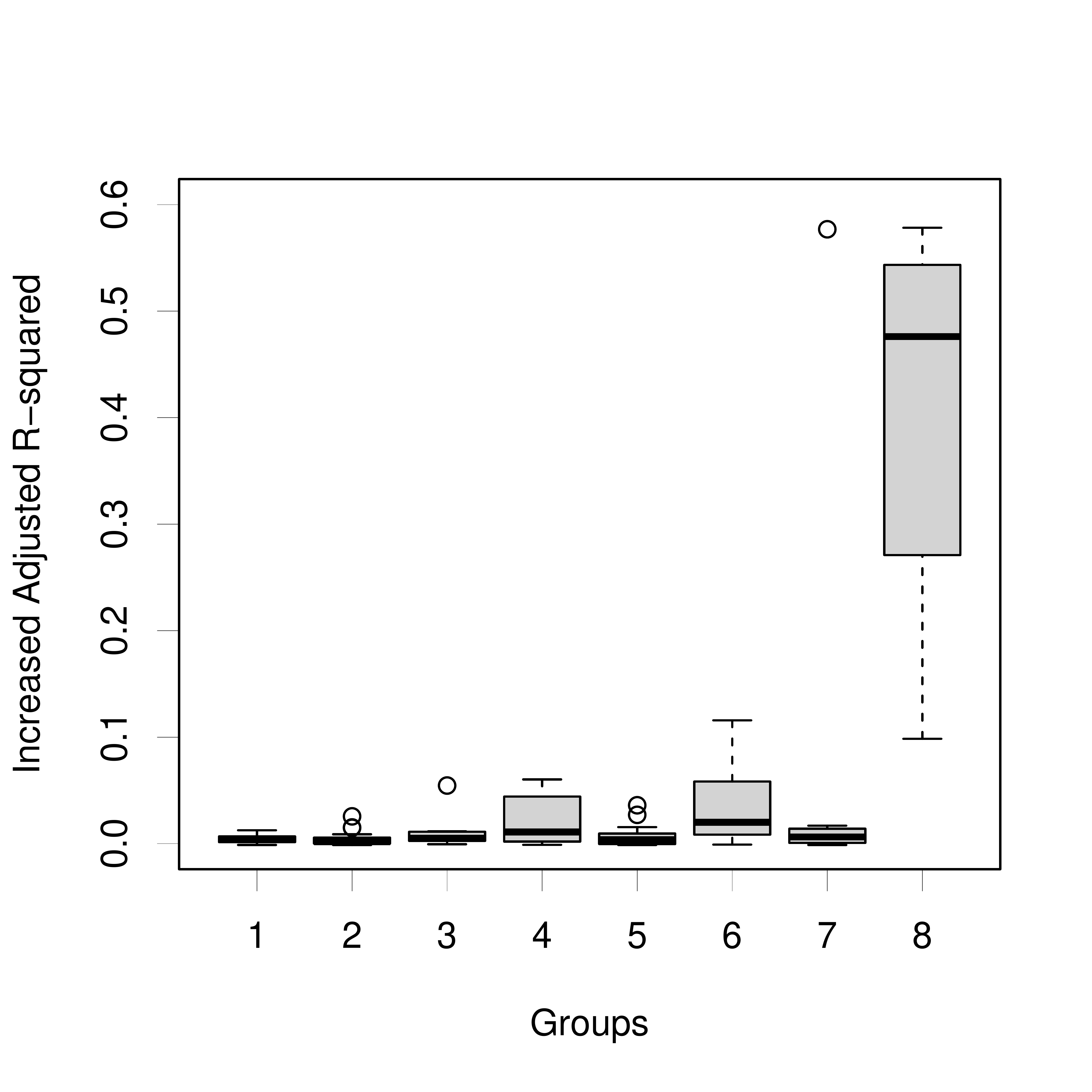}
			(b)
		\end{minipage}
		\begin{minipage}{0.24\textwidth}
			\centering
			\includegraphics[width=4cm,height=4cm]{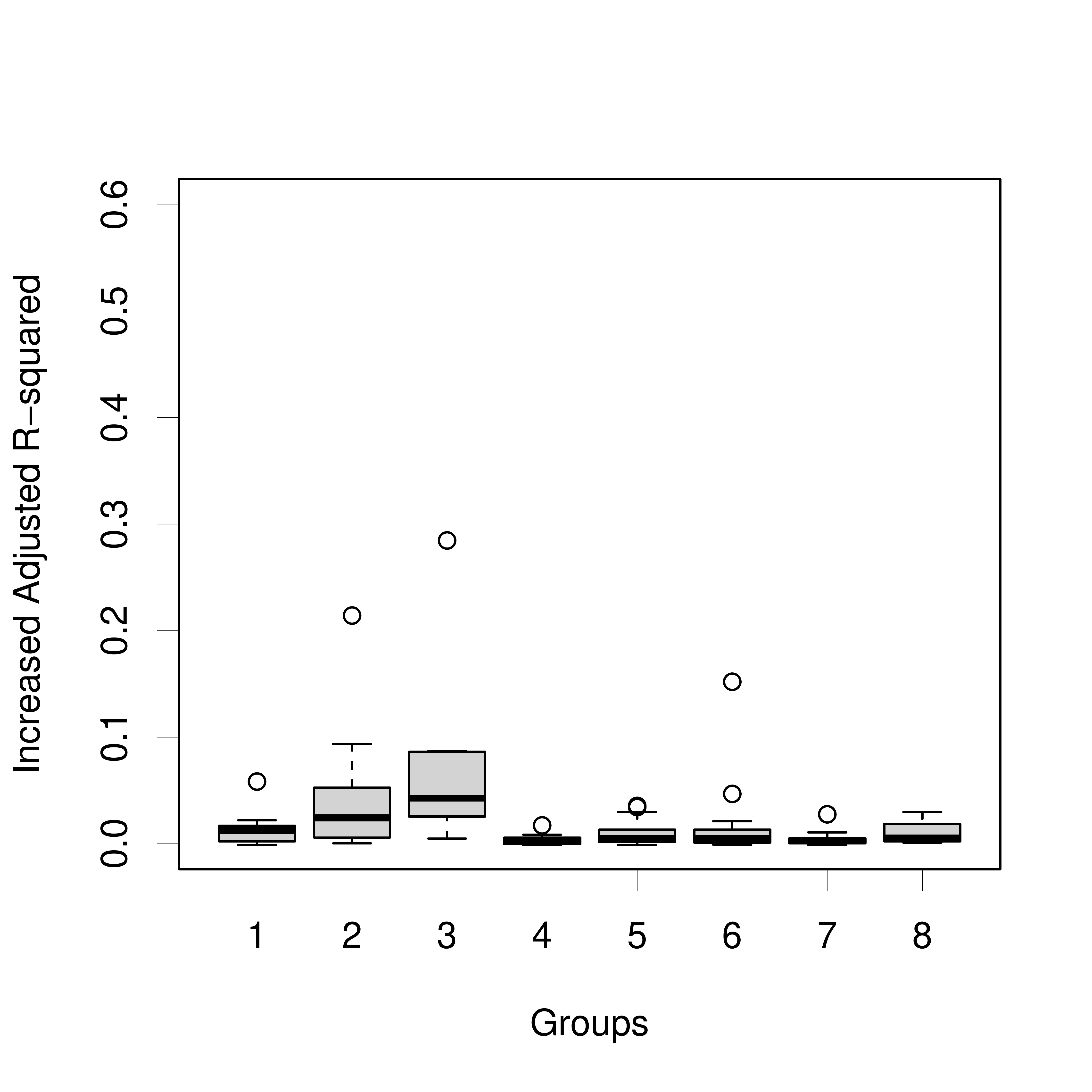}
			(c)
		\end{minipage}
		\begin{minipage}{0.24\textwidth}
			\centering
			\includegraphics[width=4cm,height=4cm]{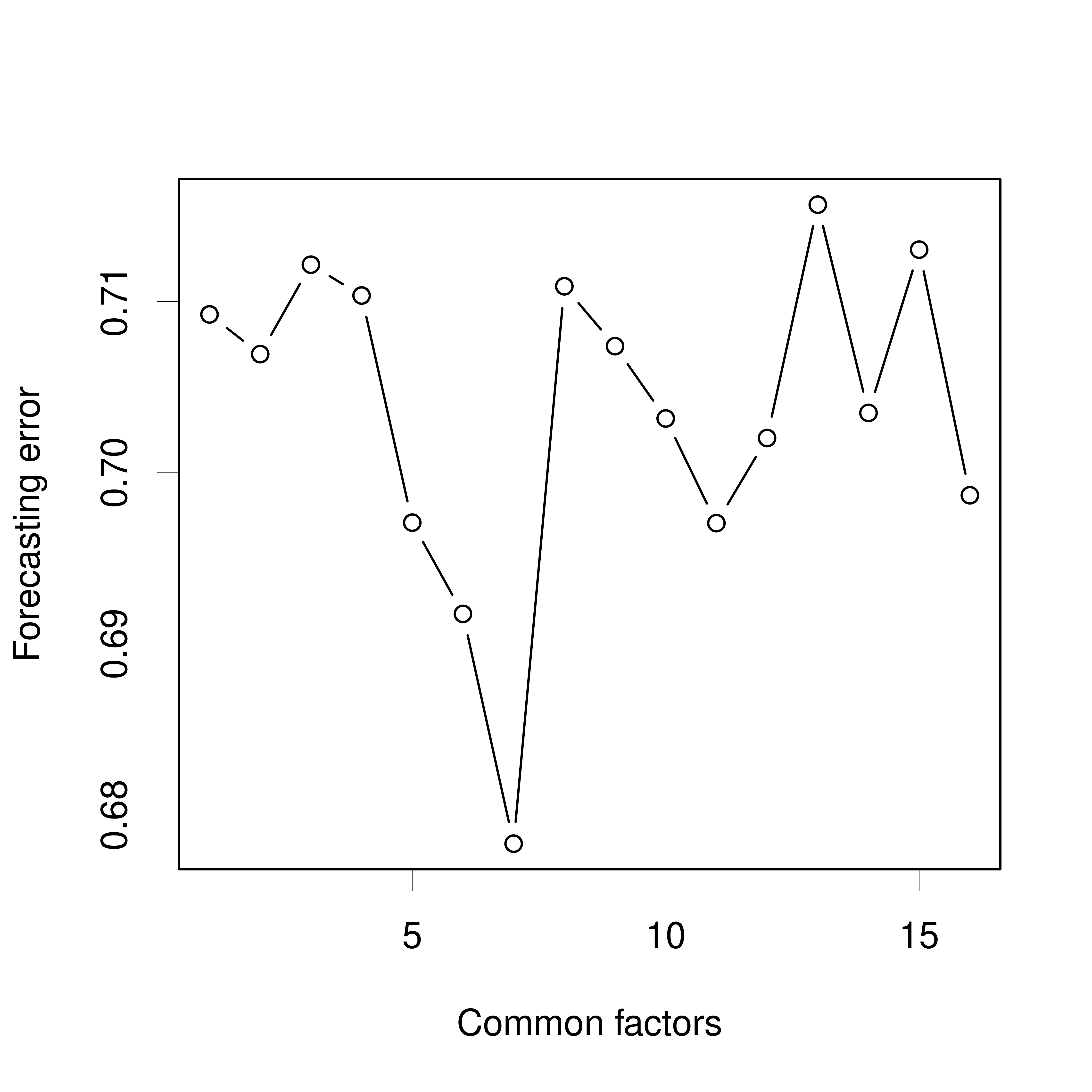}
			(d)
		\end{minipage}
		
		\caption{Figures for the FRED-MD data set: (a) eigenvalues of sample covariance matrix. (b) boxplot of the increased explained variation (adjusted R-squared) of the 122 macroeconomic in 8 groups when the 7th factor comes into system. (c) replication of (b) when the 8th factor comes into system (d) averaged diffusion-indexes forecasting error (RMSE) for the macroeconomic indices in group 8 when $\hat r$ grows.}\label{fig:fred}
	\end{figure}
	
	\begin{table}[htbt]
		\addtolength{\tabcolsep}{2pt}
		\caption{Estimated number of common factors and the computational cost (in seconds) for the FRED-MD data set by different methods.\label{tab: fred}}
		\renewcommand{\arraystretch}{0.6}
		\scalebox{0.75}{ 
			\begin{tabular}{lllllllllllllllll}
				\toprule[1.2pt]
				&$\hat r_{SMD}$&$\hat r_{SSD}$&$\hat r_{ETMD}$&$\hat r_{IC}$&$\hat r_{ABC}$&$\hat r_{ER}$&$\hat r_{TRAP}$&$\hat r_{ON}$&$\hat r_{ED}$&$\hat r_{ETC}$&$\hat r_{ETZ}$&$\hat r_{DDPA}$\\\midrule[1.2pt]
				$\hat r$&7&7&7&1&1&1&0&7&5&16&2&7\\
				Cost (s)&1.597&1.580&6.056&0.034&0.492&0.030&0.035&1.425&0.031&0.051&16.858&0.053
				
				\\\bottomrule[1.2pt]
		\end{tabular}}
	\end{table}
	
	To show why $\hat r=7$ is a reasonable result, we calculate the explanatory power of the factors, in terms of the incremental of adjusted R-squared, similarly to the financial example. The 122 variables
	are categorized into 8 groups in \cite{McCracken2015FRED} according to economic implication. Figure \ref{fig:fred} (b) and (c) are boxplots of  the explanatory power of the 7th and 8th factors, respectively to the variables in 8 groups. It's seen that the 7th common factor contributes significantly to explaining the variation of macroeconomic variables in group 8, while the gain from the 8th factor is  minor to all the 8 groups.
	
	Lastly, we investigate how new factors contribute to forecasting  macroeconomic indices. Motivated by the diffusion-indexes forecasting in \cite{stock2002macroeconomic}, we forecast $x_i$ based on $x_{i,t+1}=\alpha_i+\bbeta_i^\top\bbf_{i,t}+\gamma_i(L)x_{i,t}+\epsilon_{i,t}$ for $1\le i\le p$, where $\{\bbf_{i,t}\}$ is the factor
	process estimated from panel data excluding $x_i$ given $\hat r$, $\gamma_i(L)$ is the lag polynomials, and $\epsilon_{i,t}$ is the noise. For simplicity, we let $\gamma_i(L)=\gamma_i$ in this experiment. Motivated by Figure \ref{fig:fred}(b), we are more interested in how new factors contribute in forecasting the variables in group 8. Then, for each macroeconomic series in group 8, we use 70\% samples to train the model and estimate $\alpha_i,\bbeta_i,\gamma_i$, and calculate the forecasting error based on the remaining samples.  Figure \ref{fig:fred}(d) shows the averaged forecasting error (RMSE) for the macroeconomic indices in group 8 when $\hat r$ grows. It's clearly seen that the forecasting error is minimized at $\hat r=7$. In conclusion, we believe that $\hat r=7$ is reasonable for this data set, coinciding with the scree-plot in Figure \ref{fig:fred}(a).
	In fact, in the literature, the FRED-MD data set is generally regarded  to contain 6 to 8 common factors, which is also consistent with the economic group structure.

	\subsection{Real data: robustness to tuning parameters}
	Tables \ref{tab: ff tuning} and \ref{tab: fred tuning} report the results of $\hat r_{SMD}$, $\hat r_{SSD}$ and $\hat r_{ETMD}$ for the two real examples when different tuning parameters are used.  Motivated by our simulation results, we are more interested in the tuning of $r_{\max}$ and $\alpha$ while fixing $B=200, R=400$. The outputs are quite stable, i.e., $\hat r =3$ for the financial data set and $\hat r=7$ for the macroeconomic data set in most settings.
	
	\begin{table}
		\addtolength{\tabcolsep}{6pt}
		\caption{Estimated number of common factors by $\hat r_{ETMD}$ for the Fama-French portfolio data set with different tuning parameters.\label{tab: ff tuning}}
		\renewcommand{\arraystretch}{0.6}
		\scalebox{1}{ 
			\begin{tabular}{llllllllllllll}
				\toprule[1.2pt]
				&\multicolumn{4}{c}{$\hat r_{SMD}$}&\multicolumn{4}{c}{$\hat r_{SSD}$}&\multicolumn{4}{c}{$\hat r_{ETMD}$}\\\cmidrule(lr){2-5}\cmidrule(lr){6-9}\cmidrule(lr){10-13}
				$r_{max}$&6&7&9&10&6&7&9&10&6&7&9&10\\\midrule[1.2pt]
				$\alpha=0.02$&4&3&4&4&4&3&4&4&3&3&3&3\\
				$\alpha=0.04$&3&3&3&3&3&3&3&3&3&3&3&3\\
				$\alpha=0.06$&3&3&3&3&3&3&3&3&3&3&3&3\\
				$\alpha=0.08$&3&3&3&3&3&3&3&3&4&4&4&3\\
				$\alpha=0.10$&3&3&3&3&3&3&3&3&4&4&4&4
				\\\bottomrule[1.2pt]
		\end{tabular}}
	\end{table}
	
	\begin{table}
		\addtolength{\tabcolsep}{5pt}
		\caption{Estimated number of common factors by $\hat r_{ETMD}$ for the FRED-MD data with different tuning parameters.\label{tab: fred tuning}}
		\renewcommand{\arraystretch}{0.6}
		\scalebox{1}{ 
			\begin{tabular}{llllllllllllll}
				\toprule[1.2pt]
				&\multicolumn{4}{c}{$\hat r_{SMD}$}&\multicolumn{4}{c}{$\hat r_{SSD}$}&\multicolumn{4}{c}{$\hat r_{ETMD}$}\\\cmidrule(lr){2-5}\cmidrule(lr){6-9}\cmidrule(lr){10-13}
				$r_{max}$&8&10&14&16&8&10&14&16&8&10&14&16\\\midrule[1.2pt]
				$\alpha=0.02$&7&7&16&7&7&7&16&7&4&4&4&4\\
				$\alpha=0.04$&7&7&7&7&7&7&7&7&7&7&4&7\\
				$\alpha=0.06$&7&7&7&7&7&7&7&7&7&7&7&7\\
				$\alpha=0.08$&7&7&7&7&7&7&7&7&7&7&7&7\\
				$\alpha=0.10$&7&7&7&7&7&7&7&7&7&7&7&7
				\\\bottomrule[1.2pt]
		\end{tabular}}
	\end{table}

	\section{Proof of results in Section \ref{sec: spiked}}\label{sec: supplement A}
	\subsection{Proof of Lemma \ref{lem: preliminary}: preliminary results on $\hat\lambda_i$}
	
	We start with the population eigenvalues $\lambda_i$. Recall that $\Ab=(\Lb,\bPsi)$. By Assumption \ref{c1} and Weyl's theorem, we have
	\[
	\frac{\lambda_i}{\lambda_i(\Lb^\top\Lb)}-1\le \frac{\|\bPsi\bPsi^\top\|}{\lambda_i(\Lb^\top\Lb)}=o(1),\quad1\le i\le r,
	\]
	and  $c\le \lambda_{[c p]}(\bPsi\bPsi^\top)\le \lambda_{[cp]}\le\lambda_{r+1}\le \lambda_1(\bPsi\bPsi^\top)\le c^{-1}$. Then, $\lambda_{r}/\lambda_{r+1}\ge 1+c$. Moreover, for $1\le i\le r-1$, 
	\[
	\frac{\lambda_i}{\lambda_{i+1}}=\frac{\lambda_i}{\lambda_i(\Lb^\top\Lb)}\frac{\lambda_i(\Lb^\top\Lb)}{\lambda_{i+1}(\Lb^\top\Lb)}\frac{\lambda_{i+1}(\Lb^\top\Lb)}{\lambda_{i+1}}=\frac{\lambda_{i}(\Lb^\top\Lb)}{\lambda_{i+1}(\Lb^\top\Lb)}[1+o_p(1)]\ge 1+c.
	\] 
	
	Now we consider $\hat\lambda_i$. Recall the decomposition
	\[
	\hat{\mathcal{S}}=n^{-1}\Wb^{1/2}\Zb^\top(\bGamma_1\bLambda_1\bGamma_1^\top+\bGamma_2\bLambda_2\bGamma_2^\top)\Zb\Wb^{1/2}:=\hat{\mathcal{S}}_1+\hat{\mathcal{S}}_2.
	\]
	We first show that  $\|\hat{\mathcal{S}}_2\|\le O_{p}[n^{-1}(n\vee p)\log n]$.
	Lemma \ref{lem: sample covariance} will indicate  that
	\[
	\|\hat{\mathcal{S}}_2\|\le C\|n^{-1}\Zb^\top\Zb\|\times \max_j w_j\le O_{p}\big(n^{-1}(n\vee p)\big)\times\max_j w_j.
	\]
	Then, it's sufficient to consider  $\max_jw_j$.  
	
	For the standard bootstrap, by Jensen's equality, for any $\alpha>0$
	\begin{equation}\label{mg}
		\begin{split}
			\exp(\alpha\mathbb{E}\max_jw_j)\le &\mathbb{E}\big(\exp(\alpha\max_jw_j)\big)
			\le\mathbb{E}\big(\sum_j\exp(\alpha w_j)\big)=n\mathbb{E}\big(\exp(\alpha w_1)\big).
		\end{split}
	\end{equation}
	Note that $w_1\sim Binomial(n,n^{-1})$, whose moment generating function is
	\begin{equation}\label{Ew1}
		\mathbb{E}\big(\exp(\alpha w_1)\big)=\big(1-n^{-1}+n^{-1}\exp(\alpha)\big)^n\le e^{e^\alpha-1}+c,
	\end{equation}
	for sufficiently large $n$. Therefore,
	\begin{equation}\label{s2}
		\mathbb{E}\max_jw_j\le \alpha^{-1}(\log n+e^{\alpha}-1+c)\le O(\log n)\Longrightarrow \|\hat{\mathcal{S}}_2\|\le  O_{p}\big(n^{-1}(n\vee p)\log n\big).
	\end{equation}
	For multiplier bootstrap, the proof is similar and omitted. 
	
	Therefore, by Weyl's theorem,
	\begin{equation}\label{hat_lambda_i}
		\bigg|\frac{\hat\lambda_i}{\lambda_i}-\frac{\lambda_i(\hat{\mathcal{S}}_1)}{\lambda_i}\bigg|\le \frac{\|\hat{\mathcal{S}}_2\|}{\lambda_i}=O_{p}\bigg(\frac{(n\vee p)\log n}{n\lambda_i}\bigg)\rightarrow 0, \quad 1\le i\le r,
	\end{equation}
	while $\hat\lambda_{r+1}\le  O_p[n^{-1}(n\vee p)\log n]$. It remains to consider $\lambda_i(\hat{\mathcal{S}}_1)$, or equivalently the $i$th largest eigenvalue of $n^{-1}\bLambda_1^{1/2}\bGamma_1^\top\Zb\Wb\Zb^\top\bGamma_1\bLambda_1^{1/2}$. Let $\hat\Kb_1(x)=x\bLambda_1^{-1}-n^{-1}\bGamma_1^\top\Zb\Wb\Zb^\top\bGamma_1$. By definition, $\det\{\hat\Kb_1[\lambda_i(\hat{\mathcal{S}}_1)]\}=0$.
	Note that
	\[
	n^{-1}\bGamma_1^\top\Zb\Wb\Zb^\top\bGamma_1-\Ib=n^{-1}\bGamma_1^\top\Zb(\Wb-\Ib)\Zb^\top\bGamma_1+n^{-1}\bGamma_1^\top\Zb\Zb^\top\bGamma_1-\Ib.
	\]
	By Lemma \ref{lem: two bounds} and the independence of $w_j$'s (or weak dependence under the standard bootstrap), one can verify that
	\begin{equation}\label{W-I}
		\|n^{-1}\bGamma_1^\top\Zb(\Wb-\Ib)\Zb^\top\bGamma_1\|=O_p(n^{-1/2}).
	\end{equation}
	On the other hand, Lemma \ref{lem: two bounds} will show that $\|n^{-1}\bGamma_1^\top\Zb\Zb^\top\bGamma_1-\Ib\|=O_p(n^{-1/2})$. Then, the matrix $	\hat\Kb_1(x)$ can be written as 
	\begin{equation}\label{rewrite}
		\left(\begin{matrix}
			\frac{x}{\lambda_1}-1+O_p(n^{-1/2})&O_p(n^{-1/2})&\cdots&O_p(n^{-1/2})\\
			O_p(n^{-1/2})&	\frac{x}{\lambda_2}-1+O_p(n^{-1/2})&\cdots&O_p(n^{-1/2})\\
			\vdots&	\vdots&		\ddots&	\vdots&\\
			O_p(n^{-1/2})&	O_p(n^{-1/2})&\cdots&	\frac{x}{\lambda_r}-1+O_p(n^{-1/2})
		\end{matrix}\right).
	\end{equation}
	Let $x=\lambda_1(1-c)$ for some sufficiently small $c>0$. By Assumption \ref{c1}(c), for any $2\le k\le r$,
	\[
	x/\lambda_k=x/\lambda_1\times \lambda_1/\lambda_k\ge(1-c)(1+2c)>1+c-2c^2>0,
	\]
	as long as $c<0.5$, which further implies that $\det\hat\Kb_1(x)<0$ with probability tending to one. On the other hand, if $x=\lambda_1(1+c)$ for some sufficiently small $c>0$, we can conclude that $\det\hat\Kb_1(x)>0$ with probability tending to one. Therefore, with probability tending to one, there must be an eigenvalue $\lambda_i(\hat{\mathcal{S}}_1)$ in the interval $[\lambda_1(1-c),\lambda_1(1+c)]$. Indeed, this is the largest one $\lambda_1(\hat{\mathcal{S}}_1)$. Further, when $x$ is in this interval, we always have $|x/\lambda_i-1|\ge c_0$ for all $i\ne 1$ and some small constant $c_0>0$ by Assumption \ref{c1}(c). Therefore, by Leibniz's formula for determinant and (\ref{rewrite}), we have
	\[
	\det\hat\Kb_1(\lambda_1(\hat{\mathcal{S}}_1))=\frac{\lambda_1(\hat{\mathcal{S}}_1)}{\lambda_1}-1+O_p(n^{-1/2})+O_p(n^{-1})=0.
	\] 
	That is, $\lambda_1(\hat{\mathcal{S}}_1)/\lambda_1=1+O_p(n^{-1/2})$. Similarly, we can conclude that $\lambda_i(\hat{\mathcal{S}}_1)/\lambda_i=1+O_p(n^{-1/2})$ for any $2\le i\le r$. Combined with (\ref{hat_lambda_i}), we  conclude the lemma.

	\subsection{Proof of Lemma \ref{lem: hat  theta}: $\theta_i$ and $\hat \zeta_i$}
	\begin{proof}
		We start with $\theta_i$. By definition, for $1\le i\le r$,
		\[
		\begin{split}
			\frac{\theta_i}{\lambda_i}-1=&\frac{\theta_i}{\lambda_i}\times\frac{1}{n\theta_i}\sum_{k=1}^p\frac{\lambda_{r+k}}{1-\lambda_i^{-1}\lambda_{r+k}}
			\le \frac{2}{n\theta_i}\sum_{k=1}^p\frac{\lambda_{r+k}}{1-\lambda_i^{-1}\lambda_{r+k}}=O\bigg(\frac{\text{tr}\bLambda_2}{n\theta_i}\bigg)=O\bigg(\frac{\text{tr}\bLambda_2}{n\lambda_i}\bigg).
		\end{split}
		\]
		
		Next, for $\hat\zeta_i$, by (\ref{mg}), (\ref{Ew1}) and (\ref{s2}), we have $\max_j w_j\le O(\log n)$ with probability tending to 1. Therefore, the existence and uniqueness of $\hat\zeta_i$ are easily verified by the mean value theorem. Moreover,
		\begin{equation}\label{hat zeta i minus}
			\begin{split}
				\hat\zeta_i-\frac{\theta_i}{\lambda_i}=&\frac{1}{n}\sum_{j=1}^n\frac{\theta_i}{\lambda_i}\frac{w_j-1-\frac{w_j}{n\theta_i}\sum_{k=1}^p[\frac{\lambda_{r+k}}{1-\theta_i^{-1}\lambda_{r+k}\hat\zeta_i}-\frac{\lambda_{r+k}}{1-\lambda_i^{-1}\lambda_{r+k}}]}{1-\frac{w_j}{n\theta_i}\sum_{k=1}^p\frac{\lambda_{r+k}}{1-\theta_i^{-1}\lambda_{r+k}\hat\zeta_i}}\\
				=&\frac{1}{n}\sum_{j=1}^n\frac{\theta_i}{\lambda_i}\frac{w_j-1}{1-\frac{w_j}{n\theta_i}\sum_{k=1}^p\frac{\lambda_{r+k}}{1-\theta_i^{-1}\lambda_{r+k}\hat\zeta_i}}\\
				&-\frac{\hat\zeta_i-\theta_i/\lambda_i}{\theta_i}\times \frac{1}{n}\sum_{j=1}^n\frac{\theta_i}{\lambda_i}\frac{\frac{w_j}{n\theta_i}\sum_{k=1}^p[\frac{\lambda_{r+k}^2}{[1-\theta_i^{-1}\lambda_{r+k}\hat\zeta_i][1-\lambda_i^{-1}\lambda_{r+k}]}]}{1-\frac{w_j}{n\theta_i}\sum_{k=1}^p\frac{\lambda_{r+k}}{1-\theta_i^{-1}\lambda_{r+k}\hat\zeta_i}}\\
				=&\frac{\theta_i}{\lambda_i}\frac{1}{n}\sum_{j=1}^n\frac{w_j-1}{1-\frac{w_j}{n\theta_i}\sum_{k=1}^p\frac{\lambda_{r+k}}{1-\theta_i^{-1}\lambda_{r+k}\hat\zeta_i}}-\frac{\hat\zeta_i-\theta_i/\lambda_i}{\theta_i}\times o_p(1).
			\end{split}
		\end{equation}
		In the following, we calculate the first term on the RHS. 	By Assumption \ref{c1}(c),  we have
		\[
		\frac{1}{n\theta_i}\sum_{k=1}^p\frac{\lambda_{r+k}}{1-\theta_i^{-1}\lambda_{r+k}\hat\zeta_i}=\frac{\text{tr}\bLambda_2}{n\theta_i}\times [1+o_p(1)]=o_p(1).
		\] 
		Therefore, by Taylor's expansion for the function $f(x)=(1-x)^{-1}$, we have 
		\[
		\begin{split}
			&\frac{1}{n}\sum_{j=1}^n\frac{w_j-1}{1-\frac{w_j}{n\theta_i}\sum_{k=1}^p\frac{\lambda_{r+k}}{1-\theta_i^{-1}\lambda_{r+k}\hat\zeta_i}}\\
			=&\frac{1}{n}\sum_{j=1}^n(w_j-1)+\frac{1}{n\theta_i}\sum_{k=1}^p\frac{\lambda_{r+k}}{1-\theta_i^{-1}\lambda_{r+k}\hat\zeta_i}\times \frac{1}{n}\sum_{j=1}^nw_j(w_j-1)\\
			&+\bigg(\frac{1}{n\theta_i}\sum_{k=1}^p\frac{\lambda_{r+k}}{1-\theta_i^{-1}\lambda_{r+k}\hat\zeta_i}\bigg)^2\times \frac{1}{n}\sum_{j=1}^nw_j^2(w_j-1)
			+o_p(1)\times \bigg(\frac{1}{n\theta_i}\sum_{k=1}^p\frac{\lambda_{r+k}}{1-\theta_i^{-1}\lambda_{r+k}\hat\zeta_i}\bigg)^2.
		\end{split}
		\]
		Under the multiplier bootstrap, 
		\begin{equation}\label{wjh}
			\frac{1}{n}\sum_{j=1}^nw_j^h(w_j-1)-\mathbb{E}w_1^h(w_1-1)=O_p(\frac{1}{\sqrt{n}}),\quad, h=0,1,2.
		\end{equation}
		Under the standard bootstrap, $w_{j_1}\mid w_{j_2}\sim Bin(n-w_{j_2},(n-1)^{-1})$ for any $j_1\ne j_2$. Then,
		\[
		\begin{split}
			&\mathbb{E}\bigg(w_{j_1}^h(w_{j_1}-1)-\mathbb{E}w_1^h(w_1-1)\bigg)\bigg(w_{j_2}^h(w_{j_2}-1)-\mathbb{E}w_1^h(w_1-1)\bigg)\\
			=&\mathbb{E}\bigg[\bigg(\mathbb{E}w_{j_1}^h(w_{j_1}-1)\mid w_{j_2}-\mathbb{E}w_1^h(w_1-1)\bigg)\bigg(w_{j_2}^h(w_{j_2}-1)-\mathbb{E}w_1^h(w_1-1)\bigg)\bigg].
		\end{split}
		\]
		Note that 
		\[
		\begin{split}
			\mathbb{E}(w_{j_1}^h\mid w_{j_2})-\mathbb{E}w_1^h=\sum_{l=0}^h\left(\begin{matrix}
				h\\
				l
			\end{matrix}\right)(n-w_{j_2})^l\frac{1}{(n-1)^l}-\sum_{l=0}^h\left(\begin{matrix}
				h\\
				l
			\end{matrix}\right)=O_{L_1}(n^{-1}).
		\end{split}
		\]
		Then, after some elementary calculations, we conclude that (\ref{wjh}) also holds under the standard bootstrap. Further, for $h=0,1$, $\mathbb{E}w_j^h(w_j-1)=O(n^{-1})$. Then, using the fact that $\theta_i/\lambda_i=1+o(1)$, (\ref{hat zeta i minus}) can be written as
		\[
		\begin{split}
			\bigg(\hat\zeta_i-\frac{\theta_i}{\lambda_i}\bigg)\times [1+o_p(1)]=&\frac{1}{n}\sum_{j=1}^n(w_j-1)+\bigg(\frac{\text{tr}\bLambda_2}{n\theta_i}\bigg)^2\times \mathbb{E}w_1^2(w_1-1)\\
			&+o_p(\frac{1}{\sqrt{n}})+o_p(1)\times \bigg(\frac{\text{tr}\bLambda_2}{n\theta_i}\bigg)^2,
		\end{split}
		\]
		which concludes the lemma because $n^{-1}\sum_j(w_j-1)=O_p(n^{-1/2})$. 
	\end{proof}
	
	\subsection{Proof of Theorem \ref{thm: representation}: limiting representation for $\hat\lambda_i, 1\le i\le r$}
	\begin{proof}
		We aim to find the limiting representation of $\hat\lambda_i/\theta_i$ for $i\le r$. The proof technique is borrowed from Theorem 2.4 of \cite{cai2020limiting}, which can be regarded as the  special case where  $\Wb=\Ib$. 
		
		It suffices to prove the result for $\hat\lambda_1$, while the others can be handled similarly. By definition, $\hat\lambda_1$ is the largest eigenvalue  satisfying
		\[
		\det\big(\hat\lambda_1-n^{-1}\Wb^{1/2}\Zb^\top(\bGamma_1\bLambda_1\bGamma_1^\top+\bGamma_2\bLambda_2\bGamma_2^\top)\Zb\Wb^{1/2}\big)=0.
		\]
		By Lemma \ref{lem: preliminary}, $\det(\hat\lambda_1-\hat{\mathcal{S}}_2)\ne 0$ with probability tending to one, where $\hat{\mathcal{S}}_2$ is defined in the proof of Lemma \ref{lem: preliminary} as $n^{-1}\Wb^{1/2}\Zb^\top\bGamma_2\bLambda_2\bGamma_2^\top\Zb\Wb^{1/2}$. Then, 
		\[
		\det\big(\bLambda_1^{-1}-n^{-1}\bGamma_1^\top\Zb\Wb^{1/2}[\hat\lambda_1\Ib-\hat{\mathcal{S}}_2]^{-1}\Wb^{1/2}\Zb^\top\bGamma_1\big)=0.
		\] 
		Write  $\delta_i=(\hat\lambda_i-\theta_i)/\theta_i$ and $\Kb(x)=[\Ib-x^{-1}\hat{\mathcal{S}}_2]^{-1}$. By the matrix inverse formula 
		\begin{equation}\label{matrix inverse}
			(\Ab-\Bb)^{-1}=\Ab^{-1}+(\Ab-\Bb)^{-1}\Bb\Ab^{-1},
		\end{equation}
		as long as the associated inverses exist,  we have $\det\Mb(\theta_1)=0$ where 
		\begin{equation}\label{det}
			\begin{split}
				\Mb(\theta_1)=&\theta_1\bLambda_1^{-1}-n^{-1}\bGamma_1^\top\Zb\Wb^{1/2}\Kb(\theta_1)\Wb^{1/2}\Zb^\top\bGamma_1
				+\delta_1n^{-1}\bGamma_1^\top\Zb\Wb^{1/2}\Kb(\hat\lambda_1)\Kb(\theta_1)\Wb^{1/2}\Zb^\top\bGamma_1.
			\end{split}
		\end{equation}
		Lemma \ref{lema2} will show that 
		\[
		\begin{split}
			&n^{-1}\bGamma_1^\top\Zb\Wb^{1/2}\Kb(\theta_1)\Wb^{1/2}\Zb^\top\bGamma_1\\
			=&\frac{1}{n}\bGamma_1^\top\Zb\Wb\Zb^\top\bGamma_1-\bigg(\frac{1}{n}\text{tr}\Wb\bigg)\Ib_r+\hat\zeta_1\Ib_r+O_p\bigg(\frac{1}{\sqrt{n}}\times \frac{(n\vee p)\log n}{n\theta_1}+\frac{1}{n}\bigg),
		\end{split}
		\]
		where the $O_p$ is under Frobenius norm.  
		Then, by (\ref{W-I}), we have
		\[
		n^{-1}\bGamma_1^\top\Zb\Wb^{1/2}\Kb(\theta_1)\Wb^{1/2}\Zb^\top\bGamma_1=\hat\zeta_1\Ib_r+O_p(n^{-1/2}).
		\]
		Return to (\ref{det}). Lemma \ref{lema3} will further show that 
		\[
		\delta_1n^{-1}\bGamma_1^\top\Zb\Wb^{1/2}\Kb(\hat\lambda_1)\Kb(\theta_1)\Wb^{1/2}\Zb^\top\bGamma_1=\delta_1[\Ib+o_p(1)]. 
		\]
		Recall that $\delta_1=o_p(1)$. Consequently, the off-diagonal entries of $\Mb(\theta_1)$ all converge to 0 with rate
		\begin{equation}\label{M off}
			O_p(n^{-1/2})+\delta_1[1+o_p(1)].
		\end{equation}
		The first diagonal entry of $\Mb(\theta_1)$ can be written as 
		\begin{equation}\label{M_11}
			-n^{-1}\bgamma_1^\top\Zb\Wb\Zb^\top\bgamma_1+\frac{1}{n}\text{tr}\Wb+\frac{\theta_1}{\lambda_1}-\hat\zeta_1+O_p\bigg(\frac{1}{\sqrt{n}}\times \frac{(n\vee p)\log n}{n\theta_1}+\frac{1}{n}\bigg)+\delta_1[1+o_p(1)].
		\end{equation}
		For the other diagonal entries, note that  
		\[
		\bigg|\frac{\theta_1}{\lambda_j}-\frac{\theta_1}{\lambda_1}\bigg|\ge c,\quad 2\le j\le r,
		\]
		for some constant $c>0$. Therefore, we have
		\begin{equation}\label{M_jj}
			[\Mb(\theta_1)]_{jj}\ge c+o_p(1),\quad j\ne 1.
		\end{equation}
		Recall that $\det\Mb(\theta_1)=0$.
		Then, by  (\ref{M off}), (\ref{M_11}), (\ref{M_jj}), the Leibniz's formula for determinant and the fact that $\delta_1=o_p(1)$, we conclude that 
		\[
		\delta_1[1+o_p(1)]=n^{-1}\bgamma_1^\top\Zb\Wb\Zb^\top\bgamma_1-\frac{1}{n}\text{tr}\Wb-\frac{\theta_1}{\lambda_1}+\hat\zeta_1+O_p\bigg(\frac{1}{\sqrt{n}}\times \frac{(n\vee p)\log n}{n\theta_1}+\frac{1}{n}\bigg).
		\]
		Note that $\hat\zeta_1-\theta_1/\lambda_1=O_p((n\lambda_1)^{-2}p^2)$ by Lemma \ref{lem: hat theta}, while
		\[
		\begin{split}
			&n^{-1}\bgamma_1^\top\Zb\Wb\Zb^\top\bgamma_1-\frac{1}{n}\text{tr}\Wb\\
			=&n^{-1}\bgamma_1^\top\Zb(\Wb-\Ib)\Zb^\top\bgamma_1+n^{-1}\bgamma_1^\top\Zb\Zb^\top\bgamma_1-1+1-\frac{1}{n}\text{tr}\Wb=O_p(\frac{1}{\sqrt{n}}). 
		\end{split}
		\]
		Combined with the rate in Lemma \ref{lem: preliminary}, we have verified the theorem for $i=1$.
		For $2\le i\le r$, it's similar and we omit details. 
	\end{proof}

	\subsection{Proof of Lemma \ref{lem: without bootstrap}: without bootstrap}
	\begin{proof}
		\textbf{Part (a): eigenvalues.}
		The proof is almost the same as that for Lemma \ref{lem: preliminary} and Theorem \ref{thm: representation}, by replacing $\Wb$ with $\Ib$; see also \cite{cai2020limiting}. Therefore, we omit the details. 
		
		\textbf{Part (b): eigenvectors $\tilde\bgamma_i$ for $1\le i\le r$.}
		Let's start with $\bgamma_i^\top\tilde\bgamma_j$ for some $1\le i\ne j\le r$. By the definition of eigenvector,
		\begin{equation}\label{expansion 1}
			\begin{split}
				&\tilde\lambda_j\bgamma_i^\top\tilde\bgamma_j=n^{-1}\bgamma_i^\top\bGamma\bLambda^{1/2}\bGamma^\top\Zb\Zb^\top\bGamma\bLambda^{1/2}\bGamma^\top\tilde\bgamma_j\\
				=&\frac{\lambda_i}{n}\bgamma_i^\top\Zb\Zb^\top\bgamma_i\bgamma_i^\top\tilde\bgamma_j+\sum_{k\ne i}^r\frac{\sqrt{\lambda_i\lambda_k}}{n}\bgamma_i^\top\Zb\Zb^\top\bgamma_k\bgamma_k^\top\tilde\bgamma_j+\frac{\sqrt{\lambda_i}}{n}\bgamma_i^\top\Zb\Zb^\top\bGamma_2\bLambda_2^{1/2}\bGamma_2^\top\tilde\bgamma_j.
			\end{split}
		\end{equation}
		We already know that $n^{-1}\bgamma_i^\top\Zb\Zb^\top\bgamma_i=1+O_p(n^{-1/2})$ and $n^{-1}\bgamma_i^\top\Zb\Zb^\top\bgamma_k=O_p(n^{-1/2})$ for $i\ne k$. Moreover,
		\[
		\mathbb{E}\|n^{-1}\bgamma_i^\top\Zb\Zb^\top\bGamma_2\|^2=\sum_{k=r+1}^{r+p}\mathbb{E}|n^{-1}\bgamma_i^\top\Zb\Zb^\top\bgamma_k|^2\le O(p/n).
		\]
		Therefore, by the Cauchy--Schwartz inequality
		\begin{equation}\label{CS rate}
			|\bgamma_k^\top\tilde\bgamma_j|\le 1,\quad \bigg|\frac{1}{n}\bgamma_i^\top\Zb\Zb^\top\bGamma_2\bLambda_2^{1/2}\bGamma_2^\top\tilde\bgamma_j\bigg|\le O_p(\sqrt{p/n}).
		\end{equation}
		If $i>j$, we can write
		\[
		\bgamma_i^\top\tilde\bgamma_j\bigg[1-\frac{\lambda_i}{\tilde\lambda_j}(1+O_p(n^{-1/2}))\bigg]\le\frac{\sqrt{\lambda_1\lambda_i}}{\tilde\lambda_j}O_p(n^{-1/2})+\frac{\sqrt{\lambda_i}}{\tilde\lambda_j}O_p(\sqrt{p/n}).
		\]
		Otherwise, when $i<j$ we have
		\[
		\bgamma_i^\top\tilde\bgamma_j\bigg[\frac{\tilde\lambda_j}{\lambda_i}-1-O_p(n^{-1/2})\bigg]\le\frac{\sqrt{\lambda_1\lambda_i}}{\lambda_i}O_p(n^{-1/2})+\frac{\sqrt{\lambda_i}}{\lambda_i}O_p(\sqrt{p/n}).
		\]
		Note that $\lambda_i/\tilde\lambda_j\le (1+c)^{-1}[1+o_p(1)]$ for $i>j$ while $\tilde\lambda_j/\lambda_i-1\le -c/(1+c)[1+o_p(1)]$ for $i<j$. Therefore, we conclude that
		\begin{equation}\label{rate 1}
			\bgamma_i^\top\tilde\bgamma_j\le O_p\bigg(\frac{1}{\sqrt{n}}\frac{\sqrt{\lambda_1\lambda_i}}{\max\{\lambda_i,\lambda_j\}}+\frac{\sqrt{\lambda_i}}{\max\{\lambda_i,\lambda_j\}}\sqrt{\frac{p}{n}}\bigg)=O_p\bigg(\frac{\sqrt{(\lambda_1+p)\lambda_i}}{\sqrt{n}\times\max\{\lambda_i,\lambda_j\}}\bigg).
		\end{equation}
		The rate in (\ref{rate 1}) can further help bound $\bgamma_k^\top\tilde\bgamma_j$ in (\ref{CS rate}) for $k\ne j$. Following this idea and letting $a_n=\sqrt{(\lambda_1+p)\lambda_i}/(\sqrt{n}\times\max\{\lambda_i,\lambda_j\})$, we will have
		\[
		\begin{split}
			\bgamma_i^\top\tilde\bgamma_j\le& O_p\bigg(\frac{\sqrt{\lambda_i\lambda_j}}{\sqrt{n}\max\{\lambda_i,\lambda_j\}}+a_n\frac{\sqrt{\lambda_i\lambda_1}}{\sqrt{n}\max\{\lambda_i,\lambda_j\}}+\frac{\sqrt{p\lambda_i}}{\sqrt{n}\max\{\lambda_i,\lambda_j\}}\bigg)\\
			\le &O_p\bigg(\frac{\sqrt{(p+\lambda_j+a_n\lambda_1)\lambda_i}}{\sqrt{n}\max\{\lambda_i,\lambda_j\}}\bigg).
		\end{split}
		\]
		Repeating the above step, eventually we have
		\begin{equation}\label{gamma_i gamma_j}
			\bgamma_i^\top\tilde\bgamma_j\le O_p\bigg(\frac{\sqrt{(p+\lambda_j)\lambda_i}}{\sqrt{n}\max\{\lambda_i,\lambda_j\}}\bigg).
		\end{equation}

		Next, we consider $\bGamma_2^\top\tilde\bgamma_j$ for $1\le j\le r$. Similarly to (\ref{expansion 1}),
		\[
		\begin{split}
			&\tilde\lambda_j\bGamma_2^\top\tilde\bgamma_j=n^{-1}\bGamma_2^\top\bGamma\bLambda^{1/2}\bGamma^\top\Zb\Zb^\top\bGamma\bLambda^{1/2}\bGamma^\top\tilde\bgamma_j\\
			=&\frac{\sqrt{\lambda_j}}{n}\bLambda_2^{1/2}\bGamma_2^\top\Zb\Zb^\top\bgamma_j\bgamma_j^\top\tilde\bgamma_j+\sum_{k\ne j}^r\frac{\sqrt{\lambda_k}}{n}\bLambda_2^{1/2}\bGamma_2^\top\Zb\Zb^\top\bgamma_k\bgamma_k^\top\tilde\bgamma_j+\frac{1}{n}\bLambda_2^{1/2}\bGamma_2^\top\Zb\Zb^\top\bGamma_2\bLambda_2^{1/2}\bGamma_2^\top\tilde\bgamma_j.
		\end{split}
		\]
		Recall that $\|n^{-1}\bGamma_2^\top\Zb\Zb^\top\bgamma_j\|^2\le O_p(p/n)$. Then, 
		\[
		\begin{split}
			&	\bigg[\Ib-\frac{1}{n\tilde\lambda_j}\bLambda_2^{1/2}\bGamma_2^\top\Zb\Zb^\top\bGamma_2\bLambda_2^{1/2}\bigg]\bGamma_2^\top\tilde\bgamma_j\\
			=&\frac{\sqrt{\lambda_j}}{n\tilde\lambda_j}\bLambda_2^{1/2}\bGamma_2^\top\Zb\Zb^\top\bgamma_j\bgamma_j^\top\tilde\bgamma_j+\sum_{k\ne j}^r\frac{\sqrt{\lambda_k}}{n\tilde\lambda_j}\bLambda_2^{1/2}\bGamma_2^\top\Zb\Zb^\top\bgamma_k\bgamma_k^\top\tilde\bgamma_j\\
			\le &O_p\bigg(\frac{\sqrt{\lambda_j}}{\tilde\lambda_j}\sqrt{\frac{p}{n}}+\sum_{k\ne j}\frac{\sqrt{\lambda_k}}{\tilde\lambda_j}\sqrt{\frac{p}{n}}\times\frac{\sqrt{(p+\lambda_j)\lambda_k}}{\sqrt{n}\max\{\lambda_k,\lambda_j\}} \bigg)\\
			\le &O_p\bigg(\sqrt{\frac{p}{n\lambda_j}}\bigg),
		\end{split}
		\]
		where the $O_p$ is under Frobenius norm. Note that 
		\[
		\bigg\|\frac{1}{n\tilde\lambda_j}\bLambda_2^{1/2}\bGamma_2^\top\Zb\Zb^\top\bGamma_2\bLambda_2^{1/2}\bigg\|=o_p(1).
		\]
		Therefore, we can conclude that
		\begin{equation}\label{Gamma_2 gamma_j}
			\|\bGamma_2^\top\tilde\bgamma_j\|\le O_p\bigg(\sqrt{\frac{p}{n\lambda_j}}\bigg).
		\end{equation}
		Based on (\ref{Gamma_2 gamma_j}), we can further improve the rate in (\ref{CS rate}) to 
		\[
		\bigg|\frac{1}{n}\bgamma_i^\top\Zb\Zb^\top\bGamma_2\bLambda_2^{1/2}\bGamma_2^\top\tilde\bgamma_j\bigg|\le O_p\bigg(\frac{p}{n\sqrt{\lambda_j}}\bigg),
		\]
		which further improves the rate in (\ref{gamma_i gamma_j}) to
		\begin{equation}\label{gamma_i gamma_j final}
			\bgamma_i^\top\tilde\bgamma_j\le O_p\bigg(\frac{p}{n\sqrt{\lambda_i\lambda_j}}\frac{\lambda_i}{\max\{\lambda_i,\lambda_j\}}+\frac{\sqrt{\lambda_i\lambda_j}}{\sqrt{n}\times \max\{\lambda_i,\lambda_j\}}\bigg),\quad 1\le i\ne j\le r.
		\end{equation}
		
		Based on (\ref{gamma_i gamma_j final}), for any $1\le j\le r$, we have
		\begin{equation}\label{gamma_j gamma_j}
			\begin{split}
				&1=\tilde\bgamma_j^\top\tilde\bgamma_j=\tilde\bgamma_j^\top\bGamma\bGamma^\top\tilde\bgamma_j=(\tilde\bgamma_j^\top\bgamma_j)^2+\sum_{k\ne j}^r(\tilde\bgamma_j^\top\bgamma_k)^2+\|\tilde\bgamma_j^\top\bGamma_2\|^2,\\
				\Longrightarrow &(\tilde\bgamma_j^\top\bgamma_j)^2=1+O_p\bigg(\frac{p}{n\lambda_j}+\frac{1}{n}\bigg).
			\end{split}
		\end{equation}

		\textbf{Part (c): eigenvectors $\tilde\bu_i$ for $1\le i\le r$.}	 By the definition of eigenvector,
		\[
		\tilde\bu_i=\frac{1}{\sqrt{n\tilde\lambda_i}}\Zb^\top\bGamma\bLambda^{1/2}\bGamma^\top\tilde\bgamma_i.
		\]
		Let $\bu_i=n^{-1/2}\Zb^\top\bgamma_i$. Then, for $1\le j\le (p+r)$,
		\begin{equation}\label{expansion 2}
			\begin{split}
				\frac{\sqrt{\tilde\lambda_i}}{\sqrt{\lambda_i}}\tilde u_{ij}=u_{ij}\bgamma_i^\top\tilde\bgamma_i+\sum_{k\ne i}^r\frac{\sqrt{\lambda_k}}{\sqrt{n\lambda_i}}\bz_j^\top\bgamma_k\bgamma_k^\top\tilde\bgamma_i+\frac{1}{\sqrt{n\lambda_i}}\bz_j^\top\bGamma_2\bLambda_2^{1/2}\bGamma_2^\top\tilde\bgamma_i.
			\end{split}
		\end{equation}
		We will show that the first term on the right hand side (RHS) will dominate.  Firstly, for any $1\le i\le r$,  Lemma \ref{lem: ui} will show that $\sum_{j=1}^nu_{ij}^4=(\mathbb{E}\sum_{j=1}^nu_{ij}^4)[1+o_p(1)]$, while 
		\[
		\mathbb{E}\sum_{j=1}^nu_{ij}^4=n^{-1}(\xi_i+1)\asymp n^{-1},
		\]
		under Assumption \ref{c2}. Then, by (\ref{gamma_j gamma_j}),
		\begin{equation}\label{first term for ui}
			\sum_{j=1}^n(u_{ij}\bgamma_i^\top\tilde\bgamma_i)^4=\bigg(\sum_{j=1}^nu_{ij}^4\bigg)(\bgamma_i^\top\tilde\bgamma_i)^4=\bigg(\mathbb{E}\sum_{j=1}^nu_{ij}^4\bigg)[1+o_p(1)]\asymp n^{-1}.
		\end{equation}

		Secondly, for each $k\in\{1,\ldots,r\}\setminus \{i\}$, we have
		\begin{equation}\label{second term for ui}
			\sum_{j=1}^n\bigg(\frac{\sqrt{\lambda_k}}{\sqrt{n\lambda_i}}\bz_j^\top\bgamma_k\bgamma_k^\top\tilde\bgamma_i\bigg)^4=(\sum_{j=1}^nu_{kj}^4)\times \bigg(\frac{\lambda_k}{n\lambda_i}\bigg)^2(\bgamma_k^\top\tilde\bgamma_i)^4\le o_p(n^{-1}),
		\end{equation}
		where we use (\ref{gamma_i gamma_j final}) and Lemma \ref{lem: ui} again. 
		
		The third term on the RHS of (\ref{expansion 2}) will be more complicated. 
		Let $\Zb_j$ be the $(p+r)\times n$ random matrix by replacing $\bz_j$ with ${\bf 0}$ in $\Zb$, 
		\[
		\begin{split}
			&\tilde\Sbb_{j,21}=n^{-1}\bLambda_2^{1/2}\bGamma_2^\top\Zb_j\Zb_j^\top\bGamma_1\bLambda_1^{1/2},\quad \tilde\Sbb_{j,22}=n^{-1}\bLambda_2^{1/2}\bGamma_2^\top\Zb_j\Zb_j^\top\bGamma_2\bLambda_2^{1/2},\\
			&\Hb_{j,22}=\Ib-\tilde\lambda_{i}^{-1}\tilde\Sbb_{j,22},\quad a_{nj}=(n\tilde\lambda_{i})^{-1}\bz_j^\top\bGamma_2\bLambda_2^{1/2}\Hb_{j,22}^{-1}\bLambda_2^{1/2}\bGamma_2^\top\bz_j.
		\end{split}
		\]
		Then, by definition,
		\[
		\begin{split}
			\bGamma_2^\top\tilde\bgamma_i
			=&\frac{1}{n\tilde\lambda_{i}}\bLambda_2^{1/2}\bGamma_2^\top\Zb\Zb^\top\bGamma\bLambda^{1/2}\bGamma^\top\tilde\bgamma_i=\tilde\lambda_{i}^{-1}\tilde\Sbb_{j,21}\bGamma_1^\top\tilde\bgamma_i+\tilde\lambda_{i}^{-1}\tilde\Sbb_{j,22}\bGamma_2^\top\tilde\bgamma_i\\
			&+\frac{1}{n\tilde\lambda_{i}}\bLambda_2^{1/2}\bGamma_2^\top\bz_j\bz_j^\top\bGamma_1\bLambda_1^{1/2}\bGamma_1^\top\tilde\bgamma_i+\frac{1}{n\tilde\lambda_{i}}\bLambda_2^{1/2}\bGamma_2^\top\bz_j\bz_j^\top\bGamma_2\bLambda_2^{1/2}\bGamma_2^\top\tilde\bgamma_i.
		\end{split}
		\]
		Therefore, 
		\begin{equation}\label{third term construct}
			\begin{split}
				&\frac{1}{\sqrt{n\lambda_i}}\bz_j^\top\bGamma_2\bLambda_2^{1/2}\Hb_{j,22}^{-1}\bGamma_2^\top\tilde\bgamma_i\\
				=&	\frac{1}{\tilde\lambda_i\sqrt{n\lambda_i}}\bz_j^\top\bGamma_2\bLambda_2^{1/2}\Hb_{j,22}^{-1}\tilde\Sbb_{j,21}\bGamma_1^\top\tilde\bgamma_i+\frac{1}{\tilde\lambda_i\sqrt{n\lambda_i}}\bz_j^\top\bGamma_2\bLambda_2^{1/2}\Hb_{j,22}^{-1}\tilde\Sbb_{j,22}\bGamma_2^\top\tilde\bgamma_i\\
				&+a_{nj}\frac{1}{\sqrt{n\lambda_i}}\bz_j^\top\bGamma_1\bLambda_1^{1/2}\bGamma_1^\top\tilde\bgamma_i+a_{nj}\frac{1}{\sqrt{n\lambda_i}}\bz_j^\top\bGamma_2\bLambda_2^{1/2}\bGamma_2^\top\tilde\bgamma_i.
			\end{split}
		\end{equation}
		However, by the definition of $\Hb_{j,22}$, we have
		\[
		\Hb_{j,22}^{-1}=\Ib+\Hb_{j,22}^{-1}\tilde\lambda_i^{-1}\tilde\Sbb_{j,22}.
		\]
		Using the above decomposition in the first line of (\ref{third term construct}) and after cancellation, we have
		\begin{equation}\label{anj}
			\begin{split}
				&(1-a_{nj})\frac{1}{\sqrt{n\lambda_i}}\bz_j^\top\bGamma_2\bLambda_2^{1/2}\bGamma_2^\top\tilde\bgamma_i
				\\=&	\frac{1}{\tilde\lambda_i\sqrt{n\lambda_i}}\bz_j^\top\bGamma_2\bLambda_2^{1/2}\Hb_{j,22}^{-1}\tilde\Sbb_{j,21}\bGamma_1^\top\tilde\bgamma_i
				+a_{nj}\frac{1}{\sqrt{n\lambda_i}}\bz_j^\top\bGamma_1\bLambda_1^{1/2}\bGamma_1^\top\tilde\bgamma_i.
			\end{split}
		\end{equation}
		
		We need to discuss the magnitude of $\max_j|a_{n,j}|$. Note that $\max_j\|\tilde\Sbb_{j,22}\|\le \|n^{-1}\Zb\Zb^\top\|\le O_p(\max\{p/n,1\})$.
		{Then, $\max_j\|\Hb_{2,jj}\|=1+o_p(1)$ and $\min_j\|\Hb_{2,jj}\|=1+o_p(1)$ because  $\tilde\lambda_j\rightarrow \infty$. Further,}
		\[
		\max_j\|\Hb_{2,jj}^{-1}\|\le \frac{1}{1-o_p(1)}\le 1+o_p(1).
		\]
		Further, we have
		\[
		\max_j|a_{n,j}|\le \lambda_j^{-1}[1+o_p(1)]\max_jn^{-1}\|\bz_j\|^2\max_j\|\Hb_{j,22}^{-1}\|\le O_p(\frac{n\vee p}{n\lambda_i}) =o_p(1),
		\]
		and $\max_j|(1-a_{n,j})^{-1}|\le 1+o_p(1)$. 
		
		Return to the RHS of (\ref{anj}). For the first term, we write
		\[
		\begin{split}
			\frac{1}{\tilde\lambda_i\sqrt{n\lambda_i}}\bz_j^\top\bGamma_2\bLambda_2^{1/2}\Hb_{j,22}^{-1}\tilde\Sbb_{j,21}\bGamma_1^\top\tilde\bgamma_i
			=&\sum_{k\ne i}^r\frac{\sqrt{\lambda_k}}{\tilde\lambda_i\sqrt{n\lambda_i}}\bz_j^\top\bGamma_2\bLambda_2^{1/2}\Hb_{j,22}^{-1}\bLambda_2^{1/2}\bGamma_2^\top n^{-1}\Zb_j\Zb_j^\top\bgamma_k\bgamma_k^\top\tilde\bgamma_i\\
			&+\frac{1}{\tilde\lambda_i\sqrt{n}}\bz_j^\top\bGamma_2\bLambda_2^{1/2}\Hb_{j,22}^{-1}\bLambda_2^{1/2}\bGamma_2^\top n^{-1}\Zb_j\Zb_j^\top\bgamma_i\bgamma_i^\top\tilde\bgamma_i.
		\end{split}
		\]
		On one hand, for any $k\ne i$, 
		\[
		\begin{split}
			&\sum_{j=1}^n\bigg(\frac{\sqrt{\lambda_k}}{\tilde\lambda_i\sqrt{n\lambda_i}}\bz_j^\top\bGamma_2\bLambda_2^{1/2}\Hb_{j,22}^{-1}\bLambda_2^{1/2}\bGamma_2^\top n^{-1}\Zb_j\Zb_j^\top\bgamma_k\bgamma_k^\top\tilde\bgamma_i\bigg)^4\\
			\le &\frac{\lambda_k^2}{\tilde\lambda_i^2}(\bgamma_k^\top\tilde\bgamma_i)^4\times \sum_{j=1}^n\bigg(\frac{1}{\sqrt{n\lambda_i}}\bz_j^\top\bGamma_2\bLambda_2^{1/2}\Hb_{j,22}^{-1}\bLambda_2^{1/2}\bGamma_2^\top n^{-1}\Zb_j\Zb_j^\top\bgamma_k\bigg)^4\\
			\le &o_p(1)\times O_p\bigg(\frac{1}{n\lambda_i^4}\max_j\|\Hb_{j,22}^{-1}\|^4\|n^{-1}\Zb\Zb^\top\|^4\bigg)=o_p(n^{-1}),
		\end{split}
		\]
		where in the third line we use the results in (\ref{gamma_i gamma_j final}) and the fact that 
		\[
		\begin{split}
			&\sum_{j=1}^n\bigg(\bz_j^\top\bGamma_2\bLambda_2^{1/2}\Hb_{j,22}^{-1}\bLambda_2^{1/2}\bGamma_2^\top n^{-1}\Zb_j\Zb_j^\top\bgamma_k\bigg)^4\\
			= &\sum_{j=1}^n\bigg([(\bbf_j^\top,{\bf 0}^\top)+({\bf 0}^\top,\bepsilon_j^\top)]\bGamma_2\bLambda_2^{1/2}\Hb_{j,22}^{-1}\bLambda_2^{1/2}\bGamma_2^\top n^{-1}\Zb_j\Zb_j^\top\bgamma_k\bigg)^4\\
			\le& O_p\bigg(\max_j\|\Hb_{j,22}^{-1}\|^4\|n^{-1}\Zb\Zb^\top\|^4\sum_{j=1}^n\|(\bbf_j^\top,{\bf 0}^\top)\bGamma_2\|^4+ \sum_{j=1}^n\mathbb{E}\|\bGamma_2\bLambda_2^{1/2}\Hb_{j,22}^{-1}\bLambda_2^{1/2}\bGamma_2^\top n^{-1}\Zb_j\Zb_j^\top\bgamma_k\|^4\bigg)\\
			\le &O_p\bigg(n\max_j\|\Hb_{j,22}^{-1}\|^4\|n^{-1}\Zb\Zb^\top\|^4\bigg),
		\end{split}
		\]
		by using the fact that $\bepsilon_j$ is independent of $\Hb_{j,22}^{-1}$ and $\Zb_j$ in the third line.
		Similarly,  by (\ref{gamma_j gamma_j}) we have
		\[
		\begin{split}
			&\sum_{j=1}^n\bigg(\frac{1}{\tilde\lambda_i\sqrt{n}}\bz_j^\top\bGamma_2\bLambda_2^{1/2}\Hb_{j,22}^{-1}\bLambda_2^{1/2}\bGamma_2^\top n^{-1}\Zb_j\Zb_j^\top\bgamma_i\bgamma_i^\top\tilde\bgamma_i\bigg)^4\\
			\le &O_p\bigg(\frac{1}{n\lambda_i^4}\max_j\|\Hb_{j,22}^{-1}\|^4\|n^{-1}\Zb\Zb^\top\|^4\bigg)=o_p(n^{-1}).
		\end{split}
		\]
		Therefore, we conclude that
		\[
		\sum_{j=1}^n\bigg(\frac{1}{\tilde\lambda_i\sqrt{n\lambda_i}}\bz_j^\top\bGamma_2\bLambda_2^{1/2}\Hb_{j,22}^{-1}\tilde\Sbb_{j,21}\bGamma_1^\top\tilde\bgamma_i\bigg)^4=o_p(n^{-1}).
		\]
		
		For the second term on the RHS of (\ref{anj}), it's similar and easier, so we conclude that
		\[
		\sum_{j=1}^n\bigg(a_{nj}\frac{1}{\sqrt{n\lambda_i}}\bz_j^\top\bGamma_1\bLambda_1^{1/2}\bGamma_1^\top\tilde\bgamma_i\bigg)^4=o_p(n^{-1}),
		\]
		without showing further details. Then, by (\ref{anj}), we have
		\begin{equation}\label{third term for ui}
			\sum_{j=1}^n\bigg(\frac{1}{\sqrt{n\lambda_i}}\bz_j^\top\bGamma_2\bLambda_2^{1/2}\bGamma_2^\top\tilde\bgamma_i\bigg)^4\le \max_j|(1-a_{nj})^{-4}|\times o_p(n^{-1})=o_p(n^{-1}).
		\end{equation}
		
		Combining (\ref{first term for ui}), (\ref{second term for ui}), (\ref{third term for ui}) and returning to (\ref{expansion 2}), we have
		\[
		\tilde\sigma_i^2=\sum_{j=1}^n\tilde u_{ij}^4=\sum_{j=1}^n\mathbb{E}(u_{ij}^4)+o_p(n^{-1})=n^{-1}(\xi_i+1)[1+o_p(1)]\asymp n^{-1},
		\]
		which concludes the lemma.
	\end{proof}

	\subsection{Proof of Theorem \ref{thm: conditional on sample}: conditional on sample}
	\begin{proof}
		Following the proof of Theorem \ref{thm: representation},  with probability tending to one we have
		\[
		\frac{\hat\lambda_i}{\theta_i}-1=\frac{1}{n}\bigg(\bgamma_i^\top\Zb\Wb\Zb^\top\bgamma_i-\text{tr}\Wb\bigg)+\hat\zeta_i-\frac{\theta_i}{\lambda_i}+o_{p^*}(\frac{1}{\sqrt{n}})+o_{p^*}\bigg(\frac{p^2}{(n\theta_i)^2}\bigg),\quad 1\le i\le r,
		\]
		conditional on $\Xb$. On the other hand, by Lemma \ref{lem: without bootstrap}(1), with probability tending to one,
		\[
		\frac{\tilde\lambda_i}{\theta_i}-1=\frac{1}{n}\bigg(\bgamma_i^\top\Zb\Zb^\top\bgamma_i-n\bigg)+o(\frac{1}{\sqrt{n}})=o(1),\quad 1\le i\le r.
		\]
		Therefore,  if $(n\lambda_i)^{-1}p=o(n^{-1/4})$, by Lemma \ref{lem: hat theta} we have 
		\begin{equation}\label{hat lambda ratio}
			\begin{split}
				\frac{\hat\lambda_i}{\tilde\lambda_i}=&\frac{\hat\lambda_i}{\theta_i}\times\frac{\theta_i}{\tilde\lambda_i}=\frac{1+\frac{1}{n}\bigg(\bgamma_i^\top\Zb\Wb\Zb^\top\bgamma_i-n\bigg)+o_{p^*}(n^{-1/2})}{1+\frac{1}{n}\bigg(\bgamma_i^\top\Zb\Zb^\top\bgamma_i-n\bigg)+o(n^{-1/2})}\\
				=&\frac{\frac{1}{n}\bgamma_i^\top\Zb(\Wb-\Ib)\Zb^\top\bgamma_i}{1+\frac{1}{n}\bigg(\bgamma_i^\top\Zb\Zb^\top\bgamma_i-n\bigg)+o(n^{-1/2})}+1+o_{p^*}(n^{-1/2})\\
				=&[1+o(1)]\times \frac{1}{n}\bgamma_i^\top\Zb(\Wb-\Ib)\Zb^\top\bgamma_i+1+o_{p^*}(n^{-1/2}),
			\end{split}
		\end{equation}
		with probability tending to one conditional on sample. 
		Therefore, it remains to find the limiting distribution of 
		\[
		n^{-1/2}\bgamma_i^\top\Zb(\Wb-\Ib)\Zb^\top\bgamma_i=\sqrt{n}\sum_{j=1}^nu_{ij}^2(w_j-1),
		\] 
		conditional on $\Xb$. To this end, we handle the multiplier and standard bootstrap separately.
		
		For the multiplier bootstrap, $w_j$'s are from i.i.d. $Exp(1)$ so that 
		\[
		\begin{split}
			&\mathbb{E}\bigg(\sqrt{n}\sum_{j=1}^nu_{ij}^2(w_j-1)\mid \Xb\bigg)=0,\\
			&\mathbb{E}\bigg([\sqrt{n}\sum_{j=1}^nu_{ij}^2(w_j-1)]^2\mid \Xb\bigg)=n\sum_{j=1}^nu_{ij}^4=n\tilde\sigma_i^2+o(1)=n\tilde\sigma_i^2[1+o(1)]\rightarrow \xi_i+1,
		\end{split}
		\] 
		with probability tending to 1. It remains to verify the  Lindeberg condition. For any $\epsilon>0$,
		\[
		\begin{split}
			&\frac{1}{n\tilde\sigma_i^2}\sum_{j=1}^n\mathbb{E}\bigg([\sqrt{n}u_{ij}^2(w_j-1)]^2I[|\sqrt{n}u_{ij}^2(w_j-1)|>\epsilon \tilde\sigma_i\sqrt{n}]\mid \Xb\bigg)\\
			\le &\frac{1}{n\tilde\sigma_i^2}\sum_{j=1}^n\frac{\mathbb{E}\big([\sqrt{n}u_{ij}^2(w_j-1)]^4\mid\Xb\big)}{[\epsilon\tilde\sigma_i\sqrt{n}]^4}\le O\bigg(\frac{1}{n^2}\sum_{j=1}^n(\sqrt{n}u_{ij})^8\bigg)\rightarrow 0,
		\end{split}
		\]
		with probability tending to 1, where we use the fact  that $\mathbb{E}(\sqrt{n}u_{ij})^8\le C$ under the bounded eighth moment condition. Then, the Lindeberg condition is satisfied and 
		\[
		\frac{\hat\lambda_i/\tilde\lambda_i-1}{\tilde\sigma}\overset{d^*}{\longrightarrow}\mathcal{N}(0,1).
		\]
		
		For the standard bootstrap, $w_j$'s are from multinomial distribution so that 
		\[
		\begin{split}
			\mathbb{E}\bigg(\sqrt{n}\sum_{j=1}^nu_{ij}^2(w_j-1)\mid \Xb\bigg)=&0,\\
			\mathbb{E}\bigg([\sqrt{n}\sum_{j=1}^nu_{ij}^2(w_j-1)]^2\mid \Xb\bigg)=&n\sum_{j=1}^nu_{ij}^4-\sum_{j_1=1}^nu_{ij_1}^2\sum_{j_2\ne j_1}^nu_{ij_2}^2\\
			=&(n-1)\tilde\sigma_i^2-(\sum_{j=1}^nu_{ij}^2)^2+o(1)=[n\tilde\sigma_i^2-1][1+o(1)]\rightarrow \xi_i,
		\end{split}
		\] 
		with probability tending to 1. Now we verify the  Lindeberg condition. We can write $\bw=(w_1,\ldots,w_n)^\top$ as  a sum of $n$ independent random vectors, i.e.,
		\[
		\bw=\bw_1+\cdots\bw_n,
		\]
		where each $\bw_l=(w_{l1},\ldots,w_{ln})^\top$ is $n$-dimensional vector following $n$-dimensional multinomial distribution with $1$ trial and event probability $(n^{-1},\ldots,n^{-1})$. Then,
		\[
		w_j=\sum_{l=1}^nw_{lj}\Longrightarrow \sqrt{n}\sum_{j=1}^nu_{ij}^2w_j=\sum_{j=1}^n\sum_{l=1}^n\sqrt{n} u_{ij}^2w_{lj}=\sum_{l=1}^n\bigg(\sum_{j=1}^n\sqrt{n} u_{ij}^2w_{lj}\bigg),
		\]
		where $w_{lj}$ is the $j$-th entry of $\bw_l$. Fix $i$ and let $Q_{nl}=\sum_{j=1}^n\sqrt{n} u_{ij}^2w_{lj}$.
		Then $Q_{n1},\ldots, Q_{nn}$ are independent and it suffices to verify  
		\[
		\sum_{l=1}^n\mathbb{E}\bigg(Q_{nl}^2 I(|Q_{nl}|>\epsilon)\mid \Xb\bigg)\rightarrow 0,
		\]
		for any $\epsilon>0$. By the definition of $Q_{nl}$,
		\[
		\mathbb{P}^*\bigg(Q_{nl}=\sqrt{n}u_{ij}^2\bigg)=n^{-1}\quad \text{for any}\quad 1\le j\le n.
		\]
		Therefore, 
		\[
		\sum_{l=1}^n\mathbb{E}\bigg(Q_{nl}^4\mid \Xb\bigg)=\sum_{l=1}^n\frac{1}{n}\sum_{j=1}^n(\sqrt{n}u_{ij}^2)^4\rightarrow 0,
		\]
		with probability tending to 1. Similarly to the proof under multiplier bootstrap, this verifies the Lindeberg condition and concludes the first part of the theorem.
		
		Next, if $n^{-1/4}=o[(n\lambda_i)^{-1}p]$, similarly to (\ref{hat lambda ratio}) we will have
		\[
		\begin{split}
			\frac{\hat\lambda_i}{\tilde\lambda_i}=\bigg(\hat\zeta_i-\frac{\theta_i}{\lambda_i}\bigg)\times [1+o_p(1)]+1.
		\end{split}
		\]
		Therefore, by Lemma \ref{lem: hat theta},
		\[
		\sqrt{n}\bigg(\frac{\hat\lambda_i}{\tilde\lambda_i}-1\bigg)=\sqrt{n}\times \bigg(\frac{\text{tr}\bLambda_2}{n\lambda_i}\bigg)^2\times\mathbb{E}[w_1^2(w_1-1)]\times [1+o_p(1)]\rightarrow \infty,
		\]
		which concludes the theorem.
	\end{proof}
	
	\subsection{Proof of Corollary \ref{cor: bootstrap bias}  : bias of bootstrap}\label{sec: bootstrap bias }
	\begin{proof}
		For the bootstrapped eigenvalues $\hat\lambda_i$, $1\le i\le r$, similarly to (\ref{hat lambda ratio}) we always have
		\[
		\begin{split}
			\frac{\hat\lambda_i}{\tilde\lambda_i}
			=&[1+o_p(1)]\times \bigg[\frac{1}{n}\bgamma_i^\top\Zb(\Wb-\Ib)\Zb^\top\bgamma_i-\frac{1}{n}\text{tr}\Wb+1+\hat\zeta_i-\frac{\theta_i}{\lambda_i}\bigg]+1+o_{p}(n^{-1/2}).
		\end{split}
		\]
		Therefore, under the standard bootstrap, by Lemma \ref{lem: hat theta}, we always have
		\[
		\begin{split}
			&\mathbb{P}^*\bigg(\sqrt{n}\times\lambda_i^{-1}(\hat\lambda_i-\tilde\lambda_i)\le s\bigg)\\
			=&\mathbb{P}^*\bigg(\frac{1}{\sqrt{n}}\bgamma_i^\top\Zb(\Wb-\Ib)\Zb^\top\bgamma_i\le s\frac{\lambda_i}{\tilde\lambda_i}[1+o_p(1)]+\frac{1}{\sqrt{n}}(\text{tr}\Wb-n)-\sqrt{n}(\hat\zeta_i-\frac{\theta_i}{\lambda_i})\bigg)\\
			=&\mathbb{P}^*\bigg(\frac{1}{\sqrt{n}}\bgamma_i^\top\Zb(\Wb-\Ib)\Zb^\top\bgamma_i\le s[1+o_p(1)]-\sqrt{n}\times \frac{\text{tr}^2\bLambda_2}{(n\lambda_i)^2}\times \mathbb{E}[w_1^2(w_1-1)]\times [1+o_p(1)]\bigg)\\
			=&F_G\bigg(s\xi_i^{-1/2}-\xi_i^{-1/2}\sqrt{n}\times \frac{\text{tr}^2\bLambda_2}{(n\lambda_i)^2}\times \mathbb{E}[w_1^2(w_1-1)]\bigg) +o_p(1).
		\end{split}
		\]
		
		On the other hand, for $\tilde\lambda_i$ (without bootstrap), by Lemma \ref{lem: without bootstrap}(a), we have
		\[
		\sqrt{n}\bigg(\frac{\tilde\lambda_i}{\theta_i}-1\bigg)=\sqrt{n}\sum_{j=1}^n(u_{ij}^2-\mathbb{E}u_{ij}^2)+o_p(1)=\frac{1}{\sqrt{n}}\sum_{j=1}^n\bigg((\sqrt{n}u_{ij})^2-\mathbb{E}(\sqrt{n}u_{ij})^2\bigg)+o_p(1).
		\]
		Moreover, by Lemma \ref{lem: ui} we will  have
		\[
		\mathbb{E}\bigg((\sqrt{n}u_{ij})^2-\mathbb{E}(\sqrt{n}u_{ij})^2\bigg)^2=\xi_i+o(1).
		\]
		Therefore, it's not hard to verify that 
		\begin{equation}\label{distribution tilde lambda}
			\sqrt{n}\bigg(\frac{\tilde\lambda_i}{\theta_i}-1\bigg)\overset{d}{\longrightarrow }\mathcal{N}(0,\xi_i).
		\end{equation}
		Further, by the definition of $\theta_i$, we have
		\[
		\begin{split}
			\frac{\theta_i}{\lambda_i}-1=\frac{\text{tr}\bLambda_2}{n\theta_i}\times [1+o(1)]\Longrightarrow 1-\frac{\lambda_i}{\theta_i}=\frac{\text{tr}\bLambda_2}{n\theta_i}\times [1+o(1)].
		\end{split}
		\]
		Therefore, 
		\[
		\begin{split}
			\mathbb{P}\bigg(\sqrt{n}\times\lambda_i^{-1}(\tilde\lambda_i-\lambda_i)\le s\bigg)=&\mathbb{P}\bigg(\sqrt{n}\times\frac{\theta_i}{\lambda_i}\big(\frac{\tilde\lambda_i}{\theta_i}-1+1-\frac{\lambda_i}{\theta_i}\big)\le s\bigg)\\
			=&\mathbb{P}\bigg(\sqrt{n}\times\big(\frac{\tilde\lambda_i}{\theta_i}-1\big)\le s\frac{\lambda_i}{\theta_i}-\sqrt{n}\frac{\text{tr}\bLambda_2}{n\theta_i}\times [1+o(1)]\bigg)\\
			=&F_G\bigg(s\xi_i^{-1/2}-\xi_i^{-1/2}\sqrt{n}\frac{\text{tr}\bLambda_2}{n\lambda_i}\bigg)+o(1),
		\end{split}
		\]
		which concludes the corollary.
	\end{proof}

	\section{Technical lemmas for the proof in Section \ref{sec: spiked}}\label{sec: technical}
	{\begin{lemma}\label{lem: two bounds}
			Under Assumptions \ref{c1} and \ref{c2}, For any $1\le i\ne k\le r$, we have
			\[
			\mathbb{E}(n^{-1}\bgamma_i^\top\Zb\Zb^\top\bgamma_i-1)^2\le O(n^{-1}),\quad \mathbb{E}(n^{-1}\bgamma_i^\top\Zb\Zb^\top\bgamma_k)^2\le O(n^{-1}).
			\]
		\end{lemma}
		\begin{proof}
			Write $\bgamma_i^\top=(\bgamma_{i1}^\top,\bgamma_{i2}^\top)$, where $\bgamma_{i1}$ is composed of the first $r$ entries. Then, $\bgamma_i^\top\Zb=\bgamma_{i1}^\top\Fb^\top+\bgamma_{i2}^\top\Eb$, and 
			\[
			\begin{split}
				n^{-1}\bgamma_i^\top\Zb\Zb^\top\bgamma_i-1=&\bgamma_{i1}^\top(n^{-1}\Fb^\top\Fb-\Ib)\bgamma_{i1}+\bgamma_{i2}^\top(n^{-1}\Eb\Eb^\top-\Ib)\bgamma_{i2}+2n^{-1}\bgamma_{i1}^\top\Fb^\top\Eb^\top\bgamma_{i2}.
			\end{split}
			\]
			By the independence of the entries in $\Eb$, 	for any $\bgamma_{i2}$, we have
			\[
			\begin{split}
				&\mathbb{E}(\bgamma_{i2}^\top\bepsilon_j)^2=\|\bgamma_{i2}\|^2,\quad \mathbb{E}(n^{-1}\bgamma_{i2}^\top\Eb\Eb^\top\bgamma_{i2})=\|\bgamma_{i2}\|^2,\\
				&\mathbb{E}(n^{-1}\bgamma_{i2}^\top\Eb\Eb^\top\bgamma_{i2}-\|\bgamma_{i2}\|^2)^2=\mathbb{E}\bigg(\frac{1}{n}\sum_{j=1}^n[(\bgamma_{i2}^\top\bepsilon_j)^2-\|\bgamma_{i2}\|^2]\bigg)^2\le O(n^{-1}).
			\end{split}
			\]
			On the other hand, because $\Fb=\Cb\Fb^0$ with the entries of $\Fb^0$ being independent and $n^{-1}\text{tr}(\Cb^\top\Cb)=1$,  by elementary moment calculations  we have
			\[
			\mathbb{E}[\bgamma_{i1}^\top(n^{-1}\Fb^\top\Fb-\Ib)\bgamma_{i1}]^2\le O(n^{-2}\|\Cb^\top\Cb\|_F^2)\le O(n^{-1}\|\Cb\|^4\le O(n^{-1}).
			\]
			Further by the independence between $\Fb$ and $\Eb$, the intersection term will also be asymptotic negligible. The proof for $i\ne k$ is similar and omitted. Then,  the lemma holds.
		\end{proof}
		
		\begin{lemma}\label{lem: large deviation}
			Under Assumption \ref{c1}, for any $1\le j\le n$ and $(r+p)\times (r+p)$ deterministic symmetric matrix $\Ab$, we have
			\[
			\mathbb{E}(\bz_j^\top\Ab\bz_j-\text{tr}\Ab)^4\le O(\|\Ab\|_F^4).
			\]
		\end{lemma}
		\begin{proof}
			Write ${\bz_j^0}^\top=({\bbf_j^0}^\top,\bepsilon_j^\top)$, where $\bbf_j^0$ is the $j$-th row vector of $\Fb^0$. Let $\Ab_{ff}$ be the left-top $r\times r$ block of $\Ab$. Then,
			\[
			\begin{split}
				\mathbb{E}(\bz_j^\top\Ab\bz_j-\text{tr}\Ab)^4\lesssim& \mathbb{E}({\bz_j^0}^\top\Ab\bz_j^0-\text{tr}\Ab)^4+\mathbb{E}({\bbf_j^0}^\top\Ab_{ff}\bbf_j^0)^4+\mathbb{E}({\bbf_j}^\top\Ab_{ff}\bbf_j)^4\\
				\le &O(\|\Ab\|_F^4)+\mathbb{E}(\|\bbf_j^0\|^8+\|\bbf_j\|^8)\times O(\|\Ab\|^4)\le O(\|\Ab\|_F^4),
			\end{split}
			\]
			where in the second line we use the facts that the entries of $\bz_j^0$ are independent with the bounded 8th moments and $r$ is fixed.
		\end{proof}
		
		\begin{lemma}\label{lem: ui}
			Let $\bu_i=n^{-1/2}\Zb^\top\bgamma_i$. Under Assumptions \ref{c1} and \ref{c2}, for $1\le i\le r$ we have 
			\begin{equation}\label{sum uij 4}
				\mathbb{E}\sum_{j=1}^{n}u_{ij}^4=n^{-1}(\xi_i+1),\quad \sum_{j=1}^n u_{ij}^4-\mathbb{E}\sum_{j=1}^{n}u_{ij}^4=o_p(n^{-1}).
			\end{equation}
		\end{lemma}
		\begin{proof}
			By definition, $\sqrt{n}u_{ij}=\bz_j^\top\bgamma_i$, where the entries of $\bz_j$ are independent given $j$. Then, 
			\[
			\begin{split}
				\sum_{j=1}^n\mathbb{E}|\sqrt{n}u_{ij}|^4=&\sum_{j=1}^n\mathbb{E}|\bz_j^\top\bgamma_i|^4=\sum_{j=1}^n\sum_{k_1,k_2,k_3,k_4}^{p+r}\mathbb{E}(\gamma_{ik_1}\gamma_{ik_2}\gamma_{ik_3}\gamma_{ik_4}z_{jk_1}z_{jk_2}z_{jk_3}z_{jk_4})\\
				=&\sum_{j=1}^n\bigg[\sum_{k=1}^{p+r}\gamma_{ik}^4\nu_{jk}+3\sum_{k_1=1}^{p+r}\sum_{k_2\ne k_1}^{p+r}\gamma_{ik_1}^2\gamma_{ik_2}^2\mathbb{E}(z_{jk_1}^2)\mathbb{E}(z_{jk_2}^2)\bigg]\\
				=&\sum_{j=1}^n\bigg\{\sum_{k=1}^{p+r}\gamma_{ik}^4[\nu_{jk}-3(\mathbb{E}z_{jk}^2)^2]+3[\sum_{k=1}^{p+r}\gamma_{ik}^2\mathbb{E}(z_{jk}^2)]^2\bigg\}=n(\xi_i+1),
			\end{split}
			\]
			which proves the first result in (\ref{sum uij 4}). For the second, define the conditional expectation $\mathbb{E}_h(\cdot)=\mathbb{E}(\cdot \mid \bbf_1^0,\bepsilon_1,\ldots,\bbf_h^0,\bepsilon_h)$. Then,
			\begin{equation}\label{martingale u}
				\begin{split}
					\sum_{j=1}^n u_{ij}^4-\mathbb{E}\sum_{j=1}^{n}u_{ij}^4=\sum_{h=1}^n(\mathbb{E}_h-\mathbb{E}_{h-1})\sum_{j=1}^n u_{ij}^4.
				\end{split}
			\end{equation}
			Write $\bgamma_i^\top=(\bgamma_{i1}^\top,\bgamma_{i2}^\top)$ as in Lemma \ref{lem: two bounds}. Then,
			\[
			\begin{split}
				&(\sqrt{n} u_{ij})^4= |\bbf_j^\top\bgamma_{i1}|^4+4(\bbf_j^\top\bgamma_{i1})^3\bepsilon_j^\top\bgamma_{i2}+6(\bbf_j^\top\bgamma_{i1})^2(\bepsilon_j^\top\bgamma_{i2})^2+4(\bbf_j^\top\bgamma_{i1})(\bepsilon_j^\top\bgamma_{i2})^3+|\bepsilon_j^\top\bgamma_{i2}|^4.
			\end{split}
			\]
			Write $\vartheta_{jh}=\sum_{l\ne h}C_{jl}\bgamma_{i1}^\top\bbf_l^0$. Then, 
			\[
			\begin{split}
				&|\bbf_j^\top\bgamma_{i1}|^4=(C_{jh}\bgamma_{i1}^\top\bbf_h^0+\vartheta_{jh})^4\\
				=&(C_{jh}\bgamma_{i1}^\top\bbf_h^0)^4+4(C_{jh}\bgamma_{i1}^\top\bbf_h^0)^3\vartheta_{jh}+6(C_{jh}\bgamma_{i1}^\top\bbf_h^0)^2\vartheta_{jh}^2+4(C_{jh}\bgamma_{i1}^\top\bbf_h^0)\vartheta_{jh}^3+\vartheta_{jh}^4
				:=\sum_{l=1}^5\pi_{l,jh}.
			\end{split}
			\]
			Note that $(\mathbb{E}_h-\mathbb{E}_{h-1})\pi_{5,jh}=0$ while 
			\[
			\begin{split}
				\mathbb{E}|(\mathbb{E}_h-\mathbb{E}_{h-1})\sum_{j=1}^n\pi_{1,jh}|^2\lesssim& (\sum_j|C_{jh}|^4)^2\times \mathbb{E}|\bgamma_{i1}^\top\bbf_h^0|^8\le O(1),\\
				\mathbb{E}|(\mathbb{E}_h-\mathbb{E}_{h-1})\sum_{j=1}^n\pi_{2,jh}|^2\lesssim& (\sum_j|C_{jh}|^2)(\sum_j|C_{jh}|^4)\times \mathbb{E}|\bgamma_{i1}^\top\bbf_h^0|^6(\max_j\mathbb{E}\vartheta_{jh}^2)\le O(1),\\
				\mathbb{E}|(\mathbb{E}_h-\mathbb{E}_{h-1})\sum_{j=1}^n\pi_{3,jh}|^2\lesssim& (\sum_j|C_{jh}|^2)(\sum_j|C_{jh}|^2)\times \mathbb{E}|\bgamma_{i1}^\top\bbf_h^0|^4(\max_j\mathbb{E}\vartheta_{jh}^4)\le O(1),\\
				\mathbb{E}|(\mathbb{E}_h-\mathbb{E}_{h-1})\sum_{j=1}^n\pi_{4,jh}|^2\lesssim& \mathbb{E}|\sum_jC_{jh}(\vartheta_{jh}^3-\mathbb{E}\vartheta_{jh}^3)|^2\le \mathbb{E}|\sum_jC_{jh}[(\bgamma_{i1}^\top\bbf_j)^3-\mathbb{E}(\bgamma_{i1}^\top\bbf_j)^3]|^2+O(1).\\
			\end{split}
			\]
			Using the technique in (\ref{martingale u}) again, we have
			\[
			\begin{split}
				&\sum_jC_{jh}[(\bgamma_{i1}^\top\bbf_j)^3-\mathbb{E}(\bgamma_{i1}^\top\bbf_j)^3=\sum_{l=1}^n(\mathbb{E}_l-\mathbb{E}_{l-1})\sum_jC_{jh}[(\bgamma_{i1}^\top\bbf_j)^3\\
				=&\sum_{l=1}^n(\mathbb{E}_l-\mathbb{E}_{l-1})\sum_jC_{jh}\bigg[(C_{jl}\bgamma_{i1}^\top\bbf_l^0)^3+3(C_{jl}\bgamma_{i1}^\top\bbf_l^0)^2\vartheta_{jl}+3(C_{jl}\bgamma_{i1}^\top\bbf_l^0)\vartheta_{jl}^2+\vartheta_{jl}^3\bigg].
			\end{split}
			\]
			Similarly, $(\mathbb{E}_l-\mathbb{E}_{l-1})\vartheta_{jl}^3=0$. Write $\Cb_2=\Cb\odot \Cb$ and $\Cb_3=\Cb_2\odot\Cb$, where $\odot$ stands for the Hadamard product.  Then, after some tedious calculation, we have
			\[
			\begin{split}
				\sum_l\mathbb{E}[\sum_jC_{jh}(C_{jl}\bgamma_{i1}^\top\bbf_l^0)^3]^2\lesssim& \|\bC_{\cdot h}^\top\Cb_3\|^2,\\
				\sum_l\mathbb{E}[\sum_jC_{jh}(C_{jl}\bgamma_{i1}^\top\bbf_l^0)^2\vartheta_{jl}]^2\lesssim& \sum_l\mathbb{E}[\sum_jC_{jh}C_{jl}^2\bbf_j^\top\bgamma_{i1}]^2+\|\bC_{\cdot h}^\top\Cb_3\|^2\\
				\lesssim &\sum_j\sum_lC_{jh}^2C_{jl}^4+\|\bC_{\cdot h}^\top\Cb_3\|^2\le O(1)+\|\bC_{\cdot h}^\top\Cb_3\|^2,\\
				\sum_l\mathbb{E}[(\mathbb{E}_l-\mathbb{E}_{l-1})\sum_jC_{jh}C_{jl}\bgamma_{i1}^\top\bbf_l^0\vartheta_{jl}^2]^2\lesssim& \sum_l\mathbb{E}[(1-\mathbb{E})\sum_jC_{jh}C_{jl}(\bbf_j^\top\bgamma_{i1})^2]^2+\|\bC_{\cdot h}^\top\Cb_3\|^2\\
				\lesssim &\sum_j\sum_lC_{jh}^2C_{jl}^2+\|\bC_{\cdot h}^\top\Cb_3\|^2\le O(1)+\|\bC_{\cdot h}^\top\Cb_3\|^2.
			\end{split}
			\]
			As a result, by Burkholder's inequality, we conclude that
			\[
			\mathbb{E}|(\mathbb{E}_h-\mathbb{E}_{h-1})\sum_{j=1}^n\pi_{4,jh}|^2\lesssim O(1)+\|\bC_{\cdot h}^\top\Cb_3\|^2,
			\]
			which further implies that
			\[
			\mathbb{E}|\sum_{h=1}^n(\mathbb{E}_h-\mathbb{E}_{h-1})|\bbf_j^\top\bgamma_{i1}|^4|^2\lesssim \sum_{h=1}^n\|\bC_{\cdot h}^\top\Cb_3\|^2+O(n)\lesssim \|\Cb_3\|_F^2+O(n)\le O(n).
			\]
			Similarly, we can also prove that
			\[
			\mathbb{E}\bigg|\sum_{h=1}^n(\mathbb{E}_h-\mathbb{E}_{h-1})\sum_{j=1}^n|\bbf_j^\top\bgamma_{i1}|^{k_1}|\bepsilon_j^\top\bgamma_{i2}|^{k_2}\bigg|^2\le O(n),\quad 0\le k_1,k_2\le 4,k_1+k_2=4.
			\]
			Therefore, by (\ref{martingale u}) and Burkholder's inequality, we claim that
			\[
			\mathbb{E}\bigg(\sum_{j=1}^n u_{ij}^4-\mathbb{E}\sum_{j=1}^{n}u_{ij}^4\bigg)^2\le o(n^{-2}),
			\]
			which concludes the lemma.
	\end{proof}}
	
	\begin{lemma}\label{lem: sample covariance}
		Under Assumption \ref{c1},  there exists constant $C>0$ such that 
		\[
		\mathbb{P}\bigg(\|(n\vee p)^{-1}\Eb\Eb^\top\|>C\bigg)\rightarrow 0.
		\]
		If further $\max_{i,j}|\epsilon_{ij}|\le (np)^{1/4-c}$ for some small constant $c>0$, we have
		\[
		\mathbb{P}\bigg(\|(n\vee p)^{-1}\Eb\Eb^\top\|>C\bigg)\le (n\vee p)^{-d},
		\]
		for any constant $d>0$.
	\end{lemma}
	\begin{proof}
		We start with the case $p=n$.  The result follows directly from random matrix theory on the largest eigenvalue of sample covariance matrix, see for example Theorem 2.7 and Theorem 3.15 in \cite{ding2018necessary}.
		For $p>n$, we can always find some $p\times (p-n)$ matrix $\Eb_+$ so that $(\Eb,\Eb_+)$ is one $p\times p$ matrix satisfying all the assumptions under  the case $n=p$. Then, the result still holds.
		For $p<n$, it's parallel by transposing $\Eb$.
	\end{proof}

	\begin{lemma}\label{lema2}
		Under Assumptions \ref{c1} and \ref{c2}, for the decomposition in (\ref{det}) we have
		\[
		\begin{split}
			&n^{-1}\bGamma_1^\top\Zb\Wb^{1/2}\Kb(\theta_1)\Wb^{1/2}\Zb^\top\bGamma_1\\
			=&\frac{1}{n}\bGamma_1^\top\Zb\Wb\Zb^\top\bGamma_1-\bigg(\frac{1}{n}\text{tr}\Wb\bigg)\Ib_r+\hat\zeta_1\Ib_r+O_p\bigg(\frac{1}{\sqrt{n}}\times \frac{(n\vee p)\log n}{n\theta_1}+\frac{1}{n}\bigg),
		\end{split}
		\]
		where the $O_p$ is under Frobenius norm.
	\end{lemma}
	\begin{proof}
		{
			To ease notation,  in the following proof, $o_p(n^{-1/2})$ stands for 
			\[
			O_p\bigg(\frac{1}{\sqrt{n}}\times \frac{(n\vee p)\log n}{n\theta_1}+\frac{1}{n}\bigg).
			\]
			\noindent\textbf{Step 1: truncation.}
			
			We need to truncate the entries of $\Eb=(\epsilon_{ij})_{p\times n}$ by defining
			\[
			\epsilon_{ij}^*=\epsilon_{ij}I(|\epsilon_{ij}|<(np)^{1/4-c}),\quad \Eb^*=( \epsilon_{ij}^*),
			\]
			for some small constant $c>0$. Then, because $\mathbb{E} \epsilon_{ij}=0$, we have
			\[
			\begin{split}
				|\mathbb{E} \epsilon_{ij}^*|=&	|\mathbb{E} \epsilon_{ij}-\mathbb{E} \epsilon_{ij}^*|=	|\mathbb{E} \epsilon_{ij}I(| \epsilon_{ij}|\ge (np)^{1/4- c})|\\
				\le& (np)^{-7/4+7 c}\mathbb{E} \epsilon_{ij}^8I(| \epsilon_{ij}|\ge (np)^{1/4- c})\le C(np)^{-3/2},
			\end{split}
			\]
			for some large constant $C>0$. Similarly, we have $|1-\mathbb{E}( \epsilon_{ij}^*)^2|\le C(np)^{-1}$. 
			Further let $\epsilon_{ij}^{**}=( \epsilon_{ij}^*-\mathbb{E} \epsilon_{ij}^*)/\sqrt{\mathbb{E}( \epsilon_{ij}^*)^2}$ and $\Eb^{**}=(\epsilon_{ij}^{**})$. Then, $ \epsilon_{ij}^{**}$'s are independent random variables with mean 0, variance 1, bounded eighth moment satisfying $| \epsilon_{ij}^{**}|\le (np)^{1/4- c}$. Define $\Zb^*,\Zb^{**},\Kb^*,\Kb^{**}$ by replacing $\Eb$ with $\Eb^*, \Eb^{**}$, respectively. Then, 
			\[
			\begin{split}
				&\mathbb{P}\bigg(	\|n^{-1}\bGamma_1^\top\Zb\Wb^{1/2}\Kb(\theta_1)\Zb^\top\Wb^{1/2}\bGamma_1-n^{-1}\bGamma_1^\top\Zb^*\Wb^{1/2}\Kb^*(\theta_1){\Zb^*}^\top\Wb^{1/2}\bGamma_1\|\ge n^{-2}c\bigg)\\
				\le &\mathbb{P}(\Eb\ne \Eb^*)\le \mathbb{P}\bigg(\max_{i,j}| \epsilon_{ij}|\ge (np)^{1/4- c}\bigg)\le \sum_{i,j}\mathbb{P}\bigg(| \epsilon_{ij}|\ge (np)^{1/4-c}\bigg)\rightarrow 0,
			\end{split}
			\]
			where we use Markov's equality and the bounded eighth moment condition. Then,
			\[
			\|n^{-1}\bGamma_1^\top\Zb\Wb^{1/2}\Kb(\theta_1)\Zb^\top\Wb^{1/2}\bGamma_1-n^{-1}\bGamma_1^\top\Zb^*\Wb^{1/2}\Kb^*(\theta_1){\Zb^*}^\top\Wb^{1/2}\bGamma_1\|=o_p(n^{-2}),
			\]
			so that the error is negligible if we replace $\Eb$ with $\Eb^*$.  More tedious calculations will show that the error is also negligible if we further replace $\Eb^*$ with $\Eb^{**}$. Similar technique has been applied in Section 12 of \cite{cai2020limiting} and we omit the details. 
			As a result, without loss of generality, we can assume that
			$| \epsilon_{ij}|\le (np)^{1/4- c}$ for some small constant $ c>0$ in the proof, which only generates an error term of order $o_p(n^{-1/2})$.

			Next, we provide an upper bound for $\|\Wb\|$. Define an event $\Xi_w=\{\|\Wb\|\le 5\log n\}$. Under standard bootstrap, by (\ref{mg}) and (\ref{Ew1}), as long as $n$ is sufficiently large, we have
			\[
			1-\mathbb{P}(\Xi_w)=\mathbb{P}(\|\Wb\|> 5\log n)= \mathbb{P}\bigg(\exp(\max_jw_j)> n^5\bigg)\le\frac{n(e^{e-1}+ c)}{n^5}=o(n^{-3}).
			\]
			Under multiplier bootstrap, it's similar to conclude that
			\begin{equation}\label{Xiw}
				1-\mathbb{P}(\Xi_w)=o(n^{-3}).
			\end{equation}
			On the other hand, let $\Xi_0=\{\|n^{-1}\Eb\Eb^\top\|\le C(n\vee p)/n\}$. Then, by Lemma \ref{lem: sample covariance}, 
			\begin{equation}\label{Xi0}
				1-\mathbb{P}(\Xi_0)	\le o((n\vee p)^{-4}),\quad n\rightarrow \infty.
			\end{equation}
			Then, it's sufficient to consider  $n^{-1}\bGamma_1^\top\Zb\Wb^{1/2}\Kb(\theta_1)\Zb^\top\Wb^{1/2}\bGamma_1I(\Xi_0)I(\Xi_w)$.  Similarly, let $ \Eb_j$ be the $(p+r)\times n$ matrix by replacing the $j$-th column of $ \Eb$ with 0, and define $ \Eb_j$, ${\mathcal{S}}_{2j}$, $\Kb_j(x)$ accordingly by replacing $ \Eb$ with $ \Eb_j$. Let $\Xi_j=\{\|n^{-1} \Eb_j \Eb_j^\top\|\le C(n\vee p)/n\}$. Then, 
			\begin{equation}\label{Xij}
				1-\mathbb{P}(\Xi_j)	\le o((n\vee p)^{-4}),\quad n\rightarrow \infty, \quad 1\le j\le n. 
			\end{equation}
			In the following, we may take the  events $\Xi_0,\Xi_w,\Xi_j$ as given without further explanation.

			\noindent\textbf{Step 2: replacing $\Fb$ with $\Fb^0$.}
			
			We aim to calculate the error if replacing $\Fb$ with $\Fb^0$.  Define $\check{\Zb}=\Zb\Wb^{1/2}$ and
			\[
			\mathcal{H}(x):=\left(\begin{matrix}
				& x\Ib&n^{-1/2}\check\Zb^\top\bGamma_2\bLambda_2^{1/2}&n^{-1/2}\check\Zb^\top\bGamma_1\\
				&n^{-1/2}\bLambda_2^{1/2}\bGamma_2^\top\check\Zb&\Ib&{\bf 0}\\
				&n^{-1/2}\bGamma_1^\top\check\Zb&{\bf 0}&\Ib_r\\
			\end{matrix}\right)=\left(\begin{matrix}
				&x\Ib&n^{-1/2}\check\Zb^\top\tilde\Ab\\
				&n^{1/2}\tilde\Ab^\top\check\Zb&\Ib_{p+r}
			\end{matrix}\right),
			\]
			where $\tilde\Ab=(\bGamma_2\bLambda_2^{1/2},\bGamma_1)$. Then, by Schur's complement formula, the $(3,3)$-block of $\mathcal{H}^{-1}(\theta_1)$ is exactly equal to the inverse of $		\Ib-\theta_1^{-1}n^{-1}\bGamma_1^\top\Zb\Wb^{1/2}\Kb(\theta_1)\Zb^\top\Wb^{1/2}\bGamma_1$. That is,
			\[
			n^{-1}\bGamma_1^\top\Zb\Wb^{1/2}\Kb(\theta_1)\Zb^\top\Wb^{1/2}\bGamma_1=\theta_1(\Ib-[\mathcal{H}^{-1}(\theta_1)]_{3,3}^{-1}),
			\]
			where the subscript $(3,3)$ indicates a block.  By Shur's complement formula again,
			\[
			\begin{split}
				&[\mathcal{H}^{-1}(\theta_1)]_{n+k,n+l}= [ (\Ib-\theta_1^{-1}n^{-1}\tilde\Ab^\top\check\Zb\check\Zb^\top\tilde\Ab )^{-1} ]_{kl},\quad 1\le k,l\le r+p.
			\end{split}
			\] 
			Write $\tilde\Gb=n^{-1}\tilde\Ab^\top\check\Zb\check\Zb^\top\tilde\Ab$. 
			Define $\check{\Zb}_0$, $\mathcal{H}_0(x)$ and $\tilde\Gb_0$ by replacing $\Fb$ with $\Fb_0$, respectively. Let $\be_k$ be the $(r+p)$-dimensional unit vector with the $k$-th element being 1. By (\ref{matrix inverse}),
			\[
			\begin{split}
				&[ (\Ib-\theta_1^{-1}\tilde\Gb )^{-1} ]_{kl}-[ (\Ib-\theta_1^{-1}\tilde\Gb_0)^{-1} ]_{kl}
				=\theta_1^{-1}\be_k^\top(\Ib-\theta_1^{-1}\tilde\Gb )^{-1}(\tilde\Gb-\tilde\Gb_0)(\Ib-\theta_1^{-1}\tilde\Gb_0)^{-1}\be_l\\
				=&\theta_1^{-1}[(\tilde\Gb-\tilde\Gb_0)_{kl}+\theta_1^{-1}\be_k^\top(\tilde\Gb-\tilde\Gb_0)(\Ib-\theta_1^{-1}\tilde\Gb_0)^{-1}\tilde\Gb_0\be_l\\
				&+\theta_1^{-1}\be_k^\top(\Ib-\theta_1^{-1}\tilde\Gb)^{-1}\tilde\Gb(\tilde\Gb-\tilde\Gb_0)(\Ib-\theta_1^{-1}\tilde\Gb_0)^{-1}\be_l].
			\end{split}
			\]
			Without loss of generality, we let $k=l=p+1$.  Then,
			\[
			\begin{split}
				&\be_k^\top(\tilde\Gb-\tilde\Gb_0)(\Ib-\theta_1^{-1}\tilde\Gb_0)^{-1}\tilde\Gb_0\be_l\\=&\bgamma_1^\top\left[\begin{matrix}
					&n^{-1}(\Fb^\top\Wb\Fb-{\Fb^0}^\top\Wb\Fb^0)&n^{-1}(\Fb-\Fb^0)^\top\Wb\Eb^\top\\
					&n^{-1}\Eb\Wb(\Fb-\Fb^0)&{\bf 0}
				\end{matrix}\right](\Ib-\theta_1^{-1}\tilde\Gb_0)^{-1}\bgamma_1\le O_p(n^{-1/2}),
			\end{split}
			\]
			where we use the facts that $\Fb$, $\Fb^0$ and $\Eb$ are mutually independent and $\theta_1^{-1}\|\tilde\Gb_0\|=o_p(1)$. By similar but more tedious calculations, we can also show that
			\[
			\be_k^\top(\Ib-\theta_1^{-1}\tilde\Gb)^{-1}\tilde\Gb(\tilde\Gb-\tilde\Gb_0)(\Ib-\theta_1^{-1}\tilde\Gb_0)^{-1}\be_l\le O_p(n^{-1/2}).
			\]
			Consequently, we have
			\[
			[ (\Ib-\theta_1^{-1}\tilde\Gb )^{-1} ]_{kl}-[ (\Ib-\theta_1^{-1}\tilde\Gb_0)^{-1} ]_{kl}=\theta_1^{-1}[(\tilde\Gb-\tilde\Gb_0)_{kl}+O_p(\theta_1^{-2}n^{-1/2}),
			\]
			which further indicates that 
			\[
			[\mathcal{H}^{-1}(\theta_1)]_{3,3}-[\mathcal{H}_0^{-1}(\theta_1)]_{3,3}=(n\theta_1)^{-1}\bGamma_1^\top(\check{\Zb}\check{\Zb}^\top-\check{\Zb}_0\check{\Zb}_0^\top)\bGamma+ O_p(\theta_1^{-2}n^{-1/2}).
			\]
			Using (\ref{matrix inverse}) again, we have
			\[
			\begin{split}
				[\mathcal{H}^{-1}(\theta_1)]_{3,3}^{-1}=&[\mathcal{H}_0^{-1}(\theta_1)]_{3,3}^{-1}-(n\theta_1)^{-1}[\mathcal{H}^{-1}(\theta_1)]_{3,3}^{-1}\bGamma_1^\top(\check{\Zb}\check{\Zb}^\top-\check{\Zb}_0\check{\Zb}_0^\top)\bGamma[\mathcal{H}_0^{-1}(\theta_1)]_{3,3}^{-1}+\theta_1^{-1}o_p(n^{-1/2})\\
				=&[\mathcal{H}_0^{-1}(\theta_1)]_{3,3}^{-1}-(n\theta_1)^{-1}\bGamma_1^\top(\check{\Zb}\check{\Zb}^\top-\check{\Zb}_0\check{\Zb}_0^\top)\bGamma+\theta_1^{-1}o_p(n^{-1/2}),
			\end{split}
			\]
			where we use the fact that $	[\mathcal{H}^{-1}(\theta_1)]_{3,3}^{-1}=\Ib+O_p((n\theta_1)^{-1}(n\vee p)\log p)$ and the same rate  for $[\mathcal{H}_0^{-1}(\theta_1)]_{3,3}^{-1}$. Finally, we can conclude that
			\[
			\begin{split}
				n^{-1}\bGamma_1^\top\Zb\Wb^{1/2}\Kb(\theta_1)\Zb^\top\Wb^{1/2}\bGamma_1=&n^{-1}\bGamma_1^\top\Zb_0\Wb^{1/2}\Kb_0(\theta_1)\Zb_0^\top\Wb^{1/2}\bGamma_1\\
				&+n^{-1}\bGamma_1^\top(\check{\Zb}\check{\Zb}^\top-\check{\Zb}_0\check{\Zb}_0^\top)\bGamma+o_p(n^{-1/2}),
			\end{split}
			\]
			where $\Zb_0$ and $\Kb_0$ are obtained by replacing $\Fb$ with $\Fb^0$, respectively. 
			
			In the following, we will replace $\Fb$ with $\Fb^0$, and still write $\Fb$ to ease notation. That is, we assume the entries of $\Fb$ to be independent. The replacement error will be considered later. Moreover, we generalize the definition of events $\Xi_0,\Xi_j$ by replacing $\Eb$ with $\Zb$.
		}
		
		\noindent\textbf{Step 3: first order approximation.}
		
		Let $\bb_1$ and $\bb_2$ be any two columns of $\bGamma_1$. Now we aim to provide asymptotic representation for $\mathcal{L}:=n^{-1}\bb_1^\top\Zb\Wb^{1/2}\Kb(\theta_1)\Wb^{1/2}\Zb^\top\bb_2$,
		where $\Kb(\theta_1)=[\Ib-\theta_1^{-1}\Wb^{1/2}\tilde{\mathcal{S}}_2\Wb^{1/2}]^{-1}$ and $\tilde{\mathcal{S}}_2=n^{-1}\Zb^\top\bSigma_2\Zb$ with $\bSigma_2=\bGamma_2\bLambda_2\bGamma_2^\top$. 
		The idea is to approximate  $\mathcal{L}I(\Xi_0)I(\Xi_w)$ by its conditional expectation $\mathbb{E}[\mathcal{L}I(\Xi_0)I(\Xi_w)\mid \Wb]$ and calculate the error.
		Define the conditional expectations $\mathbb{E}_j=\mathbb{E}(\cdot\mid \bz_1,\ldots,\bz_j,\Wb), 1\le j\le n$, and $\Zb=(\bz_1,\ldots,\bz_n)$. Let $\be_j$ be the $n$-dimensional vector with the $j$-th entry being 1 and the others being 0. Then,
		\begin{equation}\label{step 1: expansion}
			\begin{split}
				&\mathcal{L}I(\Xi_0)I(\Xi_w)-\mathbb{E}[\mathcal{L}I(\Xi_0)I(\Xi_w)\mid \Wb]=\sum_{j=1}^{n}(\mathbb{E}_j-\mathbb{E}_{j-1})\mathcal{L}I(\Xi_0)I(\Xi_w)\\
				=&\frac{1}{n}\sum_{j=1}^{n}(\mathbb{E}_j-\mathbb{E}_{j-1})\bb_1^\top(\Zb_j+ \bz_j\be_j^\top)\Wb^{1/2}\Kb(\theta_1)\Wb^{1/2}(\Zb_j+ \bz_j\be_j^\top)^\top\bb_2 I(\Xi_0)I(\Xi_w)\\
				=&\frac{1}{n}\sum_{j=1}^{n}(\mathbb{E}_j-\mathbb{E}_{j-1})\bb_1^\top\Zb_j\Wb^{1/2}\Kb(\theta_1)\Wb^{1/2}\Zb_j^\top\bb_2 I(\Xi_0)I(\Xi_w)
				\\
				&+\frac{1}{n}\sum_{j=1}^{n}(\mathbb{E}_j-\mathbb{E}_{j-1})\bb_1^\top\Zb_j\Wb^{1/2}\Kb(\theta_1)\Wb^{1/2}\be_j \bz_j^\top\bb_2 I(\Xi_0)I(\Xi_w)
				\\
				&+\frac{1}{n}\sum_{j=1}^{n}(\mathbb{E}_j-\mathbb{E}_{j-1})\bb_1^\top \bz_j\be_j^\top\Wb^{1/2}\Kb(\theta_1)\Wb^{1/2}\Zb_j^\top\bb_2 I(\Xi_0)I(\Xi_w)\\
				&+\frac{1}{n}\sum_{j=1}^{n}(\mathbb{E}_j-\mathbb{E}_{j-1})\bb_1^\top \bz_j\be_j^\top\Wb^{1/2}\Kb(\theta_1)\Wb^{1/2}\be_j \bz_j^\top\bb_2 I(\Xi_0)I(\Xi_w)\\
				:=&\mathcal{I}_1+\mathcal{I}_2+\mathcal{I}_3+\mathcal{I}_4.
			\end{split}
		\end{equation}
		We aim to prove that the error is negligible if we replace $\Kb(\theta_1)$ with $\Ib_n$ in $\mathcal{I}_k$, $k=1,2,3,4$.
		
		We start with $\mathcal{I}_1$.   Under the events $\Xi_w$ and $\Xi_0$ (or $\Xi_j$), we always have $\|\bgamma_j^\top\Zb_j\|^2\le C(n\vee p)$ for any $1\le j\le p$. Therefore, combining with (\ref{Xi0}) and (\ref{Xij}), we have 
		\[
		\begin{split}
			&\bigg|\frac{1}{n}\sum_{j=1}^{n}(\mathbb{E}_j-\mathbb{E}_{j-1})\bb_1^\top\Zb_j\Wb^{1/2}\Kb(\theta_1)\Wb^{1/2}\Zb_j^\top\bb_2 [I(\Xi_0)-I(\Xi_j)]I(\Xi_w)\bigg|\\
			\le& \|\Kb(\theta_1)\|\|\Wb\|\times \frac{1}{n}\sum_{j=1}^n(n\vee p)\bigg(|I(\Xi_0)-1|+|I(\Xi_j)-1|\bigg)\le o_p(n^{-1}),
		\end{split}
		\]
		where we use the bounds in (\ref{Xi0}) and (\ref{Xij}). By the same reason, we will repeatedly exchange $\Xi_0$ and $\Xi_j$ in the proof without further explanation. Such a replacement will add at most a negligible error of order $o_p(n^{-1})$ to $\mathcal{L}$. 
		
		Now return to the definition of $\mathcal{I}_1$ in (\ref{step 1: expansion}). We write
		\[
		\begin{split}
			\mathcal{I}_1=&\frac{1}{n}\sum_{j=1}^{n}(\mathbb{E}_j-\mathbb{E}_{j-1})\bb_1^\top\Zb_j\Wb^{1/2}[\Kb(\theta_1)-\Kb_j(\theta_1)]\Wb^{1/2}\Zb_j^\top\bb_2 I(\Xi_0)I(\Xi_w)\\
			&+\frac{1}{n}\sum_{j=1}^{n}(\mathbb{E}_j-\mathbb{E}_{j-1})\bb_1^\top\Zb_j\Wb^{1/2}\Kb_j(\theta_1)\Wb^{1/2}\Zb_j^\top\bb_2 I(\Xi_j)I(\Xi_w)+o_p(n^{-1}).
		\end{split}
		\]
		The leading term in the second line is actually $0$ because the expectations under $\mathbb{E}_j$ and $\mathbb{E}_{j-1}$ are equal.  
		For the first line, recall that $\tilde{\mathcal{S}}_2=n^{-1}\Zb^\top\bSigma_2\Zb$ while $\Zb=\Zb_j+\bz_j\be_j^\top$. Then, 
		\begin{equation}\label{Rk}
			\begin{split}
				\tilde{\mathcal{S}}_2
				=&\tilde{\mathcal{S}}_{2j}+\frac{1}{n}\Zb_j^\top\bSigma_2\bz_j\be_j^\top+(\frac{1}{n}\Zb_j^\top\bSigma_2\bz_j\be_j^\top)^\top+\frac{1}{n}\be_j\bz_j^\top\bSigma_2\bz_j\be_j^\top\\
				:=&\tilde{\mathcal{S}}_{2j}+\Rb_1+\Rb_2+\Rb_3.
			\end{split}
		\end{equation}
		Then, by the matrix inverse  formula in (\ref{matrix inverse}), we have
		\begin{equation}\label{G0x}
			\Kb(x)=\Kb_j(x)+x^{-1}\Kb(x)\Wb^{1/2}(\Rb_1+\Rb_2+\Rb_3)\Wb^{1/2}\Kb_j(x),
		\end{equation}
		which implies that
		\[
		\begin{split}
			&\mathcal{I}_1-o_p(n^{-1})=\sum_{k=1}^3\mathcal{I}_{1k}\\
			:=&\sum_{k=1}^3\frac{1}{\theta_1}\sum_{j=1}^{n}(\mathbb{E}_j-\mathbb{E}_{j-1})\bb_1^\top\Zb_j\Wb^{1/2}\Kb(\theta_1)\Wb^{1/2}\Rb_k\Wb^{1/2}\Kb_j(\theta_1)\Wb^{1/2}\Zb_j^\top\bb_2 I(\Xi_0)I(\Xi_w).
		\end{split}
		\]
		It suffices to calculate $\mathcal{I}_{1k}$ for $k=1,2,3$. 
		
		By definition, the $j$th row and column entries of $\tilde{\mathcal{S}}_{2j}$ are all equal to $0$. Then, the $j$th row and column vectors of $\Kb_j(\theta_j)$ are equal to $\be_j$. As a result,
		\begin{equation}\label{ej Zj}
			\be_j^\top\Wb^{1/2}\Kb_j(\theta_1)\Wb^{1/2}=w_j\be_j^\top\Longrightarrow \be_j^\top\Wb^{1/2}\Kb_j(\theta_1)\Wb^{1/2}\Zb_j^\top=w_j\be_j^\top\Zb_j^\top={\bf 0},
		\end{equation}
		which implies that $\mathcal{I}_{11}=\mathcal{I}_{13}=0$. 
		
		It suffices to calculate $\mathcal{I}_{12}$.  To this end, we introduce some notation. Define
		\begin{equation}\label{notation a}
			\begin{split}
				&a_{\bb_1\be_j}=\bb_1^\top\Zb_j\Wb^{1/2}\Kb(\theta_1)\Wb^{1/2}\be_j,\quad a_{\bz_j\bz_j}=\bz_j^\top\bSigma_2\Zb_j\Wb^{1/2}\Kb_j(\theta_1)\Wb^{1/2}\Zb_j^\top\bSigma_2\bz_j,\\
				&\bar a_{\bz_j\bz_j}=\text{tr}[\bSigma_2\Zb_j\Wb^{1/2}\Kb_j(\theta_1)\Wb^{1/2}\Zb_j^\top\bSigma_2],\\
				&a_{\bb_1\bz_j}=\bb_1^\top\Zb_j\Wb^{1/2}\Kb(\theta_1)\Wb^{1/2}\Zb_j^\top\bSigma_2\bz_j,\quad  a_{\bb_1\bz_j}^*=\bb_1^\top\Zb_j\Wb^{1/2}\Kb_j(\theta_1)\Wb^{1/2}\Zb_j^\top\bSigma_2\bz_j.
			\end{split}
		\end{equation}
		Use (\ref{Rk}) and (\ref{ej Zj}) again so that
		\[
		\begin{split}
			a_{\bb_1\be_j}=&\frac{w_j}{n\theta_1}(a_{\bb_1\bz_j}+\bz_j^\top\bSigma_2\bz_ja_{\bb_1\be_j}),\quad a_{\bb_1\bz_j}=a_{\bb_1\bz_j}^*+\frac{1}{n\theta_1}a_{\bb_1\be_j}a_{\bz_j\bz_j}.
		\end{split}
		\]
		Therefore, 
		\begin{equation}\label{alpha_k}
			a_{\bb_1\be_j}=\alpha_j^{-1}\times\frac{w_j}{n\theta_1}a_{\bb_1\bz_j}^*,\text{ with }
			\alpha_j=1-\frac{w_j}{(n\theta_1)^2}a_{\bz_j\bz_j}-\frac{w_j}{n\theta_1}\bz_j^\top\bSigma_2\bz_j.
		\end{equation}
		Under the events $\Xi_0$ and $\Xi_w$,  we have $\|\Kb_j(\theta_1)\|\le O(1)$, and further
		\[
		|a_{\bz_j\bz_j}|\le C\|\bz_j^\top\bSigma_2\Zb_j\|^2\log n\le C(n\vee p)^2\log n,\quad |\bz_j^\top\bSigma_2\bz_j|\le C(n\vee p).
		\]
		Therefore, $|\alpha_j^{-1}|I(\Xi_0)I(\Xi_w)\le O(1)$ uniformly over $1\le j\le n$. By Burkholder's equality,
		\[
		\begin{split}
			\mathbb{E}|\mathcal{I}_{12}|^2\le &C\sum_{j=1}^n\mathbb{E}\bigg|\frac{1}{n^2\theta_1}a_{\bb_1\be_j}a_{\bb_2\bz_j}^*I(\Xi_0)I(\Xi_w)\bigg|^2=C\sum_{j=1}^n\mathbb{E}\bigg|\frac{\alpha_j^{-1}w_j}{n^3\theta_1^2}a_{\bb_1\bz_j}^*a_{\bb_2\bz_j}^*I(\Xi_0)I(\Xi_w)\bigg|^2\\
			\le & \frac{C\log^4 n}{n^6\theta_1^4}\times \sum_{j=1}^n\mathbb{E}w_j^2\|\bb_1^\top\Zb_j\|^2\|\Zb_j\|^2\|\bb_2^\top\Zb_j\|^2I(\Xi_0)I(\Xi_w)\le \frac{C}{n}\bigg(\frac{(n\vee p)\log n}{n\theta_1}\bigg)^4,
		\end{split}
		\]
		where the second line is by Lemma \ref{lem: large deviation} and the dependence between $\bz_j$ and $\Zb_j$.  Therefore, we conclude that $	\mathcal{I}_{1}=o_p(n^{-1/2})$. In other words,
		\[
		\mathcal{I}_{1}-\frac{1}{n}\sum_{j=1}^{n}(\mathbb{E}_j-\mathbb{E}_{j-1})\bb_1^\top\Zb_j\Wb\Zb_j^\top\bb_2 I(\Xi_0)I(\Xi_w)=o_p(n^{-1/2}).
		\]
		because the error is negligible to replace $\Xi_j$ with $\Xi_0$ while
		\[
		\frac{1}{n}\sum_{j=1}^{n}(\mathbb{E}_j-\mathbb{E}_{j-1})\bb_1^\top\Zb_j\Wb\Zb_j^\top\bb_2 I(\Xi_j)I(\Xi_w)=0.
		\] 
		Consequently, the error is negligible after  replacing $\Kb(\theta_1)$ in $\mathcal{I}_1$ with $\Ib_n$. 
		
		For $\mathcal{I}_2$, $\mathcal{I}_3$ and $\mathcal{I}_4$, the proof strategy is similar and omitted here. We refer to the proof of Theorem 2.4 in \cite{cai2020limiting} for the details. So we conclude directly that
		\[
		\begin{split}
			&\mathcal{I}_2-\frac{1}{n}\sum_{j=1}^{n}(\mathbb{E}_j-\mathbb{E}_{j-1})\bb_1^\top\Zb_j\Wb\be_j\bz_j^\top\bb_2 I(\Xi_0)I(\Xi_w)\le o_p(n^{-1/2}),\\
			&\mathcal{I}_3-\frac{1}{n}\sum_{j=1}^{n}(\mathbb{E}_j-\mathbb{E}_{j-1})\bb_1^\top\bz_j\be_j^\top\Wb\Zb_j^\top\bb_2 I(\Xi_0)I(\Xi_w)\le  o_p(n^{-1/2}),\\
			&\mathcal{I}_4-\frac{1}{n}\sum_{j=1}^{n}(\mathbb{E}_j-\mathbb{E}_{j-1})\bb_1^\top\bz_j\be_j^\top\Wb\be_j\bz_j^\top\bb_2 I(\Xi_0)I(\Xi_w)\le o_p(n^{-1/2}).
		\end{split}
		\]
		Therefore,  we have
		\[
		\begin{split}
			&\mathcal{L}I(\Xi_0)I(\Xi_w)-\mathbb{E}[\mathcal{L}I(\Xi_0)I(\Xi_w)\mid \Wb]\\
			=&\frac{1}{n}\sum_{j=1}^{n}(\mathbb{E}_j-\mathbb{E}_{j-1})\bb_1^\top\Zb\Wb\Zb^\top\bb_2I(\Xi_0)I(\Xi_w)+o_p(n^{-1/2})\\
			=&\frac{1}{n}\bb_1^\top\Zb\Wb\Zb^\top\bb_2-\frac{1}{n}\mathbb{E}(\bb_1^\top\Zb\Wb\Zb^\top\bb_2\mid \Wb)+o_p(n^{-1/2})\\
			=&\frac{1}{n}\bb_1^\top\Zb\Wb\Zb^\top\bb_2-\frac{1}{n}\bb_1^\top\bb_2\text{tr}(\Wb)+o_p(n^{-1/2}).
		\end{split}
		\]
		The limiting distributions are now much easier to derive because $\Kb(\theta_1)$ has been removed. However, 
		we still need to calculate $\mathbb{E}[\mathcal{L}I(\Xi_0)I(\Xi_w)\mid\Wb]$, which is organized below.
		
		\noindent\textbf{Step 4: replacing $\Zb$ with Gaussian variables.}
		
		In this step, we aim to show that
		\begin{equation}\label{step 2: target}
			\mathbb{E}[\mathcal{L}I(\Xi_0)I(\Xi_w)\mid \Wb]=\mathbb{E}[\mathcal{L}^0I(\Xi_0^0)I(\Xi_w)\mid \Wb]+o_p(n^{-1/2}),
		\end{equation}
		where $\mathcal{L}^0$, $\Xi_0^0$ are defined similarly to $\mathcal{L}$, $\Xi_0$ by replacing the entries in $\Zb$ with i.i.d. standard Gaussian variables, respectively.
		
		The key technique is the Lindeberg's replacement strategy. Let $\Zb^0$ be a $p\times n$ random matrix independent of $\Zb$ and composed of i.i.d. standard Gaussian variables. Further let 
		\[
		\Zb^k=(\bz_1,\ldots,\bz_k,\bz_{k+1}^0,\ldots,\bz_n^0),
		\] 
		which is composed of the leading $k$ columns of $\Zb$ and the last $(n-k)$ columns of $\Zb^0$.  Define $\Zb_j^k$ as the $p\times n$ matrix by replacing the $j$-th column of $\Zb^k$ with ${\bf 0}$. Similarly, define $\Zb^k$, $\Kb^k(x)$, ${\mathcal{S}}_2^k$ and $\mathcal{L}^k$ by replacing $\Zb$ with $\Zb^k$, and $\Zb_j^k$, $\Kb_j^k(x)$, ${\mathcal{S}}_{2j}^k$ by replacing  $\Zb_j$ with $\Zb_j^k$ accordingly.  Define the events $\Xi_0^0=\{\|n^{-1}\Zb^0(\Zb^0)^\top\|\le C(n\vee p)/n\}$ and $\Xi_j^0=\{\|n^{-1}\Zb_j^0(\Zb_j^0)^\top\|\le C(n\vee p)/n\}$.
		
		Similarly to \textbf{Step 3}, it suffices to consider
		\[
		\begin{split}
			&\mathbb{E}_0[\mathcal{L}I(\Xi_0)I(\Xi_w)]-\mathbb{E}_0[\mathcal{L}^0I(\Xi_0^0)I(\Xi_w)]\\
			=&\sum_{k=1}^n\bigg(\mathbb{E}_0[\mathcal{L}^kI(\Xi_0)I(\Xi_0^0)I(\Xi_w)]-\mathbb{E}_0[\mathcal{L}^{k-1}I(\Xi_0)I(\Xi_0^0)I(\Xi_w)]\bigg)+o_p(n^{-1/2}).
		\end{split}
		\]
		To ease notation, in this step we take $\Wb$ as given and assume $I(\Xi_w)=1$. 
		Expanding $\mathcal{L}^k$ using the same technique as in (\ref{step 1: expansion}) by writing $\Zb^k=\Zb_k^k+\bz_k\be_k^\top$, we have
		\[
		\begin{split}
			&\mathbb{E}_0[\mathcal{L}^kI(\Xi_0)I(\Xi_0^0)]\\
			=&\frac{1}{n}\mathbb{E}_0\bb_1^\top\Zb_k^k\Wb^{1/2}\Kb_0^k(\theta_1)\Wb^{1/2}(\Zb_k^k)^\top\bb_2 I(\Xi_0)I(\Xi_0^0)
			\\
			&+\frac{1}{n}\mathbb{E}_0\bb_1^\top\Zb_k^k\Wb^{1/2}\Kb_0^k(\theta_1)\Wb^{1/2}\be_k\bz_k^\top\bb_2 I(\Xi_0)I(\Xi_0^0)
			\\
			&+\frac{1}{n}\mathbb{E}_0\bb_1^\top\bz_k\be_k^\top\Wb^{1/2}\Kb_0^k(\theta_1)\Wb^{1/2}(\Zb_k^k)^\top\bb_2 I(\Xi_0)I(\Xi_0^0)\\
			&+\frac{1}{n}\mathbb{E}_0\bb_1^\top\bz_k\be_k^\top\Wb^{1/2}\Kb_0^k(\theta_1)\Wb^{1/2}\be_k\bz_k^\top\bb_2 I(\Xi_0)I(\Xi_0^0)\\
			:=&\mathcal{J}_1^k+\mathcal{J}_2^k+\mathcal{J}_3^k+\mathcal{J}_4^k.
		\end{split}
		\]
		Meanwhile, $\Zb^{k-1}=\Zb_k^k+\bz_k^0\be_k^\top$ so that similarly we have
		\[
		\begin{split}
			&\mathbb{E}_0[\mathcal{L}^{k-1}I(\Xi_0)I(\Xi_0^0)]\\
			=&\frac{1}{n}\mathbb{E}_0\bb_1^\top\Zb_k^k\Wb^{1/2}\Kb_0^{k-1}(\theta_1)\Wb^{1/2}(\Zb_k^k)^\top\bb_2 I(\Xi_0)I(\Xi_0^0)
			\\
			&+\frac{1}{n}\mathbb{E}_0\bb_1^\top\Zb_k^k\Wb^{1/2}\Kb_0^{k-1}(\theta_1)\Wb^{1/2}\be_k(\bz_k^0)^\top\bb_2 I(\Xi_0)I(\Xi_0^0)
			\\
			&+\frac{1}{n}\mathbb{E}_0\bb_1^\top\bz_k^0\be_k^\top\Wb^{1/2}\Kb_0^{k-1}(\theta_1)\Wb^{1/2}(\Zb_k^k)^\top\bb_2 I(\Xi_0)I(\Xi_0^0)\\
			&+\frac{1}{n}\mathbb{E}_0\bb_1^\top\bz_k^0\be_k^\top\Wb^{1/2}\Kb_0^{k-1}(\theta_1)\Wb^{1/2}\be_k(\bz_k^0)^\top\bb_2 I(\Xi_0)I(\Xi_0^0)\\
			:=&\mathcal{A}_1^k+\mathcal{A}_2^k+\mathcal{A}_3^k+\mathcal{A}_4^k.
		\end{split}
		\]
		In the following, we aim to show that 
		\begin{equation}\label{step 2 error}
			\bigg|\sum_{k=1}^n(\mathcal{J}_i^k-\mathcal{A}_i^k)\bigg|=o_p(n^{-1/2}),\quad i=1,2,3,4.
		\end{equation}
		
		Let's abuse the notation in (\ref{notation a}) slightly by replacing $\Zb_j$  with $\Zb_k^k$. Then, similarly to the proof of \textbf{Step 3}, we have
		\begin{equation}\label{J1 minus}
			\begin{split}
				&\sum_{k=1}^n\mathcal{J}_1^k-\sum_{k=1}^n\mathbb{E}_0n^{-1}\bb_1^\top\Zb_k^k\Wb^{1/2}\Kb_k^k(\theta_1)\Wb^{1/2}\Zb_k^k\bb_2 I(\Xi_0)I(\Xi_0^0)-o_p(n^{-1/2})\\
				=&\sum_{k=1}^n\frac{1}{n^2\theta_1}\mathbb{E}_0a_{\bb_1\be_k}a_{\bb_2\bz_k}^*I(\Xi_0)I(\Xi_0^0)=\sum_{k=1}^n\frac{w_k}{n^3\theta_1^2}\mathbb{E}_0\alpha_k^{-1}a_{\bb_1\bz_k}^*a_{\bb_2\bz_k}^*I(\Xi_0)I(\Xi_0^0),
			\end{split}
		\end{equation}
		where $\alpha_k$ is defined similarly to (\ref{alpha_k}). Define 
		\[
		\bar \alpha_k=1-\frac{w_k}{(n\theta_1)^2}\bar a_{\bz_k\bz_k}-\frac{w_k}{n\theta_1}\text{tr}(\bSigma_2),
		\]
		so that 
		\[
		\alpha_k^{-1}-\bar\alpha_k^{-1}=-\alpha_k^{-1}\bar\alpha_k^{-1}\bigg[\frac{w_k}{(n\theta_1)^2}(a_{\bz_k\bz_k}-\bar a_{\bz_k\bz_k})-\frac{w_k}{n\theta_1}(\bz_k^\top\bSigma_2\bz_k-\text{tr} \bSigma_2)\bigg],
		\]
		and 
		\[
		\mathbb{E}_0\alpha_k^{-1}a_{\bb_1\bz_k}^*a_{\bb_2\bz_k}^*I(\Xi_0)I(\Xi_0^0)=\mathbb{E}_0(\bar\alpha_k^{-1}+\alpha_k^{-1}-\bar\alpha_k^{-1})a_{\bb_1\bz_k}^*a_{\bb_2\bz_k}^*I(\Xi_0)I(\Xi_0^0).
		\]
		Since $|\alpha_k^{-1}|\le C$ and $|\bar\alpha_k^{-1}|\le C$ under the events $\Xi_0$, $\Xi_0^0$ and $\Xi_w$, we have
		\[
		\begin{split}
			&\mathbb{E}_0(\alpha_k^{-1}-\bar\alpha_k^{-1})a_{\bb_1\bz_k}^*a_{\bb_2\bz_k}^*I(\Xi_0)I(\Xi_0^0)\\
			\le& C\sqrt{\mathbb{E}_0|\alpha_k^{-1}-\bar\alpha_k^{-1}|^2I(\Xi_0)I(\Xi_0^0)\times \mathbb{E}_0|a_{\bb_1\bz_k}^*a_{\bb_2\bz_k}^*|^2I(\Xi_0)I(\Xi_0^0)}\\
			\le &C\sqrt{(n\theta_1)^{-1}(n\vee p)^4\log^4n}\le \frac{C}{\sqrt{n\theta_1}}\times (n\vee p)^2\log^2n.
		\end{split}
		\]
		Therefore, 
		\[
		\begin{split}
			&\sum_{k=1}^n\frac{w_k}{n^3\theta_1^2}\mathbb{E}_0\alpha_k^{-1}a_{\bb_1\bz_k}^*a_{\bb_2\bz_k}^*I(\Xi_0)I(\Xi_0^0)
			=	\sum_{k=1}^n\frac{w_k}{n^3\theta_1^2}\mathbb{E}_0\bar\alpha_k^{-1}a_{\bb_1\bz_k}^*a_{\bb_2\bz_k}^*I(\Xi_k)I(\Xi_k^0)+o_p(n^{-1/2}).
		\end{split}
		\]
		Return to (\ref{J1 minus}) so that
		\[
		\begin{split}
			&\sum_{k=1}^n\mathcal{J}_1^k-\sum_{k=1}^n\mathbb{E}_0n^{-1}\bb_1^\top\Zb_k^k\Wb^{1/2}\Kb_k^k(\theta_1)\Wb^{1/2}\Zb_k^k\bb_2 I(\Xi_0)I(\Xi_0^0)\\
			=&\sum_{k=1}^n\frac{w_k}{n^3\theta_1^2}\mathbb{E}_0\bar\alpha_k^{-1}a_{\bb_1\bz_k}^*a_{\bb_2\bz_k}^*I(\Xi_k)I(\Xi_k^0)+o_p(n^{-1/2}).
		\end{split}
		\]
		Similarly, for $\mathcal{A}_{1}^k$, we will have
		\[
		\begin{split}
			&\sum_{k=1}^n\mathcal{A}_1^k-\sum_{k=1}^n\mathbb{E}_0n^{-1}\bb_1^\top\Zb_k^k\Wb^{1/2}\Kb_k^k(\theta_1)\Wb^{1/2}\Zb_k^k\bb_2 I(\Xi_0)I(\Xi_0^0)\\
			=&\sum_{k=1}^n\frac{w_k}{n^3\theta_1^2}\mathbb{E}_0\bar\alpha_k^{-1}a_{\bb_1\bz_k^0}^*a_{\bb_2\bz_k^0}^*I(\Xi_k)I(\Xi_k^0)+o_p(n^{-1/2}).
		\end{split}
		\]
		Recall that $\mathbb{E}_0\bz_k^\top\Ab\bz_k=\mathbb{E}_0(\bz_k^0)^\top\Ab\bz_k^0$ for any matrix $\Ab$ independent of $\bz_k$ and $\bz_k^0$. Then, 
		\[
		\sum_{k=1}^n\frac{w_k}{n^3\theta_1^2}\mathbb{E}_0\bar\alpha_k^{-1}a_{\bb_1\bz_k}^*a_{\bb_2\bz_k}^*I(\Xi_k)I(\Xi_k^0)=\sum_{k=1}^n\frac{w_k}{n^3\theta_1^2}\mathbb{E}_0\bar\alpha_k^{-1}a_{\bb_1\bz_k^0}^*a_{\bb_2\bz_k^0}^*I(\Xi_k)I(\Xi_k^0),
		\]
		which further concludes (\ref{J1 minus}) when $i=1$. 
		For $\sum_{k=1}^n(\mathcal{J}_i^k-\mathcal{A}_i^k), i=2,3,4$, the proof is similar and omitted here. We conclude directly (\ref{step 2 error}) and refer to \cite{cai2020limiting} for further details.
		
		It remains to consider $\mathbb{E}_0[\mathcal{L}^0I(\Xi_0^0)I(\Xi_w)]$, or $\mathbb{E}_0[\mathcal{L}^0I(\hat\Xi_0^0)I(\Xi_w)]$ where $\hat\Xi_0^0$ is defined as the event $\{(n\vee p)^{-1}\|(\Zb^0)^\top\bSigma_2\Zb^0\|\le C\}$.  
		Instead of considering specific $\bb_1$ and $\bb_2$, in the following we calculate the whole matrix 
		\[
		\mathcal{M}:=\mathbb{E}_0\bigg(\frac{1}{n}\bGamma_1^\top\Zb^0\Wb^{1/2}[\Ib-\frac{1}{n\theta_1}\Wb^{1/2}(\Zb^0)^\top\bSigma_2\Zb^0\Wb^{1/2}]^{-1}\Wb^{1/2}(\Zb^0)^\top\bGamma_1I(\hat\Xi_0^0)I(\Xi_w)\bigg).
		\]
		Note that the entries of $\Zb^0$ are  i.i.d. from $\mathcal{N}(0,1)$. Therefore, the entries of $\bGamma^\top\Zb^0$ are also i.i.d. from $\mathcal{N}(0,1)$. In other words, $\bGamma_1^\top\Zb^0$ is independent of $\bGamma_2^\top\Zb^0$ and $\hat\Xi_0^0$. Further note the fact that $\bGamma_1^\top\bGamma_1=\Ib$. Then, 
		\[
		\mathcal{M}=\Ib_r\times \frac{1}{n}\text{tr}\mathbb{E}_0\bigg(\Wb^{1/2}[\Ib-\frac{1}{n\theta_1}\Wb^{1/2}(\Zb^0)^\top\bSigma_2\Zb^0\Wb^{1/2}]^{-1}\Wb^{1/2} I(\Xi_0^0)I(\Xi_w)\bigg).
		\]
		Since the entries of $\Zb^0$ are Gaussian, without loss of generality we can regard  $\bSigma_2$ as diagonal matrix. Then, it suffices to consider
		\[
		\zeta(\theta_1):=n^{-1}\text{tr}\mathbb{E}_0[\Wb^{1/2}\mathcal{K}(\theta_1)\Wb^{1/2} I(\Xi_0^0)I(\Xi_w)]=n^{-1}\sum_{j=1}^n w_j\mathbb{E}_0[\mathcal{K}(\theta_1)]_{jj} I( \Xi_0^0)I( \Xi_w^0),
		\]
		where $\mathcal{K}(x):=[\Ib-(nx)^{-1}\Wb^{1/2}(\Zb^0)^\top\bLambda_2\Zb^0\Wb^{1/2}]^{-1}$.
		
		\noindent\textbf{Step 5: calculating $\zeta(\theta_1)$.}
		
		Let $\mathcal{H}(x)=[\Ib-(nx)^{-1}\bLambda_2^{1/2}\Zb^0\Wb(\Zb^0)^\top\bLambda_2^{1/2}]^{-1}$. By Schur's complement formula, for $1\le j\le n$, 
		\[
		\begin{split}
			[\mathcal{K}(\theta_1)]_{jj}=&\bigg(\bigg[\begin{matrix}
				\Ib&\frac{1}{\sqrt{n\theta_1}}\Wb^{1/2}(\Zb^0)^\top\bLambda_2^{1/2}\\
				\frac{1}{\sqrt{n\theta_1}}\bLambda_2^{1/2}\Zb^0\Wb^{1/2}&\Ib
			\end{matrix}\bigg]^{-1}\bigg)_{jj}\\
			=&1+\bigg(\frac{1}{n\theta_1}\Wb^{1/2}(\Zb^0)^\top\bLambda_2^{1/2}\mathcal{H}(\theta_1)\bLambda_2^{1/2}\Zb^0\Wb^{1/2}\bigg)_{jj}\\
			=&1+\frac{w_j}{n\theta_1}(\bz_j^0)^\top\bLambda_2^{1/2}\mathcal{H}(\theta_1)\bLambda_2^{1/2}\bz_j^0.
		\end{split}
		\]
		Define $\mathcal{H}_j(x)$ by replacing $\Zb^0$ with $\Zb_j^0$. Since $\Zb^0\Wb(\Zb^0)^\top=\sum_{j=1}^nw_j\bz_j^0(\bz_j^0)^\top$, we have
		\[
		\begin{split}
			&\frac{1}{n}(\bz_j^0)^\top\bLambda_2^{1/2}\mathcal{H}(\theta_1)\bLambda_2^{1/2}\bz_j^0\\
			=&\frac{1}{n}(\bz_j^0)^\top\bLambda_2^{1/2}\mathcal{H}_j(\theta_1)\bLambda_2^{1/2}\bz_j^0+\frac{1}{n}(\bz_j^0)^\top\bLambda_2^{1/2}\mathcal{H}(\theta_1)\frac{w_j}{n\theta_1}\bLambda_2^{1/2}\bz_j^0(\bz_j^0)^\top\mathcal{H}_j(\theta_1)\bLambda_2\bz_j^0\\
			=&\frac{1}{1-\frac{w_j}{n\theta_1}(\bz_j^0)^\top\mathcal{H}_j(\theta_1)\bLambda_2\bz_j^0}\frac{1}{n}(\bz_j^0)^\top\bLambda_2^{1/2}\mathcal{H}_j(\theta_1)\bLambda_2^{1/2}\bz_j^0,
		\end{split}
		\]
		where the second line is by (\ref{matrix inverse}). As $n,p\rightarrow\infty$,
		\[
		\frac{w_j}{n\theta_1}(\bz_j^0)^\top\mathcal{H}_j(\theta_1)\bLambda_2\bz_j^0I(\Xi_0^0)I(\Xi_w)\rightarrow 0,
		\]
		which implies that 
		\begin{equation}\label{numer1}
			\begin{split}
				&\mathbb{E}_0\bigg|\frac{\frac{w_j}{n\theta_1}\bigg[(\bz_j^0)^\top\bLambda_2^{1/2}\mathcal{H}_j(\theta_1)\bLambda_2^{1/2}\bz_j^0-\text{tr}[\bLambda_2^{1/2}\mathcal{H}_j(\theta_1)\bLambda_2^{1/2}]\bigg]}{1-\frac{w_j}{n\theta_1}(\bz_j^0)^\top\mathcal{H}_j(\theta_1)\bLambda_2\bz_j^0}I(\Xi_0^0)I(\Xi_w)\bigg|\\
				\le &C\mathbb{E}_0\bigg|\frac{w_j}{n\theta_1}\bigg[(\bz_j^0)^\top\bLambda_2^{1/2}\mathcal{H}_j(\theta_1)\bLambda_2^{1/2}\bz_j^0-\text{tr}[\bLambda_2^{1/2}\mathcal{H}_j(\theta_1)\bLambda_2^{1/2}]\bigg]I(\Xi_0^0)I(\Xi_w)\bigg|\\
				\le &\frac{C\log n}{n\theta_1}\mathbb{E}_0\|\bLambda_2^{1/2}\mathcal{H}_j(\theta_1)\bLambda_2^{1/2}I(\Xi_0^0)\|_F\le o_p(n^{-1/2}),
			\end{split}
		\end{equation}
		where we use the fact that $\|\mathcal{H}_j(\theta_1)\|I(\Xi_0^0)I(\Xi_w)\le C$. On the other hand,
		\[
		\begin{split}
			&\mathbb{E}_0\bigg|\frac{1}{1-\frac{w_j}{n\theta_1}(\bz_j^0)^\top\mathcal{H}_j(\theta_1)\bLambda_2\bz_j^0}-\frac{1}{1-\frac{w_j}{n\theta_1}\text{tr}[\mathcal{H}_j(\theta_1)\bLambda_2]}\bigg|\frac{w_j}{n\theta_1}\text{tr}(\bLambda_2^{1/2}\mathcal{H}_j(\theta_1)\bLambda_2^{1/2})I(\Xi_0^0)I(\Xi_w)\\
			=&\mathbb{E}_0\bigg|\frac{\frac{w_j}{n\theta_1}[(\bz_j^0)^\top\bLambda_2^{1/2}\mathcal{H}_j(\theta_1)\bLambda_2^{1/2}\bz_j^0-\text{tr}\bLambda_2^{1/2}\mathcal{H}_j(\theta_1)\bLambda_2^{1/2}]}{[1-\frac{w_j}{n\theta_1}(\bz_j^0)^\top\mathcal{H}_j(\theta_1)\bLambda_2\bz_j^0][1-\frac{w_j}{n\theta_1}\text{tr}\mathcal{H}_j(\theta_1)\bLambda_2]}\bigg|\frac{w_j}{n\theta_1}\text{tr}(\bLambda_2^{1/2}\mathcal{H}_j(\theta_1)\bLambda_2^{1/2})I(\Xi_0^0)I(\Xi_w)\\
			\le &\frac{Cp\log n}{n\theta_1}\times\mathbb{E}_0\bigg|\frac{w_j}{n\theta_1}[(\bz_j^0)^\top\bLambda_2^{1/2}\mathcal{H}_j(\theta_1)\bLambda_2^{1/2}\bz_j^0-\text{tr}\bLambda_2^{1/2}\mathcal{H}_j(\theta_1)\bLambda_2^{1/2}]\bigg|I(\Xi_0^0)I(\Xi_w)\le o_p(n^{-1/2}). 
		\end{split}
		\]
		Consequently, we conclude that
		\begin{equation}\label{zeta step 1}
			\begin{split}
				\mathbb{E}_0[\mathcal{K}^0(\theta_1)]_{jj}I(\Xi_0^0)I(\Xi_w)
				=&1+\mathbb{E}_0\frac{\frac{w_j}{n\theta_1}\text{tr}\mathcal{H}_j(\theta_1)\bLambda_2}{1-\frac{w_j}{n\theta_1}\text{tr}\mathcal{H}_j(\theta_1)\bLambda_2}I(\Xi_0^0)I(\Xi_w)+o_p(n^{-1/2})\\
				=&\mathbb{E}_0\frac{1}{1-\frac{w_j}{n\theta_1}\text{tr}\mathcal{H}_j(\theta_1)\bLambda_2}I(\Xi_0^0)I(\Xi_w)+o_p(n^{-1/2}),
			\end{split}
		\end{equation}
		and the above convergence rate is actually uniform over $j$.  Further note that 
		\[
		\begin{split}
			&\mathbb{E}_0\bigg|\frac{1}{1-\frac{w_j}{n\theta_1}\text{tr}\mathcal{H}_j(\theta_1)\bLambda_2}-\frac{1}{1-\frac{w_j}{n\theta_1}\text{tr}\mathcal{H}(\theta_1)\bLambda_2}\bigg|I(\Xi_0^0)I(\Xi_w)\\
			=&\mathbb{E}_0\bigg|\frac{\frac{w_j}{n\theta_1}\text{tr}[\mathcal{H}(\theta_1)-\mathcal{H}_j(\theta_1)]\bLambda_2}{[1-\frac{w_j}{n\theta_1}\text{tr}\mathcal{H}_j(\theta_1)\bLambda_2][1-\frac{w_j}{n\theta_1}\text{tr}\mathcal{H}(\theta_1)\bLambda_2]}\bigg|I(\Xi_0^0)I(\Xi_w)\\
			=&\mathbb{E}_0\bigg|\frac{(\frac{w_j}{n\theta_1})^2(\bz_j^0)^\top\mathcal{H}_j(\theta_1)\bLambda_2\mathcal{H}(\theta_1)\bz_j^0}{[1-\frac{w_j}{n\theta_1}\text{tr}\mathcal{H}_j(\theta_1)\bLambda_2][1-\frac{w_j}{n\theta_1}\text{tr}\mathcal{H}(\theta_1)\bLambda_2]}\bigg|I(\Xi_0^0)I(\Xi_w)\le o_p(n^{-1/2}).
		\end{split}
		\]
		As a result,
		\begin{equation}\label{zeta step 2}
			\zeta(\theta_1)=\frac{1}{n}\sum_{j=1}^n\mathbb{E}_0\frac{w_j}{1-\frac{w_j}{n\theta_1}\text{tr}[\mathcal{H}(\theta_1)\bLambda_2]}I(\Xi_0^0)I(\Xi_w)+o_p(n^{-1/2}). 
		\end{equation}
		Furthermore, by (\ref{Xiw}) and (\ref{Xi0}),
		\[
		\begin{split}
			&\frac{1}{n}\sum_{j=1}^n\mathbb{E}_0\bigg|\frac{w_j}{1-\frac{w_j}{n\theta_1}\text{tr}\mathcal{H}(\theta_1)\bLambda_2}I(\Xi_0^0)I(\Xi_w)-\frac{w_j}{1-\frac{w_j}{n\theta_1}\mathbb{E}_0\text{tr}\mathcal{H}(\theta_1)\bLambda_2I(\Xi_0^0)I(\Xi_w)}\bigg|\\
			=&\frac{1}{n}\sum_{j=1}^n\mathbb{E}_0\bigg|\frac{w_jI(\Xi_0^0)I(\Xi_w)}{1-\frac{w_j}{n\theta_1}\text{tr}\mathcal{H}(\theta_1)\bLambda_2I(\Xi_0^0)I(\Xi_w)}-\frac{w_jI(\Xi_0^0)I(\Xi_w)}{1-\frac{w_j}{n\theta_1}\mathbb{E}_0\text{tr}\mathcal{H}(\theta_1)\bLambda_2I(\Xi_0^0)I(\Xi_w)}\bigg|+o_p(n^{-2})\\
			=&\frac{1}{n}\sum_{j=1}^n\mathbb{E}_0\bigg|\frac{\frac{w_j^2}{n\theta_1}[\text{tr}\mathcal{H}(\theta_1)\bLambda_2-\mathbb{E}_0\text{tr}\mathcal{H}(\theta_1)\bLambda_2]I(\Xi_0^0)I(\Xi_w)]}{[1-\frac{w_j}{n\theta_1}\text{tr}\mathcal{H}(\theta_1)\bLambda_2I(\Xi_0^0)I(\Xi_w)][1-\frac{w_j}{n\theta_1}\mathbb{E}_0\text{tr}\mathcal{H}(\theta_1)\bLambda_2I(\Xi_0^0)I(\Xi_w)]}\bigg|+o_p(n^{-2})\\
			\le &C\sqrt{\mathbb{E}_0\bigg|\frac{1}{n\theta_1}[\text{tr}\mathcal{H}(\theta_1)\bLambda_2-\mathbb{E}_0\text{tr}\mathcal{H}(\theta_1)\bLambda_2]I(\Xi_0^0)I(\Xi_w)\bigg|^2}+o_p(n^{-2}).
		\end{split}
		\]
		Let $\mathbb{E}_j^0=\mathbb{E}(\cdot\mid \bz_1^0,\ldots,\bz_j^0,\Wb)$. Then,
		\[
		\begin{split}
			&\frac{1}{n\theta_1}[\text{tr}\mathcal{H}(\theta_1)\bLambda_2-\mathbb{E}_0\text{tr}\mathcal{H}(\theta_1)\bLambda_2]I(\Xi_0^0)I(\Xi_w)=\sum_{j=1}^n(\mathbb{E}_j^0-\mathbb{E}_{j-1}^0)\frac{1}{n\theta_1}\text{tr}\mathcal{H}(\theta_1)\bLambda_2I(\Xi_0^0)I(\Xi_w)\\
			=&\sum_{j=1}^n(\mathbb{E}_j^0-\mathbb{E}_{j-1}^0)\bigg[\frac{1}{n\theta_1}\text{tr}\mathcal{H}_j(\theta_1)\bLambda_2+\frac{w_j}{(n\theta_1)^2}(\bz_j^0)^\top\mathcal{H}_j(\theta_1)\bLambda_2\mathcal{H}(\theta_1)\bz_j^0\bigg]I(\Xi_0^0)I(\Xi_w)\\
			=&\sum_{j=1}^n(\mathbb{E}_j^0-\mathbb{E}_{j-1}^0)\frac{w_j}{(n\theta_1)^2}(\bz_j^0)^\top\mathcal{H}_j(\theta_1)\bLambda_2\mathcal{H}(\theta_1)\bz_j^0I(\Xi_0^0)I(\Xi_w)+o_p(n^{-1}).
		\end{split}
		\]
		Therefore, by Burkholder's inequality,
		\[
		\begin{split}
			&\mathbb{E}_0\bigg|\frac{1}{n\theta_1}[\text{tr}\mathcal{H}(\theta_1)\bLambda_2-\mathbb{E}_0\text{tr}\mathcal{H}(\theta_1)\bLambda_2]I(\Xi_0^0)I(\Xi_w)\bigg|^2\\
			\le& \frac{C}{(n\theta_1)^2}\sum_{j=1}^n\mathbb{E}_j^0\bigg|\frac{1}{n\theta_1}(\bz_j^0)^\top\mathcal{H}_j(\theta_1)\bLambda_2\mathcal{H}(\theta_1)\bz_j^0I(\Xi_0^0)I(\Xi_w)\bigg|^2
			\le \frac{C}{n\theta_1^2}\times \bigg(\frac{p}{n\theta_1}\bigg)^2.
		\end{split}
		\]
		Then, we conclude that
		\begin{equation}\label{zeta step 3}
			\zeta(\theta_1)=\frac{1}{n}\sum_{j=1}^n\frac{w_j}{1-\frac{w_j}{n\theta_1}\mathbb{E}\text{tr}\mathcal{H}(\theta_1)\bLambda_2I(\Xi_0^0)I(\Xi_w)}+o_p(n^{-1/2}).
		\end{equation}
		
		It remains to calculate $n^{-1}\mathbb{E}\text{tr}\mathcal{H}(\theta_1)\bLambda_2I(\Xi_0^0)I(\Xi_w)$. In fact, this is totally parallel to $\zeta(\theta_1)$ by exchanging $\Wb$  and $\bLambda_2$ and transposing $\Zb^0$. Then, we conclude  that 
		\begin{equation}\label{zeta step 4}
			\frac{1}{n}\mathbb{E}\text{tr}\mathcal{H}(\theta_1)\bLambda_2I(\Xi_0^0)I(\Xi_w)=\frac{1}{n}\sum_{i=1}^p\frac{\lambda_{r+i}}{1-\theta_1^{-1}\lambda_{r+i}\zeta(\theta_1)}+o_p(n^{-1/2}).
		\end{equation}
		In other words, 
		\[
		\zeta(\theta_1)=\frac{1}{n}\sum_{j=1}^n\frac{w_j}{1-\frac{w_j}{n\theta_1}\sum_{i=1}^p\frac{\lambda_{r+i}}{1-\theta_1^{-1}\lambda_{r+i}\zeta(\theta_1)}}+o_p(n^{-1/2}). 
		\]
		Recall the definition of $\hat\zeta_1$,
		\[
		\hat\zeta_1=\frac{1}{n}\sum_{j=1}^n\frac{w_j}{1-\frac{w_j}{n\theta_1}\sum_{i=1}^p\frac{\lambda_{r+i}}{1-\theta_1^{-1}\lambda_{r+i}\hat\zeta_1}}.
		\]
		Then, we have
		\begin{equation}\label{zeta1 minus}
			\begin{split}
				&\zeta(\theta_1)-\hat\zeta_1-o_p(n^{-1/2})
				=\frac{1}{n}\sum_{j=1}^n\frac{\frac{w_j}{n\theta_1}\sum_{i=1}^p[\frac{\lambda_{r+i}}{1-\theta_1^{-1}\lambda_{r+i}\zeta(\theta_1)}-\frac{\lambda_{r+i}}{1-\theta_1^{-1}\lambda_{r+i}\hat\zeta_1}]}{[1-\frac{1}{n\theta_1}\sum_{i=1}^p\frac{\lambda_{r+i}}{1-\lambda_1^{-1}\lambda_{r+i}}][1-\frac{w_j}{n\theta_1}\sum_{i=1}^p\frac{\lambda_{r+i}}{1-\theta_1^{-1}\lambda_{r+i}\zeta(\theta_1)}]}\\
				=&\frac{1}{n}\sum_{j=1}^n\frac{\frac{w_j}{n\theta_1}\sum_{i=1}^p\frac{\lambda_{r+i}^2}{[1-\theta_1^{-1}\lambda_{r+i}\zeta(\theta_1)][1-\theta_1^{-1}\lambda_{r+i}\hat\zeta_1]}}{[1-\frac{1}{n\theta_1}\sum_{i=1}^p\frac{\lambda_{r+i}}{1-\lambda_1^{-1}\lambda_{r+i}}][1-\frac{w_j}{n\theta_1}\sum_{i=1}^p\frac{\lambda_{r+i}}{1-\theta_1^{-1}\lambda_{r+i}\zeta(\theta_1)}]}\times \frac{1}{\theta_1}[\zeta(\theta_1)-\hat\zeta_1]\\
				=&o_p(1)[\zeta(\theta_1)-\hat\zeta_1].
			\end{split}
		\end{equation}
		Therefore, we conclude that $\zeta(\theta_1)-\hat\zeta_1-o_p(n^{-1/2})$. 
		
		{	\noindent\textbf{Step 6: adding back replacement error.}
			
			The lemma follows \textbf{Steps 3} to \textbf{5} if $\Fb=\Fb^0$, i.e., the entries of $\Fb$ are independent with mean 0, variance 1 and bounded eighth moment. For general $\Fb$, adding back the replacement error in \textbf{Step 2}, we have
			\[
			\begin{split}
				&n^{-1}\bGamma_1^\top\Zb\Wb^{1/2}\Kb(\theta_1)\Zb^\top\Wb^{1/2}\bGamma_1\\
				=&n^{-1}\bGamma_1^\top(\check{\Zb}\check\Zb^\top-\check{\Zb}_0\check{\Zb}_0^\top)\bGamma_1+n^{-1}\bGamma_1^\top\Zb_0\Wb^{1/2}\Kb_0(\theta_1)\Zb_0^\top\Wb^{1/2}\bGamma_1+o_p(n^{-1/2})\\
				=&n^{-1}\bGamma_1^\top(\check{\Zb}\check\Zb^\top-\check{\Zb}_0\check{\Zb}_0^\top)\bGamma_1+n^{-1}\bGamma_1^\top\Zb_0\Wb\Zb_0^\top\bGamma_1-(n^{-1}\text{tr}\Wb)\times\Ib_r+\hat\zeta_1\Ib_r+o_p(n^{-1/2})\\
				=&n^{-1}\bGamma_1^\top\Zb\Wb\Zb^\top\bGamma_1-(n^{-1}\text{tr}\Wb)\times\Ib_r+\hat\zeta_1\Ib_r+o_p(n^{-1/2}),
			\end{split}
			\]
			which concludes the lemma.
		}
	\end{proof}

	\begin{lemma}\label{lema3}
		Under the conditions of Theorem \ref{thm: representation},  we have 
		\[
		n^{-1}\bGamma_1^\top\Zb\Wb^{1/2}\Kb(\hat\lambda_1)\Kb(\theta_1)\Wb^{1/2}\Zb^\top\bGamma_1=\Ib+o_{p^*}(1).
		\]
	\end{lemma}
	\begin{proof}
		By the matrix inverse formula $(\Ab+\Bb)^{-1}=\Ab^{-1}-(\Ab+\Bb)^{-1}\Bb\Ab^{-1}$, 
		\[
		\begin{split}
			&n^{-1}\bGamma_1^\top\Zb\Wb^{1/2}\Kb(\hat\lambda_1)\Kb(\theta_1)\Wb^{1/2}\Zb^\top\bGamma_1\\
			=&n^{-1}\bGamma_1^\top\Zb\Wb^{1/2}\Kb(\theta_1)\Kb(\theta_1)\Wb^{1/2}\Zb^\top\bGamma_1-\delta_1n^{-1}\bGamma_1^\top\Zb\Wb^{1/2}\Kb(\hat\lambda_1)\Kb^2(\theta_1)\Wb^{1/2}\Zb^\top\bGamma_1.
		\end{split}
		\]
		By Lemma \ref{lem: preliminary}, $\delta_1=o_p((\log n)^{-1})$, thus 
		\[
		\delta_1n^{-1}\bGamma_1^\top\Zb\Wb^{1/2}\Kb(\hat\lambda_1)\Kb^2(\theta_1)\Wb^{1/2}\Zb^\top\bGamma_1=o_p(1).
		\]
		On the other hand, similarly to the proof of Lemma \ref{lema2},
		\[
		n^{-1}\bGamma_1^\top\Zb\Wb^{1/2}\Kb^2(\theta_1)\Wb^{1/2}\Zb^\top\bGamma_1=\Ib_r+o_p(1),
		\]
		where we use the fact that $\hat\zeta_1=1+o_p(1)$. 
		Then, the lemma follows.
	\end{proof}

	\section{Proof of results in Section \ref{sec:non-spiked}: preliminaries}\label{secc}
	\subsection{Outline of the proof}
	Our major target is to prove Lemma \ref{lemma: ratio} in the main paper, which provides a sufficiently fast convergence rate for the ratio $\hat\lambda_{r+1}/\lambda_0$. Therefore,  we start with $r=0$, i.e., there are no spiked eigenvalues. Under such cases, $\Xb=\bPsi\Eb=\bPsi\Zb$, where the entries of $\Zb$ are independent with mean 0, variance 1 and bounded moments. To ease notation, without loss of generality, we assume $\bPsi=\bSigma^{1/2}$, where $\bSigma$ is the population covariance matrix. Then, the bootstrapped sample covariance matrix is
	\[
	\hat{\Sbb}=n^{-1}\Xb\Wb\Xb^\top=n^{-1}\bSigma^{1/2}\Zb\Wb\Zb^\top\bSigma^{1/2},
	\]
	where $\bSigma$ satisfies Assumptions \ref{c1} and \ref{c3}. Let $\lambda_1\ge\cdots\ge \lambda_p$ be the eigenvalues of $\bSigma$.
	
	Recall the definition of orders $\{t_1,\ldots,t_n\}$. 
	We define a series of events $\Omega_n$ satisfying:
	\begin{enumerate}
		\item $|w_{t_1}-w_{t_2}|\ge (\log n)^{-c}$;
		\item $C^{-1}\log n\le|w_{t_1}|\le C\log n$;
		\item $|w_{t_1}-w_{t_{[\sqrt{n}]}}|\ge C^{-1}\log n$;
		\item $n^{-1}\sum_{j=1}^nw_j^h\le C$, for $h=1,2$,
	\end{enumerate}
	for some constants $c, C>0$. Since $w_j$'s are i.i.d. from $Exp(1)$, we have
	\[
	w_{t_j}\overset{d}{=}\sum_{i=1}^{n-j+1}\frac{i}{n}\tilde w_i,
	\]
	where $\tilde w_i$'s are also i.i.d. from $Exp(1)$. Therefore, 
	it's not hard to verify $\Omega_n$ holds with probability tending to one. Our proof will be conditional on the events $\Omega_n$, which has negligible effects on the limiting distributions of $\hat\lambda_{r+1}$. 
	
	Motivated by the gap between $w_{t_1}$ and $w_{t_2}$, we will regard $w_{t_1}$ as an outlier from the spectrum of $\Wb$. Therefore, by definition, when $r=0$, $\hat\lambda_{r+1}$ satisfies 
	\[
	\det(\hat\lambda_{r+1}\Ib-\hat\Sbb)=0\Longrightarrow 1+n^{-1}w_{t_1}\bx_{t_1}^\top[\hat\Sbb^{(1)}-\hat\lambda_{r+1}\Ib]^{-1}\bx_{t_1}=0,
	\]
	where $\hat\Sbb^{(1)}=n^{-1}\sum_{j\ne 1}w_{t_j}\bx_{t_j}\bx_{t_j}^\top$ and  $\hat\lambda_1$ is not an eigenvalue of $\hat\Sbb^{(1)}$ for simplicity because $w_j$'s are from continuous distribution.
	Therefore, to investigate the properties of $\hat\lambda_{r+1}$, one needs to find some approximation to $[\hat\Sbb^{(1)}-\hat\lambda_{r+1}\Ib]^{-1}$.   Since $\hat\lambda_{r+1}$ is random, usually uniform convergence of $[\hat\Sbb^{(1)}-z\Ib]^{-1}$ is required for $z$ in some region of $\mathbb{C}^+$. This is  referred to the local law in random matrix theory. However, the scenario considered in the current paper has at least three differences from those considered in the literature, such as in \cite{yang2019edge}. Firstly, the support of $w_{t_1}$ is unbounded, so we don't have regular edge for the limiting spectral distribution of $\hat\Sbb^{(1)}$. In other words, $\hat\lambda_{r+1}$ tends to infinity  rather than some constant as $n\rightarrow \infty$. Secondly, the ``square-root'' type regularity conditions (see (2.18) in \cite{yang2019edge}) will not always hold. This requires us to use a larger imaginary part of $z$ (of order $\log^{1+c}n$) in the proof. Thirdly, the stability lemma (such as Lemma 5.11 in \cite{yang2019edge}) is not guaranteed. Instead, we will use the technique introduced in \cite{lee2016extremal} and \cite{kwak2021extremal}. After finding the approximation to $[\hat\Sbb^{(1)}-\hat\lambda_{r+1}\Ib]^{-1}$, Lemma \ref{lemma: ratio} can be verified similarly to the proof of Lemma \ref{thm: representation}. 
	
	Given Lemma \ref{lemma: ratio}, it will be easy to prove Lemma \ref{lemma: lambda_0} by a detailed calculation of the fluctuations of $\lambda_0$. It turns out that the limiting distribution of $\hat\lambda_{r+1}$ is mainly determined by $\lambda_0$, and further by $w_{t_1}$, as shown in Lemma \ref{lemma: lambda_0}. To extend the results to the case of $r>0$, we use the technique introduced in \cite{cai2020limiting}. Specifically, for $r>0$, $\hat\lambda_{r+1}$ is an eigenvalue of $\hat\Sbb$ if and only if
	\[
	\det\bigg(\hat\lambda_{r+1}\Ib-\frac{1}{n}\Wb^{1/2}\Zb^\top(\bGamma_1\bLambda_1\bGamma_1^\top+\bGamma_2\bLambda_2\bGamma_2^\top)\Zb\Wb^{1/2}\bigg)=0.
	\]
	Then, it suffices to verify that the determinant can take the value of 0 when $\hat\lambda_{r+1}$ is in a neighborhood of $\hat\varphi_1$, and show that this is exactly the largest non-spiked eigenvalue. The details are given in Section \ref{sec:c6}.

	\subsection{Definitions}
	Before the formal proof, we need to introduce some definitions commonly used in the literature of random matrix theory. 
	\begin{definition}[High probability event]
		We say that an $n$-dependent event $\mathcal{E}_n$ holds with high probability if for any constant $d>0$, 
		\[
		\mathbb{P}(\mathcal{E}_n)\geqslant 1-n^{-d},
		\]
		for all sufficiently large $n$. For a high probability, we may take it as given in the proof, which only brings in negligible errors.
	\end{definition}
	
	\begin{definition}[Stochastic domination]
		(a). For two families of nonnegative random variables
		\[
		A=\{A_{n}(t):n\in\mathbb{Z}_{+},t\in T_{n}\},\qquad B=\{B_{n}(t):n\in\mathbb{Z}_{+},t\in T_{n}\},
		\]
		where $T_{n}$ is
		a possibly $n$-dependent parameter set, we say that $A$ is stochastically
		dominated by $B$, uniformly on $t$ if for all (small) $\varepsilon>0$
		and (large) $d>0$ there exists $n_{0}(\varepsilon,d)\in\mathbb{Z}_{+}$
		such that 	as $n\ge n_{0}(\varepsilon,d)$,
		\[
		\sup_{t\in T_{n}}\mathbb{P}\big(A_{n}(t)>n^{\varepsilon}B_{n}(t)\big)\le n^{-d}.
		\]
		If $A$ is stochastically dominated
		by $B$, uniformly on $t$, we use notation $A\prec B$ or $A=O_{\prec}(B)$.
		Moreover, for some complex family $A$ if $|A|\prec B$ we also write $A=O_{\prec}(B)$. \\
		(b). Let $\Ab$ be a family of random matrices and $\zeta$ be a family of nonnegative random variables. Then, we denote $\Ab=O_{\prec}(\zeta)$ if $\Ab$ is dominated by $\zeta$ under weak operator norm sense, i.e. $|\langle\mathbf{v},\Ab\mathbf{w}\rangle|\prec\zeta\|\mathbf{v}\|\|\mathbf{w}\|$  for any deterministic vectors $\mathbf{v}$ and $\mathbf{w}$.\\
		(c). For two sequences of numbers $\{b_{n}\}_{n=1}^{\infty}$,
		$\{c_{n}\}_{n=1}^{\infty}$, $b_{n}\prec c_{n}$ if for all $c>0$, $b_n\leq n^{c}c_n$ for sufficiently large $n$.
	\end{definition}
	Next, we introduce the definition of Stieltjes transform. 
	Note that $\hat\Sbb_1$ has a separable structure. Motivated by \cite{yang2019edge}, we  define 
	\[
	\begin{split}
		m_{1c}(z)=&c\int\frac{t}{-z(1+m_{2c}(z)t)}dF_{\bSigma}(t),\quad m_{2c}(z)=\int\frac{t}{-z(1+m_{1c}(z)t)}\exp(-t)dt,\\
		m_c(z)=&\int\frac{1}{-z(1+tm_{2c}(z))}dF_{\bSigma}(t),\quad z\in\mathbb{C}^+.
	\end{split}
	\]
	Indeed, $m_{1c}(z)$, $m_{2c}(z)$ and $m_{c}(z)$ are the limits of some Stieltjes transforms, corresponding to some deterministic probability functions, shown in \cite{yang2019edge}.  Then, $m_{1c}(z)$ and $m_{2c}(z)$ have unique solutions in $\mathbb{C}^+$ according to \cite{couillet2014analysis} and \cite{yang2019edge}. Remember $r=0$ so $\lambda_i$'s are bounded. Define the finite sample versions as 
	\[
	\begin{split}
		m_{1n}(z)=&\frac{1}{n}\sum_{i=1}^p\frac{\lambda_i}{-z(1+m_{2n}(z)\lambda_i)},\quad m_{2n}(z)=\frac{1}{n}\sum_{j= 1}^n\frac{w_j}{-z(1+w_jm_{1n}(z))},\\
		m_n(z)=&\frac{1}{p}\sum_{i=1}^p\frac{1}{-z(1+\lambda_im_{2n}(z))},\quad z\in\mathbb{C}^+.
	\end{split}
	\]
	Note the different notation for complex number $z$ and the entries in $\Zb$, i.e., $z_{ij}$. The latter  always has double subscript index. Then, for any $z\in\mathbb{C}^+$, $m_{1n}(z)$ converges to $m_{1c}(z)$ as $n\rightarrow\infty$, and similar results hold for $m_{2n}(z)$, $m_n(z)$.  Let $\lambda_{(1)}$ be the largest solution satisfying
	\[
	1+(w_{t_1}+n^{-1/2+c})m_{1n}(\lambda_{(1)})=0.
	\]
	Note the difference between $\lambda_{(1)}$ and $\lambda_0$. We first show that there is a solution to the above equation. By definition, 
	\begin{equation}\label{lambda1}
		\begin{split}
			&-\frac{1}{w_{t_1}+n^{-1/2+c}}=-\frac{1}{n}\sum_{i=1}^p\frac{\lambda_i}{\lambda_{(1)}-\frac{\lambda_i}{n}\sum_{j=1}^n\frac{w_j}{1-(w_{t_1}+n^{-1/2+c})^{-1}w_j}}\\
			\Rightarrow& 1=\frac{1}{n}\sum_{i=1}^p\frac{\lambda_i}{\frac{\lambda_{(1)}}{w_{t_1}+n^{-1/2+c}}-\frac{\lambda_i}{n}\sum_{j}\frac{w_j}{w_{t_1}+n^{-1/2+c}-w_j}}.
		\end{split}
	\end{equation}
	Therefore, by continuity and monotonicity on $\lambda_{(1)}$, the equation always has only one  solution in the interval 
	\[
	\bigg(\frac{\lambda_1}{n}\sum_{j}\frac{w_j(w_{t_1}+n^{-1/2+c})}{w_{t_1}+n^{-1/2+c}-w_j},+\infty\bigg).
	\]
	Under $\Omega_n$, for sufficiently large $n$,
	\begin{equation}\label{sum}
		\begin{split}
			&\frac{1}{n}\sum_{j}\frac{w_j}{w_{t_1}+n^{-1/2+c}-w_j}
			= \frac{1}{n}(\sum_{j=1}^{[\sqrt{n}]}+\sum_{j=[\sqrt{n}]+1}^{n})\frac{w_j}{w_{t_1}+n^{-1/2+c}-w_j}\\
			\le &\frac{\sqrt{n}}{n}\log^{1+c}n+\frac{1}{n}\sum_{j=\sqrt{n}+1}^n\frac{w_j}{C^{-1}\log n}
			\le C(\log n)^{-1}.
		\end{split}
	\end{equation}
	Then, combining (\ref{lambda1}), we conclude that
	\begin{equation}\label{ratio}
		\frac{\lambda_{(1)}}{w_{t_1}+n^{-1/2+c}}\sim \frac{1}{n}\sum_i\lambda_i:=\phi_n\bar\lambda,
	\end{equation}
	where $\phi_n:=p/n$ and $\bar\lambda:=p^{-1}\sum_{i=1}^p\lambda_i$. Actually, the definition of $m_{1n}(z)$ can be extended to $z\in\mathbb{R}$ for $z>w_{t_1}$ by letting the imaginary part $\operatorname{Im}z\downarrow 0$. 
	In the following, we aim to prove that the largest eigenvalue of $\hat\Sbb$ will not exceed $\lambda_{(1)}$ with high probability. The result is shown in Lemma \ref{upper bound}.
	
	Return to the sample covariance matrix. Recall the companion matrix defined by
	\[
	\hat{\mathcal{S}}=n^{-1}\Wb^{1/2}\Xb^\top\Xb\Wb^{1/2}.
	\]
	Define the corresponding Green functions by
	\[
	\mathcal{G}(z)=(\hat{\mathcal{S}}-z\Ib)^{-1},\quad \Gb(z)=(\hat{\Sbb}-z\Ib)^{-1}.
	\]
	In the following, we may suppress the dependence on $z$ and write $\mathcal{G}$, $\Gb$ directly. Define the Stieltjes transform corresponding to $\hat{\Sbb}$ and $\hat{\mathcal{S}}$ as
	\[
	\tilde m(z)=\frac{1}{n}\text{tr}\mathcal{G}(z),\quad m(z)=\frac{1}{p}\text{tr}\Gb(z),
	\] 
	and two related quantities 
	\[
	m_1(z)=\frac{1}{n}\text{tr}\Gb(z)\bSigma,\quad m_2(z)=\frac{1}{n}\sum_{i=1}^nw_i[\mathcal{G}(z)]_{ii}.
	\]
	Since $\hat{\Sbb}$ and $\hat{\mathcal{S}}$ have at most $|n-p|$ zero non-identical eigenvalues, we have
	\[
	n\tilde  m(z)=pm(z)-\frac{n-p}{z}.
	\]
	Now we introduce the definition of minors  in \cite{pillai2014universality}.
	\begin{definition}[Minors]
		For any index set $\mathbb{T}\subset\{1,\ldots,n\}$, define $\Xb^{(\mathbb{T})}$ as the $p\times (n-|\mathbb{T}|)$ subset of $\Xb$ by removing the columns of $\Xb$ indexed by $\mathbb{T}$. However, we keep the names of indices of $\Xb$, i.e.,
		\[
		(\Xb^{(\mathbb{T})})_{ij}=\mathbf{1}(j\notin \mathbb{T})X_{ij}.
		\]
		Define $\hat{\Sbb}^{(\mathbb{T})}$, $\hat{\mathcal{S}}^{(\mathbb{T})}$, $\mathcal{G}^{(\mathbb{T})}$ and $\Gb^{(\mathbb{T})}$ by replacing $\Xb$ with $\Xb^{(\mathbb{T})}$. Further, define $m^{(\mathbb{T})}(z)$, $\tilde m^{(\mathbb{T})}(z)$, $m_1^{(\mathbb{T})}(z)$, $m_2^{(\mathbb{T})}(z)$ using $\mathcal{G}^{(\mathbb{T})}$, $\Gb^{(\mathbb{T})}$. Abbreviate $(\{i\})$ as $(i)$ and $\{i\}\cup\mathbb{T}$ as $(i\mathbb{T})$. To ease notation, we may suppress the dependence on $z$ in the proof.
	\end{definition}
	Then, we have the next lemma.
	\begin{lemma}[Resolvent identity]\label{resolvent}
		Write $\by_i=\sqrt{w_i}\bx_i$. Then,
		\[
		\begin{split}
			\mathcal{\mathcal{G}}_{ii}(z)=&\frac{1}{-z-zn^{-1}\by_i^\top\Gb^{(i)}\by_i},\\
			\mathcal{\mathcal{G}}_{ij}(z)=&zn^{-1}\mathcal{G}_{ii}(z)\mathcal{G}_{jj}^{(i)}(z)\by_i^\top\Gb^{(ij)}\by_j,\quad i\ne j,\\
			\mathcal{\mathcal{G}}_{ij}(z)=&\mathcal{G}_{ij}^{(k)}(z)+\frac{\mathcal{G}_{ik}(z)\mathcal{G}_{kj}(z)}{\mathcal{G}_{kk}(z)},\quad i,j\ne k.
		\end{split}
		\]
		The results also hold after replacing $\mathcal{G}$ with $\mathcal{G}^{(\mathbb{T})}$. 
	\end{lemma}
	\begin{proof}
		See Lemma 2.3 in \cite{pillai2014universality}.
	\end{proof}

	\subsection{Some useful lemmas}
	In the following, we present some useful lemmas for the proof related to $\hat \lambda_{r+1}$, such as the local law and eigenvalue rigidity properties. These lemma commonly appear in the literature of random matrix theory, such as \cite{ERDOS20121435}, \cite{pillai2014universality}, \cite{ding2018necessary} and \cite{yang2019edge}  to characterize the fluctuations of a non-spiked sample eigenvalue. 
	Following the definitions above, actually it suffices  to consider  $z$ in  the region
	\[
	D:=\{z=\lambda_{(1)}+\tau+i\eta: 0< \tau\le C\log n,n^{-2/3}\le \eta\le (\log n)^{1+c}\},
	\]
	for some small constant $c>0$. The following lemma holds.
	\begin{lemma}\label{m_{1n}}
		For the multiplier bootstrap, if Assumptions \ref{c1}, \ref{c3} and the events $\Omega_n$ hold, as $n\rightarrow\infty$ we have
		\[
		m_{1n}(\lambda_{(1)}+\tau)\in\bigg[-\frac{1}{w_{t_1}+n^{-1/2+c}},-\frac{1}{1+c}\frac{\phi_n\bar\lambda}{ \lambda_{(1)}+\tau}\bigg],
		\]
		for any $0\le \tau\le C\log n$, and 
		\[
		\begin{split}
			&-\frac{1}{w_{t_1}+n^{-1/2+c}}<\operatorname{Re}m_{1n}(z)<-\frac{\phi_n\bar\lambda}{1+c}\frac{E}{(E^2+\eta^2)},\\
			&\frac{\phi_n\bar\lambda}{1+c}\frac{\eta}{E^2+\eta^2}<\operatorname{Im}m_{1n}(z)<\eta|\operatorname{Re}m_{1n}(z)|<O(\eta\log ^{-1}n),\\
			&|m_{2n}(z)|=o(1),\quad |m_n(z)|=o(1),\\
			&\operatorname{Im}m_{2n}(z)=o(\eta),\quad \operatorname{Im}m_{n}(z)=o(\eta),
		\end{split}
		\]
		for any $z=E+i\eta\in D$ and small constant $c>0$, where $E:=\lambda_{(1)}+\tau$.
	\end{lemma}
	\begin{proof}
		We first calculate $m_{1n}(\lambda_{(1)}+\tau)$. Let $z=\lambda_{(1)}+\tau$. Define
		\begin{equation}\label{f}
			f_{1n}:=f_{1n}(z,m_{1n}(z)):=-m_{1n}(z)-\frac{1}{n}\sum_{i=1}^p\frac{\lambda_i}{z-\frac{\lambda_i}{n}\sum_{j}\frac{w_j}{1+m_{1n}(z)w_j}}=0.
		\end{equation}
		On one hand, if $m_{1n}(\lambda_{(1)}+\tau)=-(w_{t_1}+n^{-1/2+c})^{-1}$, then for $\tau\ge0$ we always have $f_{1n}\ge0$. On the other hand, if
		\[
		m_{1n}(\lambda_{(1)}+\tau)=-\frac{1}{1+c}\frac{\phi_n\bar\lambda}{ \lambda_{(1)}+\tau}, 
		\]
		for some constant $c>0$, then similarly to (\ref{sum}) we have 
		\[
		f_{1n}=\frac{1}{1+c}\frac{\phi_n\bar\lambda}{ \lambda_{(1)}+\tau}-\frac{1}{n}\sum_{i=1}^p\frac{\lambda_i}{(\lambda_{(1)}+\tau)[1+o(1)]}<0,
		\]
		for sufficiently large $n$, where the $o(1)$ in denominator is uniform on $i\in[1,p]$.
		Then, by continuity, for sufficiently large $n$ there is always a solution satisfying
		\[
		m_{1n}(\lambda_{(1)}+\tau)\in\bigg[-\frac{1}{w_{t_1}+n^{-1/2+c}},-\frac{1}{1+c}\frac{\phi_n\bar\lambda}{ \lambda_{(1)}+\tau}\bigg],
		\]
		for any $0\le \tau\le C\log n$, because $f_{1n}$ takes opposite signs at the two end points.
		
		Now we add the imaginary part into the equation. Let $z=E+i\eta\in D$ and $E=\lambda_{(1)}+\tau$.  Taking real part in (\ref{f}) and writing $m_{1n}$ for $m_{1n}(z)$, we have
		\[
		\begin{split}
			\operatorname{Re}f_{1n}=-\operatorname{Re}m_{1n}-\frac{1}{n}\sum_{i=1}^p\frac{\lambda_i\operatorname{Re}(z-\frac{\lambda_i}{n}\sum_{j}\frac{w_j}{1+m_{1n}w_j})}{\operatorname{Re}^2(z-\frac{\lambda_i}{n}\sum_{j}\frac{w_j}{1+m_{1n}w_j})+\operatorname{Im}^2(z-\frac{\lambda_i}{n}\sum_{j}\frac{w_j}{1+m_{1n}w_j})}=0,
		\end{split}
		\]
		while
		\[
		\begin{split}
			\operatorname{Re}(z-\frac{\lambda_i}{n}\sum_{j}\frac{w_j}{1+m_{1n}w_j})=&E-\frac{\lambda_i}{n}\sum_{j}\frac{w_j[1+w_j\operatorname{Re}m_{1n}]}{[1+w_j\operatorname{Re}m_{1n}]^2+w_j^2\operatorname{Im}^2m_{1n}},\\
			\operatorname{Im}(z-\frac{\lambda_i}{n}\sum_{j}\frac{w_j}{1+m_{1n}w_j})=&\eta+\frac{\lambda_i}{n}\sum_{j}\frac{w_j^2\operatorname{Im} m_{1n}}{[1+w_j\operatorname{Re}m_{1n}]^2+w_j^2\operatorname{Im}^2m_{1n}}.
		\end{split}
		\]
		If $\operatorname{Re}m_{1n}=m_{1n}(E)$, we have
		\[
		\bigg|\frac{1}{n}\sum_{j}\frac{w_j[1+w_j\operatorname{Re}m_{1n}]}{[1+w_j\operatorname{Re}m_{1n}]^2+w_j^2\operatorname{Im}^2m_{1n}}\bigg|\le \bigg|\frac{1}{n}\sum_{j}\frac{w_j(w_{t_1}+n^{-1/2+c})}{w_{t_1}+n^{-1/2+c}-w_j}\bigg|\le C,
		\]
		which further indicates that 
		\[
		\operatorname{Re}(z-\frac{\lambda_i}{n}\sum_{j}\frac{w_j}{1+m_{1n}w_j})>0,
		\]
		for sufficiently large $n$ and
		\[
		\begin{split}
			\operatorname{Re}f_{1n}\ge&-m_{1n}(E)-\frac{1}{n}\sum_{i=1}^p\frac{\lambda_i}{E-\frac{\lambda_i}{n}\sum_{j}\frac{w_j[1+w_j\operatorname{Re}m_{1n}]}{[1+w_j\operatorname{Re}m_{1n}]^2+w_j^2\operatorname{Im}^2m_{1n}}}\\
			> &-m_{1n}(E)-\frac{1}{n}\sum_{i=1}^p\frac{\lambda_i}{E-\frac{\lambda_i}{n}\sum_{j}\frac{w_j}{[1+w_jm_{1n}(E)]}}=0.
		\end{split}
		\]
		On the other hand, if $\operatorname{Re}m_{1n}=-\phi_n\bar\lambda(1+c)^{-1}E/(E^2+\eta^2)>-\phi_n\bar\lambda \lambda_{(1)}^{-1}$, we will have
		\[
		\begin{split}
			\frac{1}{n}\sum_{j=1}^n\frac{w_j^k}{[1+w_j\operatorname{Re}m_{1n}]^2}=&\frac{1}{n}\sum_{j=1}^{\sqrt{n}}\frac{w_j^k}{[1+w_j\operatorname{Re}m_{1n}]^2}+\frac{1}{n}\sum_{j=\sqrt{n}+1}^n\frac{w_j^k}{[1+w_j\operatorname{Re}m_{1n}]^2}\\
			\le &C\operatorname{Re}^{-2}m_{1n}(\log n)^{-2}=E\times o(1),
		\end{split}
		\]
		with high probability for any constant $k>0$. 
		Then, one can verify that
		\begin{equation}\label{ReIm}
			\operatorname{Re}(z-\frac{\lambda_i}{n}\sum_{j}\frac{w_j}{1+m_{1n}w_j})=E[1+o(1)],\quad \operatorname{Im}(z-\frac{\lambda_i}{n}\sum_{j}\frac{w_j}{1+m_{1n}w_j})=\eta+O\bigg(\frac{E^2+\eta^2}{E\log n}\bigg).
		\end{equation}
		Therefore, by definition,
		\[
		\operatorname{Re}f_{1n}=\frac{-c\phi_n\bar\lambda}{1+c}E/(E^2+\eta^2)[1+o(1)]<0.
		\]
		Then, for sufficiently large $n$, there is a solution 
		\begin{equation}\label{re m1n}
			\operatorname{Re}m_{1n}(z)\in\bigg(-\frac{1}{w_{t_1}+n^{-1/2+c}},-\frac{\phi_n\bar\lambda}{1+c}\frac{E}{(E^2+\eta^2)}\bigg).
		\end{equation}
		
		Now we focus on $\operatorname{Im} m_{1n}$. Similarly,
		\[
		\operatorname{Im}f_{1n}(z)=-\operatorname{Im}m_{1n}+\frac{1}{n}\sum_{i=1}^p\frac{\lambda_i\operatorname{Im}(z-\frac{\lambda_i}{n}\sum_{j}\frac{w_j}{1+m_{1n}w_j})}{\operatorname{Re}^2(z-\frac{\lambda_i}{n}\sum_{j}\frac{w_j}{1+m_{1n}w_j})+\operatorname{Im}^2(z-\frac{\lambda_i}{n}\sum_{j}\frac{w_j}{1+m_{1n}w_j})}=0.
		\]
		When $\operatorname{Re}m_{1n}$ satisfies (\ref{re m1n}), the results in (\ref{ReIm}) still hold. Then, if $\operatorname{Im}m_{1n}=-\eta\operatorname{Re}m_{1n}$, 
		\[
		\operatorname{Im}(z-\frac{\lambda_i}{n}\sum_{j}\frac{w_j}{1+m_{1n}w_j})=\eta+O(\eta(\log^n)^{-2}\operatorname{Re}^{-1}m_{1n})=\eta[1+o(1)],
		\]
		which further implies  that
		\[
		\begin{split}
			\operatorname{Im}f_{1n}=&\eta\operatorname{Re}m_{1n}+\frac{\phi_n\bar\lambda\eta}{E^2+\eta^2}[1+o(1)]=\eta\operatorname{Re}m_{1n}[1+o(1)]<0.
		\end{split}
		\]
		On the other hand, if $\operatorname{Im}m_{1n}=(1+c)^{-1}\phi_n\bar\lambda\eta/(E^2+\eta^2)$, we still have
		\[
		\operatorname{Im}(z-\frac{\lambda_i}{n}\sum_{j}\frac{w_j}{1+m_{1n}w_j})=\eta[1+o(1)],
		\]
		which further implies 
		\[
		\begin{split}
			\operatorname{Im}f_{1n}=&-\frac{\phi_n\bar\lambda}{1+c}\frac{\eta}{E^2+\eta^2}+\phi_n\bar\lambda\frac{\eta}{E^2+\eta^2}[1+o(1)]>0.
		\end{split}
		\]
		Therefore, for $\operatorname{Re}m_{1n}$ satisfying (\ref{re m1n}), we always have a solution $m_{1n}(z)$ satisfying
		\begin{equation}\label{Imm1n}
			\frac{\phi_n\bar\lambda}{1+c}\frac{\eta}{E^2+\eta^2}<\operatorname{Im}m_{1n}(z)<\eta|\operatorname{Re}m_{1n}|<C\eta\log ^{-1}n,\quad z\in D.
		\end{equation}
		Recall that $m_{1n}$ converges to $m_{1c}$  while $m_{1c}$ has a unique solution in $\mathbb{C}^+$. Since $\operatorname{Im}m_{1n}>0$, we claim that for sufficiently large $n$, the solutions for $\operatorname{Re}m_{1n}$ and $\operatorname{Im}m_{1n}$ are unique in $\mathbb{C}^+$, which are given by (\ref{re m1n}) and (\ref{Imm1n}).  
		
		Now we calculate $m_{2n}$ and $m_n$. By definition,
		\[
		m_{2n}=\frac{1}{n}\sum_j\frac{w_j[-E-w_j(E\operatorname{Re}m_{1n}-\eta\operatorname{Im}m_{1n})+i\eta+iw_j(E\operatorname{Im}m_{1n}+\eta\operatorname{Re}m_{1n})]}{[-E-w_j(E\operatorname{Re}m_{1n}-\eta\operatorname{Im}m_{1n})]^2+[-\eta-w_j(E\operatorname{Im}m_{1n}+\eta\operatorname{Re}m_{1n})]^2}.
		\]
		Then, when $\eta\le O(1)$, by (\ref{re m1n}) and (\ref{Imm1n}), we have
		\[
		|m_{2n}|\le \frac{C}{E|\operatorname{Re}m_{1n}|\log n}=o(1).
		\]
		On the other hand, if $\eta\rightarrow \infty$,
		\[
		|m_{2n}|\le C\frac{1}{n}\sum_j\frac{w_j}{\eta(1+w_j\operatorname{Re}m_{1n})}\le O\bigg(\frac{1}{\eta|\operatorname{Re}m_{1n}|\log n}\bigg)\le O\bigg(\frac{1}{\eta}+\frac{\eta}{E\log n}\bigg)=o(1).
		\]
		Therefore, we always have $|m_{2n}(z)|=o(1)$ for $z\in D$. By a similar procedure, we can also prove $|m_n|=o(1)$. A more careful but elementary calculation will lead to
		\[
		\operatorname{Im}m_{2n}= o(\eta),\quad \operatorname{Im}m_{n}= O(\eta(\log n)^{-1}),
		\]
		which concludes the lemma.
	\end{proof}

	Now we provide local law for large $\eta$. 
	\begin{lemma}[Average local law for large $\eta$]\label{large eta}
		For the multiplier bootstrap, if Assumptions \ref{c1}, \ref{c3} and the events $\Omega_n$ hold, then uniformly on $z\in D$ with $\eta=(\log n)^{1+c}$, it holds that
		\[
		\begin{split}
			&m_1-m_{1n}\prec n^{-1/2},\quad m_2-m_{2n}\prec n^{-1/2},\quad m-m_n\prec n^{-1/2},\\
			&\max_{i,j}\bigg(\mathcal{G}+z^{-1}(\Ib+m_{1n}\Wb)^{-1}\bigg)_{ij}\prec n^{-1/2}.
		\end{split}
		\]
	\end{lemma}
	\begin{proof}
		When $\eta=\log^{1+c}n$, directly we have $\max\{\|\Gb^{(\mathbb{T})}\|, \|\mathcal{G}^{(\mathbb{T})}\|\}\le \eta^{-1}\le (\log n)^{-1-c}$,  for any $\mathbb{T}\subset \{1,\ldots,n\}$.  To prove the lemma, we need to find the relationship between $m_1$ and $m_2$. By Lemma \ref{resolvent} and the definition of $m_2$, we have
		\begin{equation}\label{m2 and Zi}
			\begin{split}
				m_2=&\frac{1}{n}\sum_{i=1}^n\frac{w_{i}}{-z[1+n^{-1}\by_i^\top\Gb^{(i)}\by_i]}=\frac{1}{n}\sum_{i=1}^n\frac{w_{i}}{-z[1+w_in^{-1}\text{tr}\Gb^{(i)}\bSigma+Z_i]},\quad\text{where}\\
				Z_i:=&n^{-1}\by_i^\top\Gb^{(i)}\by_i-w_in^{-1}\text{tr}\Gb^{(i)}\bSigma.
			\end{split}
		\end{equation}
		Note that $\by_i$ is independent of $\Gb^{(i)}$. Then, by large deviation bounds,
		\begin{equation}\label{Zi}
			Z_i\prec \frac{w_i}{n}\bigg(\|\Gb^{(i)}\bSigma\|_F^2\bigg)^{1/2}\prec n^{-1/2}.
		\end{equation}
		On the other hand, 
		\begin{equation}\label{m1 minor}
			\frac{1}{n}\text{tr}\Gb^{(i)}\bSigma-m_1(z)=\frac{1}{n^2}\by_i^\top\Gb^{(i)}\Gb\by_i\prec n^{-1}.
		\end{equation}
		As a result, we write 
		\begin{equation}\label{m2=m1}
			m_2=\frac{1}{n}\sum_{i=1}^n\frac{w_{i}}{-z[1+w_im_1+O_{\prec}(n^{-1/2})]}=\frac{1}{n}\sum_{i=1}^n\frac{w_{i}}{-z[1+w_im_1]}+O_{\prec}(n^{-1/2}),
		\end{equation}
		where we use the fact $|1+w_im_1|\ge |1-C\log n\times \eta^{-1}|\ge C^{-1}$. 
		
		Conversely, we can also use $m_2$ to represent $m_1$. Below we show the details. By definition,
		\[
		\hat{\Sbb}-z\Ib=\frac{1}{n}\sum_{i}\by_i\by_i^\top+zm_2(z)\bSigma-z[\Ib+m_2(z)\bSigma].
		\]
		Taking inverse on both sides, 
		\[
		\Gb=-z^{-1}[\Ib+m_2(z)\bSigma]^{-1}+z^{-1}\Gb\bigg[\frac{1}{n}\sum_{i}\by_i\by_i^\top+zm_2(z)\bSigma\bigg][\Ib+m_2(z)\bSigma]^{-1}.
		\]
		By elementary matrix inverse formulas,
		\[
		\begin{split}
			\Gb\by_i=\Gb^{(i)}\by_i-n^{-1}\Gb\by_i\by_i^\top\Gb^{(i)}\by_i=\frac{1}{1+n^{-1}\by_i^\top\Gb^{(i)}\by_i}\Gb^{(i)}\by_i.
		\end{split}
		\]
		Therefore, we can write
		\begin{equation}\label{R1+R2}
			\begin{split}
				\Gb=&-z^{-1}[\Ib+m_2(z)\bSigma]^{-1}+z^{-1}\frac{1}{n}\sum_i\frac{\Gb^{(i)}[\by_i\by_i^\top-w_i\bSigma]}{1+n^{-1}\by_i^\top\Gb^{(i)}\by_i}[\Ib+m_2(z)\bSigma]^{-1}\\
				&+z^{-1}\frac{1}{n}\sum_i\frac{w_i[\Gb^{(i)}-\Gb]\bSigma}{1+n^{-1}\by_i^\top\Gb^{(i)}\by_i}[\Ib+m_2(z)\bSigma]^{-1}\\
				:=&-z^{-1}[\Ib+m_2(z)\bSigma]^{-1}+\Rb_1+\Rb_2.
			\end{split}
		\end{equation}
		In the following, we bound the error terms. For $\Rb_1$, 
		\begin{equation}\label{R11+R12}
			\begin{split}
				&\frac{z}{n}\text{tr}\Rb_1\bSigma=\frac{1}{n^2}\sum_i\text{tr}\bigg(\frac{\Gb^{(i)}[\by_i\by_i^\top-w_i\bSigma]}{1+n^{-1}\by_i^\top\Gb^{(i)}\by_i}[\Ib+m_2^{(i)}(z)\bSigma]^{-1}\bSigma\bigg)\\
				&+\frac{1}{n^2}\sum_i\text{tr}\bigg(\frac{\Gb^{(i)}[\by_i\by_i^\top-w_i\bSigma]}{1+n^{-1}\by_i^\top\Gb^{(i)}\by_i}[\Ib+m_2(z)\bSigma]^{-1}[m_2^{(i)}(z)-m_2(z)]\bSigma[\Ib+m_2^{(i)}(z)\bSigma]^{-1}\bSigma\bigg)\\
				:=&R_{11}+R_{12}.
			\end{split}
		\end{equation}
		Since $\eta=\log^{1+c}n$, we have $\|\Gb^{(i)}\|\rightarrow 0$, $|m_2^{(i)}|\rightarrow 0$, and with high probability $|1+n^{-1}\by_i^\top\Gb^{(i)}\by_i|\ge c$. Then, for each $i$,
		\[
		\frac{1}{n}\text{tr}\bigg(\frac{\Gb^{(i)}[\by_i\by_i^\top-w_i\bSigma]}{1+n^{-1}\by_i^\top\Gb^{(i)}\by_i}[\Ib+m_2^{(i)}(z)\bSigma]^{-1}\bSigma\bigg)\prec \frac{1}{n}\|\bSigma\|_F\prec n^{-1/2},
		\]
		which further indicates $R_{11}\prec n^{-1/2}$. For $R_{12}$, note that
		\[
		m_2^{(i)}(z)-m_2(z)=\frac{1}{n}\sum_{\mu\ne i}w_\mu[\mathcal{G}_{\mu\mu}-\mathcal{G}_{\mu\mu}^{(i)}]=\frac{1}{n}\sum_{\mu\ne i}\frac{\mathcal{G}_{i\mu}\mathcal{G}_{\mu i}}{\mathcal{G}_{ii}},
		\]
		while by Lemma \ref{resolvent},
		\[
		\frac{1}{\mathcal{G}_{ii}}=-z-zn^{-1}\by_i^\top\Gb^{(i)}\by_i\le C|z|,\quad \mathcal{G}_{ij}(z)\le C|z|n^{-1}\|\Gb^{(ij)}\|_F\prec n^{-1/2}, \quad i\ne j.
		\]
		Therefore, we can conclude that
		\[
		m_2^{(i)}(z)-m_2(z)\prec n^{-1}\Rightarrow R_{12}\prec n^{-1}\Rightarrow \frac{z}{n}\text{tr}\Rb_1\bSigma\prec n^{-1/2}.
		\]
		
		For $\Rb_2$, note that
		\[
		\frac{1}{n}\text{tr}(\Gb^{(i)}-\Gb)\bSigma[\Ib+m_2(z)\bSigma]^{-1}\bSigma\le\bigg| \frac{1}{n^2}\by_i^\top\Gb^{(i)}\bSigma[\Ib+m_2(z)\bSigma]^{-1}\bSigma\Gb\by_i\bigg|\prec n^{-1}.
		\]
		Then, directly we have 
		\[
		\frac{1}{n}\text{tr}\Rb_2\bSigma\prec n^{-1}.
		\]
		Consequently, we have
		\begin{equation}\label{m1=m2}
			m_1=\frac{1}{n}\text{tr}\Gb\bSigma=-z^{-1}\frac{1}{n}\text{tr}[\Ib+m_2\bSigma]^{-1}\bSigma+O_{\prec}(n^{-1/2})=-\frac{1}{n}\sum_{i=1}^p\frac{\lambda_i}{z[1+\lambda_im_2]}+O_{\prec}(n^{-1/2}).
		\end{equation}
		
		Combine (\ref{m2=m1}) and (\ref{m1=m2}) to get
		\[
		\begin{split}
			&m_2-m_{2n}=\frac{1}{n}\sum_{i=1}^n\frac{w_i^2[m_1-m_{1n}]}{-z[1+w_im_1][1+w_im_{1n}]}+O_{\prec}(n^{-1/2})\\
			=&\bigg(\frac{1}{n}\sum_{i=1}^n\frac{w_i^2}{-z[1+w_im_1][1+w_im_{1n}]}\bigg)\bigg(\frac{1}{n}\sum_{i=1}^p\frac{\lambda_i^2[m_2-m_{2n}]}{-z[1+\lambda_im_2][1+\lambda_im_{2n}]}\bigg)+O_{\prec}(n^{-1/2}).
		\end{split}
		\]
		Since $|z(1+w_im_{1n})|\ge c$, $|1+w_im_1|\ge c$, $|1+\lambda_im_2|\ge c$,  $|1+\lambda_im_{2n}|\ge c$, and $|z|\rightarrow \infty$, we conclude that
		\begin{equation}\label{m2-m2n}
			m_2-m_{2n}=O(|z|^{-1})(m_2-m_{2n})+O_{\prec}(n^{-1/2}),
		\end{equation}
		which further implies $m_2-m_{2n}\prec n^{-1/2}$. By a parallel procedure, we have $m_1-m_{1n}\prec n^{-1/2}$. 
		
		Next, we show the result for $m(z)$ (abbreviated as $m$). Indeed, similarly to (\ref{m1=m2}) and (\ref{m2-m2n}), we can easily conclude that
		\[
		m=\frac{1}{p}\text{tr}\Gb=-\frac{1}{n}\sum_{i=1}^p\frac{1}{z[1+\lambda_im_2]}+O_{\prec}(n^{-1/2})=m_{n}+O_{\prec}(n^{-1/2}).
		\]
		Since $m, m_1, m_2$ are  Lipschitz on $z$ with Lipschitz coefficient $n^2$, the results hold uniformly on $z$ by a standard lattice technique. For example, see the argument below (5.51) in \cite{kwak2021extremal}.
		
		The last step is to prove the result for $\mathcal{G}$. For the diagonal entries, by Lemma \ref{resolvent},
		\[
		\begin{split}
			\mathcal{G}_{ii}=&-\frac{1}{z[1+n^{-1}\by_i^\top\Gb^{(i)}\by_i]}=-\frac{1}{z[1+w_in^{-1}\text{tr}\Gb^{(i)}\bSigma+O_{\prec}(n^{-1/2})]}\\
			=&-\frac{1}{z[1+w_im_1+O_{\prec}(n^{-1/2})]}=-\frac{1}{z[1+w_im_{1n}]}+O_{\prec}(n^{-1/2}).
		\end{split}
		\]
		For the off-diagonal entries,  we have
		\[
		|\mathcal{G}_{ij}|\le z|\mathcal{G}_{ii}||\mathcal{G}_{jj}^{(i)}||n^{-1}\by_i^\top\Gb^{(ij)}\by_j|\prec n^{-1}\|\Gb^{(ij)}\|_F\prec n^{-1/2}.
		\]
		The lemma is then verified.
	\end{proof}

	The next step is to show that the results in Lemma \ref{large eta} also hold for small $\eta$. We need the following self-improvement lemma.
	\begin{lemma}[Self-improvement]\label{self-improvement}
		Under multiplier bootstrap, Assumptions \ref{c1}, \ref{c3} and the events $\Omega_n$,  for any $z\in D$, if 
		\[
		\begin{split}
			\psi_1:=&\max_{i,j}\bigg(\mathcal{G}+z^{-1}(\Ib+m_{1n}\Wb)^{-1}\bigg)_{ij}\prec n^{-1/2+c},\\
			\psi_2:=&|m_1-m_{1n}|+|m_2-m_{2n}|+|m-m_n|\prec n^{-1/2+c},
		\end{split}
		\]
		for some constant $0<c<1/5$, then we have
		\[
		\psi_1\prec n^{-1/2},\quad \psi_2\prec n^{-1/2}.
		\]
	\end{lemma}
	\begin{proof}
		We essentially follow the same strategy as that in the proof of Lemma \ref{large eta}. The major difference is that $\|\mathcal{G}\|\prec 1$ no longer holds because $\eta$ can be very small.  To overcome this challenge, we will mainly rely on the preliminary bounds of $\psi_1$ and $\psi_2$ to control all the error terms. Below we show the details.
		
		We first use $m_1$ to represent $m_2$. By the priori bound  of $\psi_1$ and Lemma \ref{m_{1n}}, we conclude that $|\mathcal{G}_{ii}|\prec 1$ and $|G_{ij}|\prec n^{-1/2+c}$ for $i\ne j$. To prove (\ref{m2=m1}) for small $\eta$, note that (\ref{m2 and Zi}) still holds but we need to reconsider (\ref{Zi}) and (\ref{m1 minor}). For $Z_i$,
		\[
		\begin{split}
			Z_i\prec&\frac{w_i}{n}\|\Gb^{(i)}\bSigma\|_F\prec n^{-1}\|\Gb^{(i)}\|_F\prec n^{-1}\|\mathcal{G}^{(i)}\|_F+n^{-1/2}=n^{-1/2}\sqrt{\frac{\operatorname{Im} \tilde m^{(i)}}{n\eta}}+n^{-1/2}\\
			\prec &n^{-1/2}\bigg(\sqrt{\frac{\operatorname{Im}\tilde m}{n\eta}}+\sqrt{\frac{|\tilde m-\tilde m^{(i)}|}{n\eta}}+1\bigg)\prec n^{-1/2}\bigg(1+\sqrt{\frac{1}{n\eta}\frac{1}{n}\sum_{j\ne i}\frac{|\mathcal{G}_{ij}|^2}{|\mathcal{G}_{ii}|}}\bigg)\\
			\prec & n^{-1/2}.
		\end{split}
		\]
		On the other hand, 	
		\[
		\begin{split}
			\frac{1}{n}\text{tr}\Gb^{(i)}\bSigma-m_1(z)=&\frac{1}{n^2}\by_i^\top\Gb^{(i)}\bSigma\Gb\by_i=\frac{1}{n^2}\by_i^\top\Gb^{(i)}\bSigma\Gb^{(i)}\by_i-\frac{1}{n^3}\by_i^\top\Gb^{(i)}\bSigma\Gb\by_i\by_i^\top\Gb^{(i)}\by_i\\
			=&\frac{\frac{1}{n^2}\by_i^\top\Gb^{(i)}\bSigma\Gb^{(i)}\by_i}{1+n^{-1}\by_i^\top\Gb^{(i)}\by_i}\prec \frac{\frac{1}{n^2}\|\Gb^{(i)}\|_F^2}{|1+n^{-1}\by_i^\top\Gb^{(i)}\by_i|}\prec \frac{\frac{1}{n^2}\|\mathcal{G}^{(i)}\|_F^2+n^{-1}}{|1+n^{-1}\by_i^\top\Gb^{(i)}\by_i|}\\
			=&\frac{\frac{1}{n\eta}\operatorname{Im}\tilde m^{(i)}+n^{-1}}{|1+n^{-1}\by_i^\top\Gb^{(i)}\by_i|}\prec\frac{\frac{1}{n\eta}\operatorname{Im}\tilde m+\frac{1}{n\eta}\times n^{-1+2c}+n^{-1}}{|1+n^{-1}\by_i^\top\Gb^{(i)}\by_i|}\\
			\le&\frac{\frac{1}{n\eta}\operatorname{Im} m+\frac{1}{n\eta}|\operatorname{Im}z^{-1}|+n^{-1}}{|1+n^{-1}\by_i^\top\Gb^{(i)}\by_i|}\prec \frac{n^{-1}+(n\eta)^{-1}|\operatorname{Im}m-\operatorname{Im}m_n|}{|1+n^{-1}\by_i^\top\Gb^{(i)}\by_i|}\\
			\prec&n^{-5/6+c} |z\mathcal{G}_{ii}|\prec n^{-1/2}.
		\end{split}
		\]
		Therefore,
		\begin{equation}\label{m2=m1 general}
			\begin{split}
				m_2=&\frac{1}{n}\sum_{i=1}^n\frac{w_{i}}{-z[1+w_im_1+O_{\prec}(n^{-1/2})]}\\=&\frac{1}{n}\sum_{i=1}^n\frac{w_{i}}{-z[1+w_im_1]}+\frac{1}{n}\sum_{i=1}^n\frac{w_{i}O_{\prec}(n^{-1/2})}{-z[1+w_im_1][1+w_im_1+O_{\prec}(n^{-1/2})]}\\
				=&\frac{1}{n}\sum_{i=1}^n\frac{w_{i}}{-z[1+w_im_1]}+\frac{1}{n}\frac{w_{t_1}O_{\prec}(n^{-1/2})}{-z[1+w_{t_1}m_1][1+w_{t_1}m_1+O_{\prec}(n^{-1/2})]}\\
				&+\frac{1}{n}\sum_{i=2}^n\frac{w_{t_i}O_{\prec}(n^{-1/2})}{-z[1+w_{t_i}m_1][1+w_{t_i}m_1+O_{\prec}(n^{-1/2})]}\\
				=&\frac{1}{n}\sum_{i=1}^n\frac{w_{i}}{-z[1+w_im_1]}+O_{\prec}(n^{-1/2}),
			\end{split}
		\end{equation}
		where we use the facts that
		\[
		\begin{split}
			1+w_{t_1}m_1=&1+w_{t_1}m_{1n}+w_{t_1}(m_1-m_{1n})\ge 1-\frac{w_{t_1}}{w_{t_1}+n^{-1/2+c}}+O_{\prec}(n^{-1/2+c})\\
			\ge &n^{-1/2-2c},	
		\end{split}
		\]
		and meanwhile for $i\ge 2$,
		\[
		\begin{split}
			1+w_{t_i}m_1=&1+w_{t_i}m_{1n}+w_{t_i}(m_1-m_{1n})\ge 1+w_{t_2}m_{1n}+O_{\prec}(n^{-1/2+c})\\
			\ge& [w_{t_1}-w_{t_2}]m_{1n}+O_{\prec}(n^{-1/2+2c})\ge c(\log n)^{-2}.
		\end{split}
		\]
		
		Next, we aim to represent $m_1$ using $m_2$, i.e., prove (\ref{m1=m2}) based on the priori bounds of $\psi_1$ and $\psi_2$. The decompositions in (\ref{R1+R2}) and (\ref{R11+R12}) still hold, but we need to reconsider how to bound $R_{11}$, $R_{12}$ and $\Rb_2$. For $R_{11}$, note that 
		\[
		\begin{split}
			m_2^{(i)}=&m_2^{(i)}-m_2+m_2-m_{2n}+m_{2n}
			\prec\frac{1}{n}\sum_{\mu\ne i}\frac{|\mathcal{G}_{i\mu}|^2}{|\mathcal{G}_{ii}|}+O_{\prec}(n^{-1/2+c})+m_{2n}\\
			=&m_{2n}+O_{\prec}(n^{-1/2+c})\rightarrow 0.
		\end{split}
		\]
		Then, for each $i$,
		\[
		\begin{split}
			&\frac{1}{n}\text{tr}\bigg(\frac{\Gb^{(i)}[\by_i\by_i^\top-w_i\bSigma]}{1+n^{-1}\by_i^\top\Gb^{(i)}\by_i}[\Ib+m_2^{(i)}(z)\bSigma]^{-1}\bSigma\bigg)\prec |z\mathcal{G}_{ii}|\frac{1}{n}\|\bSigma[\Ib+m_2^{(i)}(z)\bSigma]^{-1}\bSigma\Gb^{(i)}\|_F\\
			\prec &\frac{1}{n}\|\Gb^{(i)}\|_F\prec n^{-1/2},
		\end{split}
		\]
		where the last step follows from the bound of $Z_i$. This further indicates that $R_{11}\prec n^{-1/2}$.  On the other hand, for $R_{12}$,  since $m_2^{(i)}-m_2\prec n^{-1+2c}$ and $m_2^{(i)}\rightarrow 0$ with high probability, we have
		\[
		\begin{split}
			R_{12}\prec& \frac{n^{-1+2c}}{n^2}\sum_i\text{tr}\bigg(\frac{\Gb^{(i)}[\by_i\by_i^\top-w_i\bSigma]}{1+n^{-1}\by_i^\top\Gb^{(i)}\by_i}[\Ib+m_2(z)\bSigma]^{-1}\bSigma[\Ib+m_2^{(i)}(z)\bSigma]^{-1}\bSigma\bigg)\\
			\prec &	\frac{n^{-1+2c}}{n^2}\sum_i\text{tr}\bigg(\frac{\Gb^{(i)}[\by_i\by_i^\top-w_i\bSigma]}{1+n^{-1}\by_i^\top\Gb^{(i)}\by_i}[\Ib+m_2^{(i)}(z)\bSigma]^{-1}\bSigma[\Ib+m_2^{(i)}(z)\bSigma]^{-1}\bSigma\bigg)\\
			&+O_{\prec}\bigg(\frac{n^{-2+4c}}{n^2}\sum_i|\by_i^\top\Gb^{(i)}\by_i|+n\|\bSigma\|\bigg)\\
			\prec &\frac{n^{-1+2c}}{n^2}\sum_i\|\Gb^{(i)}\|_F+O_{\prec}\bigg(n^{-2+4c}(\eta^{-1}+1)\bigg)\prec n^{-1/2}.
		\end{split}
		\]
		Lastly, for $\Rb_2$,
		\[
		\begin{split}
			&\frac{1}{n}\text{tr}(\Gb^{(i)}-\Gb)\bSigma[\Ib+m_2(z)\bSigma]^{-1}\bSigma\le\bigg| \frac{1}{n^2}\by_i^\top\Gb^{(i)}\bSigma[\Ib+m_2(z)\bSigma]^{-1}\bSigma\Gb\by_i\bigg|\\
			= &\bigg|\frac{1}{1+n^{-1}\by_i^\top\Gb^{(i)}\by_i} \frac{1}{n^2}\by_i^\top\Gb^{(i)}\bSigma[\Ib+m_2(z)\bSigma]^{-1}\bSigma\Gb^{(i)}\by_i\bigg|\\
			\prec &|z\mathcal{G}_{ii}|n^{-2}\|\Gb^{(i)}\|_F^2\prec n^{-1}.
		\end{split}
		\]
		Consequently, (\ref{m1=m2}) still holds. The remaining  proof is almost the same as that of Lemma \ref{large eta}, and we omit the details.
	\end{proof}

	Now we can show the local law for $z\in D$.
	\begin{lemma}[Average local law]\label{average local law}
		Let Assumptions \ref{c1}, \ref{c3} and the events $\Omega_n$ hold. Under the multiplier bootstrap, it holds uniformly on $z\in D$ that $\psi_1\prec n^{-1/2}$ and $\psi_2\prec n^{-1/2}$, where $\psi_1$ and $\psi_2$ are defined in Lemma \ref{self-improvement}.
	\end{lemma}
	\begin{proof}
		We use a standard  discrete continuity argument to prove the result. For each $z=E+i\eta\in D$, fix $E$ and consider a sequence $\{\eta_j\}$ defined by $\eta_j=\log^{1+c}n-jn^{-2}$. Then, $\eta$ must fall in an interval $[\eta_{j-1},\eta_j]$ for some $1\le j\le Cn^3$.
		
		We start with $\eta_0$ and use induction to complete the proof. For $\eta_0$, Lemma \ref{large eta} already indicates the results. Now assume $\psi_1\prec n^{-1/2}$ and $\psi_2\prec n^{-1/2}$ for some $j=K$. For any $\eta^\prime$ satisfying $\eta_{j-1}\le \eta^\prime\le \eta_j$, write $z^\prime=E+i\eta^\prime$ and $z_j=E+i\eta_j$. We then always have
		\begin{equation}\label{difference of G}
			\|\mathcal{G}(z^\prime)-\mathcal{G}(z_j)\|\le \frac{|\eta^\prime-\eta_j|}{|\eta_j\eta^\prime|}\le n^{-2/3},
		\end{equation}
		where we use the fact that $\eta\ge n^{-2/3}$ for $z\in D$. On the other hand, 
		\[
		\begin{split}
			m_{1n}(z^\prime)-m_{1n}(z_j)=&\frac{1}{n}\sum_i\bigg(\frac{\lambda_i}{-z^\prime[1+\lambda_im_{2n}(z^\prime)]}-\frac{\lambda_i}{-z_j[1+\lambda_im_{2n}(z_j)]}\bigg)\\
			=&\frac{1}{n}\sum_i\bigg(\frac{\lambda_i}{-z^\prime[1+\lambda_im_{2n}(z^\prime)]}-\frac{\lambda_i}{-z^\prime[1+\lambda_im_{2n}(z_j)]}\bigg)\\
			&+\frac{1}{n}\sum_i\bigg(\frac{\lambda_i}{-z^\prime[1+\lambda_im_{2n}(z_j)]}-\frac{\lambda_i}{-z_j[1+\lambda_im_{2n}(z_j)]}\bigg)\\
			:=&\mathcal{J}_{1}+\mathcal{J}_2.
		\end{split}
		\]
		For $\mathcal{J}_1$, we have
		\[
		\begin{split}
			\mathcal{J}_1=&\frac{1}{n}\sum_i\frac{\lambda_i^2}{-z^\prime[1+\lambda_im_{2n}(z^\prime)][1+\lambda_im_{2n}(z_j)]}[m_{2n}(z_j)-m_{2n}(z^\prime)]\\
			=&o(1)\times \frac{1}{n}\sum_{i}\bigg(\frac{w_i}{-z_j[1+w_im_{1n}(z_j)]}-\frac{w_i}{-z^\prime[1+w_im_{1n}(z^\prime)]}\bigg)\\
			=&o(1)\times \frac{1}{n}\sum_{i}\bigg(\frac{w_i}{-z_j[1+w_im_{1n}(z_j)]}-\frac{w_i}{-z_j[1+w_im_{1n}(z^\prime)]}\bigg)\\
			&+o(1)\times \frac{1}{n}\sum_{i}\bigg(\frac{w_i}{-z_j[1+w_im_{1n}(z^\prime)]}-\frac{w_i}{-z^\prime[1+w_im_{1n}(z^\prime)]}\bigg)\\
			=&o(1)\times [m_{1n}(z^\prime)-m_{1n}(z_j)]+O_{\prec}(z_j-z^\prime)\\
			=&o(1)\times [m_{1n}(z^\prime)-m_{1n}(z_j)]+O_{\prec}(n^{-2}).
		\end{split}
		\]
		Similarly, for $\mathcal{J}_2$, we have
		\[
		\mathcal{J}_2\prec O(n^{-2}).
		\]
		Therefore, 
		\begin{equation}\label{z prime}
			|m_{1n}(z^\prime)-m_{1n}(z_j)|\prec n^{-2}.
		\end{equation}
		By similar procedures, we have
		\[
		\begin{split}
			&|m_n(z^\prime)-m_{n}(z_j)|\prec n^{-2},\quad |m_{2n}(z^\prime)-m_{2n}(z_j)|\prec n^{-2}, \\
			& \|(z^\prime)^{-1}(\Ib+m_{1n}(z^\prime)\Wb)^{-1}-z_j^{-1}(\Ib+m_{1n}(z_j)\Wb)^{-1}\|\prec n^{-2}.
		\end{split}
		\]
		Therefore, combining (\ref{difference of G}), we  conclude that
		\[
		\begin{split}
			&\|\mathcal{G}(z^\prime)-(z^\prime)^{-1}(\Ib+m_{1n}(z^\prime)\Wb)^{-1}\|\\
			\le& \|\mathcal{G}(z^\prime)-\mathcal{G}(z_j)\|+\|\mathcal{G}(z_j)-z_j^{-1}(\Ib+m_{1n}(z_j)\Wb)^{-1}\|\\
			&+\|(z^\prime)^{-1}(\Ib+m_{1n}(z^\prime)\Wb)^{-1}-z_j^{-1}(\Ib+m_{1n}(z_j)\Wb)^{-1}\|\prec n^{-1/2}.
		\end{split}
		\]
		Similarly, we have
		\[
		|m_1(z^\prime)-m_{1n}(z^\prime)|\prec n^{-1/2},\quad |m_2(z^\prime)-m_{2n}(z^\prime)|\prec n^{-1/2},\quad |m(z^\prime)-m_{n}(z^\prime)|\prec n^{-1/2}.
		\]
		Then, by induction and Lemma \ref{self-improvement}, we conclude that $\psi_1\prec n^{-1/2}$ and $\psi_2\prec n^{-1/2}$ for any $z\in D$. The uniform bound is by standard lattice argument and we omit details.
	\end{proof}

	\begin{lemma}[upper bound for eigenvalues]\label{upper bound}
		Let $r=0$, and Assumptions \ref{c1}, \ref{c3} and the events $\Omega_n$ hold. Under the multiplier bootstrap, with high probability there is no eigenvalue of $\hat{\Sbb}$ in the interval $(\lambda_{(1)},C\log n)$ for some large constant $C>0$ as $n\rightarrow\infty$.
	\end{lemma}
	\begin{proof}
		We prove the lemma by indirect argument. Assume that there is an eigenvalue of $\hat{\Sbb}$ in the interval, denoted as $\hat\lambda$. Then, we let $z=\hat\lambda+in^{-2/3}$. Since $z\in D$, by Lemma \ref{m_{1n}} we have $\operatorname{Im}m_n(z)=o(\eta)=o(n^{-2/3})$. Therefore,
		\begin{equation}\label{Im mz}
			\begin{split}
				\operatorname{Im}m(z)=&\operatorname{Im}m_n(z)+\operatorname{Im}[m(z)-m_n(z)]\\
				\le&o(n^{-2/3})+\max_{z\in D}|m(z)-m_n(z)|=O_{\prec}(n^{-1/2}).
			\end{split}
		\end{equation}
		However, by the definition of $m(z)$, we know that
		\[
		\operatorname{Im}m(z)\ge \frac{1}{n}\operatorname{Im}\frac{1}{\hat\lambda-z}=n^{-1}\eta^{-1}=n^{-1/3},
		\]
		which is a contradiction to (\ref{Im mz}). Therefore, there is no eigenvalue in this interval.
	\end{proof}
	
	Now we can consider the limiting properties of the largest eigenvalue of $\hat{\Sbb}$, i.e., $\hat \lambda_1$ for $r=0$. We have the next lemma.
	\begin{lemma}[Eigenvalue rigidity]\label{Largest eigenvalue bound}
		Let $r=0$, and Assumptions \ref{c1}, \ref{c3} and the events $\Omega_n$ hold. Under the multiplier bootstrap,  it holds that $|\hat\lambda_1-\lambda_0|\prec n^{-1/2+2c}$.
	\end{lemma}
	\begin{proof}
		By definition,  $\hat\lambda_1$ is a non-zero eigenvalue of $\hat{\Sbb}$ if
		\[
		\det(\hat\lambda_1\Ib-\hat{\Sbb})=0.
		\]
		Define $\hat{\Sbb}^{(1)}$ as
		\[
		\hat{\Sbb}^{(1)}:=\hat{\Sbb}-n^{-1}\by_{t_1}\by_{t_1}^\top=n^{-1}\Yb^{(1)}\Wb^{(1)}(\Yb^{(1)})^\top,
		\]
		where $\by_{t_1}$ is corresponding  to $w_{t_1}$, $\Yb^{(1)}$ is obtained by removing $\by_{t_1}$ from $\Yb$, and $\Wb^{(1)}$ is obtained by removing the row and column corresponding to $w_{t_1}$ from $\Wb$. 
		We can assume that $\hat\lambda_1$ is not an eigenvalue of $\hat{\Sbb}^{(1)}$ because $w_j$'s follow continuous distribution. Then,
		\[
		\begin{split}
			&\det(\hat\lambda_1\Ib-n^{-1}w_{t_1}\bx_{t_1}\bx_{t_1}^\top-\hat{\Sbb}^{(1)})=0
			\Longrightarrow 1+n^{-1}w_{t_1}\bx_{t_1}^\top[\hat{\Sbb}^{(1)}-\hat\lambda_1\Ib]^{-1}\bx_{t_1}=0.
		\end{split}
		\]
		Define 
		\[
		h(\lambda):=1+n^{-1}w_{t_1}\bx_{t_1}^\top[\hat{\Sbb}^{(1)}-\lambda\Ib]^{-1}\bx_{t_1}.
		\]
		We aim to show that $h(\lambda)$ will change sign when $\lambda$ grows from $\lambda_0-n^{-1/2+2c}$ to $\lambda_0+n^{-1/2+2c}$. 
		
		The first step is to find approximation to $n^{-1}\bx_{t_1}^\top[\hat{\Sbb}^{(1)}-\lambda\Ib]^{-1}\bx_{t_1}$. 
		Given the order of $w_j$'s, we provide an upper bound for the eigenvalues of $\hat{\Sbb}^{(1)}$. Actually, following the proof of Lemma \ref{large eta} to Lemma \ref{upper bound},  a direct upper bound for $\lambda_1(\hat{\Sbb}^{(1)})$ is $\lambda_{(2)}$, which is defined as
		\[
		1+(w_{t_2}+n^{-1/2+c})m_{1n}^{(1)}(\lambda_{(2)})=0,
		\]
		where $m_{1n}^{(1)}(z)$ is the solution in $\mathbb{C}^+$ to
		\[
		m_{1n}^{(1)}(z)=\frac{1}{n}\sum_{i=1}^p\frac{\lambda_i}{-z\bigg[1+\frac{\lambda_i}{n}\sum_{j=2}^n\frac{w_{t_j}}{-z(1+w_{t_j}m_{1n}^{(1)}(z))}\bigg]}, \quad z\in \mathbb{C}^+,
		\]
		and $m_{1n}^{(1)}(\lambda)=\lim_{\eta\downarrow0}m_{1n}^{(1)}(\lambda+i\eta)$. With $\lambda_{(2)}$, we can define a new region $D_{(2)}\subset \mathbb{C}^+$ replacing $\lambda_{(1)}$ with $\lambda_{(2)}$ in the definition of $D$.
		
		Before moving forward, we need to calculate the gap between $\lambda_0$ and $\lambda_{(2)}$.
		By definition, 
		\[
		\begin{split}
			\frac{\lambda_0}{w_{t_1}}=&\frac{1}{n}\sum_{i=1}^p\frac{\lambda_i}{1-\frac{w_{t_1}}{\lambda_0}\frac{\lambda_i}{n}\sum_{j=2}^n\frac{w_{t_j}}{w_{t_1}-w_{t_j}}},\\
			\frac{\lambda_{(2)}}{w_{t_2}+n^{-1/2+c}}=&\frac{1}{n}\sum_{i=1}^p\frac{\lambda_i}{1-\frac{w_{t_2}+n^{-1/2+c}}{\lambda_{(2)}}\frac{\lambda_i}{n}\sum_{j=2}^n\frac{w_{t_j}}{w_{t_2}+n^{-1/2+c}-w_{t_j}}}.
		\end{split}
		\]
		Similarly to (\ref{sum}) and (\ref{ratio}), one can conclude that
		\[
		\begin{split}
			&\frac{1}{n}\sum_{j=2}^n\frac{w_{t_j}}{w_{t_2}+n^{-1/2+c}-w_{t_j}}\le C(\log n)^{-1},\quad 	\frac{\lambda_{(2)}}{w_{t_2}+n^{-1/2+c}}\asymp C,\\
			&\frac{1}{n}\sum_{j=1}^n\frac{w_{t_j}}{w_{t_1}-w_{t_j}}\le C(\log n)^{-1},\quad  \frac{\lambda_0}{w_{t_1}}\asymp C.
		\end{split}
		\]
		Then, 
		\[
		\begin{split}
			&\frac{\lambda_0}{w_{t_1}}-	\frac{\lambda_{(2)}}{w_{t_2}+n^{-1/2+c}}=\bigg(\frac{1}{n}\sum_i\lambda_i^2+o(1)\bigg)\\
			&\times \bigg(\frac{w_{t_1}}{\lambda_0}\frac{1}{n}\sum_{j=2}^n\frac{w_{t_j}}{w_{t_1}-w_{t_j}}-\frac{w_{t_2}+n^{-1/2+c}}{\lambda_{(2)}}\frac{1}{n}\sum_{j=2}^n\frac{w_{t_j}}{w_{t_2}+n^{-1/2+c}-w_{t_j}}\bigg)\\
			=&\bigg(\frac{w_{t_1}}{\lambda_0}-\frac{w_{t_2}+n^{-1/2+c}}{\lambda_{(2)}}\bigg)O((\log n)^{-1})+\frac{w_{t_2}+n^{-1/2+c}}{\lambda_{(2)}}\times \frac{C}{n}\frac{w_{t_1}}{n^{-1/2+c}}\\
			=&O((\log n)^{-1})\bigg(\frac{\lambda_0}{w_{t_1}}-	\frac{\lambda_{(2)}}{w_{t_2}+n^{-1/2+c}}\bigg)+O(n^{-3/4}).
		\end{split}
		\]
		As a result, we write
		\[
		\frac{\lambda_0}{w_{t_1}}-	\frac{\lambda_{(2)}}{w_{t_2}+n^{-1/2+c}}=O(n^{-3/4}).
		\]
		This further indicates that
		\[
		\begin{split}
			\lambda_0-\lambda_{(2)}=&w_{t_1}\bigg(\frac{\lambda_{0}}{w_{t_1}}-\frac{\lambda_{(2)}}{w_{t_2}+n^{-1/2+c}}\bigg)+\lambda_{(2)}\bigg(\frac{w_{t_1}}{w_{t_2}+n^{-1/2+c}}-1\bigg)\\
			=&O(n^{-3/4+c})+\frac{\lambda_{(2)}}{w_{t_2}+n^{-1/2+c}}(w_{t_1}-w_{t_2}-n^{-1/2+c}).
		\end{split}
		\]
		Therefore, we conclude that
		\begin{equation}\label{lambda 02 gap}
			\lambda_0-\lambda_{(2)}\ge c(\log n)^{-c}.
		\end{equation}
		
		Now we calculate the sign of $h(\lambda)$ when $\lambda$ takes values of the two end points. Let $\lambda=\lambda_0-n^{-1/2+2c}$. Then, $\lambda$ is larger than $\lambda_{(2)}$ for sufficiently large $n$. Conditional on the order of $\{w_i\}$, we always have $\bx_{t_1}$ is independent of $\hat{\Sbb}^{(1)}$ and $\lambda_{0}$. Then,
		\[
		n^{-1}\bx_{t_1}^\top[\hat{\Sbb}^{(1)}-\lambda\Ib]^{-1}\bx_{t_1}-m_1^{(1)}(\lambda)\prec n^{-1}\|\hat{\Sbb}^{(1)}-\lambda\Ib\|_F\prec n^{-1/2},
		\]
		where we use the fact $|\lambda-\lambda_1(\hat{\Sbb}^{(1)})|\ge c(\log n)^{-1}$ because $\lambda_1(\hat{\Sbb}^{(1)})\le \lambda_{(2)}$ similarly to Lemma \ref{upper bound}. Meanwhile, let $z=\lambda+in^{-1/2-c}$, so $z\in D_{(2)}$ and
		\[
		m_1^{(1)}(\lambda)=m_1^{(1)}(\lambda)-m_1^{(1)}(z)+m_1^{(1)}(z)-m_{1n}^{(1)}(z)+m_{1n}^{(1)}(z)-m_{1n}^{(1)}(\lambda)+m_{1n}^{(1)}(\lambda).
		\]
		Similarly to Lemma \ref{average local law}, we have $|m_1^{(1)}(z)-m_{1n}^{(1)}(z)|\prec n^{-1/2}$. For the fist term, let $\hat\lambda_j^{(1)}$ and $\hat\bgamma_j^{(1)}$ be the eigenvalues and eigenvectors of $\hat{\Sbb}^{(1)}$, so
		\[
		m_1^{(1)}(\lambda)-m_1^{(1)}(z)=\frac{1}{n}\sum_{j=1}^p|\bSigma^{1/2}\hat\bgamma_j^{(1)}|^2\bigg(\frac{1}{\hat\lambda_j^{(1)}-\lambda}-\frac{1}{\hat\lambda_j^{(1)}-z}\bigg).
		\]
		Recall that $\lambda-\hat\lambda_j^{(1)}\ge c(\log n)^{-c}-n^{-1/2+2c}\ge c\log^{-c}n$ because $c$ is arbitrary. Then,
		\[
		\begin{split}
			&|m_1^{(1)}(\lambda)-m_1^{(1)}(z)|\le \frac{1}{n}\sum_{j=1}^p|\bSigma^{1/2}\hat\bgamma_j^{(1)}|^2\bigg|\frac{1}{\hat\lambda_j^{(1)}-\lambda}-\frac{1}{\hat\lambda_j^{(1)}-z}\bigg|\\
			= &\frac{1}{n}\sum_{j=1}^p|\bSigma^{1/2}\hat\bgamma_j^{(1)}|^2\bigg|\frac{in^{-1/2-c}}{(\hat\lambda_j^{(1)}-\lambda)(\hat\lambda_j^{(1)}-z)}\bigg|
			\le\frac{1}{n}\sum_{j=1}^p|\bSigma^{1/2}\hat\bgamma_j^{(1)}|^2\bigg|\frac{in^{-1/2-c}+O_{\prec}(n^{-1})}{|\hat\lambda_j^{(1)}-z|^2}\bigg|\\
			\prec &\operatorname{Im}m_1^{(1)}(z)+O_{\prec}(n^{-1})
			\prec \operatorname{Im} m_{1n}^{(1)}(z)+n^{-1/2}
			\prec o(\eta)+n^{-1/2}\prec n^{-1/2}.
		\end{split}
		\]
		Next, by definition,
		\[
		\begin{split}
			&m_{1n}^{(1)}(z)-m_{1n}^{(1)}(\lambda)\\
			=&\frac{1}{n}\sum_{i=1}^p\frac{\lambda_i}{-z\bigg[1+\frac{\lambda_i}{n}\sum_{j=2}^n\frac{w_{t_j}}{-z[1+m_{1n}(z)w_{t_j}]}\bigg]}-\frac{1}{n}\sum_{i=1}^p\frac{\lambda_i}{-\lambda\bigg[1+\frac{\lambda_i}{n}\sum_{j=2}^n\frac{w_{t_j}}{-\lambda[1+m_{1n}(\lambda)w_{t_j}]}\bigg]}\\
			=&\frac{1}{n}\sum_{i=1}^p\frac{\lambda_i}{-z\bigg[1+\frac{\lambda_i}{n}\sum_{j=2}^n\frac{w_{t_j}}{-z[1+m_{1n}(z)w_{t_j}]}\bigg]}-\frac{1}{n}\sum_{i=1}^p\frac{\lambda_i}{-z\bigg[1+\frac{\lambda_i}{n}\sum_{j=2}^n\frac{w_{t_j}}{-z[1+m_{1n}(\lambda)w_{t_j}]}\bigg]}\\
			&+\frac{1}{n}\sum_{i=1}^p\frac{\lambda_i}{-z\bigg[1+\frac{\lambda_i}{n}\sum_{j=2}^n\frac{w_{t_j}}{-z[1+m_{1n}(\lambda)w_{t_j}]}\bigg]}-\frac{1}{n}\sum_{i=1}^p\frac{\lambda_i}{-\lambda\bigg[1+\frac{\lambda_i}{n}\sum_{j=2}^n\frac{w_{t_j}}{-\lambda[1+m_{1n}(\lambda)w_{t_j}]}\bigg]}.
		\end{split}
		\]
		Similarly to (\ref{z prime}), we have
		\[
		|m_{1n}^{(1)}(z)-m_{1n}^{(1)}(\lambda)|\prec n^{-1/2}.
		\]
		As a result, we conclude that
		\[
		m_1^{(1)}(\lambda)=m_{1n}^{(1)}(\lambda)+O_{\prec}(n^{-1/2}).
		\]
		Then, for $\lambda=\lambda_0-n^{-1/2+2c}$,
		\begin{equation}\label{determination left}
			1+n^{-1}w_{t_1}\bx_{t_1}^\top[\hat{\Sbb}^{(1)}-\lambda\Ib]^{-1}\bx_{t_1}=1+w_{t_1}m_{1n}^{(1)}(\lambda)+O_{\prec}(n^{-1/2}).
		\end{equation}
		It's then sufficient to consider $1+w_{t_1}m_{1n}^{(1)}(\lambda)$. We write
		\[
		\begin{split}
			1+w_{t_1}m_{1n}^{(1)}(\lambda)=&1+w_{t_1}m_{1n}^{(1)}(\lambda_0)-w_{t_1}[m_{1n}^{(1)}(\lambda_0)-m_{1n}^{(1)}(\lambda)].
		\end{split}
		\]
		By definition, $1+w_{t_1}m_{1n}^{(1)}(\lambda_0)=0$, while
		\[
		\begin{split}
			&m_{1n}^{(1)}(\lambda_0)-m_{1n}^{(1)}(\lambda)\\	=&\frac{1}{n}\sum_{i=1}^p\frac{\lambda_i}{-\lambda_0\bigg[1+\frac{\lambda_i}{n}\sum_{j=2}^n\frac{w_{t_j}}{-\lambda_0[1+m_{1n}^{(1)}(\lambda_0)w_{t_j}]}\bigg]}-\frac{1}{n}\sum_{i=1}^p\frac{\lambda_i}{-\lambda\bigg[1+\frac{\lambda_i}{n}\sum_{j=2}^n\frac{w_{t_j}}{-\lambda[1+m_{1n}^{(1)}(\lambda)w_{t_j}]}\bigg]}\\
			=&\frac{1}{n}\sum_{i=1}^p\frac{\lambda_i}{-\lambda_0\bigg[1+\frac{\lambda_i}{n}\sum_{j=2}^n\frac{w_{t_j}}{-\lambda_0[1+m_{1n}^{(1)}(\lambda_0)w_{t_j}]}\bigg]}-\frac{1}{n}\sum_{i=1}^p\frac{\lambda_i}{-\lambda_0\bigg[1+\frac{\lambda_i}{n}\sum_{j=2}^n\frac{w_{t_j}}{-\lambda_0[1+m_{1n}^{(1)}(\lambda)w_{t_j}]}\bigg]}\\
			&+\frac{1}{n}\sum_{i=1}^p\frac{\lambda_i}{-\lambda_0\bigg[1+\frac{\lambda_i}{n}\sum_{j=2}^n\frac{w_{t_j}}{-\lambda_0[1+m_{1n}^{(1)}(\lambda)w_{t_j}]}\bigg]}-\frac{1}{n}\sum_{i=1}^p\frac{\lambda_i}{-\lambda\bigg[1+\frac{\lambda_i}{n}\sum_{j=2}^n\frac{w_{t_j}}{-\lambda[1+m_{1n}^{(1)}(\lambda)w_{t_j}]}\bigg]}\\
			:=&\mathcal{L}_1+\mathcal{L}_2.
		\end{split}
		\]
		For $\mathcal{L}_1$ we have,
		\[
		\begin{split}
			\mathcal{L}_1=&\frac{1}{n}\sum_{i=1}^p\frac{\lambda_i^2\frac{1}{n}\sum_{j=2}^n\bigg[\frac{w_{t_j}}{-\lambda_0[1+m_{1n}^{(1)}(\lambda_0)w_{t_j}]}-\frac{w_{t_j}}{-\lambda_0[1+m_{1n}^{(1)}(\lambda)w_{t_j}]}\bigg]}{-\lambda_0\bigg[1+\frac{\lambda_i}{n}\sum_{j=2}^n\frac{w_{t_j}}{-\lambda_0[1+m_{1n}^{(1)}(\lambda_0)w_{t_j}]}\bigg]\bigg[1+\frac{\lambda_i}{n}\sum_{j=2}^n\frac{w_{t_j}}{-\lambda_0[1+m_{1n}^{(1)}(\lambda)w_{t_j}]}\bigg]}\\
			=&O((\log n)^{-1})\times \frac{1}{n}\sum_{j=2}^n\bigg[\frac{w_{t_j}}{-\lambda_0[1+m_{1n}^{(1)}(\lambda_0)w_{t_j}]}-\frac{w_{t_j}}{-\lambda_0[1+m_{1n}^{(1)}(\lambda)w_{t_j}]}\bigg]\\
			=&o(1)\times [m_{1n}^{(1)}(\lambda_0)-m_{1n}^{(1)}(\lambda)].
		\end{split}
		\]
		On the other hand, for $\mathcal{L}_2$, we have
		\[
		\begin{split}
			\mathcal{L}_2=&\frac{1}{n}\sum_{i=1}^p\frac{\lambda_i[-\lambda+\lambda_0]}{\lambda\lambda_0\bigg[1+\frac{\lambda_i}{n}\sum_{j=2}^n\frac{w_{t_j}}{-\lambda_0[1+m_{1n}^{(1)}(\lambda)w_{t_j}]}\bigg]\bigg[1+\frac{\lambda_i}{n}\sum_{j=2}^n\frac{w_{t_j}}{-\lambda[1+m_{1n}^{(1)}(\lambda)w_{t_j}]}\bigg]}\\
			=&\frac{1}{n}\sum_{i=1}^p\frac{\lambda_i}{\lambda\lambda_0}[1+o(1)]\times n^{-1/2+2c}.
		\end{split}
		\]
		Therefore, we have
		\[
		[1+o(1)][m_{1n}^{(1)}(\lambda_0)-m_{1n}^{(1)}(\lambda)]=\frac{\phi_n \bar\lambda}{\lambda_0^2}[1+o(1)]n^{-1/2+2c},
		\]
		which further indicates that
		\begin{equation}\label{sign}
			1+w_{t_1}m_{1n}^{(1)}(\lambda)=-\frac{\phi_n \bar\lambda w_{t_1}}{\lambda_0^2}[1+o(1)]n^{-1/2+2c}<0.
		\end{equation}
		That  is to say, when $\lambda=\lambda_0-n^{-1/2+2c}$, with high probability,
		\[
		1+n^{-1}w_{t_1}\bx_{t_1}^\top[\hat{\Sbb}^{(1)}-\lambda\Ib]^{-1}\bx_{t_1}<0.
		\]
		
		Now let $\lambda=\lambda_0+n^{-1/2+2c}$, and we aim to show that with high probability
		\begin{equation}\label{right edge}
			1+n^{-1}w_{t_1}\bx_{t_1}^\top[\hat{\Sbb}^{(1)}-\lambda\Ib]^{-1}\bx_{t_1}>0.
		\end{equation}
		By an almost parallel procedure  in the proof of (\ref{determination left}), 
		\[
		1+n^{-1}w_{t_1}\bx_{t_1}^\top[\hat{\Sbb}^{(1)}-\lambda\Ib]^{-1}\bx_{t_1}=1+w_{t_1}m_{1n}^{(1)}(\lambda)+O_{\prec}(n^{-1/2}).
		\]
		Similarly to (\ref{sign}), we can show that
		\[
		1+w_{t_1}m_{1n}^{(1)}(\lambda)=\frac{\phi_n \bar\lambda w_{t_1}}{\lambda_0^2}[1+o(1)]n^{-1/2+2c}>0.
		\]
		Then, for $\lambda=\lambda_0+n^{-1/2+2c}$, with high probability
		\[
		1+n^{-1}w_{t_1}\bx_{t_1}^\top[\hat{\Sbb}^{(1)}-\lambda\Ib]^{-1}\bx_{t_1}>0.
		\]
		By continuity, with high probability there always exists an eigenvalue of $\hat{\Sbb}$ in the interval $[\lambda_0-n^{-1/2+2c},\lambda_0+n^{-1/2+2c}]$ for any constant $c>0$, which is also the largest eigenvalue. 
		
		Note that we are conditional on the order of $\{w_i\}$ in the above proof.  Recall that the orders are independent of $\Yb$ and follow uniform distribution. Write $\bt=\{t_1,\ldots,t_n\}$ and $\mathcal{T}_n=\{\bt_1,\ldots,\bt_{2^n}\}$. Then,
		\[
		\begin{split}
			&\mathbb{P}\bigg(|\hat\lambda_1-\lambda_0|\le n^{-1/2+2c}\bigg)=\sum_{j=1}^{2^n}\mathbb{P}\bigg(|\hat\lambda_1-\lambda_0|\le n^{-1/2+2c}, \bt=\bt_j\bigg)\\
			=&\sum_{j=1}^{2^n}\mathbb{P}\bigg(|\hat\lambda_1-\lambda_0|\le n^{-1/2+2c}\mid \bt=\bt_j\bigg)\mathbb{P}(\bt=\bt_j)\\
			\ge &\frac{1}{2^n}\sum_j\bigg(1-n^{-d}\bigg)
			\ge 1-n^{-d},
		\end{split}
		\]
		for any $d>0$. The lemma is then verified.
	\end{proof}

	\section{Proof of the results in Section \ref{sec:non-spiked}}\label{sec: proof non spike}
	Based on the preliminary results in the last section, now we are ready to prove Lemmas \ref{lemma: ratio}, \ref{lemma: lambda_0} and Theorem \ref{thm:non-spike}.
	\subsection{Proof of Lemma \ref{lemma: ratio}: ratio}
	\begin{proof}
		Note that $r=0$. Following Lemma \ref{Largest eigenvalue bound}, $\hat\lambda_1$ satisfies \[
		h(\hat\lambda_1)=1+n^{-1}w_{t_1}\bx_{t_1}^\top[\hat\Sbb_1^{(1)}-\hat\lambda_1\Ib]^{-1}\bx_{t_1}=0,
		\]
		which is equivalent to
		\[
		\begin{split}
			0=&\frac{\lambda_0}{w_{t_1}}+n^{-1}\bx_{t_1}^\top\bigg[\lambda_0^{-1}\hat\Sbb_1^{(1)}-\frac{\hat\lambda_1}{\lambda_0}\Ib\bigg]^{-1}\bx_{t_1},\\
			\Rightarrow	-\frac{\lambda_0}{w_{t_1}}=&\frac{1}{n}\bx_{t_1}^\top\bigg[\frac{1}{\lambda_0}\hat\Sbb_1^{(1)}-\Ib\bigg]^{-1}\bx_{t_1}+\frac{1}{n}\bx_{t_1}^\top\bigg[\frac{1}{\lambda_0}\hat\Sbb_1^{(1)}-\Ib\bigg]^{-1}\frac{\hat\lambda_1-\lambda_0}{\lambda_0}\bigg[\frac{1}{\lambda_0}\hat\Sbb_1^{(1)}-\frac{\hat\lambda_1}{\lambda_0}\Ib\bigg]^{-1}\bx_{t_1}\\
			:=&L_1+L_2.
		\end{split}
		\]
		
		We take the order of $\{w_j\}$ as given and start with an approximation to $L_1$. Write
		\[
		\begin{split}
			L_1=&\bigg\{\frac{1}{n}\bx_{t_1}^\top\bigg[\frac{1}{\lambda_0}\hat\Sbb_1^{(1)}-\Ib\bigg]^{-1}\bx_{t_1}-\frac{1}{n}\text{tr}\bigg[\frac{1}{\lambda_0}\hat\Sbb_1^{(1)}-\Ib\bigg]^{-1}\bSigma\bigg\}+\frac{1}{n}\text{tr}\bigg[\frac{1}{\lambda_0}\hat\Sbb_1^{(1)}-\Ib\bigg]^{-1}\bSigma\\
			:=&L_{11}+L_{12}.
		\end{split}
		\]
		By elementary matrix inverse formula in (\ref{matrix inverse}), 
		\[
		\begin{split}
			-L_{11}=& n^{-1}(\bx_{t_1}^\top\bx_{t_1}-\text{tr}\bSigma)+\frac{1}{n}\bx_{t_1}^\top\bigg[\hat\Sbb_1^{(1)}-\lambda_0\Ib\bigg]^{-1}\hat\Sbb_1^{(1)}\bx_{t_1}-\frac{1}{n}\text{tr}\bigg[\hat\Sbb_1^{(1)}-\lambda_0\Ib\bigg]^{-1}\hat\Sbb_1^{(1)}\bSigma.
		\end{split}
		\]
		Given the order of $w_j$'s, $\bx_{t_1}$ is independent of $\hat\Sbb_1^{(1)}$ and $\lambda_0$. Then,
		\[
		\begin{split}
			&\mathbb{E}\bigg|\frac{1}{n}\bx_{t_1}^\top\bigg[\hat\Sbb_1^{(1)}-\lambda_0\Ib\bigg]^{-1}\hat\Sbb_1^{(1)}\bx_{t_1}-\frac{1}{n}\text{tr}\bigg[\hat\Sbb_1^{(1)}-\lambda_0\Ib\bigg]^{-1}\hat\Sbb_1^{(1)}\bSigma\bigg|^2\\
			\le & \mathbb{E}\frac{C}{n^2}\bigg\|\bigg[\hat\Sbb_1^{(1)}-\lambda_0\Ib\bigg]^{-1}\hat\Sbb_1^{(1)}\bigg\|_F^2=\mathbb{E}\frac{C}{n^2}\sum_j\bigg(\frac{\hat\lambda_j^{(1)}}{\lambda_0-\hat\lambda_j^{(1)}}\bigg)^2\\
			\le& \frac{c}{n}\sqrt{\mathbb{E}\frac{1}{n}\sum_j(\hat\lambda_j^{(1)})^4\mathbb{E}\frac{1}{n}\sum_j\frac{1}{(\lambda_0-\hat\lambda_j^{(1)})^4}}\\&
			\le \frac{C}{n}\sqrt{\mathbb{E}\frac{1}{n}\|(\hat\Sbb^{(1)})^2\|_F^2\times\mathbb{E} \frac{C^2\log^{2c}n}{n}\sum_j\frac{1}{(\lambda_0-\hat\lambda_j^{(1)})^2}},
		\end{split}
		\]
		where we use the fact $\lambda_0-\hat\lambda_1^{(1)}\ge c\log^{-c}n$ in the last line.
		By definition,
		\[
		\|(\hat\Sbb_1^{(1)})^2\|_F^2\le \|\tilde\Sbb^{(1)}\|^4\sum_{j=2}^n w_{t_j}^4,
		\]
		and we conclude $\mathbb{E}n^{-1}\|(\hat\Sbb^{(1)})^2\|_F^2\le C$. On the other hand, let $z=\lambda_0+in^{-1/3}$, so
		\[
		\frac{1}{n}\sum_j\frac{1}{(\lambda_0-\hat\lambda_j^{(1)})^2}-	\frac{1}{n}\sum_j\frac{1}{|z-\hat\lambda_j^{(1)}|^2}=\frac{1}{n}\sum_j\frac{n^{-2/3}}{(\lambda_0-\hat\lambda_j^{(1)})^2|z-\hat\lambda_j^{(1)}|^2}\le \frac{C\log^{4c}n}{n^{2/3}}.
		\]
		However, we already know that
		\[
		\frac{1}{n}\sum_j\frac{1}{|z-\hat\lambda_j^{(1)}|^2}=\frac{1}{\eta}\operatorname{Im}m^{(1)}(z)=O((\log n)^{-1}).
		\]
		Therefore, 
		\[
		-L_{11}=n^{-1}(\bx_{t_1}^\top\bx_{t_1}-\text{tr}\bSigma)+O_p\bigg(\frac{\log^{\frac{c}{2}}n}{\sqrt{n}\log^{\frac{1}{4}}n}\bigg)=o_p(n^{-1/2}).
		\]
		
		Next, we consider $L_{12}$. Let $z=\lambda_0+in^{-2/3}$, and recall that $\hat\bgamma_j^{(1)}$'s are the eigenvectors of $\hat\Sbb^{(1)}$. Then, by elementary calculation,
		\begin{equation}\label{L12}
			\begin{split}
				|\lambda_0^{-1}L_{12}-m_1^{(1)}(z)|\le&\frac{1}{n}\sum_{j=1}^p|\bSigma^{1/2}\hat\bgamma_j|^2\bigg|\frac{1}{\hat\lambda_j^{(1)}-\lambda_0}-\frac{1}{\hat\lambda_j^{(1)}-z}\bigg|\le C\operatorname{Im}m^{(1)}(z)\\
				\le &C\operatorname{Im}m_n^{(1)}(z)+C|m^{(1)}(z)-m_n^{(1)}(z)|.
			\end{split}
		\end{equation}
		We already know that $m_1^{(1)}(z)$ and $m^{(1)}(z)$ converge with rate $n^{-1/2+c}$. However, here we need a sharper bound. For simplicity, we follow the proof of Lemmas \ref{large eta} and \ref{self-improvement} and only show the key steps. That is, we focus on $m(z)$, $m_1(z)$, while the proof for  $m_1^{(1)}(z)$, $m^{(1)}(z)$ is essentially similar. 
		
		We first improve the rate in Lemma \ref{large eta}, i.e., $\eta=\log^{1+c}n$. We rewrite (\ref{m2=m1}) as
		\begin{equation}\label{m2=m1 detail}
			\begin{split}
				m_2=\frac{1}{n}\sum_{i=1}^n\frac{w_i}{-z[1+w_im_1]}+\frac{1}{n}\sum_i\frac{w_iZ_i}{z[1+w_im_1]^2}+O_{\prec}(n^{-1}).
			\end{split}
		\end{equation}
		Then, following a standard fluctuation averaging argument (e.g., see Lemma 5.13 in \cite{yang2019edge}), we claim that
		\[
		\frac{1}{n}\sum_i\frac{w_iZ_i}{z[1+w_im_1]^2}\prec n^{-1}.
		\]
		Therefore, we can improve the rate in (\ref{m2=m1}) to 
		\begin{equation}\label{m2=m1 improve}
			m_2=\frac{1}{n}\sum_{i=1}^n\frac{w_i}{-z[1+w_im_1]}+O_{\prec}(n^{-1}).
		\end{equation}
		Next, to improve the rate in (\ref{m1=m2}), it suffices to reconsider $R_{11}$. In fact, since $\by_i$ is independent of $\Gb^{(i)}$ and $m_2^{(i)}$, similarly to the fluctuation averaging argument,  we can conclude that $R_{11}\prec n^{-1}$. That is to say,
		\begin{equation}\label{m1=m2 improve}
			m_1=\frac{1}{n}\sum_{i=1}^p\frac{\lambda_i}{-z[1+\lambda_im_2]}+O_{\prec}(n^{-1}).
		\end{equation}
		With (\ref{m2=m1 improve}) and (\ref{m1=m2 improve}), following the proof of Lemma \ref{large eta}, we conclude that $	\psi_2\prec n^{-1} $ for $\eta=\log^{1+c}n$, where $\psi_2$ is defined in Lemma \ref{self-improvement}. However, the rate for $\psi_1$ is still $n^{-1/2}$.
		
		It remains to show the local law also holds for small $\eta$. Following the argument in Lemmas \ref{self-improvement} and \ref{average local law}, it suffices to reconsider the self-improvement step in Lemma \ref{self-improvement}. We revise the statement in Lemma \ref{self-improvement} as follows.\\
		\begin{claim}\label{claim self improve}
			Under the same condition  in Lemma \ref{self-improvement}, for any $z\in D$, if $\psi_1\prec n^{-1/2+c}$ and $\psi_2\prec n^{-2/3+c}$ for some $0<c<1/6$, we then have
			\[
			\psi_1\prec n^{-1/2},\quad \psi_2\prec n^{-2/3}.
			\]
		\end{claim}
		To prove the claim, we need to improve the rate for $m_1$ and $m_2$ in Lemma \ref{self-improvement}. Indeed, following our proof of Lemma \ref{self-improvement} and (\ref{m2=m1 detail}), once again we can rewrite (\ref{m2=m1 general}) as
		\[
		m_2=\frac{1}{n}\sum_{i=1}^n\frac{w_i}{-z[1+w_im_1]}+\frac{1}{n}\sum_i\frac{w_iZ_i}{z[1+w_im_1]^2}+O_{\prec}(n^{-2/3}).
		\]
		By the fluctuation averaging argument, 
		\[
		\frac{1}{n}\sum_i\frac{w_iZ_i}{z[1+w_im_1]^2}\prec \psi_1^2\prec n^{-2/3}.
		\]
		Therefore, 
		\[
		m_2=\frac{1}{n}\sum_{i=1}^n\frac{w_i}{-z[1+w_im_1]}+O_{\prec}(n^{-2/3}).
		\]
		On the other hand, to obtain sharper rate for $m_1$, it's still sufficient to reconsider the bound of $R_{11}$, since we have already shown in the proof of Lemma \ref{self-improvement} that
		\[
		R_{12}\prec n^{-1+2c}+n^{-2+4c}\eta^{-1}\prec n^{-2/3},\quad \frac{1}{n}\text{tr}\Rb_2\bSigma\prec n^{-1}.
		\]
		For $R_{11}$, we can apply the fluctuation averaging argument to obtain that $R_{11}\prec \psi_1^2\prec n^{-2/3}$ for $0<c<1/6$. Therefore, 
		\[
		m_1=\frac{1}{n}\sum_{i=1}^p\frac{\lambda_i}{-z[1+\lambda_im_2]}+O_{\prec}(n^{-2/3}),
		\]
		which concludes our Claim \ref{claim self improve}. Then, following the induction step in the proof of Lemma \ref{average local law}, we have
		\[
		m_1(z)-m_{1n}(z)\prec n^{-2/3},\quad m(z)-m_n(z)\prec n^{-2/3},
		\]
		uniformly on $z\in D$. Similar procedures lead to 
		\[
		m_1^{(1)}(z)-m_{1n}^{(1)}(z)\prec n^{-2/3},\quad m^{(1)}(z)-m_n^{(1)}(z)\prec n^{-2/3},
		\]
		uniformly on $z=E+i\eta$ with $\lambda_{(2)}<E\le C\log n$ and $n^{-2/3}\le \eta\le \log^{1+c}n$.  Return to (\ref{L12}), so
		\[
		|\lambda_0^{-1}L_{12}-m_{1n}^{(1)}(z)|\le O_p(n^{-2/3+c}),\quad z=\lambda_0+in^{-2/3}.
		\]
		Moreover, a similar technique in proving (\ref{z prime}) leads to $|m_{1n}^{(1)}(z)-m_{1n}^{(1)}(\lambda_0)|\prec n^{-2/3}$. Consequently,
		\[
		L_1= -n^{-1}(\bx_{t_1}^\top\bx_{t_1}-\text{tr}\bSigma)+\lambda_0m_{1n}^{(1)}(\lambda_0)+o_p(n^{-1/2}).
		\]
		
		Next, we calculate $L_2$. By definition and (\ref{matrix inverse}),
		\[
		\begin{split}
			L_2=&\frac{\hat\lambda_1-\lambda_0}{\lambda_0}\frac{1}{n}\bx_{t_1}^\top\bigg[\frac{1}{\lambda_0}\hat\Sbb_1^{(1)}-\Ib\bigg]^{-2}\bx_{t_1}\\
			&+\bigg(\frac{\hat\lambda_1-\lambda_0}{\lambda_0}\bigg)^2\frac{1}{n}\bx_{t_1}^\top\bigg[\frac{1}{\lambda_0}\hat\Sbb_1^{(1)}-\Ib\bigg]^{-1}\bigg[\frac{1}{\lambda_0}\hat\Sbb_1^{(1)}-\frac{\hat\lambda_1}{\lambda_0}\Ib\bigg]^{-1}\bx_{t_1}\\
			=&\frac{\hat\lambda_1-\lambda_0}{\lambda_0}\bigg(\frac{1}{n}\text{tr}\bigg[\frac{1}{\lambda_0}\hat\Sbb_1^{(1)}-\Ib\bigg]^{-2}\bSigma+O_{\prec}(n^{-1/2})+\frac{\hat\lambda_1-\lambda_0}{\lambda_0}\times O_{\prec}(1)\bigg).
		\end{split}
		\]
		Note that
		\[
		\begin{split}
			&\frac{1}{n}\text{tr}\bigg[\frac{1}{\lambda_0}\hat\Sbb_1^{(1)}-\Ib\bigg]^{-2}\bSigma-\frac{1}{n}\text{tr}\bSigma\\
			=&\frac{2}{n\lambda_0}\text{tr}\bigg[\Ib-\frac{1}{\lambda_0}\hat\Sbb_1^{(1)}\bigg]^{-1}\hat\Sbb_1^{(1)}\bSigma+\frac{1}{n\lambda_0^2}\text{tr}\bigg[\Ib-\frac{1}{\lambda_0}\hat\Sbb_1^{(1)}\bigg]^{-1}\hat\Sbb_1^{(1)}\bigg[\Ib-\frac{1}{\lambda_0}\hat\Sbb_1^{(1)}\bigg]^{-1}\hat\Sbb_1^{(1)}\bSigma\\
			=&\frac{2}{n\lambda_0}\text{tr}\bigg[\Ib-\frac{1}{\lambda_0}\hat\Sbb_1^{(1)}\bigg]^{-1}\hat\Sbb_1^{(1)}\bSigma+o_p(1),
		\end{split}
		\]
		where the third line is by a similar technique as  bounding $L_{11}$. On the other hand, 
		\[
		\begin{split}
			\frac{1}{n\lambda_0}\text{tr}\bigg[\Ib-\frac{1}{\lambda_0}\hat\Sbb_1^{(1)}\bigg]^{-1}\hat\Sbb_1^{(1)}\bSigma=&\frac{1}{n}\sum_{j=1}^p|\bSigma^{1/2}\hat\bgamma_j^{(1)}|^2\frac{\hat\lambda_j^{(1)}}{\lambda_0-\hat\lambda_j^{(1)}}\\
			\le&C\sqrt{\frac{1}{n}\sum_j(\hat\lambda_j^{(1)})^2\times \frac{1}{n}\sum_j\frac{1}{(\lambda_0-\hat\lambda_j^{(1)})^2}}.
		\end{split}
		\]
		It's easy to see that $n^{-1}\|\hat\Sbb^{(1)}\|_F^2\le C$. Let $z=\lambda_0+in^{-1/3}$. Then,
		\[
		\begin{split}
			\frac{1}{n}\sum_j\frac{1}{(\lambda_0-\hat\lambda_j^{(1)})^2}=&\frac{1}{n}\sum_j\frac{1}{|z-\hat\lambda_j^{(1)}|^2}+O_{\prec}(n^{-2/3})
			=\frac{1}{\eta}\operatorname{Im}m^{(1)}(z)+O_{\prec}(n^{-2/3})\\
			=&\frac{1}{\eta}\operatorname{Im}m_{n}^{(1)}(z)+O_{\prec}(\eta^{-1}n^{-1/2})+O_{\prec}(n^{-2/3})
			=o_p(1).
		\end{split}
		\]
		Therefore, we conclude that
		\[
		L_2=\frac{\hat\lambda_1-\lambda_0}{\lambda_0}\times [\phi_n\bar\lambda+o_p(1)],
		\]
		which further implies that
		\[
		-\frac{\lambda_0}{w_{t_1}}= -n^{-1}(\bx_{t_1}^\top\bx_{t_1}-\text{tr}\bSigma)+\lambda_0m_{1n}^{(1)}(\lambda_0)+o_p(n^{-1/2})+\frac{\hat\lambda_1-\lambda_0}{\lambda_0}\times [\phi_n\bar\lambda+o_p(1)].
		\]
		By definition, we have $1+w_{t_1}m_{1n}^{(1)}(\lambda_0)=0$. Then,
		\[
		\sqrt{n}\frac{\hat\lambda_1-\lambda_0}{\lambda_0}=\frac{1}{\phi_n\bar\lambda}\frac{1}{\sqrt{n}}(\bx_{t_1}^\top\bx_{t_1}-p\bar\lambda)+o_p(1).
		\]
		The above argument is conditional on the order of $w_j$'s. It can be extended to unconditional results by similar technique  in the proof of Lemma \ref{Largest eigenvalue bound}. We omit details.
	\end{proof}

	\subsection{Proof of Lemma \ref{lemma: lambda_0}: $\lambda_0$}
	\begin{proof}
		We prove the lemma based on the definition of $\lambda_0$, i.e.,
		\[
		m_{1n}^{(1)}(\lambda_0)=-\frac{1}{w_{t_1}}=\frac{1}{n}\sum_{i=1}^p\frac{\lambda_i}{-\lambda_0\bigg[1+\frac{\lambda_i}{n}\sum_{j=2}^n\frac{w_{t_j}}{-\lambda_0[1-w_{t_j}/w_{t_{1}}]}\bigg]}.
		\]
		After some calculations,
		\[
		\begin{split}
			\frac{\lambda_0}{w_{t_1}}=&\frac{1}{n}\sum_i\lambda_i+\frac{1}{n}\sum_i\frac{\frac{\lambda_i^2}{n}\sum_{j=2}^n\frac{w_{t_j}}{\lambda_0[1-w_{t_j}/w_{t_1}]}}{1-\frac{\lambda_i}{n}\sum_{j=2}^n\frac{w_{t_j}}{\lambda_0[1-w_{t_j}/w_{t_{1}}]}}\\
			=&\phi_n\bar\lambda+\frac{1}{n}\sum_i\frac{\lambda_i^2}{n}\sum_{j=2}^n\frac{w_{t_j}}{\lambda_0[1-w_{t_j}/w_{t_1}]}+\frac{1}{n}\sum_i\frac{\lambda_i^3\bigg[\frac{1}{n}\sum_{j=2}^n\frac{w_{t_j}}{\lambda_0[1-w_{t_j}/w_{t_1}]}\bigg]^2}{1-\frac{\lambda_i}{n}\sum_{j=2}^n\frac{w_{t_j}}{\lambda_0[1-w_{t_j}/w_{t_{1}}]}}.
		\end{split}
		\]
		On the other hand, under the events $\Omega_n$ we can prove that
		\[
		\begin{split}
			\frac{1}{n}\sum_{j=2}^n\frac{w_{t_j}}{\lambda_0[1-w_{t_j}/w_{t_1}]}=&\frac{w_{t_1}}{\lambda_0}\frac{1}{n}\sum_{j=2}^n\frac{w_{t_j}}{w_{t_1}-w_{t_j}}\le O((\log n)^{-1}),
		\end{split}
		\]
		where we use the fact $\lambda_0>\lambda_{(2)}\ge C^{-1}\log n$	for large $n$ under the events $\Omega_n$. Therefore,
		\begin{equation}\label{lambda_0/w_a1}
			\frac{\lambda_0}{w_{t_1}}-\phi_n\bar\lambda=\frac{w_{t_1}}{\lambda_0}\times\frac{1}{n}\sum_i\lambda_i^2\times\frac{1}{n}\sum_{j=2}^n\frac{w_{t_j}}{w_{t_1}-w_{t_j}}+O((\log n)^{-2}).
		\end{equation}
		Then,  $C^{-1}\le \lambda_0/w_{t_1}\le C$ for sufficiently large $n$. We rewrite (\ref{lambda_0/w_a1}) as 
		\[
		\begin{split}
			\frac{\lambda_0}{w_{t_1}}-\phi_n\bar\lambda=&\bigg[\frac{1}{\phi_n\bar\lambda}+\frac{\phi_n\bar\lambda-\frac{\lambda_0}{w_{t_1}}}{\phi_n\bar\lambda\frac{\lambda_0}{w_{t_1}}}\bigg]O((\log n)^{-1})+O(\log^{-2}n),
		\end{split}
		\]
		which further leads to
		\[
		[1+o(1)]\bigg[\frac{\lambda_0}{w_{t_1}}-\phi_n\bar\lambda\bigg]=O((\log n)^{-1}).
		\]
		Return to (\ref{lambda_0/w_a1}), so
		\[
		\frac{\lambda_0}{w_{t_1}}-\phi_n\bar\lambda=\frac{1}{\phi_n\bar\lambda}\times\frac{1}{n}\sum_i\lambda_i^2\times\frac{1}{n}\sum_{j=2}^n\frac{w_{t_j}}{w_{t_1}-w_{t_j}}+O(\log^{-2}n).
		\]
		Further, a more detailed calculation shows that
		\[
		\begin{split}
			\frac{1}{n}\sum_{j=2}^n\frac{w_{t_j}}{w_{t_1}-w_{t_j}}=&\frac{1}{nw_{t_1}}\sum_{j=2}^n\frac{w_{t_j}}{1-w_{t_j}/w_{t_1}}
			=\frac{1}{nw_{t_1}}\sum_{j=2}^nw_{t_j}+\frac{1}{nw_{t_1}}\sum_{j=2}^n\frac{w_{t_j}^2/w_{t_1}}{1-w_{t_j}/w_{t_1}}\\
			=&\frac{1}{w_{t_1}}+O_p(\log^{-2}n).
		\end{split}
		\]
		Therefore, we conclude that
		\[
		\lambda_0=\phi_n\bar\lambda w_{t_1}+\frac{1}{\phi_n\bar\lambda}\frac{1}{n}\sum_i\lambda_i^2+O_p((\log n)^{-1}).
		\]
		Lemma \ref{lemma: ratio} has already shown that
		\[
		\hat\lambda_1-\lambda_0=\lambda_0\times O_p(n^{-1/2})=O_p(n^{-1/2}\log n),\quad r=0.
		\]
		Therefore, directly we have
		\[
		\hat\lambda_1=\phi_n\bar\lambda w_{t_1}+\frac{1}{\phi_n\bar\lambda}\frac{1}{n}\sum_i\lambda_i^2+O_p\bigg(\frac{1}{\log n}+\frac{\log n}{\sqrt{n}}\bigg),
		\]
		which concludes the unconditional results. When conditional on $\Xb$, the results follow a similar strategy because the limiting distribution is mainly determined by the resampling weights $w_j$'s. 
	\end{proof}

	\subsection{Proof of Theorem \ref{thm:non-spike}: $r\ge  0$}\label{sec:c6}
	\begin{proof}
		Now we consider the case $r>0$. Under such cases,  
		\[
		\Ab^\top\Ab=\bGamma\bLambda\bGamma^\top=\bGamma_1\bLambda_1\bGamma_1^\top+\bGamma_2\bLambda_2\bGamma_2^\top,
		\] 
		where $\bGamma_1$ and $\bLambda_1$ are the eigenvector and eigenvalue matrices corresponding to the $r$ spikes. Recall the definition $\bSigma_2=\bGamma_2\bLambda_2\bGamma_2^\top$. For simplicity, let $\varphi_j$ be the eigenvalues of $\bSigma_2$ (descending). Let $\hat\varphi_j$ be the $j$-th largest eigenvalue of $\hat{\mathcal{Q}}:=n^{-1}\Wb^{1/2}\Zb^\top\bSigma_2\Zb\Wb^{1/2}$, while $\hat\bbeta_j$ is the corresponding eigenvector.  Define $\hat\varphi_j^0$ and $\hat{\mathcal{Q}}^0$ accordingly by replacing $\Fb$ with $\Fb^0$, where $\Fb=\Cb\Fb^0$ by Assumption \ref{c1}.

		By Weyl's theorem, $\hat\lambda_{r+1}\le \hat\varphi_1$. On the other hand, we already know that $\hat\lambda_r\gg \hat\varphi_1\asymp \log n$.   Therefore, it suffices to prove that with probability tending to one, there is at least one eigenvalue of $\hat\Sbb$ in the interval $[\hat\varphi_1-n^{-2/3+c},\hat\varphi_1]$ for arbitrary small constant $c>0$. 
		
		By definition, $\hat\lambda_{r+1}$ satisfies
		\[
		\det\bigg(\hat\lambda_{r+1}\Ib-\frac{1}{n}\Wb^{1/2}\Zb^\top(\bGamma_1\bLambda_1\bGamma_1^\top+\bGamma_2\bLambda_2\bGamma_2^\top)\Zb\Wb^{1/2}\bigg)=0.
		\]
		We assume that $\hat\lambda_{r+1}$ is not an eigenvalue of $\hat{\mathcal{Q}}$ because $w_j$'s follow continuous distribution. For simplicity,  write $\check\Zb=\Zb\Wb^{1/2}=(\check\bz_1,\ldots,\check\bz_n)$. Then,
		\[
		\det(\Ib-n^{-1}\bLambda_1^{1/2}\bGamma_1^\top\check\Zb[\hat\lambda_{r+1}\Ib-\hat{\mathcal{Q}}]^{-1}\check\Zb^\top\bGamma_1\bLambda^{1/2})=0.
		\]
		Define 
		\[
		g(\lambda)=\det(\bLambda_1^{-1}-n^{-1}\bGamma_1^\top\check\Zb[\lambda\Ib-\hat{\mathcal{Q}}]^{-1}\check\Zb^\top\bGamma_1):=\det(\bLambda_1^{-1}-\bPi(\lambda)).
		\]
		We aim to prove that the sign of $g(\lambda)$ will change if $\lambda$ grows from $\hat\varphi_1-n^{-2/3+c}$ to $\hat\varphi_1$.
		
		Decompose $\bPi$ according  to
		\[
		\begin{split}
			\frac{1}{n}\bGamma_1^\top\check\Zb[\lambda\Ib-\hat{\mathcal{Q}}]^{-1}\check\Zb^\top\bGamma_1
			=&\frac{1}{\lambda-\hat\varphi_1}\frac{1}{n}\bGamma_1^\top\check\Zb\hat\bbeta_1\hat\bbeta_1^\top\check\Zb^\top\bGamma_1+\sum_{j=2}^n\frac{1}{\lambda-\hat\varphi_j}\frac{1}{n}\bGamma_1^\top\check\Zb\hat\bbeta_j\hat\bbeta_j^\top\check\Zb^\top\bGamma_1\\
			:=&\bPi_1(\lambda)+ \bPi_2(\lambda).
		\end{split}
		\]
		The key step of the proof is to find approximations to $\bPi_1$ and $\bPi_2$. 
		
		\noindent\textbf{Step 1: calculate $\bPi_1$.}
		
		We start with $\bPi_1$ and  show that 
		\begin{equation}\label{delocalization}
			n^{-1}\|\bGamma_1^\top\check\Zb\hat\bbeta_1\|^2\le O_p( n^{-3/4+c}).
		\end{equation}
		Since $r$ is fixed, it suffices to consider $n^{-1}|\bgamma_1^\top\check\Zb\hat\bbeta_1|^2$. (\ref{delocalization}) can be verified if we show that
		\begin{equation}\label{delocalization priori}
			\|\bPi(z)\|\le O_p(n^c),\quad \text{for }z=\hat\varphi_1+i\eta,\text{ where } \eta=n^{-3/4}.
		\end{equation}
		This is because as long as (\ref{delocalization priori}) holds, 
		\[
		\eta^{-1}\frac{1}{n}\bgamma_1^\top\check\Zb\hat\bbeta_1\hat\bbeta_1^\top\check\Zb^\top\bgamma_1=\operatorname{Im}\bigg(\frac{1}{\hat\varphi_1-z}\frac{1}{n}\bgamma_1^\top\check\Zb\hat\bbeta_1\hat\bbeta_1^\top\check\Zb^\top\bgamma_1\bigg)\le|\operatorname{Im} (\bPi(z))_{11} |\le O_p(n^c).
		\]
		Therefore, we aim to prove (\ref{delocalization priori}) in the following.
		Without loss of generality, we assume that $w_1\ge\cdots \ge w_n$ and all the arguments below are conditional on this event. The extension to unconditional results is similar to the argument in the proof of Lemma \ref{Largest eigenvalue bound}. Write $\check\Zb=\check\bz_1\be_1^\top+\check\Zb_{(1)}$. Then, 
		\[
		\begin{split}
			\bPi(z)=&\frac{1}{n}\bGamma_1^\top(\check\bz_1\be_1^\top+\check\Zb_ {(1)})\bigg[z\Ib-n^{-1}(\check\Zb_{(1)}^\top+\be_1\check\bz_1^\top)\bSigma_2(\check\bz_1\be_1^\top+\check\Zb_ {(1)})\bigg]^{-1}(\check\Zb_{(1)}^\top+\be_1\check\bz_1^\top)\bGamma_1\\
			=&\frac{1}{n}\bGamma_1^\top\check\bz_1\be_1^\top\bigg[z\Ib-n^{-1}(\check\Zb_{(1)}^\top+\be_1\check\bz_1^\top)\bSigma_2(\check\bz_1\be_1^\top+\check\Zb_ {(1)})\bigg]^{-1}\be_1\check\bz_1^\top\bGamma_1\\
			&+\frac{1}{n}\bGamma_1^\top\check\Zb_ {(1)}\bigg[z\Ib-n^{-1}(\check\Zb_{(1)}^\top+\be_1\check\bz_1^\top)\bSigma_2(\check\bz_1\be_1^\top+\check\Zb_ {(1)})\bigg]^{-1}\be_1\check\bz_1^\top\bGamma_1\\
			&+\frac{1}{n}\bGamma_1^\top\check\bz_1\be_1^\top\bigg[z\Ib-n^{-1}(\check\Zb_{(1)}^\top+\be_1\check\bz_1^\top)\bSigma_2(\check\bz_1\be_1^\top+\check\Zb_ {(1)})\bigg]^{-1}\check\Zb_{(1)}^\top\bGamma_1\\
			&+\frac{1}{n}\bGamma_1^\top\check\Zb_ {(1)}\bigg[z\Ib-n^{-1}(\check\Zb_{(1)}^\top+\be_1\check\bz_1^\top)\bSigma_2(\check\bz_1\be_1^\top+\check\Zb_ {(1)})\bigg]^{-1}\check\Zb_{(1)}^\top\bGamma_1\\
			:=&\bPi_{11}(z)+\bPi_{12}(z)+\bPi_{13}(z)+\bPi_{14} (z).
		\end{split}
		\] 
		We deal with the four terms one by one. Write $\hat{\mathcal{Q}}_{(1)}=n^{-1}\check\Zb_{(1)}^\top\bSigma_2\check\Zb_{(1)}$. By (\ref{matrix inverse}),
		\[
		\begin{split}
			&\be_1^\top\bigg[z\Ib-n^{-1}(\check\Zb_{(1)}^\top+\be_1\check\bz_1^\top)\bSigma_2(\check\bz_1\be_1^\top+\check\Zb_ {(1)})\bigg]^{-1}\be_1
			=\be_1^\top[z\Ib-\hat{\mathcal{Q}}_{(1)}]^{-1}\be_1\\
			&+\be_1^\top[z\Ib-\hat{\mathcal{Q}}]^{-1}\frac{1}{n}(\check\Zb_{(1)}^\top\bSigma_2\check\bz_1\be_1^\top+\be_1\check\bz_1^\top\bSigma_2\check\bz_1\be_1^\top+\be_1\check\bz_1^\top\bSigma_2\check\Zb_{(1)})[z\Ib-\hat{\mathcal{Q}}_{(1)}]^{-1}\be_1.
		\end{split}
		\]
		It's easy to see that $\be_1^\top[z\Ib-\hat{\mathcal{Q}}_{(1)}]^{-1}=z^{-1}\be_1^\top$ and $\be_1^\top\check\Zb_{(1)}={\bf 0}$. Therefore,
		\[
		\begin{split}
			&\be_1^\top\bigg[z\Ib-n^{-1}(\check\Zb_{(1)}^\top+\be_1\check\bz_1^\top)\bSigma_2(\check\bz_1\be_1^\top+\check\Zb_ {(1)})\bigg]^{-1}\be_1\\
			=&z^{-1}+\frac{1}{nz}\be_1^\top[z\Ib-\hat{\mathcal{Q}}]^{-1}\check\Zb_{(1)}^\top\bSigma_2\check\bz_1+\be_1^\top[z\Ib-\hat{\mathcal{Q}}]^{-1}\be_1\times\frac{1}{nz}\check\bz_1^\top\bSigma_2\check\bz_1.
		\end{split}
		\]
		Using (\ref{matrix inverse}) once again, we have
		\[
		\begin{split}
			&\frac{1}{n}\be_1^\top[z\Ib-\hat{\mathcal{Q}}]^{-1}\check\Zb_{(1)}^\top\bSigma_2\check\bz_1=\frac{1}{n}\be_1^\top[z\Ib-\hat{\mathcal{Q}}_{(1)}]^{-1}\check\Zb_{(1)}^\top\bSigma_2\check\bz_1\\
			&+\frac{1}{n}\be_1^\top[z\Ib-\hat{\mathcal{Q}}]^{-1}\frac{1}{n}(\check\Zb_{(1)}^\top\bSigma_2\check\bz_1\be_1^\top+\be_1\check\bz_1^\top\bSigma_2\check\bz_1\be_1^\top+\be_1\check\bz_1^\top\bSigma_2\check\Zb_{(1)})[z\Ib-\hat{\mathcal{Q}}_{(1)}]^{-1}\check\Zb_{(1)}^\top\bSigma_2\check\bz_1\\
			=&\frac{1}{n}\be_1^\top[z\Ib-\hat{\mathcal{Q}}]^{-1}\be_1\times \frac{1}{n}\check\bz_1^\top\bSigma_2\check\Zb_{(1)}[z\Ib-\hat{\mathcal{Q}}_{(1)}]^{-1}\check\Zb_{(1)}^\top\bSigma_2\check\bz_1.
		\end{split}
		\]
		Consequently, we claim that
		\[
		|\be_1^\top[z\Ib-\hat{\mathcal{Q}}]^{-1}\be_1|=\frac{1}{|z-n^{-1}\check\bz_1^\top\bSigma_2\check\Zb_{(1)}[z\Ib-\hat{\mathcal{Q}}_{(1)}]^{-1}\check\Zb_{(1)}^\top\bSigma_2\check\bz_1-n^{-1}\check\bz_1^\top\bSigma_2\check\bz_1|}\le (\operatorname{Im}z)^{-1}\le n^{3/4}.
		\]
		Further, we have $\|\bGamma_1^\top\check{\bz}_1\|\le O_p(n^c)$ by the independence of the entries in $\bz_1$. Then, 
		\[
		\|\bPi_{11}(z)\|\le O_p(n^{-1/4+c}).
		\]
		Next, for $\bPi_{12}(z)$, by (\ref{matrix inverse}) we have
		\[
		\begin{split}
			&\frac{1}{n}\bGamma_1^\top\check\Zb_{(1)}[z\Ib-\hat{\mathcal{Q}}]^{-1}\be_1=\frac{1}{n}\bGamma_1^\top\check\Zb_{(1)}[z\Ib-\hat{\mathcal{Q}}_{(1)}]^{-1}\be_1\\
			&+\frac{1}{n}\bGamma_1^\top\check\Zb_{(1)}[z\Ib-\hat{\mathcal{Q}}]^{-1}\frac{1}{n}(\check\Zb_{(1)}^\top\bSigma_2\check\bz_1\be_1^\top+\be_1\check\bz_1^\top\bSigma_2\check\bz_1\be_1^\top+\be_1\check\bz_1^\top\bSigma_2\check\Zb_{(1)})[z\Ib-\hat{\mathcal{Q}}_{(1)}]^{-1}\be_1\\
			=&\frac{1}{nz}\bGamma_1^\top\check\Zb_{(1)}[z\Ib-\hat{\mathcal{Q}}]^{-1}\frac{1}{n}(\check\Zb_{(1)}^\top\bSigma_2\check\bz_1+\be_1\check\bz_1^\top\bSigma_2\check\bz_1).
		\end{split}
		\]
		Further by (\ref{matrix inverse}),
		\[
		\begin{split}
			&\frac{1}{n^2z}\bGamma_1^\top\check\Zb_{(1)}[z\Ib-\hat{\mathcal{Q}}]^{-1}\check\Zb_{(1)}^\top\bSigma_2\check\bz_1=\frac{1}{n^2z}\bGamma_1^\top\check\Zb_{(1)}[z\Ib-\hat{\mathcal{Q}}_{(1)}]^{-1}\check\Zb_{(1)}^\top\bSigma_2\check\bz_1\\
			&+\frac{1}{n^2z}\bGamma_1^\top\check\Zb_{(1)}[z\Ib-\hat{\mathcal{Q}}]^{-1}\frac{1}{n}(\check\Zb_{(1)}^\top\bSigma_2\check\bz_1\be_1^\top+\be_1\check\bz_1^\top\bSigma_2\check\bz_1\be_1^\top+\be_1\check\bz_1^\top\bSigma_2\check\Zb_{(1)})[z\Ib-\hat{\mathcal{Q}}_{(1)}]^{-1}\check\Zb_{(1)}^\top\bSigma_2\check\bz_1\\
			=&\frac{1}{n^2z}\bGamma_1^\top\check\Zb_{(1)}[z\Ib-\hat{\mathcal{Q}}_{(1)}]^{-1}\check\Zb_{(1)}^\top\bSigma_2\check\bz_1+\frac{1}{n^2z}\bGamma_1^\top\check\Zb_{(1)}[z\Ib-\hat{\mathcal{Q}}]^{-1}\be_1\frac{1}{n}\check\bz_1^\top\bSigma_2\check\Zb_{(1)})[z\Ib-\hat{\mathcal{Q}}_{(1)}]^{-1}\check\Zb_{(1)}^\top\bSigma_2\check\bz_1.
		\end{split}
		\]
		Consequently,
		\begin{equation}\label{Pi 12 part}
			\begin{split}
				\|\frac{1}{n}\bGamma_1^\top\check\Zb_{(1)}[z\Ib-\hat{\mathcal{Q}}]^{-1}\be_1\|=&\frac{\|\frac{1}{n^2}\bGamma_1^\top\check\Zb_{(1)}[z\Ib-\hat{\mathcal{Q}}_{(1)}]^{-1}\check\Zb_{(1)}^\top\bSigma_2\check\bz_1\|}{|z-n^{-1}\check\bz^\top\bSigma_2\check\bz_1-n^{-2}\check\bz_1^\top\bSigma_2\check\Zb_{(1)})[z\Ib-\hat{\mathcal{Q}}_{(1)}]^{-1}\check\Zb_{(1)}^\top\bSigma_2\check\bz_1|}\\
				\le &\|n^{-2}\bGamma_1^\top\check\Zb_{(1)}[z\Ib-\hat{\mathcal{Q}}_{(1)}]^{-1}\check\Zb_{(1)}^\top\bSigma_2\check\bz_1\|(\operatorname{Im} z)^{-1}.
			\end{split}
		\end{equation}
		Define $\mathcal{Q}^0,\hat{\mathcal{Q}}_{(1)}^0$ by replacing $\Fb$ with $\Fb^0$, respectively. By the assumption that $n^{-1}\|\Cb\|_F^2=1$, we have $\|\mathcal{Q}-\mathcal{Q}^0\|=O_p(n^{-1/2+c})$. Then, $|\operatorname{Re}z-\lambda_{\max}(\mathcal{Q}^0)|\le O_p(n^{-1/2+c})$.  Further,  similarly to Lemma \ref{Largest eigenvalue bound} and (\ref{lambda 02 gap}), we can show that $\lambda_{\max}(\mathcal{Q}^0)-\lambda_{\max}(\hat{\mathcal{Q}}_{(1)}^0)\ge (\log n)^{-c}$ with probability tending to 1 for any $c>0$. That is, $\|[z\Ib-\hat{\mathcal{Q}}_{(1)}^0]^{-1}\|\le O_p(n^c)$. The same bound holds for $\hat{\mathcal{Q}}_{(1)}$ because $\|\hat{\mathcal{Q}}_{(1)}-\hat{\mathcal{Q}}_{(1)}^0\|\le O_p(n^{-1/2+c})$. Return to (\ref{Pi 12 part}), and write $\check\bz_1^\top=w_1(\bC_{1\cdot}^\top\Fb^0,\bepsilon_1^\top)$, where $\bC_{1\cdot}$ is the first row vector of $\Cb$. Then
		\[
		\begin{split}
			&\|n^{-2}\bGamma_1^\top\check\Zb_{(1)}[z\Ib-\hat{\mathcal{Q}}_{(1)}]^{-1}\check\Zb_{(1)}^\top\bSigma_2\check\bz_1\|\\
			\le& w_1\|n^{-2}\bGamma_1^\top\check\Zb_{(1)}[z\Ib-\hat{\mathcal{Q}}_{(1)}]^{-1}\check\Zb_{(1)}^\top\bSigma_2\|\|\bC_{1\cdot}^\top\Fb^0\|+w_1\|n^{-2}\bGamma_1^\top\check\Zb_{(1)}[z\Ib-\hat{\mathcal{Q}}_{(1)}]^{-1}\check\Zb_{(1)}^\top\bSigma_2({\bf 0}^\top,\bepsilon_1^\top)^\top\|\\
			\le &O_p(n^{-1+c})+O_p(w_1\|n^{-2}\bGamma_1^\top\check\Zb_{(1)}[z\Ib-\hat{\mathcal{Q}}_{(1)}]^{-1}\check\Zb_{(1)}^\top\bSigma_2\|)=O_p(n^{-1+c}),
		\end{split}
		\]
		where the third line is by the independence of between $\bepsilon_1$ and $\check{\Zb}_{(1)}$. 
		As a result, 
		\[
		\|n^{-2}\bGamma_1^\top\check\Zb_{(1)}[z\Ib-\hat{\mathcal{Q}}_{(1)}]^{-1}\check\Zb_{(1)}^\top\bSigma_2\check\bz_1\|=O_p(n^{-1+c}),\quad \|n^{-1}\bGamma_1^\top\check\Zb_{(1)}[z\Ib-\hat{\mathcal{Q}}]^{-1}\be_1\|=O_p(n^{-1/4+c}),
		\]
		which further indicates that
		\[
		\|\bPi_{12}(z)\|\le O_p( n^{-1/4+c}),\quad \|\bPi_{13}(z)\|\le O_p( n^{-1/4+c}).
		\]
		Lastly, for $\bPi_{14}(z)$, we write
		\[
		\begin{split}
			&\bPi_{14}(z)=\frac{1}{n}\bGamma_1^\top\check\Zb_ {(1)}[z\Ib-\hat{\mathcal{Q}}_{(1)}]^{-1}\check\Zb_{(1)}^\top\bGamma_1\\
			&+\frac{1}{n}\bGamma_1^\top\check\Zb_ {(1)}[z\Ib-\hat{\mathcal{Q}}]^{-1}\frac{1}{n}(\check\Zb_{(1)}^\top\bSigma_2\check\bz_1\be_1^\top+\be_1\check\bz_1^\top\bSigma_2\check\bz_1\be_1^\top+\be_1\check\bz_1^\top\bSigma_2\check\Zb_{(1)})[z\Ib-\hat{\mathcal{Q}}_{(1)}]^{-1}\check\Zb_{(1)}^\top\bGamma_1\\
			=&\frac{1}{n}\bGamma_1^\top\check\Zb_ {(1)}[z\Ib-\hat{\mathcal{Q}}_{(1)}]^{-1}\check\Zb_{(1)}^\top\bGamma_1+\frac{1}{n^2}\bGamma_1^\top\check\Zb_ {(1)}[z\Ib-\hat{\mathcal{Q}}]^{-1}\be_1\check\bz_1^\top\bSigma_2\check\Zb_{(1)}[z\Ib-\hat{\mathcal{Q}}_{(1)}]^{-1}\check\Zb_{(1)}^\top\bGamma_1\le O_p(n^c).
		\end{split}
		\]
		Then, (\ref{delocalization priori}) holds, and (\ref{delocalization}) follows.

		\noindent\textbf{Step 2: calculate $\bPi_2$.}
		
		Now we turn to $\bPi_2(\lambda)$. Let $z_1=\lambda+in^{-2/3}$, $z_2=\lambda+n^{-1/6}+in^{-2/3}$ for $\lambda\in[\hat\varphi_1-n^{-2/3+c},\hat\varphi_1]$. According to the definitions of $\lambda_0$ and $ \lambda_{(1)}$ in Section \ref{secc}, similarly define $\varphi_0$ and $\varphi_{(1)}$ by replacing $\bSigma$ with $\bSigma_2$. Then, if the entries of $\Fb$ are independent,  by Lemmas \ref{upper bound}, \ref{Largest eigenvalue bound} and the fact that $|\varphi_0-\varphi_{(1)}|\le O(n^{-1/2+c})$, we have $|\hat\varphi_1-\varphi_{(1)}|\le O_p(n^{-1/2+c})$.  Moreover, by Weyl's theorem, $\hat\varphi_2\le \lambda_{\max}(\hat{\mathcal{Q}}_{(1)})\le \hat\varphi_1-(\log n)^{-c}$.  Therefore, with probability tending to 1 we have
		\begin{equation}\label{z1 z2}
			\operatorname{Re}z_2>\varphi_{(1)}+n^{-1/6-c},\quad \operatorname{Re}z_1-\hat\varphi_2\ge (\log n)^{-c}.
		\end{equation}
		The results also holds for $\Fb=\Cb\Fb^0$ because the replacement error to the eigenvalues is upper-bounded by $O_p(n^{-1/2+c})$. All the convergence rates hereafter are uniform on $\lambda$.
		
		By definition, we write
		\[
		\begin{split}
			\bPi_2(\lambda)=&\frac{1}{n}\sum_{j=1}^n\frac{1}{z_2-\hat\varphi_j}\bGamma_1^\top\check\Zb\hat\bbeta_j\hat\bbeta_j^\top\check\Zb^\top\bGamma_1+\bPi_2(\lambda)-\frac{1}{n}\sum_{j=2}^n\frac{1}{z_1-\hat\varphi_j}\bGamma_1^\top\check\Zb\hat\bbeta_j\hat\bbeta_j^\top\check\Zb^\top\bGamma_1\\
			&+\frac{1}{n}\sum_{j=2}^n\frac{1}{z_1-\hat\varphi_j}\bGamma_1^\top\check\Zb\hat\bbeta_j\hat\bbeta_j^\top\check\Zb^\top\bGamma_1-\frac{1}{n}\sum_{j=1}^n\frac{1}{z_2-\hat\varphi_j}\bGamma_1^\top\check\Zb\hat\bbeta_j\hat\bbeta_j^\top\check\Zb^\top\bGamma_1.
		\end{split}
		\]
		Elementary calculation leads to
		\[
		\begin{split}
			&\bigg|\bigg(\bPi_2(\lambda)-\frac{1}{n}\sum_{j=2}^n\frac{1}{z_1-\hat\varphi_j}\bGamma_1^\top\check\Zb\hat\bbeta_j\hat\bbeta_j^\top\check\Zb^\top\bGamma_1\bigg)_{ii}\bigg|\\
			\le&\frac{1}{n}\sum_{j=2}^n|\bgamma_i^\top\check\Zb\hat\bbeta_j|^2\bigg|\frac{1}{z_1-\hat\varphi_j}-\frac{1}{\lambda-\hat\varphi_j}\bigg|
			\lesssim\frac{1}{n}\sum_{j=2}^n|\bgamma_i^\top\check\Zb\hat\bbeta_j|^2\frac{\operatorname{Im}z_1}{|z_1-\hat\varphi_j|^2}\\
			\le& O_p(n^{-2/3+c})\bigg(n^{-1}\bgamma_i^\top\check\Zb\check\Zb^\top\bgamma_i+O_p(n^{-3/4+c})\bigg)\le O_p(n^{-2/3+c}),
		\end{split}
		\]
		where we use (\ref{delocalization}) in the last line. 
		Similarly,
		\[
		\begin{split}
			&\bigg|\bigg(\frac{1}{n}\sum_{j=2}^n\frac{1}{z_1-\hat\varphi_j}\bGamma_1^\top\check\Zb\hat\bbeta_j\hat\bbeta_j^\top\check\Zb^\top\bGamma_1-\frac{1}{n}\sum_{j=1}^n\frac{1}{z_2-\hat\varphi_j}\bGamma_1^\top\check\Zb\hat\bbeta_j\hat\bbeta_j^\top\check\Zb^\top\bGamma_1\bigg)_{ik}\bigg|\\
			\le &O_p(n^{-3/4+c+2/3})+\frac{1}{n}\sum_{j=2}^n|\bgamma_i^\top\check\Zb\hat\bbeta_j||\bgamma_k^\top\check\Zb\hat\bbeta_j|\bigg|\frac{1}{z_1-\hat\varphi_j}-\frac{1}{z_2-\hat\varphi_j}\bigg|\le O_p(n^{-1/12+c}).
		\end{split}
		\]
		Therefore,
		\begin{equation}\label{PI 21}
			\bPi_2(\lambda)=\frac{1}{n}\bGamma_1^\top\check\Zb[z_2\Ib-\hat{\mathcal{Q}}]^{-1}\check\Zb^\top\bGamma_1+O_p(n^{-1/12+c}).
		\end{equation}

		\noindent\textbf{Step 3: finish the proof.}
		
		Similarly to Step 2 in the proof of Lemma \ref{lema2}, write
		\[
		\mathcal{H}(z):=\left(\begin{matrix}
			& z\Ib&n^{-1/2}\check\Zb^\top\bGamma_2\bLambda_2^{1/2}&n^{-1/2}\check\Zb^\top\bGamma_1\\
			&n^{-1/2}\bLambda_2^{1/2}\bGamma_2^\top\check\Zb&\Ib&{\bf 0}\\
			&n^{-1/2}\bGamma_1^\top\check\Zb&{\bf 0}&\Ib_r\\
		\end{matrix}\right)=\left(\begin{matrix}
			&z\Ib&n^{-1/2}\check\Zb^\top\tilde\Ab\\
			&n^{1/2}\tilde\Ab^\top\check\Zb&\Ib_{p+r}
		\end{matrix}\right),
		\]
		where $\tilde\Ab=(\bGamma_2\bLambda_2^{1/2},\bGamma_1)$. Then, the lowest-rightest $r\times r$ block of $\mathcal{H}^{-1}(z_2)$ is exactly equal to the inverse of $\Ib-n^{-1}\bGamma_1^\top\check\Zb[z_2\Ib-\hat{\mathcal{Q}}]^{-1}\check\Zb^\top\bGamma_1$. By Shur's complement formula, for $1\le k,l\le r+p$,
		\[
		\begin{split}
			&[\mathcal{H}^{-1}(z_2)]_{n+k,n+l}= [ (\Ib-z_2^{-1}n^{-1}\tilde\Ab^\top\check\Zb\check\Zb^\top\tilde\Ab )^{-1} ]_{kl}.
		\end{split}
		\] 
		Similarly to (\ref{z1 z2}), we can show that $\operatorname{Re}z_2\ge \lambda_{\max}(n^{-1}\tilde\Ab^\top\check{\Zb}\check{\Zb}^\top\tilde\Ab)+n^{-1/6-c}\ge (\log n)^{1-c}$  with probability tending to 1 for any $c>0$.  Therefore, $\|\Ib-z_2^{-1}n^{-1}\tilde\Ab^\top\check{\Zb}\check{\Zb}^\top\tilde\Ab\|\le O_p(n^{1/6+2c})$. Define $\check{\Zb}_0$ by replacing $\Fb$ with $\Fb^0$.  Then, we also have $\|\Ib-z_2^{-1}n^{-1}\tilde\Ab^\top\check{\Zb}_0\check{\Zb}_0^\top\tilde\Ab\|\le O_p(n^{1/6+2c})$. Moreover, by (\ref{matrix inverse}),
		\[
		\begin{split}
			&\|(\Ib-z_2^{-1}n^{-1}\tilde\Ab^\top\check\Zb\check\Zb^\top\tilde\Ab )^{-1} -(\Ib-z_2^{-1}n^{-1}\tilde\Ab^\top\check\Zb_0\check\Zb_0^\top\tilde\Ab )^{-1} \|\\
			\le &\|(\Ib-z_2^{-1}n^{-1}\tilde\Ab^\top\check\Zb\check\Zb^\top\tilde\Ab )^{-1} \frac{1}{nz_2}\tilde\Ab^\top(\check{\Zb}\check{\Zb}^\top-\check{\Zb}_0\check{\Zb}_0^\top)\tilde\Ab(\Ib-z_2^{-1}n^{-1}\tilde\Ab^\top\check\Zb_0\check\Zb_0^\top\tilde\Ab )^{-1}\|\\
			\le &O_p(n^{1/3+4c})\|n^{-1}\tilde\Ab^\top(\check{\Zb}\check{\Zb}^\top-\check{\Zb}_0\check{\Zb}_0^\top)\tilde\Ab\|\le O_p(n^{-1/6+5c}),
		\end{split}
		\]
		where we use the fact that $n^{-1}\|\tilde\Ab^\top(\check{\Zb}\check{\Zb}^\top-\check{\Zb}_0\check{\Zb}_0^\top)\tilde\Ab\|\le O_p(n^{-1/2+c})$ for arbitrary $c>0$.  Consequently, we can assume $\Cb=\Ib$ in the following because $c$ can be arbitrary small.

		Note that  $(\Ib-z_2^{-1}n^{-1}\tilde\Ab^\top\check\Zb\check\Zb^\top\tilde\Ab )^{-1}$ is exactly the matrix $-z_2\Gb$ defined in Section \ref{secc} by taking $\bSigma=\tilde\bSigma:=\tilde\Ab^\top\tilde\Ab$.  Because $\operatorname{Re}z_2\ge \hat\varphi_1-n^{-2/3+c}+n^{-1/6}$ while $\hat\varphi_1=\lambda_{(1)}+O_{\prec}(n^{-1/2+c})$, we have $z_2\in D$ with high probability and $\|\Gb\|\prec n^{1/6+c}, \|\Gb^{(i)}\|\prec n^{1/6+c}$ for $1\le i\le n$. According to (\ref{R1+R2}), it suffices to consider $|z(\Rb_1)_{kl}|$ and $|z(\Rb_2)_{kl}|$ for $z\in D$. 
		
		We start with $|z(\Rb_1)_{kk}|$. By Lemma \ref{self-improvement}, $|m_2(z)-m_{2n}(z)|\prec n^{-1/2+c}$, while Lemma \ref{m_{1n}} indicates that $|m_{2n}(z)|=o(1)$. Moreover, we already know that  $(1+n^{-1}\by_i^\top\Gb^{(i)}\by_i)^{-1}=-z\mathcal{G}_{ii}=1+m_{1n}(z)w_i+O_{\prec}(n^{-1/2+c})$. Then, we conclude that
		\[
		|z(\Rb_1)_{kk}|=\bigg|\frac{1}{n}\sum_i(1-\mathbb{E}_i)(1+m_{1n}(z)w_i)\be_k^\top\Gb^{(i)}\by_i\by_i^\top[\Ib+m_{2n}(z)\tilde\bSigma]\be_k\bigg|+O_{\prec}(n^{-1/6+c}),
		\]
		where $\mathbb{E}_i$ indicates the conditional expectation given $\bz_i$ and $\Wb$. Note that $\by_i$ and $\Gb^{(i)}$ is independent conditional on $\Wb$. Therefore, by Burkholder's inequality we will have 
		\[
		|z(\Rb_1)_{kk}|\prec n^{-1/2}\|\be_k^\top\Gb^{(i)}\|\|[\Ib+m_{2n}(z)\tilde\bSigma]\be_k\|+O_p(n^{-1/6+c})\prec n^{-1/6+c}.
		\]
		The same bound holds for $k\ne l$ by similar arguments.  Next, for $|z(\Rb_2)_{kk}|$, similarly we have
		\[
		|z(\Rb_2)_{kk}|=\bigg|\frac{1}{n}\sum_i(1+m_{1n}(z)w_i)w_i\be_k^\top(\Gb^{(i)}-\Gb)\tilde\bSigma[\Ib+m_{2n}(z)\tilde\bSigma]\be_k\bigg|+O_{\prec}(n^{-1/6+c}).
		\]
		Given $i$, 
		\[
		\begin{split}
			&|\be_k^\top(\Gb^{(i)}-\Gb)\tilde\bSigma[\Ib+m_{2n}(z)\tilde\bSigma]\be_k|=|n^{-1}\be_k^\top\Gb\by_i\by_i^\top\Gb^{(i)}\tilde\bSigma[\Ib+m_{2n}(z)\tilde\bSigma]\be_k|\\
			=&|n^{-1}z\mathcal{G}_{ii}\be_k^\top\Gb^{(i)}\by_i\by_i^\top\Gb^{(i)}\tilde\bSigma[\Ib+m_{2n}(z)\tilde\bSigma]\be_k|\prec n^{-1}\|\Gb^{(i)}\|^2\prec n^{-1/6+c}.
		\end{split}
		\]
		Therefore, $|z(\Rb_2)_{kk}|\prec n^{-1/6+c}$ and similarly $|z(\Rb_2)_{kl}|\prec n^{-1/6+c}$ for $k\ne l$. Consequently, by (\ref{R1+R2}) we have 
		\[
		[-z\Gb]_{kl}=\big([\Ib+m_{2n}(z)\tilde\bSigma]^{-1}\big)_{kl}+O_{\prec}(n^{-1/6+c}),\quad z\in D,
		\]
		which further indicates that 
		\[
		[(\Ib-z_2^{-1}n^{-1}\tilde\Ab^\top\check\Zb\check\Zb^\top\tilde\Ab )^{-1}]_{kl}=([\Ib+m_{2n}(z)\tilde\bSigma]^{-1})_{kl}+O_p(n^{-1/6+c}),
		\]
		where the $O_p$ is  uniform on $\lambda$.
		Note that the lower-right $r\times r$ block of $\tilde\bSigma$ is  the identity matrix. Then,
		\[
		\bigg(\Ib-n^{-1}\bGamma_1^\top\check\Zb[z\Ib-\hat{\mathcal{Q}}]^{-1}\check\Zb^\top\bGamma_1\bigg)^{-1}=\bigg[(\Ib+m_{2n}(z)\Ib)\bigg]^{-1}+O_p(n^{-1/6+c}).
		\]
		Following the proof of Lemma \ref{m_{1n}}, one can get that
		\[
		-C(\log n)^{-1}\le\operatorname{Re}m_{2n}(z_2)\le- c(\log n)^{-1},\quad \operatorname{Im}m_{2n}(z_2)=o(n^{-2/3}).
		\]
		Eventually, we have uniformly on $\lambda$ that
		\begin{equation}\label{PI 22}
			n^{-1}\bGamma_1^\top\check\Zb[z_2\Ib-\hat{\mathcal{Q}}]^{-1}\check\Zb^\top\bGamma_1-m_{2n}(z_2)\Ib\le O_p(n^{-1/6+c}).
		\end{equation}
		
		Now, combining (\ref{delocalization}), (\ref{PI 21}) and (\ref{PI 22}), we can write
		\[
		\bLambda_1^{-1}-\bPi(\lambda)=\bLambda_1^{-1}-\frac{1}{\lambda-\hat\varphi_1}\frac{1}{n}\bGamma_1^\top\check\Zb\hat\bbeta_1\hat\bbeta_1^\top\check\Zb^\top\bGamma_1+m_{2n}(z_2)\Ib+O_p(n^{-1/12+c}),
		\]
		and  $n^{-1}\|\bGamma_1^\top\check{\Zb}\hat\bbeta_1\|^2\le O_p(n^{-3/4+c})$ uniformly on $\lambda$ for arbitrary $c$. 
		Therefore, when $\lambda=\hat\varphi_1-n^{-2/3+c}$,  $m_{2n}(z_2)\Ib $ dominates and all the eigenvalues of $\bLambda_1^{-1}-\bPi(\lambda)$ are negative. However, as $\lambda$ approaches  $\hat\varphi_1$, the second term $(\lambda-\hat\varphi_1)^{-1}n^{-1}\bGamma_1^\top\check\Zb\hat\bbeta_1\hat\bbeta_1^\top\check\Zb^\top\bGamma_1$ will dominate and the largest eigenvalue of $\bLambda_1^{-1}-\bPi(\lambda)$ will become positive. Then, by continuity, there must be some $\lambda\in[\hat\varphi_1-n^{-2/3+c},\hat\varphi_1]$ such that the largest eigenvalue of $\bLambda_1^{-1}-\bPi(\lambda)$ is equal to 0, i.e., $g(\lambda)=0$. The theorem is  verified.

	\end{proof}

\end{appendices}

\bibliographystyle{model2-names} 
\bibliography{Ref}   

\newpage
\begin{appendices}
	\begin{center}
		\begin{Huge}
	Errata on ``Testing the number of common factors by bootstrapped sample covariance matrix  in high-dimensional factor models''
		\end{Huge}
	\end{center}
	
	This is an errata on the published version ``Yu, L., Zhao, P., \& Zhou, W. (2025). Testing the Number of Common Factors by Bootstrapped Sample Covariance Matrix in High-Dimensional Factor Models. {\it Journal of the American Statistical Association}, 120(549), 448–459.''
	
	In the proof of Lemma 3.2, we mistakenly used $\mathbb{E}[w_j^h(w_j-1)]=O(n^{-1})$ for $h=0,1$, while in fact $\mathbb{E}[w_j(w_j-1)]=1+O(n^{-1})$ for $h=1$; see the first paragraph on page 12 of the supplementary materials. Then, following the proof procedure of the original paper, the representation of $\hat\zeta_i$ in Lemma 3.2 should be revised to
	\[
	\begin{split}
		\hat\zeta_i-\frac{\theta_i}{\lambda_i}=&\frac{1}{n}\sum_{j=1}^n(w_j-1)+\frac{\mathrm{tr}\bLambda_2}{n\lambda_i}\times\mathbb{E}[w_1(w_1-1)]
		+o_p\bigg(\frac{\mathrm{tr}\bLambda_2}{n\lambda_i}\bigg)+o_p(n^{-1/2})\\
		=&\frac{1}{n}\sum_{j=1}^n(w_j-1)+\frac{\mathrm{tr}\bLambda_2}{n\lambda_i}
		+o_p\bigg(\frac{\mathrm{tr}\bLambda_2}{n\lambda_i}\bigg)+o_p(n^{-1/2}),
	\end{split}
	\]
	to include the first-order Taylor expansion rather than the second-order expansion as in the original Lemma. Correspondingly, following the proof of Theorem 3.3 in the original paper, the last $o_{\rm p}$ term in (4) of Theorem 3.3 should be revised to $o_{\rm p}((n\lambda_i)^{-1}p)$ instead of $o_{\rm p}((n\lambda_i)^{-2}p^2)$, since it arises from an error term of the form $o_{\rm p}(\hat\zeta_i-\lambda_i^{-1}\theta_i)$ in the proof. Similarly, (5) in Theorem 3.5 now holds if $(n\lambda_i)^{-1}p=o(n^{-1/2})$ instead of $(n\lambda_i)^{-1}p=o(n^{-1/4})$, while (6) holds if $(n\lambda_i)^{-1}p\gg n^{-1/2}$. A more fundamental revision stemming from this mistake concerns Corollary 1 on the bootstrap bias. The bias claimed in the original paper was based on a discrepancy between the tail probabilities with and without bootstrap when approximating with a first-order expansion. With the revised Lemma 3.2, the leading bias terms with and without bootstrap agree; the remaining discrepancy appears only at the second-order level, which is negligible relative to the first-order term. We rewrite Corollary 1 as follows and prove it in this Errata.

	\begin{cor}[Bias of bootstrap]\label{rev cor: bootstrap bias}
		Suppose Assumptions 1, 2 hold as $\min\{n,p\}\rightarrow\infty$, and 
		\begin{equation}\label{bias condition}		\sqrt{n}\bigg(\frac{p^3}{n^3\lambda_i^3}+\frac{p^2}{n^2\lambda_i^3}+\frac{p}{n\lambda_i^3}\bigg)=o(1),\quad i\le r.
		\end{equation}  
		Then, for any constant $s\in\mathbb{R}$ and $1\le i\le r$, under the standard bootstrap, we have
		\begin{equation}\label{rev bias hat}
			\begin{split}
				&\mathbb{P}^*\bigg(\sqrt{n}\times\lambda_i^{-1}(\hat\lambda_i-\tilde\lambda_i)\le s\bigg)\\
				=&F_G\Bigg(
				\frac{s}{\sqrt{\xi_i}}
				-\sqrt{\frac{n}{\xi_i}}\Big(
				\frac{\mathrm{tr}\bLambda_2}{n\lambda_i}
				+\frac{\mathrm{tr}(\bLambda_2^2)}{n\lambda_i^2}
				+\big(\frac{\mathrm{tr}\bLambda_2}{n\lambda_i}\big)^2
				\big[\mathbb{E}\{w_1^2(w_1-1)\}-1\big]
				\Big)
				\Bigg)+o_p(1),
			\end{split}
		\end{equation}
		where $F_G(\cdot)$ is the cumulative distribution function (CDF) of a standard normal random variable. (\ref{rev bias hat}) also holds under the multiplier bootstrap by replacing $\xi_i$ with $\xi_i+1$. On the other hand, without bootstrap, we have
		\begin{equation}\label{rev bias tilde}
			\begin{split}
				&\mathbb{P}\bigg(\sqrt{n}\times\lambda_i^{-1}(\tilde\lambda_i-\lambda_i)\le s\bigg)
				=F_G\Bigg(
				\frac{s}{\sqrt{\xi_i}}
				-\sqrt{\frac{n}{\xi_i}}\Bigg\{
				\frac{\mathrm{tr}\bLambda_2}{n\lambda_i}
				+\frac{\mathrm{tr}(\bLambda_2^2)}{n\lambda_i^2}
				\Bigg\}
				\Bigg)+o(1).
			\end{split}
		\end{equation}
		If \eqref{bias condition} does not hold,  the corresponding CDFs in \eqref{rev bias hat} and \eqref{rev bias tilde} tend to $0$ for any fixed $s$.
	\end{cor}
	Thus, under the standard bootstrap, the two tail probabilities match asymptotically whenever $(n\lambda_i)^{-1}p=O(n^{-1/2})$, because $\sqrt{n}\,\mathrm{tr}(\bLambda_2^2)/(n\lambda_i^2)=o(1)$ and $\sqrt{n}((n\lambda_i)^{-1}\mathrm{tr}(\bLambda_2))^2=o(1)$ under this condition. When $(n\lambda_i)^{-1}p\gg n^{-1/2}$, the limits of \eqref{rev bias hat} and \eqref{rev bias tilde} are both $0$, but this does not imply that the bootstrap accurately approximates the limiting distribution of $\tilde\lambda_i$ in this case. For the multiplier bootstrap, the bias arises mainly because the asymptotic variances of the two limiting distributions do not match, as claimed in the original paper. 
	
	Lastly, there is a typo in (12) in Lemma 4.2 of the original paper, where ``$>x$'' should be revised to ``$\le x$''.

	We now prove the revised Corollary \ref{rev cor: bootstrap bias}.

	\begin{proof}
		Following the derivation of the original proof, we still have
		\begin{equation}\label{eq: ratio for bias}
			\begin{split}
				\frac{\hat\lambda_i}{\tilde\lambda_i}
				=&[1+o_p(1)]\times \bigg[\frac{1}{n}\bgamma_i^\top\Zb(\Wb-\Ib)\Zb^\top\bgamma_i
				-\frac{1}{n}\mathrm{tr}\Wb+1+\hat\zeta_i-\frac{\theta_i}{\lambda_i}\bigg]
				+1+o_p(n^{-1/2}).
			\end{split}
		\end{equation}
		The difference is that we construct a revised representation for $\hat\zeta_i-\lambda_i^{-1}\theta_i$. Following the argument leading to (B.8) of the original supplementary materials, we have
		\begin{equation*}
			\begin{split}
				\frac1n\sum_{j=1}^n\frac{w_j-1}{1-w_j D_n}
				=&\frac1n\sum_{j=1}^n(w_j-1)
				+D_n\,\frac1n\sum_{j=1}^nw_j(w_j-1)
				+D_n^2\,\frac1n\sum_{j=1}^nw_j^2(w_j-1)
				+O_p(D_n^3),
			\end{split}
		\end{equation*}
		where
		\[
		D_n:=\frac1{n\theta_i}\sum_{k=1}^p\frac{\lambda_{r+k}}{1-\theta_i^{-1}\lambda_{r+k}\hat\zeta_i}.
		\]
		(B.8) of the original supplementary materials shows that 
		\[
		\frac{1}{n}\sum_{j=1}^nw_j^h(w_j-1)-\mathbb{E}[w_1^h(w_1-1)]=O_p(n^{-1/2}),\quad h=0,1,2.
		\]
		Since $\mathbb{E}[w_1(w_1-1)]=1+O(n^{-1})$, we obtain
		\begin{equation}\label{eq: zeta tilde mu}
			\begin{split}
				\hat\zeta_i-\frac{\theta_i}{\lambda_i}
				=\frac1n\sum_{j=1}^n(w_j-1)
				+D_n
				+D_n^2\,\mathbb{E}[w_1^2(w_1-1)]
				+O_p(D_n^3)+o_p(n^{-1/2}).
			\end{split}
		\end{equation}
		Next expand $D_n$. Since $\hat\zeta_i=1+O_p\big(\mathrm{tr}(\bLambda_2)/(n\lambda_i)+n^{-1/2}\big)$, we have
		\[
		\frac{\lambda_{r+k}}{1-\theta_i^{-1}\lambda_{r+k}\hat\zeta_i}
		=\lambda_{r+k}+\theta_i^{-1}\lambda_{r+k}^2\hat\zeta_i+O(\theta_i^{-2}\lambda_{r+k}^3).
		\]
		Thus, 
		\begin{equation}\label{Dn}
			D_n=\frac{\mathrm{tr}(\bLambda_2)}{n\theta_i}+\frac{\mathrm{tr}(\bLambda_2^2)}{n\theta_i^2}\hat\zeta_i+O\Big(\frac{p}{n\theta_i^3}\Big).
		\end{equation}
		Further note that
		\begin{equation}\label{eq: theta identity}
			\frac{\theta_i}{\lambda_i}-1
			=\frac1{n\lambda_i}\sum_{k=1}^p\frac{\lambda_{r+k}}{1-\lambda_i^{-1}\lambda_{r+k}}=\frac{\mathrm{tr}(\bLambda_2)}{n\lambda_i}+\frac{\mathrm{tr}(\bLambda_2^2)}{n\lambda_i^2}+O\Big(\frac{p}{n\lambda_i^3}\Big).
		\end{equation}
		Therefore, 
		\[
		\begin{split}
			\frac{\lambda_i}{\theta_i}=&\{1+(\lambda_i^{-1}\theta_i-1)\}^{-1}=1-(\lambda_i^{-1}\theta_i-1)+(\lambda_i^{-1}\theta_i-1)^2+O((\lambda_i^{-1}\theta_i-1)^3)\\
			=&1-\frac{\mathrm{tr}(\bLambda_2)}{n\lambda_i}-\frac{\mathrm{tr}(\bLambda_2^2)}{n\lambda_i^2}+\Big(\frac{\mathrm{tr}(\bLambda_2)}{n\lambda_i}\Big)^2+O\bigg(\frac{p}{n\lambda_i^3}+\frac{p}{n\lambda_i}\frac{p}{n\lambda_i^2}+\frac{p^3}{n^3\lambda_i^3}\bigg).
		\end{split}
		\]
		Returning to \eqref{Dn}, we obtain
		\begin{equation}\label{eq: tilde delta expansion}
			\begin{split}
				D_n=&\frac{\mathrm{tr}(\bLambda_2)}{n\lambda_i}\frac{\lambda_i}{\theta_i}+\frac{\mathrm{tr}(\bLambda_2^2)}{n\lambda_i^2}\Big(\frac{\lambda_i}{\theta_i}\Big)^2\hat\zeta_i+O\Big(\frac{p}{n\lambda_i^3}\Big)\\
				=&\frac{\mathrm{tr}(\bLambda_2)}{n\lambda_i}
				+\frac{\mathrm{tr}(\bLambda_2^2)}{n\lambda_i^2}
				-\bigg(\frac{\mathrm{tr}(\bLambda_2)}{n\lambda_i}\bigg)^2+O_p\Bigg(
				\Big(\frac{p}{n\lambda_i}\Big)^3
				+\frac{p}{n\lambda_i}\cdot\frac{p}{n\lambda_i^2}
				+\frac{1}{\sqrt{n}}\frac{p}{n\lambda_i^2}+\frac{p}{n\lambda_i^3}
				\Bigg).
			\end{split}
		\end{equation}
		Substituting \eqref{eq: tilde delta expansion} into \eqref{eq: zeta tilde mu} yields
		\begin{equation}\label{eq: zeta in delta}
			\begin{split}
				\hat\zeta_i-\frac{\theta_i}{\lambda_i}
				=&\frac1n\sum_{j=1}^n(w_j-1)
				+\frac{\mathrm{tr}(\bLambda_2)}{n\lambda_i}
				+\frac{\mathrm{tr}(\bLambda_2^2)}{n\lambda_i^2}
				+\bigg(\frac{\mathrm{tr}(\bLambda_2)}{n\lambda_i}\bigg)^2
				\big(\mathbb{E}[w_1^2(w_1-1)]-1\big)\\
				&+O_p\bigg(
				\Big(\frac{p}{n\lambda_i}\Big)^3
				+\frac{p}{n\lambda_i}\cdot\frac{p}{n\lambda_i^2}
				+\frac{p}{n\lambda_i^3}
				\bigg)+o_p(n^{-1/2}).
			\end{split}
		\end{equation}
		Substituting \eqref{eq: zeta in delta} into \eqref{eq: ratio for bias} and using $\tilde\lambda_i/\lambda_i=1+o_p(1)$ yields
		\begin{equation}\label{eq: boot fluctuation}
			\begin{split}
				\sqrt{n}\,\lambda_i^{-1}(\hat\lambda_i-\tilde\lambda_i)
				=&[1+o_p(1)]\Bigg\{
				\frac1{\sqrt n}\bgamma_i^\top\Zb(\Wb-\Ib)\Zb^\top\bgamma_i\\
				&\qquad+\sqrt n\Bigg(
				\frac{\mathrm{tr}(\bLambda_2)}{n\lambda_i}
				+\frac{\mathrm{tr}(\bLambda_2^2)}{n\lambda_i^2}
				+\bigg(\frac{\mathrm{tr}(\bLambda_2)}{n\lambda_i}\bigg)^2
				\big[\mathbb{E}\{w_1^2(w_1-1)\}-1\big]
				\Bigg)
				\Bigg\}\\
				&+O_p\Bigg\{\sqrt n\Big(
				\frac{p^3}{n^3\lambda_i^3}
				+\frac{p}{n\lambda_i}\cdot\frac{p}{n\lambda_i^2}+\frac{p}{n\lambda_i^3}
				\Big)\Bigg\}+o_p(1).
			\end{split}
		\end{equation}
		Under \eqref{bias condition}, the $O_p$ remainder is $o_p(1)$. By the proof in the original supplementary materials, conditional on $\Xb$ and under the standard bootstrap, we have
		\[
		\frac1{\sqrt n}\bgamma_i^\top\Zb(\Wb-\Ib)\Zb^\top\bgamma_i
		\overset{d^*}{\longrightarrow}\mathcal{N}(0,\xi_i),
		\]
		while under the multiplier bootstrap the limiting variance is $\xi_i+1$. Combining this limit with \eqref{eq: boot fluctuation} yields \eqref{rev bias hat} and the multiplier analogue.
		
		On the other hand, for $\tilde\lambda_i$ (without bootstrap), by Lemma 3.4-(a) and Lemma C.3 in the original paper,
		\begin{equation*}
			\sqrt{n}\Big(\frac{\tilde\lambda_i}{\theta_i}-1\Big)\overset{d}{\longrightarrow}\mathcal{N}(0,\xi_i).
		\end{equation*}
		Using \eqref{eq: theta identity}, we have
		\begin{align*}
			\sqrt{n}\,\lambda_i^{-1}(\tilde\lambda_i-\lambda_i)
			&=\sqrt{n}\,\frac{\theta_i}{\lambda_i}\Big(\frac{\tilde\lambda_i}{\theta_i}-1\Big)
			+\sqrt{n}\Big(\frac{\theta_i}{\lambda_i}-1\Big)\\
			&=\sqrt{n}\,\frac{\theta_i}{\lambda_i}\Big(\frac{\tilde\lambda_i}{\theta_i}-1\Big)
			+\sqrt{n}\Bigg(
			\frac{\mathrm{tr}(\bLambda_2)}{n\lambda_i}
			+\frac{\mathrm{tr}(\bLambda_2^2)}{n\lambda_i^2}
			\Bigg)
			+O\Bigg(\sqrt{n}\times\frac{p}{n\lambda_i^3}\Bigg).
		\end{align*}
		Then, \eqref{rev bias tilde} holds by the fact that $\theta_i/\lambda_i\rightarrow 1$ and the condition \eqref{bias condition}. If \eqref{bias condition} does not hold, the CDFs tend to $0$ for any fixed $s$, because 
		\[
		\sqrt{n}\times \min\bigg\{\frac{\mathrm{tr}(\bLambda_2)}{n\lambda_i},\frac{\mathrm{tr}(\bLambda_2^2)}{n\lambda_i^2}, \bigg(\frac{\mathrm{tr}(\bLambda_2)}{n\lambda_i}\bigg)^2\bigg\}\rightarrow \infty,
		\]
		which concludes the corollary.
	\end{proof}

\end{appendices}

\end{document}